%% file: doubly-periodic.tex
\documentclass[a4,10pt]{article}
\usepackage{amsmath,amscd,amssymb}

\input{notation}
\input{new_theorem}

\begin{document}

\title{Doubly periodic monopoles
and $q$-difference modules}

\author{Takuro Mochizuki}
\date{}
\maketitle

\begin{abstract}

An interesting theme in complex differential geometry
is to find a correspondence between 
algebraic objects and differential geometric objects.
One of the most attractive is 
the non-abelian Hodge theory of Simpson.
In this paper, 
pursuing an analogue of the non-abelian Hodge theory
in the context of $q$-difference modules,
we study Kobayashi-Hitchin correspondences between 
doubly periodic monopoles and 
parabolic $q$-difference modules,
depending on twistor parameters.

\vspace{.1in}
\end{abstract}

\input{1}

\input{reference}
\noindent
{\em Address\\
Research Institute for Mathematical Sciences,
Kyoto University,
Kyoto 606-8502, Japan\\
takuro@kurims.kyoto-u.ac.jp
}

\end{document}

%% file: notation.tex
\newcommand{\nbiga}{\mathcal{A}}
\newcommand{\nbigb}{\mathcal{B}}
\newcommand{\nbigc}{\mathcal{C}}

\newcommand{\nbige}{\mathcal{E}}
\newcommand{\nbigf}{\mathcal{F}}
\newcommand{\nbigg}{\mathcal{G}}
\newcommand{\nbigh}{\mathcal{H}}

\newcommand{\nbigk}{\mathcal{K}}
\newcommand{\nbigl}{\mathcal{L}}
\newcommand{\nbigm}{\mathcal{M}}
\newcommand{\nbign}{\mathcal{N}}
\newcommand{\nbigo}{\mathcal{O}}
\newcommand{\nbigp}{\mathcal{P}}

\newcommand{\nbigr}{\mathcal{R}}
\newcommand{\nbigs}{\mathcal{S}}

\newcommand{\nbigu}{\mathcal{U}}
\newcommand{\nbigv}{\mathcal{V}}
\newcommand{\nbigw}{\mathcal{W}}

\newcommand{\nbigz}{\mathcal{Z}}

\newcommand{\proj}{\mathbb{P}}
\newcommand{\seisuu}{{\mathbb Z}}
\newcommand{\rnum}{{\mathbb Q}}

\newcommand{\cnum}{{\mathbb C}}
\newcommand{\real}{{\mathbb R}}

\newcommand{\EE}{\mathbb{E}}

\newcommand{\LL}{\mathbb{L}}
\newcommand{\VV}{\mathbb{V}}

\newcommand{\gbiga}{\mathfrak A}

\newcommand{\gbigc}{\mathfrak C}

\newcommand{\gbigf}{\mathfrak F}

\newcommand{\gbigi}{\mathfrak I}

\newcommand{\gbigu}{\mathfrak U}
\newcommand{\gbigv}{\mathfrak V}

\newcommand{\gminia}{\mathfrak a}

\newcommand{\gminie}{\mathfrak e}

\newcommand{\gminig}{\mathfrak g}

\newcommand{\gminio}{\mathfrak o}
\newcommand{\gminip}{\mathfrak p}
\newcommand{\gminiq}{\mathfrak q}

\newcommand{\gminit}{\mathfrak t}

\newcommand{\gminiv}{\mathfrak v}

\newcommand{\vece}{{\boldsymbol e}}

\newcommand{\vecv}{{\boldsymbol v}}
\newcommand{\vecu}{{\boldsymbol u}}
\newcommand{\vecw}{{\boldsymbol w}}

\newcommand{\veca}{{\boldsymbol a}}
\newcommand{\vecb}{{\boldsymbol b}}

\newcommand{\vecs}{{\boldsymbol s}}
\newcommand{\vect}{{\boldsymbol t}}

\newcommand{\vecn}{{\boldsymbol n}}

\newcommand{\vecV}{{\boldsymbol V}}
\newcommand{\vecA}{{\boldsymbol A}}

\newcommand{\lrarr}{\longrightarrow}

\newcommand{\pf}{{\bf Proof}\hspace{.1in}}
\newcommand{\qed}{\mbox{\rule{1.2mm}{3mm}}}

\def\Hom{\mathop{\rm Hom}\nolimits}

\def\End{\mathop{\rm End}\nolimits}

\def\Image{\mathop{\rm Im}\nolimits}

\def\Re{\mathop{\rm Re}\nolimits}

\def\Gr{\mathop{\rm Gr}\nolimits}

\def\GL{\mathop{\rm GL}\nolimits}

\def\rank{\mathop{\rm rank}\nolimits}
\def\Spec{\mathop{\rm Spec}\nolimits}

\def\length{\mathop{\rm length}\nolimits}
\def\Gr{\mathop{\rm Gr}\nolimits}
\def\Sym{\mathop{\rm Sym}\nolimits}

\def\ad{\mathop{\rm ad}\nolimits}

\def\ord{\mathop{\rm ord}\nolimits}
\def\degpar{\mathop{\rm par\textrm{-}deg}\nolimits}

\def\Tr{\mathop{\rm Tr}\nolimits}
\def\vol{\mathop{\rm vol}\nolimits}
\def\dvol{\mathop{\rm dvol}\nolimits}
\def\Diff{\mathop{\rm Diff}\nolimits}

\def\id{\mathop{\rm id}\nolimits}

\def\gcd{\mathop{\rm g.c.d.}\nolimits}

\newcommand{\del}{\partial}
\newcommand{\delbar}{\overline{\del}}

\newcommand{\pardeg}{\degpar}

\newcommand{\vecP}{\boldsymbol P}

\newcommand{\ttt}{{\tt t}}

\newcommand{\ttK}{{\tt K}}
\newcommand{\ttG}{{\tt G}}
\newcommand{\ttH}{{\tt H}}
\newcommand{\ttE}{{\tt E}}
\newcommand{\ttD}{{\tt D}}

\newcommand{\ttP}{{\tt P}}

\newcommand{\ttF}{{\tt F}}

\newcommand{\barz}{\overline{z}}
\newcommand{\zbar}{\barz}

\newcommand{\baralpha}{\overline{\alpha}}
\newcommand{\alphabar}{\baralpha}
\newcommand{\barlambda}{\overline{\lambda}}
\newcommand{\lambdabar}{\barlambda}

\newcommand{\etabar}{\overline{\eta}}
\newcommand{\xibar}{\overline{\xi}}

\newcommand{\Par}{{\mathcal Par}}
\newcommand{\Sp}{{\mathcal Sp}}

\newcommand{\paramap}{\gminip}

\def\reg{\mathop{\rm reg}\nolimits}

\newcommand{\vecnbigl}{{\boldsymbol{\mathcal L}}}

\newcommand{\vecnbigk}{{\boldsymbol{\mathcal K}}}

\newcommand{\closedopen}[2]{[#1,#2[}
\newcommand{\openclosed}[2]{]#1,#2]}
\newcommand{\openopen}[2]{]#1,#2[}
\newcommand{\closedclosed}[2]{[#1,#2]}

\newcommand{\vecvhat}{\widehat{\vecv}}

\newcommand{\rhotilde}{\widetilde{\rho}}

\newcommand{\Etilde}{\widetilde{E}}

\newcommand{\thetatilde}{\widetilde{\theta}}

\newcommand{\Vtilde}{\widetilde{V}}
\newcommand{\nablatilde}{\widetilde{\nabla}}
\newcommand{\vecvtilde}{\widetilde{\vecv}}

\newcommand{\vecutilde}{\widetilde{\vecu}}
\newcommand{\zerohat}{\widehat{0}}

\newcommand{\fhat}{\widehat{f}}
\newcommand{\vecVhat}{\widehat{\vecV}}

\newcommand{\wbar}{\overline{w}}

\newcommand{\nbigetilde}{\widetilde{\nbige}}

\newcommand{\htilde}{\widetilde{h}}

\newcommand{\vtilde}{\widetilde{v}}

\newcommand{\ftilde}{\widetilde{f}}
\newcommand{\utilde}{\widetilde{u}}
\newcommand{\Ftilde}{\widetilde{F}}

\newcommand{\stilde}{\widetilde{s}}
\newcommand{\mubar}{\overline{\mu}}

\newcommand{\nbigltilde}{\widetilde{\nbigl}}

\newcommand{\pihat}{\widehat{\pi}}

\newcommand{\Fhat}{\widehat{F}}

\newcommand{\nbigvtilde}{\widetilde{\nbigv}}
\newcommand{\nbigvhat}{\widehat{\nbigv}}

\newcommand{\Utilde}{\widetilde{U}}

\newcommand{\Ltilde}{\widetilde{L}}

\def\ord{\mathop{\rm ord}\nolimits}

\def\Gal{\mathop{\rm Gal}\nolimits}

\def\Mod{\mathop{\rm Mod}\nolimits}

\def\Vect{\mathop{\rm Vect}\nolimits}

\def\ad{\mathop{\rm ad}\nolimits}

\def\Aut{\mathop{\rm Aut}\nolimits}

\def\sgn{\mathop{\rm sign}\nolimits}
\def\cov{\mathop{\rm cov}\nolimits}
\def\Vol{\mathop{\rm Vol}\nolimits}
\def\Slope{\mathop{\rm Slope}\nolimits}

\def\LFM{\mathop{\rm LFM}\nolimits}
\def\an{\mathop{\rm an}\nolimits}
\def\AHN{\mathop{\rm AHN}\nolimits}

\newcommand{\Ubar}{\overline{U}}

\newcommand{\ubar}{\overline{u}}

\newcommand{\nbigutilde}{\widetilde{\nbigu}}
\newcommand{\nbigstilde}{\widetilde{\nbigs}}

\newcommand{\Phat}{\widehat{P}}

\newcommand{\Atilde}{\widetilde{A}}

\newcommand{\Ybar}{\overline{Y}}

\newcommand{\betabar}{\overline{\beta}}

\newcommand{\nbigmlambda}{\nbigm^{\lambda}}

\newcommand{\Hhat}{\widehat{H}}

\newcommand{\phitilde}{\widetilde{\phi}}

\newcommand{\nbigmbar}{\overline{\nbigm}}

\newcommand{\vecVtilde}{\widetilde{\vecV}}

\newcommand{\vecetilde}{\widetilde{\vece}}
\newcommand{\etilde}{\widetilde{e}}

\newcommand{\inftyhat}{\widehat{\infty}}
\newcommand{\nbigmbarlambda}{\nbigmbar^{\lambda}}

\newcommand{\chitilde}{\widetilde{\chi}}

\newcommand{\LLtilde}{\widetilde{\LL}}
\newcommand{\VVtilde}{\widetilde{\VV}}

\newcommand{\nuhat}{\widehat{\nu}}

\newcommand{\Ihat}{\widehat{I}}

\newcommand{\sfv}{{\sf v}}
\newcommand{\sff}{{\sf f}}

\newcommand{\sfp}{{\sf p}}

\newcommand{\sfQ}{{\sf Q}}
\newcommand{\slQ}{{\textsl Q}}

\newcommand{\sfA}{{\sf A}}

\newcommand{\slq}{{\textsl q}}

\newcommand{\ttu}{{\tt u}}
\newcommand{\ttv}{{\tt v}}

\newcommand{\slw}{\textsl w}

\newcommand{\tts}{{\tt s}}

\newcommand{\ttg}{{\tt g}}

\newcommand{\sle}{{\textsl e}}
\newcommand{\tte}{{\tt e}}
\newcommand{\ttx}{{\tt x}}
\newcommand{\tty}{{\tt y}}
\newcommand{\ttp}{{\tt p}}
\newcommand{\tti}{{\tt i}}

\newcommand{\ttU}{{\tt U}}
\newcommand{\ttA}{{\tt A}}
\newcommand{\ttB}{{\tt B}}
\newcommand{\ttC}{{\tt C}}
\newcommand{\ttR}{{\tt R}}
\newcommand{\ttM}{{\tt M}}
\newcommand{\ttm}{{\tt m}}

\newcommand{\sfL}{{\sf L}}
\newcommand{\vecsfL}{\boldsymbol{\sf L}}
\newcommand{\ttL}{{\tt L}}

\newcommand{\ttl}{{\tt l}}

\newcommand{\zetatilde}{\widetilde{\zeta}}
\newcommand{\sfz}{{\sf z}}
\newcommand{\sfzbar}{\overline{\sfz}}
\newcommand{\sfw}{{\sf w}}
\newcommand{\sfwbar}{\overline{{\sf w}}}
\newcommand{\nbigubar}{\overline{\nbigu}}
\newcommand{\sfa}{{\sf a}}
\newcommand{\sfb}{{\sf b}}

\newcommand{\sfLtilde}{\widetilde{\sfL}}

\newcommand{\sfvtilde}{\widetilde{\sfv}}

\newcommand{\ttgbar}{\overline{\ttg}}
\newcommand{\ttubar}{\overline{\ttu}}
\newcommand{\ttvbar}{\overline{\ttv}}
\newcommand{\ttxbar}{\overline{\ttx}}
\newcommand{\ttybar}{\overline{\tty}}
\newcommand{\vecLL}{\pmb\LL}
\newcommand{\gminietilde}{\widetilde{\gminie}}

\newcommand{\ttUbar}{\overline{\ttU}}
\newcommand{\slwbar}{\overline{\slw}}

\newcommand{\vecVV}{\pmb\VV}

\newcommand{\sfbbar}{\overline{\sfb}}
\newcommand{\su}{\mathfrak{su}}
\newcommand{\sfatilde}{\widetilde{\sfa}}
\newcommand{\sfbtilde}{\widetilde{\sfb}}

\newcommand{\gbigvhat}{\widehat{\gbigv}}
\newcommand{\ttGhat}{\widehat{\ttG}}
\newcommand{\ttKtilde}{\widetilde{\ttK}}

\newcommand{\ttEtilde}{\widetilde{\ttE}}

%% file: new_theorem.tex
\newtheorem{thm}{Theorem}[section]
\newtheorem{cor}[thm]{Corollary}

\newtheorem{rem}[thm]{Remark}
\newtheorem{lem}[thm]{Lemma}
\newtheorem{prop}[thm]{Proposition}
\newtheorem{df}[thm]{Definition}

\newtheorem{example}[thm]{Example}
\newtheorem{condition}[thm]{Condition}
\newtheorem{assumption}[thm]{Assumption}

\newtheorem{notation}[thm]{Notation}

%% file: 1.tex
\section{Introduction}

In \cite{Mochizuki-difference-modules},
we studied Kobayashi-Hitchin correspondences
between periodic monopoles and 
difference modules with parabolic structure、
depending on the twistor parameters.
It is an interesting variant of Kobayashi-Hitchin correspondences
for harmonic bundles
pioneered by 
Corlette \cite{corlette},
Donaldson \cite{don2},
Hitchin \cite{Hitchin-self-duality}
and particularly Simpson 
\cite{Simpson88, Simpson90, s5, s3, s4}.
See \cite[\S1]{Mochizuki-difference-modules}
for more background.

In this paper, as another interesting variant,
we shall study Kobayashi-Hitchin correspondences
between doubly-periodic monopoles
and $\gminiq$-difference modules,
depending on the twistor parameters.

\subsection{Meromorphic doubly periodic monopoles}

Let $\Gamma$ be any lattice in $\real^2$.
It naturally acts on $\real^2$ by the addition.
We obtain the induced action of $\Gamma$
on $\real\times\real^2$.
Let $\nbigm$ denote the quotient space.
It is naturally equipped
with the metric $g_{\nbigm}$
induced by the Euclidean metric of $\real^3$.
Let $Z$ be a finite subset in $\nbigm$.

Let $E$ be a complex vector bundle on $\nbigm\setminus Z$
equipped with a Hermitian metric $h$,
a unitary connection $\nabla$,
and an anti-self-adjoint endomorphism $\phi$
satisfying the Bogomolny equation
\[
 F(\nabla)=\ast\nabla\phi.
\]
Here, $F(\nabla)$ denotes the curvature of $\nabla$,
and $\ast$ denotes the Hodge star operator with respect to $g_{\nbigm}$.
Such a tuple $(E,h,\nabla,\phi)$
is called a doubly periodic monopole
because it can be regarded as
a singular monopole on $\real^3$ with periodicity
in two directions.
It is called {\em meromorphic} in this paper
if the following is satisfied:
\begin{itemize}
\item
Each point of $Z$ is Dirac type singularity
of the monopole.
\item
There exists a compact subset $C$ which contains $Z$
such that
$F(\nabla)$ is bounded
with respect to $h$ and $g_{\nbigm}$
on $\nbigm\setminus C$.
\end{itemize}

\subsubsection{Examples}

We use the coordinate system
$(y_0,y_1,y_2)$ on $\real\times\real^2$.
We may regard $\real^2_{(y_1,y_2)}$ as $\cnum$
by the complex coordinate $z=y_1+\sqrt{-1}y_2$,
and we regard $T^0:=\cnum/\Gamma$ as an elliptic curve.
It is equipped with the Euclidean metric $dz\,d\zbar$.
The Riemannian manifold $\nbigm$
is naturally identified with the product $\real\times T^0$.

Take a holomorphic line bundle $L_m$
of degree $-m$,
i.e., $\int_{T^0}c_1(L_m)=-m$.
There exists a Hermitian metric $h_{L_{m}}$
such that the curvature of the Chern connection
$\nabla_{h_{L_m}}$ is equal to
$\frac{-m\pi}{\Vol(T^0)}\,dz\,d\zbar$.
Let $p:\real\times T^0\lrarr T^0$
denote the projection.
We obtain
$(E_m,h_m,\nabla_m)$
as the pull back of
$(L_m,h_{L_m},\nabla_{L_m})$.
Set $\phi_m:=-\sqrt{-1}\frac{2\pi m}{\Vol(T^0)}y_0$.
Then, 
$(E_m,h_m,\nabla_m,\phi_m)$
is a meromorphic doubly periodic monopole.

Let $\Gamma'\subset\Gamma$ be a sub-lattice
such that $|\Gamma/\Gamma'|=k$.
We set $T':=\cnum/\Gamma'$.
Let $T'\lrarr T^0$ be the induced covering of degree $k$.
Take a holomorphic line bundle
$L'_m$ of degree $m$ on $T'$.
Let $h_{L'_m}$ and $\nabla_{L_m'}$ be as above.
We obtain a monopole
$(E'_m,h_m',\nabla_{m}',\phi'_m)$
on $\real\times T'$.
Set $\omega=m/k$.
By taking the push-forward with respect to
the induced covering
$\real\times T'\lrarr\real\times T^0$,
we obtain a monopole
$(E_{\omega},h_{\omega},\nabla_{\omega},\phi_{\omega})$
of rank $k$ on $\real\times T^0$.

\vspace{.1in}

Let $\veca=(a_0,a_1,a_2)\in\real^3$.
Let $\underline{\cnum}\,e$ be the product line bundle
on $\nbigm$ with a global frame $e$.
Let $h$ be the metric determined by $h(e,e)=1$.
Let $\nabla_{\veca}$ and $\phi_{\veca}$ be 
determined by
\[
 \nabla_{\veca} e=e\sqrt{-1}(a_1\,dy_1+a_2\,dy_2),
\quad
 \phi_{\veca}=\sqrt{-1}a_0.
\]
Then, 
$(\underline{\cnum},h,\nabla_{\veca},\phi_{\veca})$
is a meromorphic monopole on $\nbigm$.

\subsection{Parabolic $\gminiq$-difference modules}
\label{subsection;19.1.26.1}

\subsubsection{$\gminiq$-difference modules}

Let $\gminiq\in\cnum^{\ast}$.
Let $\Phi^{\ast}$ be the automorphism of the algebra
$\cnum[y,y^{-1}]$ determined by
$\Phi^{\ast}(f)=f(\gminiq y)$.
A $\gminiq$-difference $\cnum[y,y^{-1}]$-module
is a $\cnum[y,y^{-1}]$-module $\vecV$
equipped with a $\cnum$-linear automorphism
$\Phi^{\ast}$
such that
$\Phi^{\ast}(fs)=\Phi^{\ast}(f)\Phi^{\ast}(s)$
for any $f\in \cnum[y,y^{-1}]$ and $s\in \vecV$.

We set $\nbiga_{\gminiq}:=
\bigoplus_{n\in\seisuu}
 \cnum[y,y^{-1}](\Phi^{\ast})^n$.
It is a non-commutative algebra
endowed with the multiplication
induced by 
$(\Phi^{\ast})^my^k=y^k\gminiq^{km}(\Phi^{\ast})^m$.
Then, $\gminiq$-difference modules
are equivalent to $\nbiga_{\gminiq}$-modules.

\begin{rem}
The automorphism $\Phi^{\ast}$
is extended to automorphisms of 
$R:=\cnum(\!(y)\!)$,
$\cnum(\!(y^{-1})\!)$ and $\cnum(y)$.
The notion of $\gminiq$-difference $R$-modules
are defined similarly.
\hfill\qed
\end{rem}

In this section,
we impose the following condition
to $\gminiq$-difference $\cnum[y,y^{-1}]$-modules $\vecV$
unless otherwise specified.
\begin{itemize}
\item
 It is torsion-free as $\cnum[y,y^{-1}]$-module.
\item
 There exists a free $\cnum[y,y^{-1}]$-submodule
 $V\subset \vecV$ of finite rank such that
 $V\otimes_{\cnum[y,y^{-1}]}\cnum(y)=\vecV\otimes_{\cnum[y,y^{-1}]}\cnum(y)$
 and 
 $\nbiga_{\gminiq}\cdot V=\vecV$.
\end{itemize}

\subsubsection{Parabolic $\gminiq$-difference $\cnum[y,y^{-1}]$-modules}

We introduce parabolic structure
on $\gminiq$-difference
$\cnum[y,y^{-1}]$-modules,
which consists of good parabolic structure at infinity
and parabolic structure at finite place.

\paragraph{Good parabolic structure at infinity}

Let $(\nbigvhat,\Phi^{\ast})$ be 
a $\gminiq$-difference $\cnum(\!(y)\!)$-module,
for which we always assume that
$\dim_{\cnum(\!(y)\!)}\nbigvhat<\infty$.
As known classically
(see \cite{van-der-Put-Singer, Sauloy2004, Soibelman-Vologodsky}),
there exists a slope decomposition of
$(\nbigvhat,\Phi^{\ast})=
 \bigoplus_{\omega\in\rnum}
 (\nbigvhat_{\omega},\Phi^{\ast})$
such that the following holds.
\begin{itemize}
\item
 Let $\omega=\ell/k$,
 where $\ell\in\seisuu$ and $k\in\seisuu_{>0}$.
 Then, 
 there exists a $\cnum[\![y]\!]$-lattice
 $\nbigl_{\omega}\subset
 \nbigvhat_{\omega}$
such that
 $y^{\ell}(\Phi^{\ast})^k\nbigl_{\omega}=\nbigl_{\omega}$.
\end{itemize}
Recall that a filtered bundle $\nbigp_{\ast}\nbigvhat$ over $\nbigvhat$
means an increasing sequence of 
$\cnum[\![y]\!]$-lattices
$\nbigp_{a}\nbigvhat\subset\nbigvhat$ $(a\in\real)$
such that
(i) $\nbigp_{a+n}\nbigvhat=y^{-n}\nbigp_a\nbigvhat$
for any $a\in\real$ and $n\in\seisuu$,
(ii) $\nbigp_a\nbigvhat=\bigcap_{b>a}\nbigp_b\nbigvhat$.
A filtered bundle 
$\nbigp_{\ast}\nbigvhat$ over $\nbigvhat$
is called good
if the following holds.
\begin{itemize}
\item
The filtration $\nbigp_{\ast}\nbigvhat$
is compatible with the slope decomposition,
i.e.,
 $\nbigp_{\ast}\nbigvhat=\bigoplus
 \nbigp_{\ast}\nbigvhat_{\omega}$.
\item
 $\Phi^{\ast}\nbigp_a(\nbigvhat_{\omega})
=\nbigp_{a+\omega}(\nbigvhat_{\omega})$ holds.
\end{itemize}

Let $\vecV$ 
be a $\gminiq$-difference $\cnum[y,y^{-1}]$-module.
We set
$\vecV_{|\zerohat}:=\vecV\otimes\cnum(\!(y)\!)$
and 
$\vecV_{|\inftyhat}:=\vecV\otimes\cnum(\!(y^{-1})\!)$.
Then, a good parabolic structure of $\vecV$ at infinity
is defined to be good filtered bundles
$\nbigp_{\ast}\vecV_{|\zerohat}$
and 
$\nbigp_{\ast}\vecV_{|\inftyhat}$
over $\vecV_{|\zerohat}$
and $\vecV_{|\inftyhat}$,
respectively.

\paragraph{Parabolic structure at finite place}

Set $y_{\alpha}:=y-\alpha$ for any $\alpha\in \cnum^{\ast}$.
For any subset $S\subset \cnum^{\ast}$,
let $\cnum[y,y^{-1}](\ast S)$
denote the localization of
$\cnum[y,y^{-1}]$ with respect to
$(y_{\alpha}\,|\,\alpha\in S)$.
For any $\cnum[y,y^{-1}]$-module $M$,
we set $M(\ast S):=M\otimes_{\cnum[y,y^{-1}]} \cnum[y,y^{-1}](\ast S)$.

A parabolic structure of $\vecV$ at  finite place is the following data:
\begin{itemize}
\item
A free $\cnum[y,y^{-1}]$-submodule $V\subset\vecV$
such that 
$V\otimes_{\cnum[y,y^{-1}]}\cnum(y)=
 \vecV\otimes_{\cnum[y,y^{-1}]}\cnum(y)$
and $\nbiga_{\gminiq}\cdot V=\vecV$.
\item
A finite subset $D\subset\cnum^{\ast}$
such that 
$V(\ast D)=(\Phi^{\ast})^{-1}(V)(\ast D)$ in $\vecV$.
\item
A sequence
$\vect_{\alpha}=(0\leq t_{\alpha,0}<t_{\alpha,1}<\cdots<t_{\alpha,m(\alpha)}<1)$
and a tuple 
$\vecnbigl_{\alpha}=(\nbigl_{\alpha,i}\,|\,i=1,\ldots,m(\alpha)-1)$
of $\cnum[\![y_{\alpha}]\!]$-lattices
$\nbigl_{\alpha,i}$ of
     $V\otimes_{\cnum[y,y^{-1}]}\cnum(\!(y_{\alpha})\!)$
are attached
to each $\alpha\in D$.
We formally set
$\nbigl_{\alpha,0}:=V\otimes\cnum[\![y_{\alpha}]\!]$
and 
$\nbigl_{\alpha,n(\alpha)}:=
 (\Phi^{\ast})^{-1}(V)\otimes\cnum[\![y_{\alpha}]\!]$.
\end{itemize}
If we fix $D$ and $\vect_{\alpha}$ $(\alpha\in D)$,
it is called a parabolic structure at 
$(D,(\vect_{\alpha})_{\alpha\in D})$
or just $(\vect_{\alpha})_{\alpha\in D}$.

\paragraph{Parabolic $\gminiq$-difference modules
and stability condition}

A parabolic $\gminiq$-difference $\cnum[y,y^{-1}]$-module
$\vecV_{\ast}$
consists of a $\gminiq$-difference $\cnum[y,y^{-1}]$-module
$\vecV$
with a good parabolic structure at infinity
$(\nbigp_{\ast}\vecV_{|\zerohat},\nbigp_{\ast}\vecV_{|\inftyhat})$
and a parabolic structure at finite place
$(D,(\vect_{\alpha},\vecnbigl_{\alpha})_{\alpha\in D})$.

We define the parabolic degree of $\vecV_{\ast}$.
Note that we obtain a parabolic vector bundle 
$\nbigp_{\ast}V$ on $\proj^1$
from $V$
and the filtered bundles 
$(\nbigp_{\ast}\vecV_{|\zerohat},\nbigp_{\ast}\vecV_{|\inftyhat})$.
For each $\alpha\in D$
and $i=0,\ldots,m(\alpha)$,
we define
\[
 \deg\bigl(\nbigl_{\alpha,i+1},\nbigl_{\alpha,i}\bigr):=
 \length\bigl(\nbigl_{\alpha,i+1}/\nbigl_{\alpha,i+1}\cap\nbigl_{\alpha,i}\bigr)
-\length\bigl(\nbigl_{\alpha,i}/\nbigl_{\alpha,i+1}\cap\nbigl_{\alpha,i}\bigr).
\]
Then, we set
\begin{multline}
 \deg(\vecV_{\ast}):=
 \deg(\nbigp_{\ast}V)
+\sum_{\alpha\in D}
 \sum_{i=0}^{m(\alpha)}
 (1-t_{\alpha,i})\deg(\nbigl_{\alpha,i+1},\nbigl_{\alpha,i})
 \\
-\sum_{\omega\in\rnum}
  \frac{\omega}{2}
 \Bigl(
 \dim_{\cnum(\!(y^{-1})\!)}\bigl((\vecV_{|\inftyhat})_{\omega}\bigr)
+\dim_{\cnum(\!(y)\!)}\bigl((\vecV_{\zerohat})_{\omega}\bigr)
 \Bigr).
\end{multline}

The stability condition is defined in a standard way.
Let $\vecVtilde'$ be 
a $\gminiq$-difference $\cnum(y)$-subspace of
$\vecVtilde:=\vecV\otimes\cnum(y)$.
We obtain a $\gminiq$-difference $\cnum[y,y^{-1}]$-submodule
$\vecV':=\vecVtilde'\cap\vecV$,
which is equipped with the induced parabolic structure.
We say that $\vecV_{\ast}$ is stable (resp. semistable)
if 
\[
 \frac{\deg(\vecV'_{\ast})}{\dim_{\cnum(y)}\vecVtilde'}
<\frac{\deg(\vecV_{\ast})}{\dim_{\cnum(y)}\vecVtilde}
\quad
\left(
\mbox{\rm resp. }
 \frac{\deg(\vecV'_{\ast})}{\dim_{\cnum(y)}\vecVtilde'}
\leq\frac{\deg(\vecV_{\ast})}{\dim_{\cnum(y)}\vecVtilde}
\right)
\]
for any $\gminiq$-difference $\cnum(y)$-subspace
$0\neq\vecVtilde'\subsetneq \vecVtilde$.
The polystability condition is also defined in the standard way.

\subsection{Geometrization of
parabolic $\gminiq$-difference
$\cnum[y,y^{-1}]$-modules}
\label{subsection;19.2.4.2}

It is the purpose in this paper to study 
the relationship between
meromorphic doubly periodic monopoles
and stable parabolic $\gminiq$-difference
$\cnum[y,y^{-1}]$-modules
of degree $0$.
As a bridge to connect them,
let us explain geometric objects 
directly corresponding to 
parabolic $\gminiq$-difference $\cnum[y,y^{-1}]$-modules.
We have already used a similar geometrization
in the context of difference modules in \cite{Mochizuki-difference-modules}.

\subsubsection{Spaces}

We consider the action of $\seisuu$ on
$\nbigm^{\cov}_{\gminiq}:=\cnum^{\ast}\times\real$
and  $\nbigmbar^{\cov}_{\gminiq}:=\proj^1\times\real$
determined by $n\bullet(y,t)=(\gminiq^n y,t+n)$.
We set $\nbigm_{\gminiq}:=\nbigm^{\cov}_{\gminiq}/\seisuu$
and $\nbigmbar_{\gminiq}:=\nbigm_{\gminiq}^{\cov}/\seisuu$.
For $\nu=0,\infty$,
we set
$H^{\cov}_{\gminiq,\nu}:=\{\nu\}\times\real$
and
$H_{\gminiq}:=H^{\cov}_{\gminiq,\nu}/\seisuu$.
We put
$H^{\cov}_{\gminiq}:=
 H^{\cov}_{\gminiq,0}\cup
 H^{\cov}_{\gminiq,\infty}$
and 
$H_{\gminiq}:=
 H_{\gminiq,0}\cup
 H_{\gminiq,\infty}$.

Let
$\nbigo_{\nbigmbar_{\gminiq}^{\cov}}(\ast H^{\cov}_{\gminiq})$
denote the sheaf of algebras on $\nbigmbar^{\cov}_{\gminiq}$
obtained as the pull back of $\nbigo_{\proj^1}(\ast\{0,\infty\})$
via the natural projection
$\nbigmbar^{\cov}_{\gminiq}\lrarr\proj^1$.
It is naturally equivariant with respect to the $\seisuu$-action.
Therefore,
we obtain a sheaf of algebras
$\nbigo_{\nbigmbar_{\gminiq}}(\ast H_{\gminiq})$
on $\nbigmbar_{\gminiq}$.
For any subset
$\nbigu\subset\nbigmbar_{\gminiq}$,
the restriction
of $\nbigo_{\nbigmbar_{\gminiq}}(\ast H_{\gminiq})$
to $\nbigu$ is denoted by $\nbigo_{\nbigu}(\ast (\nbigu\cap H_{\gminiq}))$.
We use a similar notation for 
the restriction of
$\nbigo_{\nbigmbar^{\cov}_{\gminiq}}(\ast H^{\cov}_{\gminiq})$
to subsets of $\nbigmbar^{\cov}_{\gminiq}$.

\subsubsection{Locally free sheaves with Dirac type singularity}

Let $Z\subset\nbigm_{\gminiq}$ be a finite subset.
Let $Z^{\cov}$ denote the subset of 
$\nbigm_{\gminiq}^{\cov}$
obtained as the pull back of $Z$.
Let $\gbigv$ be a locally free
$\nbigo_{\nbigmbar_{\gminiq}\setminus Z}(\ast H_{\gminiq})$-module.
Let $\gbigv^{\cov}$ denote
the $\seisuu$-equivariant locally free
$\nbigo_{\nbigmbar^{\cov}_{\gminiq}\setminus Z^{\cov}}
 (H^{\cov}_{\gminiq})$-module
obtained as the pull back of $\gbigv$.

Let $U$ be an open subset in $\proj^1$.
If $U\times\{t\}\subset \nbigmbar^{\cov}\setminus Z^{\cov}$,
the restriction
$\gbigv^{\cov}_{|U\times\{t\}}$ is 
naturally a locally free $\nbigo_U(\ast(U\cap\{0,\infty\}))$-module.
Note that any local sections of
$\nbigo_{\nbigmbar^{\cov}_{\gminiq}}(\ast H^{\cov}_{\gminiq})$
are locally constant in the $t$-direction.
Therefore, 
if $(U\times \closedclosed{t_1}{t_2})\cap Z^{\cov}=\emptyset$,
then there exists a naturally induced isomorphism
$\gbigv^{\cov}_{|U\times\{t_1\}}\simeq
 \gbigv^{\cov}_{|U\times\{t_2\}}$.
We call it  the scattering map
by following \cite{Charbonneau-Hurtubise}.

Let $(\alpha_0,t_0)\in Z^{\cov}$.
Take a neighbourhood $U$ of $\alpha_0$ in $\cnum^{\ast}$
and small $\epsilon>0$.
Set $U^{\ast}:=U\setminus\{\alpha_0\}$.
We have the isomorphism of $\nbigo_{U^{\ast}}$-modules
$\gbigv^{\cov}_{|U^{\ast}\times\{t_0-\epsilon\}}
\simeq
\gbigv^{\cov}_{|U^{\ast}\times\{t_0+\epsilon\}}$
induced by the scattering map.
We say that $(\alpha_0,t_0)$ is Dirac type singularity
if it is extended to an isomorphism of
$\nbigo_U(\ast\alpha_0)$-modules
$\gbigv^{\cov}_{|U\times\{t_0-\epsilon\}}(\ast\alpha_0)
\simeq
\gbigv^{\cov}_{|U\times\{t_0+\epsilon\}}(\ast\alpha_0)$.

If any $(\alpha_0,t_0)\in Z^{\cov}$ is Dirac type singularity,
we say that 
$\gbigv$ is a locally free 
$\nbigo_{\nbigmbar_{\gminiq}\setminus Z}(\ast H_{\gminiq})$-module
with Dirac type singularity.

\subsubsection{$\gminiq$-difference
$\cnum[y,y^{-1}]$-modules with parabolic structure at finite place}
\label{subsection;19.2.4.1}

Let $\gbigv$ be a locally free
$\nbigo_{\nbigmbar_{\gminiq}\setminus Z}(\ast H_{\gminiq})$-module
with Dirac type singularity.
Let $D$ denote the image of 
$Z^{\cov}\cap(\proj^1\times\closedopen{0}{1})$
by the projection
$\proj^1\times \closedopen{0}{1}\lrarr\proj^1$.
For $\alpha\in D$,
the sequence
$0\leq t_{\alpha,0}<t_{\alpha,1}<\cdots<t_{\alpha,m(\alpha)}<1$
is determined by
$\{(\alpha,t_{\alpha,i})\}=
 Z^{\cov}\cap(\{\alpha\}\times\closedopen{0}{1})$.
Let us observe that
$\gbigv$ naturally induces
a $\gminiq$-difference
$\cnum[y,y^{-1}]$-module 
with parabolic structure at 
$(D,\{\vect_{\alpha}\}_{\alpha\in D})$.

Take a sufficiently small $\epsilon>0$
such that
$(\proj^1\times\closedopen{-\epsilon}{0})
 \cap Z^{\cov}=\emptyset$.
The restriction of $\gbigv^{\cov}$
to $\proj^1\times\{-\epsilon\}$
induces 
a locally free $\nbigo_{\proj^1}(\ast\{0,\infty\})$-module
$\gbigv^{\cov}_{-\epsilon}$.
We obtain a $\cnum[y,y^{-1}]$-module
$V:=H^0(\proj^1,\gbigv^{\cov}_{-\epsilon})$.
It is independent of a choice of $\epsilon$
up to canonical isomorphisms.
Similarly, the restriction of 
$\gbigv^{\cov}$
to $\proj^1\times\{1-\epsilon\}$
induces 
a locally free $\nbigo_{\proj^1}(\ast\{0,\infty\})$-module
$\gbigv^{\cov}_{1-\epsilon}$.
We obtain a $\cnum[y,y^{-1}]$-module
$V':=H^0(\proj^1,\gbigv^{\cov}_{1-\epsilon})$.

Let $\Phi:\proj^1\lrarr\proj^1$
be the morphism defined by
$\Phi(y)=\gminiq y$.
We have the natural isomorphism
$\Phi^{\ast}:\gbigv^{\cov}_{1-\epsilon}
\simeq
 \gbigv^{\cov}_{-\epsilon}$,
which induces a $\cnum$-linear isomorphism
\begin{equation}
\label{eq;19.2.3.10}
\Phi^{\ast}:V'\simeq V
\end{equation}
such that $\Phi^{\ast}(fs)=\Phi^{\ast}(f)\Phi^{\ast}(s)$
for any $f\in\cnum[y,y^{-1}]$ and $s\in V'$.
The scattering map induces an isomorphism
\begin{equation}
\label{eq;19.2.3.11}
 V(\ast D)\simeq
 V'(\ast D).
\end{equation}
The isomorphisms (\ref{eq;19.2.3.10})
and (\ref{eq;19.2.3.11})
induce
a $\cnum$-linear automorphism $\Phi^{\ast}$
on $\vecVtilde:=V\otimes\cnum(y)$
such that
$\Phi^{\ast}(fs)=\Phi^{\ast}(f)\Phi^{\ast}(s)$
for any $f\in\cnum(y)$ and $s\in\vecVtilde$.
We set
$\vecV:=\nbiga_{\gminiq}\cdot V$
in $\vecVtilde$.

For $\alpha\in D$
and $t_{\alpha,i}$ $(1\leq i\leq m(\alpha))$,
we obtain the $\cnum[\![y_{\alpha}]\!]$-lattices
$\nbigl_{\alpha,i}$ of 
$V\otimes\cnum(\!(y_{\alpha})\!)$
induced by the formal completion of
the stalks of $\gbigv$ at $(\alpha,t_{\alpha,i}-\epsilon)$
for any sufficiently small $\epsilon>0$.
They induce a parabolic structure 
$\bigl\{
 (\vect_{\alpha},\vecnbigl_{\alpha})
 \bigr\}_{\alpha\in D}$
of $\vecV$ at $(D,(\vect_{\alpha})_{\alpha\in D})$.
The following lemma is easy to observe.

\begin{lem}
The above construction induces
an equivalence between the following objects:
\begin{itemize}
\item
Locally free 
 $\nbigo_{\nbigmbar_{\gminiq}\setminus Z}(\ast H_{\gminiq})$-modules
 with Dirac type singularity.
\item
$\gminiq$-difference $\cnum[y,y^{-1}]$-modules
with parabolic structure at
 $(D,(\vect_{\alpha})_{\alpha\in D})$.
\hfill\qed
\end{itemize}
\end{lem}

\subsubsection{Good filtered bundles
over equivariant
$\nbigo_{\Hhat^{\cov}_{\gminiq,\nu}}(\ast H^{\cov}_{\gminiq,\nu})$-modules}

We set $y_{0}:=y$ and $y_{\infty}:=y^{-1}$.
We also set $\gminiq_{0}:=\gminiq$
and $\gminiq_{\infty}:=\gminiq^{-1}$.
For $\nu=0,\infty$,
let
$\nbigo_{\Hhat^{\cov}_{\gminiq,\nu}}(\ast H^{\cov}_{\gminiq,\nu})$
denote the sheaf of locally constant
$\cnum(\!(y_{\nu})\!)$-valued functions
on $H^{\cov}_{\gminiq,\nu}$.
It is $\seisuu$-equivariant
by the action
$n^{\ast}(f)(y_{\nu})=f(\gminiq_{\nu}^ny_{\nu})$.

For any $\seisuu$-equivariant locally free
$\nbigo_{\Hhat^{\cov}_{\gminiq,\nu}}(\ast H^{\cov}_{\gminiq,\nu})$-module
$\gbigvhat^{\cov}$,
let $\gbigvhat^{\cov}_{|t}$
denote the restriction of $\gbigvhat^{\cov}$
to $t\in\real$
which is naturally a $\cnum(\!(y_{\nu})\!)$-vector space.
For any $t_1,t_2\in\real$,
we have the isomorphism called the scattering map:
\begin{equation}
\label{eq;19.2.3.20}
 \gbigvhat^{\cov}_{|t_1}\simeq
 \gbigvhat^{\cov}_{|t_2}.
\end{equation}
By the $\seisuu$-action,
we have the isomorphism
\begin{equation}
\label{eq;19.2.3.21}
 \Phi^{\ast}:
 \gbigvhat^{\cov}_{|t+1}
\simeq
 \gbigvhat^{\cov}_{|t}.
\end{equation}
Therefore,
$\gbigvhat^{\cov}_{|0}$
is naturally a $\gminiq$-difference
$\cnum(\!(y_{\nu})\!)$-module.
It is easy to observe that
this procedure induces
an equivalence between
$\seisuu$-equivariant
locally free
$\nbigo_{\Hhat^{\cov}_{\gminiq,\nu}}(\ast H^{\cov}_{\gminiq,\nu})$-modules
and
$\gminiq_{\nu}$-difference $\cnum(\!(y_{\nu})\!)$-modules.

Let $\gbigvhat^{\cov}$ be a $\seisuu$-equivariant locally free
$\nbigo_{\Hhat^{\cov}_{\gminiq,\nu}}
 (\ast H^{\cov}_{\gminiq,\nu})$-module.
There exists a decomposition
$\gbigvhat^{\cov}=
 \bigoplus_{\omega\in\rnum}
 \gbigvhat^{\cov}_{\omega}$
corresponding to the slope decomposition of
the $\gminiq$-difference $\cnum(\!(y_{\nu})\!)$-module 
$\gbigvhat^{\cov}_{|0}$.
A good filtered bundle $\nbigp_{\ast}\gbigvhat^{\cov}$
over $\gbigvhat^{\cov}$
is defined to be a family of filtered bundles
$(\nbigp_{\ast}(\gbigvhat^{\cov}_{|t})\,|\,t\in\real)$
such that the following holds.
\begin{itemize}
\item
$\nbigp_{\ast}(\gbigvhat^{\cov}_{|t})
=\bigoplus_{\omega}\nbigp_{\ast}(\gbigvhat^{\cov}_{\omega|t})$.
\item
The isomorphism (\ref{eq;19.2.3.20}) induces
$\nbigp_{a}(\gbigvhat^{\cov}_{\omega|t_1})
\simeq
 \nbigp_{a+\omega(t_2-t_1)}(\gbigvhat^{\cov}_{\omega|t_2})$
for any $a\in\real$ and $t_1,t_2\in\real$.
\item
The isomorphism (\ref{eq;19.2.3.21}) induces
$\nbigp_{a}(\gbigvhat^{\cov}_{\omega|t+1})
\simeq
 \nbigp_a(\gbigvhat^{\cov}_{\omega|t})$
for any $t\in\real$ and $a\in\real$.
\end{itemize}
Clearly,
good filtered bundles over 
a $\seisuu$-equivariant
$\nbigo_{\Hhat^{\cov}_{\gminiq,\nu}}(\ast H^{\cov}_{\gminiq,\nu})$-module
$\gbigvhat^{\cov}$
are equivalent to 
good filtered bundles over
$\gminiq_{\nu}$-difference
$\cnum(\!(y_{\nu})\!)$-module
$\gbigvhat^{\cov}_{|0}$.

\subsubsection{Good parabolic structure at infinity}
\label{subsection;19.2.4.10}

Let $\gbigv$ be a locally free
$\nbigo_{\nbigmbar_{\gminiq}\setminus Z}(\ast H_{\gminiq})$-module
with Dirac type singularity.
Let $\gbigv^{\cov}$ be the $\seisuu$-equivariant
$\nbigo_{\nbigmbar^{\cov}_{\gminiq}\setminus Z^{\cov}}
 (\ast H^{\cov}_{\gminiq})$-module
obtained as the pull back of $\gbigv$.
For any $t\in \real$ and $\nu=0,\infty$,
we obtain the formal completions
$\gbigvhat^{\cov}_{\nu|t}$
of $\gbigvhat_{|(\proj^1\times\{t\})\setminus Z^{\cov}}$
at $(\nu,t)$.
They induce $\seisuu$-equivariant
locally free 
$\nbigo_{\Hhat^{\cov}_{\gminiq,\nu}}(\ast H^{\cov}_{\gminiq,\nu})$-modules
$\gbigvhat^{\cov}_{\nu}$ $(\nu=0,\infty)$.

Let $\vecV$  be the $\gminiq$-difference $\cnum[y,y^{-1}]$-module
with a parabolic structure 
$(\vect_{\alpha},\vecnbigl_{\alpha})_{\alpha\in D}$
at finite place
corresponding to $\gbigv$
as in \S\ref{subsection;19.2.4.1}.
Note that
$\vecV_{|\widehat{\nu}}$
is naturally identified with $\gbigvhat^{\cov}_{\nu|0}$.
Under the identification,
good filtered bundles
$\nbigp_{\ast}\gbigvhat^{\cov}_{\nu}=
 (\nbigp_{\ast}\gbigvhat^{\cov}_{\nu|t}\,|\,t\in\real)$
over $\gbigvhat^{\cov}_{\nu}$
are equivalent to good filtered bundles
$\nbigp_{\ast}\vecV_{|\widehat{\nu}}$
over $\vecV_{|\nuhat}$.

\subsubsection{Geometrization of parabolic $\gminiq$-difference
$\cnum[y,y^{-1}]$-modules}

By the considerations in 
\S\ref{subsection;19.2.4.1} and \S\ref{subsection;19.2.4.10},
we obtain the following.
\begin{prop}
\label{prop;19.2.4.20}
The following objects are equivalent.
\begin{itemize}
\item
 $\gminiq$-difference $\cnum[y,y^{-1}]$-modules
 with a good parabolic structure at infinity
 and a parabolic structure at $(D,(\vect_{\alpha})_{\alpha\in D})$.
\item
Good filtered bundles with Dirac type singularity
 over $(\nbigmbar_{\gminiq};H_{\gminiq},Z)$,
i.e.,
locally free
$\nbigo_{\nbigmbar_{\gminiq}\setminus Z}(\ast H_{\gminiq})$-modules
$\gbigv$ with Dirac type singularity
enhanced by good filtered bundles
$\nbigp_{\ast}(\gbigvhat^{\cov}_{\nu})$ over
$\gbigvhat^{\cov}_{\nu}$.
\end{itemize}
Here, $Z$ and $(D,(\vect_{\alpha})_{\alpha\in D})$
are related as in {\rm \S\ref{subsection;19.2.4.1}}.
\hfill\qed
\end{prop}

Let $\gbigv$ be a locally free
$\nbigo_{\nbigmbar_{\gminiq}\setminus Z}(\ast H_{\gminiq})$-module
with Dirac type singularity
enhanced with good filtered bundles
$\nbigp_{\ast}\gbigvhat^{\cov}_{\nu}=
 (\nbigp_{\ast}\gbigvhat^{\cov}_{\nu|t}\,|\,t\in\real)$
$(\nu=0,\infty)$.
Let $\pi^{\cov}:\nbigmbar^{\cov}_{\gminiq}\lrarr\real$
denote the projection.
For any $t\in\closedopen{0}{1}\setminus \pi^{\cov}(Z^{\cov})$,
let $\gbigv^{\cov}_{\proj^1\times\{t\}}$ denote 
the $\nbigo_{\proj^1}(\ast\{0,\infty\})$-module
obtained as the restriction of
$\gbigv^{\cov}$ to $\proj^1\times\{t\}$.
We obtain a filtered bundle
$\nbigp_{\ast}\gbigv^{\cov}_{\proj^1\times\{t\}}$
from $\gbigv^{\cov}_{\proj^1\times\{t\}}$
and $(\nbigp_{\ast}\gbigvhat^{\cov}_{0|t},\nbigp_{\ast}\gbigv^{\cov}_{\infty|t})$.
We set
\[
 \deg\bigl(
 \gbigv,
 \bigl(
 \nbigp_{\ast}\gbigvhat^{\cov}_{\nu}
 \bigr)_{\nu=0,\infty}
 \bigr)
 :=
 \int_{0}^1
 \deg(\nbigp_{\ast}\gbigv^{\cov}_{\proj^1\times\{t\}})\,dt.
\]
We define the stability condition
in the standard way.
The following is easy to see by the construction.
\begin{lem}
The degree is preserved by the equivalence
in Proposition {\rm\ref{prop;19.2.4.20}}.
Therefore, the stability condition is also 
preserved by the equivalence.
\hfill\qed
\end{lem}

\subsection{From monopoles to $\gminiq$-difference modules}

Let us explain how a meromorphic monopole
on $\nbigm$ induces geometric objects
as in \S\ref{subsection;19.2.4.2},
and hence $\gminiq$-difference modules.
More detailed explanation will be used later.

\subsubsection{Space}

Take $\mu_i\in\cnum$ $(i=1,2)$
such that
(i) $\mu_1$ and $\mu_2$ are linearly independent over $\real$,
(ii) $\Image(\mu_2/\mu_1)>0$.
Let $\Gamma$ denote the lattice of $\cnum$
generated by $\mu_1$ and $\mu_2$.
Let $\Vol(\Gamma)$ denote the volume 
of the quotient $\cnum/\Gamma$
with respect to the volume form $\frac{\sqrt{-1}}{2}\,dz\,d\zbar$.

We set $X:=\cnum_z\times\cnum_w$
with the Euclidean metric $dz\,d\zbar+dw\,d\wbar$.
Let us consider the action of
$\real\,\tte_0\oplus\seisuu\tte_1\oplus\seisuu\tte_2$
on $X$ by
$\tte_0(z,w)=(z,w+1)$
and
$\tte_i(z,w)=(z+\mu_i,w)$ $(i=1,2)$.

Let $\nbigm^{\cov}$ be the quotient space of
$X$ by the action of $\real\tte_0\oplus\seisuu\tte_1$.
It is equipped with an induced action of $\seisuu\tte_2$.
The quotient space $\nbigm^{\cov}/\seisuu\tte_2$
is naturally identified with $\nbigm$.

\subsubsection{Mini-complex coordinate system}
\label{subsection;19.1.26.10}

Let $\lambda$ be a complex number
such that $\lambda\neq \pm\sqrt{-1}\mu_1|\mu_1|^{-1}$.
As in Lemma \ref{lem;18.3.16.2} below,
there exist $\tts_1\in\real$
and $\ttg_1\in\cnum$ with $|\ttg_1|=1$
such that 
\[
 -\lambda\mubar_1+\tts_1=
 \ttg_1(\mu_1+\lambda \tts_1)\neq 0.
\]
If $|\lambda|\neq 1$,
there are two such choices.
If $|\lambda|=1$ and $\lambda\neq \pm\sqrt{-1}\mu_1|\mu_1|^{-1}$
there is a unique choice.
We consider the complex coordinate system
$(\ttu,\ttv)$ given as follows:
\[
 \ttu=\frac{1}{1-\ttg_1\lambda}\bigl(z+\lambda^2\zbar+\lambda(\wbar-w)\bigr),
\quad\quad
 \ttv=\frac{1}{1-\ttg_1\lambda}
 \bigl(-\ttg_1z-\lambda\zbar+w-\lambda\ttg_1\wbar\bigr).
\]
Note that
\[
 \tte_0(\ttu,\ttv)=(\ttu,\ttv)+(0,1),
\quad\quad
 (\tte_1+\tts_1\tte_0)(\ttu,\ttv)
=(\ttu,\ttv)+(\mu_1+\lambda\tts_1,0).
\]
We define
\[
 \ttU:=\exp\Bigl(
 \frac{2\pi\sqrt{-1}}{\mu_1+\lambda\tts_1}\ttu
 \Bigr),
\quad\quad
 \ttt:=\Image(\ttv).
\]
Then, $(\ttU,\ttt)$ induces
an isomorphism
$\nbigm^{\cov}\simeq\cnum_{\ttU}^{\ast}\times \real_{\ttt}$.
We set
\begin{equation}
\label{eq;19.1.26.11}
 \gminiq^{\lambda}:=
 \exp\Bigl(
 2\pi\sqrt{-1}
 \frac{\mu_2+\lambda^2\mubar_2}{\mu_1+\lambda^2\mubar_1}
 \Bigr),
\quad\quad
 \gminit^{\lambda}:=
 -\frac{\Vol{\Gamma}}{\Re(\ttg_1\mu_1)}.
\end{equation}
The following holds:
\[
 \tte_2(\ttU,\ttt)=(\gminiq^{\lambda}\ttU,\ttt+\gminit^{\lambda}).
\]
Note that $\gminit^{\lambda}$  is non-zero,
but that $\gminit^{\lambda}$ is not necessarily positive.
We also remark that
$|\gminiq^{\lambda}|=1$ 
if and only if $|\lambda|=1$.

When we consider the above coordinate system $(\ttU,\ttt)$,
$\nbigm^{\cov}$ and $\nbigm$
are also denoted by
$\nbigm^{\lambda\cov}$ and $\nbigm^{\lambda}$,
respectively.

\subsubsection{Compactification}

We set 
$\nbigmbar^{\lambda\cov}:=
 \proj_{\ttU}^1\times\real_{\ttt}$
which we regard a partial compactification of
$\nbigm^{\lambda\cov}\simeq \cnum_{\ttU}^{\ast}\times\real_{\ttt}$.
It is equipped with the naturally induced
$\seisuu\tte_2$-action.
We put
$\nbigmbar^{\lambda}:=\nbigmbar^{\lambda\cov}/\seisuu\tte_2$
and 
$\nbigm^{\lambda}:=\nbigm^{\lambda\cov}/\seisuu\tte_2$.
Set
$H^{\lambda\cov}:=
\nbigmbar^{\lambda\cov}
\setminus\nbigm^{\lambda\cov}$.
We obtain
$H^{\lambda}\subset
 \nbigmbar^{\lambda}$
as the quotient of $H^{\lambda\cov}$
by the $\seisuu\tte_2$-action.

We have the 
$\seisuu$-equivariant isomorphisms
$\nbigmbar^{\lambda\cov}
\simeq
 \nbigmbar^{\cov}_{\gminiq^{\lambda}}$,
$\nbigm^{\lambda\cov}
\simeq
 \nbigm^{\cov}_{\gminiq^{\lambda}}$
and 
$H^{\lambda\cov}
 \simeq 
 H^{\cov}_{\gminiq^{\lambda}}$
given by $\ttU=y$ and $\ttt=\gminit^{\lambda}t$.
It induces an isomorphism
$\nbigmbar^{\lambda}
\simeq
\nbigmbar_{\gminiq^{\lambda}}$,
$\nbigm^{\lambda}
\simeq
\nbigm_{\gminiq^{\lambda}}$,
and
$H^{\lambda}
 \simeq
H_{\gminiq^{\lambda}}$.

\subsubsection{Mini-holomorphic bundles associated to monopoles}

Let us explain
how a meromorphic monopole induces
$\nbigo_{\nbigmbar^{\lambda}\setminus Z}(\ast H^{\lambda})$-modules
with Dirac type singularity
enhanced with good filtered bundles at infinity.
It depends on the choice of
$(\lambda,\tte_1,\tts_1)$.

Let $(E,h,\nabla,\phi)$ be a meromorphic monopole on 
$\nbigm\setminus Z$.
We have the naturally defined  operators
$\del_{E,\ttUbar}$
and 
$\del_{E,\ttt}$ on $E$
such  that
$\bigl[
 \del_{E,\ttUbar},\del_{E,\ttt}
\bigr]=0$,
which is a consequence of the Bogomolny equation.
(Note that the vector fields
$\del_{\ttUbar}$ and $\del_{\ttt}$
are not necessarily orthogonal.)

Let $Z^{\cov}$ denote the subset of $\nbigm^{\cov}$
obtained as the pull back of $Z$.
Let $E^{\cov}$ denote the vector bundle
on $\nbigm^{\cov}\setminus Z^{\cov}$
obtained as the pull back of $E$.
It is equipped with the induced operators
$\del_{E^{\cov},\ttUbar}$
and
$\del_{E^{\cov},\ttt}$.

We obtain 
a $\seisuu$-equivariant
locally free $\nbigo_{\nbigm^{\cov}\setminus Z^{\cov}}$-module
$\nbige^{\cov}$
as the sheaf of $C^{\infty}$-sections $s$ of $E^{\cov}$
such that $\del_{E^{\cov},\ttUbar}s=\del_{E^{\cov},\ttt}s=0$.
Each point of $Z^{\cov}$
is Dirac type singularity of $\nbige^{\cov}$
under the assumption that
each point of $Z$ is Dirac type singularity
of the monopole
$(E,h,\nabla,\phi)$.

For $\ttt\in\real$,
let $E^{\cov}(\ttt)$
denote the restriction of
$E^{\cov}$
to $(\cnum_{\ttU}^{\ast}\times\{\ttt\})\setminus Z^{\cov}$.
Together with the operator 
$\del_{E^{\cov},\ttUbar}$,
it is naturally a holomorphic vector bundle.
The sheaf $\nbige^{\cov}(\ttt)$
of holomorphic sections of $E^{\cov}(\ttt)$
is identified with the restriction of
$\nbige^{\cov}$
to $(\cnum_{\ttU}^{\ast}\times\{\ttt\})\setminus Z^{\cov}$.

Let $h(\ttt)$ be the restriction of the metric $h$
to $E^{\cov}(\ttt)$.
Because the monopole is meromorphic,
it turns out that
$(E^{\cov}(\ttt),\delbar_{E^{\cov},\ttUbar},h(\ttt))$
is acceptable around $\ttU=0,\infty$,
i.e.,
the curvature of the Chern connection
is bounded with respect to $h(\ttt)$
and the metric
$|\ttU|^{-2}(\log|\ttU|)^{-2}d\ttU\,d\ttUbar$.
(See 
Proposition \ref{prop;18.9.2.1},
Lemma \ref{lem;19.2.7.10}
and Corollary \ref{cor;18.12.17.3}.)
Therefore,
$\nbige^{\cov}(\ttt)$ is extended to
a locally free 
$\nbigo_{(\proj^1\times\{\ttt\}\setminus Z^{\cov})}(\ast\{0,\infty\})$-module
$\nbigp\nbige^{\cov}(\ttt)$.
Moreover, we obtain filtered bundles
$\nbigp_{\ast}\nbige^{\cov}(\ttt)_{|\widehat{\nu}}$ $(\nu=0,\infty)$
over the formal completions
$\nbigp\nbige^{\cov}(\ttt)_{|\widehat{\nu}}$
by considering the growth orders of 
the norms of sections with respect to $h(\ttt)$.

It is easy to see that 
the scattering map induces an isomorphism
$\nbigp\nbige^{\cov}(\ttt_1)
\simeq
\nbigp\nbige^{\cov}(\ttt_1)$
for $\ttt_1\in\ttt_2$
on neighbourhoods of $\ttU=0,\infty$.
Therefore, 
$\nbigp\nbige^{\cov}(\ttt)$ $(\ttt\in\real)$
induce a $\seisuu\tte_2$-equivariant
locally free
$\nbigo_{\nbigmbar^{\lambda\cov}\setminus Z^{\cov}}
 (\ast H^{\lambda\cov})$-module
$\nbigp\nbige^{\cov}$
with Dirac type singularity.
We obtain a locally free
$\nbigo_{\nbigmbar^{\lambda}\setminus Z}(\ast H^{\lambda})$-module
$\nbigp\nbige$
with Dirac type singularity
as the descent of 
$\nbigp\nbige^{\cov}$.
Moreover,
the families of filtrations
$\nbigp_{\ast}\nbige^{\cov}(\ttt)_{|\nuhat}$ $(\ttt\in\real)$
are good filtered bundles
(Theorem \ref{thm;18.12.17.110}).
In this way,
a meromorphic monopole on $\nbigm\setminus Z$
induces a good filtered bundle with Dirac type singularity
over $(\nbigmbar^{\lambda};H^{\lambda},Z)$,
and hence
a parabolic $\gminiq^{\lambda}$-difference
$\cnum[y,y^{-1}]$-module.

Then, the following theorem is the main result  of this paper.
\begin{thm}[Theorem 
 \ref{thm;18.11.21.20}]
The above construction induces 
an equivalence between
meromorphic doubly periodic monopoles
and polystable parabolic $\gminiq^{\lambda}$-difference modules
of degree $0$.
\end{thm}

\subsection{Filtered objects on elliptic curves}

As the ``Betti'' side,
we shall also give a minor complement
on the parabolic version of the Riemann-Hilbert correspondence
of $\gminiq$-difference modules ($|\gminiq|\neq 1$)
and its relation with the Kobayashi-Hitchin correspondence
in \S\ref{section;19.1.27.1}.

\subsubsection{Riemann-Hilbert correspondence for 
$\gminiq$-difference modules with $|\gminiq|\neq 1$}
\label{subsection;19.2.3.1}

Suppose that $|\gminiq|\neq 1$.
The Riemann-Hilbert correspondence for 
germs of analytic $\gminiq$-difference modules
was established by van der Put and Reversat \cite{van-der-Put-Reversat},
and Ramis, Sauloy and Zhang \cite{Ramis-Sauloy-Zhang}.
The global Riemann-Hilbert correspondence
for $\gminiq$-difference $\cnum[y,y^{-1}]$-modules
is due to Kontsevich and Soibelman.

Set $\gminiq^{\seisuu}:=\{\gminiq^n\,|\,n\in\seisuu\}$.
Let $\Phi:\cnum^{\ast}\lrarr \cnum^{\ast}$
be defined by $\Phi(y)=\gminiq y$.
It induces a $\gminiq^{\seisuu}$-action on $\cnum^{\ast}$.
We set $T:=\cnum^{\ast}/\gminiq^{\seisuu}$.
Clearly,
$\gminiq^{\seisuu}$-equivariant 
coherent $\nbigo_{\cnum^{\ast}}$-modules
are equivalent to coherent $\nbigo_{T}$-modules.

For a locally free $\nbigo_T$-module $\ttE$,
an anti-Harder-Narasimhan filtration of $\ttE$
is a filtration $\gbigf$ indexed by $(\rnum\cup\{\infty\},\leq)$
such that
(i) $\ttE_{\mu}:=\Gr^{\gbigf}_{\mu}(\ttE)$ is semistable
with $\deg(\ttE_{\mu})/\rank(\ttE_{\mu})=\mu$ if $\mu\neq \infty$,
(ii) $\Gr^{\gbigf}_{\infty}(\ttE)$ is torsion.
If $\Gr^{\gbigf}_{\infty}(\ttE)=0$
then we call it an anti-Harder-Narasimhan filtration
indexed by $(\rnum,\leq)$.

According to \cite{van-der-Put-Reversat}
and \cite{Ramis-Sauloy-Zhang},
$\gminiq^{\seisuu}$-equivariant
locally free $\nbigo_{\cnum_y}(\ast 0)$-modules
are equivalent to 
locally free $\nbigo_{T}$-modules
equipped with an anti-Harder-Narasimhan filtration
indexed by $(\rnum,\leq)$.
For $(\ttE,\gbigf)$,
let $\ttK_0(\ttE,\gbigf)$ denote the corresponding
$\gminiq^{\seisuu}$-equivariant
locally free $\nbigo_{\cnum_y}(\ast 0)$-module.
It is equipped with the filtration
induced by $\gbigf$
so that
$\Gr^{\gbigf}_{\mu}\ttK_0(\ttE,\gbigf)_{|\zerohat}$
has pure slope $\varrho(\gminiq)\mu$,
where $\varrho(\gminiq)\in\{\pm1\}$ 
is the signature of $\log|\gminiq|\neq 0$.
Similarly, $\gminiq^{\seisuu}$-equivariant
locally free $\nbigo_{\cnum_{y^{-1}}}(\ast \infty)$-modules
are also equivalent to locally free $\nbigo_{T}$-modules
equipped with an anti-Harder-Narasimhan filtration
indexed by $(\rnum,\leq)$.
For $(\ttE,\gbigf)$,
we have the $\gminiq^{\seisuu}$-equivariant
locally free $\nbigo_{\cnum_{y^{-1}}}(\ast\infty)$-module
$\ttK_{\infty}(\ttE,\gbigf)$.
For the induced filtration,
$\Gr^{\gbigf}_{\mu}(\ttK_{\infty}(\ttE,\gbigf)_{|\inftyhat})$
has pure slope $-\varrho(\gminiq)\mu$.

According to Kontsevich-Soibelman,
$\gminiq$-difference $\cnum[y,y^{-1}]$-modules
are equivalent to locally free $\nbigo_T$-modules $\ttE$
equipped with two anti-Harder-Narasimhan filtrations $\gbigf_{\pm}$
indexed by $(\rnum\cup\{\infty\},\leq)$.

\subsubsection{Filtered objects on elliptic curves}

The Riemann-Hilbert correspondence
for $\gminiq$-difference modules in 
\S\ref{subsection;19.2.3.1}
is enhanced to the correspondence
for filtered objects.
Let us explain the filtered counterpart 
on the side of elliptic curves.

Let $\ttD\subset T$ be a finite subset.
Let $\ttEtilde$ be a locally free
$\nbigo_T(\ast\ttD)$-module.
For each $P\in T$,
let $\ttEtilde_{|\Phat}$ denote the formal completion
of the stalk of $\ttEtilde$ at $P$.
A $\gminiq$-difference parabolic structure
on $\ttEtilde$ consists of the following data:
\begin{itemize}
\item
A finite sequence
$\vecs_P=\bigl(s_{P,1}<s_{P,2}<\cdots<s_{P,m(P)}\bigr)$
in $\real$ for each $P\in \ttD$.

We formally set
$s_{P,0}:=-\infty$
and $s_{P,m(P)+1}:=\infty$.
\item
A tuple of lattices
$\vecnbigk_P=\bigl(\nbigk_{P,i}\,|\,i=0,\ldots,m(P)\bigr)$
of $\ttEtilde_{|\Phat}$.

Note that we obtain 
the lattice $\ttE_{-}\subset\ttEtilde$
determined by  $\nbigk_{P,0}$ $(P\in \ttD)$
and the lattice
$\ttE_{+}\subset\ttEtilde$
determined by 
$\nbigk_{P,m(P)}$ $(P\in \ttD)$.
\item
Let $\gbigf_{\pm}$ be anti-Harder-Narasimhan
filtrations of $\ttE_{\pm}$
indexed by $(\rnum,\leq)$.
\item
 Filtrations $\nbigf_{\pm}$
 on $\Gr^{\gbigf_+}_{\mu}(\ttE_{\pm})$ $(\mu\in\rnum)$
 indexed by $(\real,\leq)$ 
 such that
 $\ttE_{a,\mu,\pm}:=
 \Gr^{\nbigf_{\pm}}_a\Gr^{\gbigf_{\pm}}_{\mu}(\ttE_{\pm})$
 are also semistable 
 with $\deg(\ttE_{a,\mu,\pm})/\rank(\ttE_{a,\mu,\pm})=\mu$.
\end{itemize}
When we fix $(\vecs_P)_{P\in\ttD}$,
it is called $\gminiq$-difference parabolic structure
at $(\vecs_{P})_{P\in\ttD}$.

We define the degree of
$\ttEtilde_{\ast}=(\ttEtilde,(\vecs_P,\vecnbigk_P)_{P\in \ttD},
 (\gbigf_{\pm},\nbigf_{\pm}))$
as follows:
\begin{multline}
 \deg(\ttEtilde_{\ast}):=
 -\sum_{P\in D}
 \sum_{i=1}^{m(P)}
 s_{P,i}\deg(\nbigk_{P,i},\nbigk_{P,i-1})
 \\
-\sum_{\mu\in\rnum}
 \sum_{b\in\real}
 b\rank\Gr^{\nbigf_-}_b\Gr^{\gbigf_-}_{\mu}(\ttE_-)
-\sum_{\mu\in\rnum}
 \sum_{b\in\real}
 b\rank\Gr^{\nbigf_+}_b\Gr^{\gbigf_+}_{\mu}(\ttE_+).
\end{multline}
By using the degree,
we define the stability, semistability and polystability conditions
for filtered objects in the standard ways.

\paragraph{Rescaling of $\gminiq$-difference parabolic structure}
There is a rescaling of
$\gminiq$-difference parabolic structure.
For $\gminit>0$,
we obtain a sequence
$\vecs^{(\gminit)}_P:=(\gminit s_{P,i})$.
We set
$\vecnbigk_P^{(\gminit)}:=\vecnbigk_P$
and 
$\gbigf^{(\gminit)}_{\pm}:=\gbigf_{\pm}$.
We also obtain filtrations
$\nbigf_{\pm}^{(\gminit)}$ by
$(\nbigf_{\pm}^{(\gminit)})_{\gminit a}\Gr^{\gbigf_{\pm}}(\ttE_{\pm})
 :=(\nbigf_{\pm})_{a}\Gr^{\gbigf_{\pm}}(\ttE_{\pm})$.
We set
\[
 \ttH^{(\gminit)}\bigl(
 \ttEtilde_{\ast}\bigr):=
 \bigl(
 \ttEtilde,(\vecs^{(\gminit)}_P,\vecnbigk^{(\gminit)}_P)_{P\in \ttD},
 (\gbigf^{(\gminit)}_{\pm},\nbigf^{(\gminit)}_{\pm})
 \bigr). 
\]
In the case $\gminit<0$,
we set
$s^{(\gminit)}_{P,i}:=
 \gminit s_{P,m(P)-i+1}$,
and 
$\vecs^{(\gminit)}_P:=
 \bigl(
 s^{(\gminit)}_{P,i}
 \bigr)$.
We set
$\nbigk^{(\gminit)}_{P,i}:=
 \nbigk_{P,m(P)-i}$
and 
$\vecnbigk^{(\gminit)}_P:=
 \bigl(
 \nbigk^{(\gminit)}_{P,i}
 \bigr)$.
We also set
$(\gbigf^{(\gminit)}_{\pm})_{\mu}:=
 \gbigf_{\mp,\mu}$
and 
$(\nbigf^{(\gminit)}_{\pm})_{|\gminit|a}:=
 \nbigf_{\mp,a}$.
Then, we define
\[
 \ttH^{(\gminit)}\bigl(
 \ttEtilde_{\ast}\bigr):=
 \bigl(
 \ttEtilde_{\ast},
 (\vecs^{(\gminit)}_P,\vecnbigk^{(\gminit)}_P)_{P\in\ttD},
 \gbigf^{(\gminit)}_{\pm},
 \nbigf^{(\gminit)}_{\pm}
 \bigr).
\]
It is easy to see
$\deg(\ttH^{(\gminit)}(\ttEtilde_{\ast}))
=|\gminit|\deg(\ttEtilde_{\ast})$.

\subsubsection{Equivalence}

The natural projection
$\nbigm^{\cov}_{\gminiq}\lrarr\cnum^{\ast}$
induces
$\sfp:\nbigm_{\gminiq}\lrarr T$.
Let $\sff:\nbigm^{\cov}_{\gminiq}\lrarr\real$
be defined by
\[
 \sff(y,t):=t-\frac{\log|y|}{\log|\gminiq|}.
\]
It induces the map
$\sff:\nbigm_{\gminiq}\lrarr \real$.
Let $Z\subset \nbigm_{\gminiq}$ be a finite subset.
We set $\ttD:=\sfp(Z)$.
For each $P\in\ttD$,
we obtain
\begin{equation}
\label{eq;19.2.4.40}
\vecs_P=(s_{P,1}<s_{P,2}<\cdots<s_{P,m(P)}):=
 \sff(\sfp^{-1}(P)\cap Z).
\end{equation}

Let $\gbigv$ be a locally free
$\nbigo_{\nbigmbar_{\gminiq}\setminus Z}(\ast H_{\gminiq})$-module
with Dirac type singularity
enhanced by
good filtered bundles
$\nbigp_{\ast}\gbigvhat_{\nu}$ $(\nu=0,\infty)$.
Due to the scattering map,
the restriction
$\gbigv_{|\nbigm_{\gminiq}\setminus \sfp^{-1}(\ttD)}$
induces a locally free $\nbigo_{T\setminus \ttD}$-module $\ttE'$.
For $P\in \ttD$,
we take $(\alpha_P,t_P)\in \nbigm^{\cov}\setminus Z^{\cov}$
which is mapped to $P$.
Take $U_{P}$ be a small neighbourhood of
$\alpha_P$ in $\cnum^{\ast}$.
Set $U_P^{\ast}:=U_P\setminus\{\alpha_P\}$.
There exists a natural isomorphism
$\ttE'_{|\sfp(U_P^{\ast}\times\{t_P\})}
\simeq
 \gbigv^{\cov}_{|U_P^{\ast}\times\{t_P\}}$.
By gluing $\ttE'$
and $\bigl(\gbigv^{\cov}_{|U_P\times\{t_P\}}\bigr)(\ast \alpha_P)$
$(P\in \ttD)$,
we obtain a locally free
$\nbigo_T(\ast\ttD)$-module $\ttEtilde$.
It is independent of a choice of $(\alpha_P,t_P)$.

For $P\in \ttD$,
choose $\alpha_P\in\cnum^{\ast}$
which is mapped to $P$ by the projection $\cnum^{\ast}\lrarr T$.
We set
\[
t_{P,i}=s_{P,i}+
 \frac{\log|\alpha_P|}{\log|\gminiq|}.
\]
Then, 
$Z^{\cov}\cap
 (\{\alpha_P\}\times\real)
=\bigl\{
 (\alpha_P,t_{P,i})\,|\,i=1,\ldots,m(P)
 \bigr\}$ holds.
We formally set
$t_{P,0}:=-\infty$
and $t_{P,m(P)+1}:=\infty$.
We choose $t_{P,i}<t'_{P,i}<t_{P,i+1}$
for $i=0,\ldots,m(P)$.
Let $\nbigk_{P,i}$ $(i=0,\ldots,m(P))$
denote the formal completion of
$\gbigv^{\cov}$ at $(\alpha_{P},t_{P,i}')$.
They induce lattices of
$\ttEtilde_{|\Phat}$.

We obtain a locally free $\nbigo_{T}$-submodule
$\ttE_-$ of $\ttEtilde$
determined by the lattices $\nbigk_{P,0}$ $(P\in \ttD)$.
Similarly,
we obtain a locally free $\nbigo_T$-submodule
$\ttE_+$ of $\ttEtilde$
determined by the lattices $\nbigk_{P,m(P)}$ $(P\in \ttD)$.
Note that we have
the $\gminiq$-difference $\cnum[y,y^{-1}]$-module $\vecV$
with a parabolic structure at finite place
corresponding to $\gbigv$.
\begin{itemize}
\item
If $\log|\gminiq|>0$,
let $\gbigf_-$ be the anti-Harder-Narasimhan filtration
indexed by $(\rnum,\leq)$ on $\ttE_-$
corresponding to the germ of $\vecV$ at $y=0$,
and 
let $\gbigf_+$ be the anti-Harder-Narasimhan filtration
indexed by $(\rnum,\leq)$ on $\ttE_+$
corresponding to the germ of $\vecV$ at $y=\infty$,
\item
If $\log|\gminiq|<0$,
we replace $y=0$ and $y=\infty$.
\end{itemize}
Moreover, good filtered bundles
over $\vecV_{|\nuhat}$
induce filtrations
on $\Gr^{\gbigf_{\pm}}(\ttE_{\pm})$
as in \S\ref{subsection;19.2.4.30}.
(See \S\ref{subsection;19.2.4.31}
for the relation with the growth order of the norms.)

In this way, 
good filtered bundles with Dirac type singularity 
on $(\nbigmbar_{\gminiq};H_{\gminiq},Z)$
induces locally free $\nbigo_{T}(\ast\ttD)$-modules
with $\gminiq$-difference parabolic structure.
The following is easy to see.
\begin{prop}
The above procedure induces an equivalence between 
good filtered bundles with Dirac type singularity 
on $(\nbigmbar_{\gminiq};H_{\gminiq},Z)$
and locally free $\nbigo_{T}(\ast\ttD)$-modules
with $\gminiq$-difference parabolic structure
at $(\vecs_{P})_{P\in\ttD}$.
Here, $Z$ and $(\vecs_P)_{P\in\ttD}$
are related as in {\rm(\ref{eq;19.2.4.40})}.
Moreover,
it preserves the degree.
\hfill\qed
\end{prop}

\subsubsection{Filtered objects associated to meromorphic monopoles}

Let $(E,h,\nabla,\phi)$ be a meromorphic monopole
on $\nbigm\setminus Z$.
We fix $\lambda\in\cnum$ such that $|\lambda|\neq 1$.
Take $(\tte_1,\tts_1)$ as in \S\ref{subsection;19.1.26.10}.
Let $\gminiq^{\lambda}(\tte_1,\tts_1)$
and 
$\gminit^{\lambda}(\tte_1,\tts_1)$ 
denote 
$\gminiq^{\lambda}$ and $\gminit^{\lambda}$
in (\ref{eq;19.1.26.11})
to emphasize the dependence on $(\tte_1,\tts_1)$.
Then, we have the associated parabolic
$\gminiq^{\lambda}(\tte_1,\tts_1)$-difference module,
and hence the associated
filtered object $\ttEtilde_{(\tte_1,\tts_1)\,\ast}$
on the elliptic curve
$T^{\lambda}(\tte_1,\tts_1)=
 \cnum^{\ast}/\gminiq^{\lambda}(\tte_1,\tts_1)^{\seisuu}$.
It is easy to observe that
$T^{\lambda}(\tte_1,\tts_1)$
is independent of $(\tte_1,\tts_1)$
by the construction.
Moreover, we obtain the following.

\begin{thm}
$\ttH^{(\gminit^{\lambda}(\tte_1,\tts_1))}(\ttEtilde_{(\tte_1,\tts_1)\,\ast})$
are independent of the choice of $(\tte_1,\tts_1)$.
\end{thm}

\subsection{Acknowledgement}

I owe much to Carlos Simpson whose ideas 
on the Kobayashi-Hitchin correspondence
are fundamental in this study.
I have been stimulated by the works
of Maxim Kontsevich and Yan Soibelman
on $\gminiq$-difference modules.
I am clearly  influenced by the works of
Benois Charbonneau and Jacques Hurtubise
\cite{Charbonneau-Hurtubise}
and Sergey Cherkis and Anton Kapustin
\cite{Cherkis-Kapustin1, Cherkis-Kapustin2}.
I am grateful to Claude Sabbah
for his kindness and discussions on many occasions.
A part of this study was done
during my visits at the Tata Institute of Fundamental Research
and the International Center for Theoretical Sciences.
I appreciate Indranil Biswas for his excellent hospitality.
I thank Yoshifumi Tsuchimoto
and Akira Ishii for their constant encouragement.
I thank
Indranil Biswas,
Sergey Cherkis, Jacques Hurtubise,
Ko-ki Ito,
Hisashi Kasuya,
Maxim Kontsevich,
Masa-Hiko Saito,
Yota Shamoto,
Carlos Simpson,
and Masaki Yoshino
for discussions.

I am partially supported by
the Grant-in-Aid for Scientific Research (S) (No. 17H06127),
the Grant-in-Aid for Scientific Research (S) (No. 16H06335),
and the Grant-in-Aid for Scientific Research (C) (No. 15K04843),
Japan Society for the Promotion of Science.

\section{Good filtered formal $\gminiq$-difference modules}

\subsection{Formal $\gminiq$-difference modules}

We review a classification of formal $\gminiq$-difference modules
to prepare notations.
See \cite{van-der-Put-Reversat, Ramis-Sauloy-Zhang, Sauloy2000, Sauloy2004}.
Some statements will be proved
though they are standard and well known.
It is just to explain 
that the statements are valid even in the case 
where $\gminiq$ is a root of unity.

\subsubsection{Preliminary}

Take any non-zero complex number $\gminiq$.
Set $\gminiq^{\seisuu}:=\{\gminiq^n\,|\,n\in\seisuu\}$,
which is a subgroup of $\cnum^{\ast}$.
If $\gminiq$ is not a root of $1$,
then $\gminiq^{\seisuu}$ is naturally isomorphic to
$\seisuu$.
If $\gminiq$ is a primitive $k$-th root of $1$,
then 
$\gminiq^{\seisuu}=
 \bigl\{\mu\in\cnum^{\ast}\,\big|\,\mu^k=1\bigr\}$.
We fix $a\in\cnum$ such that $\exp(a)=\gminiq$,
and we put $\gminiq_m:=\exp(a/m)$
for any positive integer $m$.

We set $\nbigk:=\cnum(\!(y)\!)$ and $\nbigr:=\cnum[\![y]\!]$
where $y$ is a variable.
We fix $m$-th roots $y_m$ of $y$
for any positive integers $m$
such that $(y_{mn})^n=y_m$
for any $(m,n)\in\seisuu_{>0}^2$.
We set $\nbigk_m:=\cnum(\!(y_m)\!)$
and $\nbigr_m:=\cnum[\![y_m]\!]$.
Let $\Phi^{\ast}$ be the automorphisms of
$\nbigk_m$ determined by
$\Phi^{\ast}(f)(y_m)=f(\gminiq_my_m)$.

A $\gminiq_m$-difference $\nbigk_m$-module
is a finite dimensional $\nbigk_m$-vector space $\nbigv$
equipped with a $\cnum$-linear isomorphism
$\Phi^{\ast}:\nbigv\lrarr\nbigv$
such that
$\Phi^{\ast}(fs)=\Phi^{\ast}(f)\Phi^{\ast}(s)$
for any $f\in\nbigk_m$ and $s\in\nbigv$.

A morphism of
$\gminiq_m$-difference $\nbigk_m$-modules
$g:(\nbigv_1,\Phi^{\ast})\lrarr (\nbigv_2,\Phi^{\ast})$
is defined to be 
a morphism of $\nbigk_m$-vector spaces
$g:\nbigv_1\lrarr\nbigv_2$
such that 
$g\circ \Phi^{\ast}=\Phi^{\ast}\circ g$.

Let $\Diff_m(\nbigk,\gminiq)$
be the category of $\gminiq_m$-difference
$\nbigk_m$-modules.
If $m=1$, it is also denoted by
$\Diff(\nbigk,\gminiq)$.

Let $(\nbigv_i,\Phi^{\ast})\in\Diff_m(\nbigk,\gminiq)$.
The operators $\Phi^{\ast}$ on $\nbigv_1\oplus\nbigv_2$
and $\nbigv_1\otimes\nbigv_2$
are defined by
$\Phi^{\ast}(v_1\oplus v_2)=\Phi^{\ast}(v_1)\oplus\Phi^{\ast}(v_2)$
and 
$\Phi^{\ast}(v_1\otimes v_2)=\Phi^{\ast}(v_1)\otimes\Phi^{\ast}(v_2)$.
Thus,
we obtain the direct sum and the tensor product
on $\Diff_m(\nbigk,\gminiq)$.
For $(\nbigv,\Phi^{\ast})\in \Diff_m(\nbigk,\gminiq)$,
let $\nbigv^{\lor}:=\Hom_{\nbigk}(\nbigv,\nbigk)$.
We define the operator $\Phi^{\ast}$
on $\nbigv^{\lor}$
by $\Phi^{\ast}(f)(v):=f(\Phi^{\ast}(v))$.
We set
$(\nbigv,\Phi^{\ast})^{\lor}
:=(\nbigv^{\lor},\Phi^{\ast})$.

\subsubsection{Pull back and push-forward}

Let 
$(\nbigv,\Phi^{\ast})\in\Diff_m(\nbigk,\gminiq)$.
For any $n\in\seisuu_{>0}$,
we define a $\cnum$-automorphism $\Phi^{\ast}$
on $\nbigv\otimes_{\nbigk_m}\nbigk_{mn}$
by $\Phi^{\ast}(s\otimes g)=\Phi^{\ast}(s)\otimes\Phi^{\ast}(g)$.
In this way,
we obtain a $\gminiq_{mn}$-difference $\nbigk_{mn}$-module
$(\nbigv\otimes_{\nbigk_m}\nbigk_{mn},\Phi^{\ast})$.
It induces a functor
$(\sfp_{m,nm})^{\ast}:
 \Diff_m(\nbigk,\gminiq)\lrarr\Diff_{nm}(\nbigk,\gminiq)$.

Let $(\nbigv,\Phi^{\ast})\in \Diff_{nm}(\nbigk,\gminiq)$
for $n,m\in\seisuu_{>0}$.
We may naturally regard $\nbigv$ 
as a $\gminiq_m$-difference $\nbigk_m$-module.
Thus, we obtain a functor
$(\sfp_{m,nm})_{\ast}:
 \Diff_{nm}(\nbigk,\gminiq)
\lrarr
 \Diff_{m}(\nbigk,\gminiq)$.

For any $(\nbigv,\Phi^{\ast})\in\Diff_m(\nbigk,\gminiq)$,
there exists a natural isomorphism
\[
 (\sfp_{nm,m})_{\ast}
  (\sfp_{nm,m})^{\ast}(\nbigv,\Phi^{\ast})
\simeq
 (\nbigv,\Phi^{\ast})
 \otimes
 (\sfp_{nm,m})_{\ast}(\nbigk_{nm},\Phi^{\ast}).
\]

Let $\Gal(nm,m)$ denote the Galois group of 
$\nbigk_{nm}/\nbigk_m$,
which is naturally identified with
$\{\mu\in \cnum^{\ast}\,|\,\mu^{n}=1\}$
by the action $(\mu\bullet f)(y_{nm})=f(\mu y_{nm})$.
Note that
$\Phi^{\ast}(\mu\bullet f)=\mu\bullet \Phi^{\ast}(f)$.

Let $(\nbigv,\Phi^{\ast})\in\Diff_{nm}(\nbigk,\gminiq)$.
We set
$\mu^{\ast}(\nbigv):=\nbigv$ as $\cnum$-vector space.
Any element $v\in\nbigv$ is denoted by $\mu^{\ast}(v)$
when we regard it as an element of $\mu^{\ast}(\nbigv)$.
We regard $\mu^{\ast}\nbigv$
as a $\nbigk_{nm}$-vector space
by $f\cdot \mu^{\ast}(v):=
 \mu^{\ast}\bigl( (\mu^{-1})^{\ast}(f)v\bigr)$.
Note that for any $(\nbigv,\Phi^{\ast})\in\Diff_{nm}(\nbigk,\gminiq)$,
there exists a natural isomorphism
\[
 (\sfp_{m,nm})^{\ast}\circ(\sfp_{m,nm})_{\ast}(\nbigv,\Phi^{\ast})
\simeq 
 \bigoplus_{\mu\in\Gal(nm,m)}
 \mu^{\ast}(\nbigv,\Phi^{\ast}).
\]

A $\nbigk_{nm}$-vector space $\nbigv$ is called
$\Gal(nm,m)$-equivariant
when a homomorphism
$\Gal(nm,m)\lrarr \Aut_{\cnum}(\nbigv)$
is given
such that 
$\mu\bullet (fv)=(\mu\bullet f)\cdot (\mu\bullet v)$
for any $\mu\in\Gal(nm,n)$, $f\in\nbigk_m$ and $v\in\nbigv$.

An object $(\nbigv,\Phi^{\ast})\in\Diff_m(\nbigk,\gminiq)$
is called $\Gal(nm,m)$-equivariant
when 
$\nbigv$ is a finite dimensional $\Gal(nm,m)$-equivariant
$\nbigk_{nm}$-vector space
such that 
$\Phi^{\ast}\circ(\mu\bullet v)=\mu\bullet \Phi^{\ast}(v)$
for any $\mu\in\Gal(nm,n)$ and $v\in\nbigv$.

For any $(\nbigv,\Phi^{\ast})\in\Diff_{m}(\nbigk,\gminiq)$,
$(\sfp_{m,nm})^{\ast}(\nbigv,\Phi^{\ast})$
is naturally $\Gal(nm,n)$-equivariant.
Conversely,
let $(\nbigv_1,\Phi^{\ast})$
be a $\Gal(nm,m)$-equivariant object in
$\Diff_{nm}(\nbigk,\gminiq)$.
We set $\nbigv_0:=
 \{v\in\nbigv_1\,|\,\mu\bullet v=v\,\,(\forall\mu\in
 \Gal(nm,n))\}$.
Thus, we obtain
$(\nbigv_0,\Phi^{\ast})\in
 \Diff_{m}(\nbigk,\gminiq)$,
which is called the descent of
$(\nbigv_1,\Phi^{\ast})$.
Then,
$(\sfp_{m,nm})^{\ast}(\nbigv_0,\Phi^{\ast})$
is $\Gal(nm,n)$-equivariantly isomorphic to
$(\nbigv_1,\Phi^{\ast})$.
In particular,
$(\nbigv_0,\Phi^{\ast})$
is isomorphic to the descent of
$(\sfp_{m,nm})^{\ast}(\nbigv_0,\Phi^{\ast})$.

\subsubsection{A splitting lemma}

Let $(\nbigv,\Phi^{\ast})\in\Diff_m(\nbigk,\gminiq)$.
For any $\nbigr_m$-lattice $\nbigl\subset\nbigv$
such that 
$y_m^{\ell}\Phi^{\ast}(\nbigl)\subset\nbigl$,
we have the induced endomorphism
$\sigma(y_m^{\ell}\Phi^{\ast};\nbigl)$ 
of $\nbigl_{|0}:=\nbigl/y_m\nbigl$
obtained as follows:
for any $s\in\nbigl_{|0}$,
we take $\stilde\in\nbigl$ which induces $s$,
and let $\sigma(y_m^{\ell}\Phi^{\ast};\nbigl)(s)\in\nbigl_{|0}$
denote the element induced by 
$y_m^{\ell}\Phi^{\ast}(\stilde)\in\nbigl$.

The following lemma is standard.
\begin{prop}
\label{prop;18.8.7.10}
Suppose that there exist an $\nbigr_m$-lattice 
$\nbigl\subset\nbigv$ and an integer $\ell$
such that the following holds.
\begin{itemize}
\item
$y_m^{\ell}\Phi^{\ast}(\nbigl)\subset\nbigl$
holds.
In particular,
we obtain the induced endomorphism
$F:=\sigma(y_m^{\ell}\Phi^{\ast};\nbigl)$
of $\nbigl_{|0}$.
\item
There exists a decomposition
$\nbigl_{|0}=L_1\oplus L_2$
such that 
$F(L_i)\subset L_i$.
\item
Let $\Sp(F,L_i)$
be the set of eigenvalues of $F_{|L_i}$.
Then, 
$\bigl(\gminiq_m^{\seisuu}\cdot \Sp(F,L_1)\bigr)
\cap \Sp(F,L_2)=\emptyset$.
\end{itemize}
Then, there exists a unique decomposition 
$\nbigl=\nbigl_1\oplus\nbigl_2$
of $\nbigr_m$-modules
such that
(i) $y_m^{\ell}\Phi^{\ast}(\nbigl_i)\subset\nbigl_i$,
(ii) $\nbigl_{i|0}=L_i$.
\end{prop}
\pf
We give only an indication.
For any ring $R$ and a positive integer $r$,
let $M_r(R)$ 
denote the space of
$r$-square matrices with $R$-coefficient.
For any ring $R$ and positive integers $r_i$ $(i=1,2)$,
let $M_{r_1,r_2}(R)$ denote the space of
$(r_1\times r_2)$-matrices with $R$-coefficient.
For a decomposition $r=r_1+r_2$ $(r_i>0)$,
any element $C$ of $M_r(R)$
is expressed as 
\[
 C=\left(
 \begin{array}{cc}
 C_{11} & C_{12}\\
 C_{21} & C_{22}
 \end{array}
 \right),
\]
where
$C_{ij}\in M_{r_i,r_j}(R)$.

\begin{lem}
\label{lem;18.8.7.11}
Let $r=r_1+r_2$ $(r_i> 0)$ be a decomposition.
Let $A\in y_m^{-\ell}M_r(\nbigr_m)$.
We obtain $A_{ij}$ $(1\leq i,j\leq 2)$ as above,
which have the expansions
$A_{ij}=\sum_{k=-\ell}^{\infty} A_{ij;k}y_m^k$.
We assume the following.
\begin{itemize}
\item
 $A_{ij,-\ell}=0$ if $i\neq j$.
\item
 $\bigl(
 q_m^{\seisuu}\Sp(A_{11,-\ell})
 \bigr)
 \cap
 \Sp(A_{22,-\ell})=\emptyset$,
where 
 $\Sp(A_{ii,-\ell})$ denote the sets of eigenvalues of
 $A_{ii,-\ell}$.
\end{itemize}
Then, there exists
$G\in \GL_r(\nbigr_m)$
such that 
(i) $G_{ii}$ are identity matrices in $M_{r_i}(\nbigr_m)$,
(ii) $G_{ij|0}=0$ $(i\neq j)$,
(iii) $\bigl(G(y_m)^{-1}AG(\gminiq_my_m)\bigr)_{ij}=0$ $(i\neq j)$.
\end{lem}
\pf
Let $\Atilde\in y_m^{-\ell}M_r(\nbigr_m)$
determined by
(i) $\Atilde_{ij}=0$ $(i\neq j)$,
(ii) $\Atilde_{ii}=A_{ii}$.
Let $U$ denote a matrix in $y_m^{-\ell}M_r(\nbigr_m)$
such that
(i) $U_{ij}=0$ $(i\neq j)$,
(ii) $U_{ii,-\ell}=0$.
We consider the following equation for $G$ and $U$:
\[
A(y_m)G(\gminiq_my_m)=
 G(y_m)\bigl(\Atilde(y_m)+U(y_m)\bigr).
\]
It is equivalent to the following equations:
\begin{equation}
\label{eq;18.11.19.1}
 A_{12}(y_m)G_{21}(\gminiq_my_m)=U_{11}(y_m),
\quad
 A_{22}(y_m)G_{21}(\gminiq_my_m)-G_{21}(y_m)A_{11}(y_m)
+A_{21}(y_m)-G_{21}(y_m)U_{11}(y_m)=0,
\end{equation}
\begin{equation}
\label{eq;18.11.19.2}
 A_{21}(y_m)G_{12}(\gminiq_my_m)=U_{22}(y_m),
\quad
 A_{11}(y_m)G_{12}(\gminiq_my_m)
-G_{12}(y_m)A_{22}(y_m)
+A_{12}(y_m)
-G_{12}(y_m)U_{22}(y_m)=0.
\end{equation}
From (\ref{eq;18.11.19.1}),
we obtain the following equation for $G_{21}$:
\[
 A_{22}(y_m)G_{21}(\gminiq_my_m)
-G_{21}(y_m)A_{11}(y_m)
+A_{21}(y_m)
-G_{21}(y_m)A_{12}(y_m)G_{21}(\gminiq_my_m)=0.
\]
It is equivalent to the following equations for $G_{21;k}$
$(k\in\seisuu_{\geq 0})$.
\begin{multline}
\label{eq;18.8.7.1}
 A_{22;-\ell}G_{21;k}\gminiq_m^{k}
-G_{21;k}A_{11;-\ell}
 \\
+\sum_{\substack{i+j=k-\ell\\ 0\leq j<k}}
 A_{22;i}G_{21;j}\gminiq_m^j
-\sum_{\substack{i+j=k-\ell\\ 0\leq j<k}}
 G_{21;j}A_{11;i}
+A_{21;-\ell+k}
-\sum_{\substack{i+j+p=k-\ell\\ j>-\ell}}
 G_{21;i}A_{12;j}G_{21;p}\gminiq_m^p=0.
\end{multline}
For $k=0$, 
we have a solution $G_{21;0}=0$.
For $k\geq 1$,
we can determine $G_{21;k}$ in an inductive way
by using (\ref{eq;18.8.7.1}).
We obtain $U_{11}$
from (\ref{eq;18.11.19.1}).
Similarly, we obtain
$G_{12}$ and $U_{22}$
from (\ref{eq;18.11.19.2}).
\hfill\qed

\vspace{.1in}
The following lemma is also standard
and easy to see by using 
the power series expansions.
\begin{lem}
\label{lem;18.8.7.12}
Let $(r_1,r_2)\in\seisuu_{>0}^2$
and $A_i=\sum_{k\geq -\ell}A_{i;k}y_m^k\in y_m^{-\ell}M_{r_i}(\nbigr_m)$.
Assume the following.
\begin{itemize}
\item
$\bigl(
 \gminiq^{\seisuu}\Sp(A_{1;-\ell})
 \bigr)
\cap
 \Sp(A_{2;-\ell})=\emptyset$,
where 
$\Sp(A_{i;-\ell})$
denote the sets of the eigenvalues of $A_{i;-\ell}$.
\end{itemize}
Let $H\in M_{r_1,r_2}(\nbigr_m)$
such that
$A_2(y_m)H(\gminiq_my_m)=H(y_m)A_1(y_m)$.
Then, $H=0$.
\hfill\qed
\end{lem}

We obtain the claim of Proposition \ref{prop;18.8.7.10}
from 
Lemma \ref{lem;18.8.7.11}
and 
Lemma \ref{lem;18.8.7.12}.
\hfill\qed

\subsubsection{Fuchsian $\gminiq$-difference modules}

We recall the Fuchsian (regular singular) condition
of $\gminiq$-difference modules
by following \cite{Ramis-Sauloy-Zhang}.

\begin{df}
\label{df;18.8.19.1}
A $\gminiq_m$-difference $\nbigk_m$-module
$(\nbigv,\Phi^{\ast})$
is called Fuchsian if there exists 
an $\nbigr_{m}$-lattice $\nbigl$
such that $\Phi^{\ast}(\nbigl)=\nbigl$.
Let $\Diff_m(\nbigk,\gminiq;0)
\subset\Diff_m(\nbigk,\gminiq)$
denote the full subcategory 
of Fuchsian $\gminiq_m$-difference $\nbigk_m$-modules.

\hfill\qed
\end{df}

Let $(\nbigv,\Phi^{\ast})\in\Diff_m(\nbigk,\gminiq;0)$.
Let $\nbigl$ be an $\nbigr_m$-lattice such that
$\Phi^{\ast}(\nbigl)=\nbigl$.
We obtain the induced automorphism
$\sigma(\Phi^{\ast};\nbigl)$ of $\nbigl_{|0}$,
and let $\Sp(\sigma(\Phi^{\ast};\nbigl))\subset\cnum^{\ast}$ 
denote the set of eigenvalues.
Let $[\Sp(\sigma(\Phi^{\ast});\nbigl)]$ denote the image of
$\Sp(\sigma(\Phi^{\ast};\nbigl))$ by
$\cnum^{\ast}\lrarr \cnum^{\ast}/\gminiq_m^{\seisuu}$.
There exists the decomposition
$\nbigl_{|0}=\bigoplus_{\gminio\in[\Sp(\sigma(\Phi^{\ast};\nbigl))]}
 L_{\gminio}$
such that
(i) $\sigma(\Phi^{\ast};\nbigl)(L_{\gminio})=L_{\gminio}$,
(ii) the eigenvalues of
$\sigma(\Phi^{\ast};\nbigl)_{|L_{\gminio}}$
are contained in $\gminio$.
We also obtain the following lemma
from 
Lemma \ref{lem;18.8.7.11}
and 
Lemma \ref{lem;18.8.7.12}.

\begin{lem}
There exists a unique decomposition
$(\nbigl,\Phi^{\ast})
=\bigoplus_{\gminio\in [\Sp(\sigma(\Phi^{\ast});\nbigl)]}
 (\nbigl_{\gminio},\Phi^{\ast})$
such that
$\nbigl_{\gminio|0}=L_{\gminio}$.
The set
$[\Sp(\sigma(\Phi^{\ast};\nbigl))]$
is independent of the choice of 
an $\nbigr_m$-lattice such that $\Phi^{\ast}(\nbigl)=\nbigl$.
\hfill\qed
\end{lem}

We set $[\Sp(\sigma(\Phi^{\ast};\nbigv))]:=
 [\Sp(\sigma(\Phi^{\ast};\nbigl))]$
for an $\nbigr_m$-lattice $\nbigl$ such that
$\Phi^{\ast}(\nbigl)=\nbigl$,
which is independent of the choice of $\nbigl$.

\begin{example}
\label{example;18.12.20.1}
Let $V$ be a finite dimensional $\cnum$-vector space.
For any $f\in\GL(V)$,
we set $\VV_m(V,f):=V\otimes\nbigk_m$,
and we define the $\gminiq_m$-difference operator $\Phi^{\ast}$
on $\VV_m(V,f)$
by  $\Phi^{\ast}(s)=f(s)$ for any $s\in V$.
Then, $(\VV_m(V,f),\Phi^{\ast})\in \Diff_m(\nbigk,\gminiq;0)$.
We have
$[\Sp(\sigma(\Phi^{\ast},\VV_m(V,f)))]=[\Sp(f)]$.

Similarly,
for any $r\in\seisuu_{>0}$ and $A\in \GL_r(\cnum)$,
let $\VV_m(A)$ denote 
the $\nbigk_m$-vector space 
with a frame $\vece=(e_1,\ldots,e_r)$
equipped with the $\gminiq_m$-difference operator
defined by
$\Phi^{\ast}(\vece)=\vece A$.
Then, $\VV_m(A)\in\Diff_m(\nbigk,\gminiq;0)$.
We have
$[\Sp(\sigma(\Phi^{\ast},\VV_m(A)))]=[\Sp(A)]$.
\hfill\qed
\end{example}

\begin{rem}
Let $S\subset\cnum^{\ast}$ be any subset 
such that the induced map
$S\lrarr \cnum^{\ast}/\gminiq^{\seisuu}$ is a bijection.
As proved in {\rm\cite{Sauloy2004}},
if $\gminiq_m$ is not a root of $1$,
for any  $(\nbigv,\Phi^{\ast})\in\Diff_m(\nbigk,\gminiq;0)$,
there exists $A\in \GL_r(\cnum)$
such that
(i) $(\nbigv,\Phi^{\ast})\simeq \VV_m(A)$,
(ii) $\Sp(A)\subset S$.
If $\gminiq_m$ is a root of $1$,
it does not hold in general.
\hfill\qed
\end{rem}

\subsubsection{Formal pure isoclinic
$\gminiq$-difference modules}
\label{subsection;18.11.22.10}

We recall the notion of pure isoclinic $\gminiq$-difference modules
\cite{Ramis-Sauloy-Zhang}.

\begin{df}
Let $\omega$ be a rational number.
We say that
$(\nbigv,\Phi^{\ast})\in \Diff_m(\nbigk,\gminiq)$ 
is pure isoclinic of slope $\omega$
if the following holds.
\begin{itemize}
\item Take any $m_1\in m\seisuu_{>0}$
such that $\omega\in \frac{1}{m_1}\seisuu$.
Then, there exists an $\nbigr_{m_1}$-lattice 
$\nbigl\subset \sfp_{m_1,m}^{\ast}\nbigv$
such that
$y_{m_1}^{m_1\omega}\Phi^{\ast}(\nbigl)=\nbigl$.
\end{itemize}
Let $\Diff_m(\nbigk,\gminiq;\omega)
\subset \Diff_m(\nbigk,\gminiq)$
denote the full subcategory of
pure isoclinic $\gminiq_m$-difference
$\nbigk_m$-modules of slope $\omega$.
\hfill\qed
\end{df}

\begin{rem}
Recall that $|\gminiq|>1$ is assumed in {\rm\cite{Ramis-Sauloy-Zhang}}.
In the case $|\gminiq|<1$,
it seems better to change the signature of the slope
in the relation with the analytic classification
of $\gminiq$-difference modules.
However, because we also study the case $|\gminiq|=1$,
we do not change the signature.
\hfill\qed
\end{rem}

\begin{lem}
Let $(\nbigv_i,\Phi^{\ast})\in\Diff_m(\nbigk,\gminiq;\omega_i)$
$(i=1,2)$.
Let 
$f:(\nbigv_1,\Phi^{\ast})\lrarr(\nbigv_2,\Phi^{\ast})$
be a morphism in $\Diff_m(\nbigk,\gminiq)$.
If $\omega_1\neq\omega_2$, then $f=0$.
\end{lem}
\pf
It follows from Lemma \ref{lem;18.8.7.12}.
\hfill\qed

\begin{lem}
Let $(\nbigv,\Phi^{\ast})\in\Diff_m(\nbigk,\gminiq)$.
Suppose that there exist
a rational number $\omega$
and a finite family of
subobjects
$(\nbigv_i,\Phi^{\ast})
\subset
 (\nbigv,\Phi^{\ast})$
in $\Diff_m(\nbigk,\gminiq)$
such that
(i) $\nbigv=\sum_{i=1}^N\nbigv_i$,
(ii) $(\nbigv_i,\Phi^{\ast})\in\Diff_m(\nbigk,\gminiq;\omega)$.
Then, 
$\nbigv\in\Diff_m(\nbigk,\gminiq;\omega)$.
\end{lem}
\pf
We may assume that $\omega\in\frac{1}{m}\seisuu$.
There exist $\nbigr_m$-lattices $\nbigl_i\subset\nbigv_i$
such that 
$y_m^{m\omega}\Phi^{\ast}(\nbigl_i)=\nbigl_i$.
We put $\nbigl:=\sum\nbigl_i$,
which is an $\nbigr_m$-lattice of $\nbigv$.
We have
$y_m^{m\omega}\Phi^{\ast}(\nbigl)=\nbigl$.
\hfill\qed

\begin{lem}
Let $\omega=\ell/k\in\rnum$,
where $\ell\in\seisuu$ and $k\in\seisuu_{>0}$.
Then, $(\nbigv,\Phi^{\ast})\in\Diff_m(\nbigk,\gminiq)$
is pure isoclinic of slope $\omega$
if and only if there exists a lattice 
$\nbigl\subset\nbigk$
such that
$(\Phi^{\ast})^k\nbigl=y_m^{-\ell}\nbigl$.
\end{lem}
\pf
Set $m_1:=km$.
Suppose that $(\nbigv,\Phi^{\ast})\in\Diff_m(\nbigk,\gminiq;\omega)$.
We obtain
$(\nbigv_1,\Phi^{\ast}):=
 \sfp_{m,m_1}^{\ast}(\nbigv,\Phi^{\ast})$.
There exists a lattice
$\nbigl_1\subset\nbigv_1$
such that 
$y_{m_1}^{\ell}\Phi^{\ast}\nbigl_1=\nbigl_1$.
We set
$\nbigl_2:=\sum_{\mu\in\Gal(m_1,m)}
 \mu^{\ast}\nbigl_1$.
Then, $\nbigl_2$ is $\Gal(m_1,m)$-equivariant,
and 
$y_{m_1}^{\ell}\Phi^{\ast}\nbigl_2=\nbigl_2$
holds.
Let $\nbigl$ be the $\Gal(m_1,m)$-invariant part
of $\nbigl_2$.
We obtain
$y_{m_1}^{k\ell}(\Phi^{\ast})^k\nbigl=\nbigl$.
Hence, we obtain a lattice with the desired property.

Suppose that a lattice $\nbigl$ of $\nbigv$
has the desired property.
We set $\nbigl':=\nbigl\otimes_{\nbigr_m}\nbigr_{m_1}$.
We have
$(y_{m_1}^{\ell}\Phi^{\ast})^k\nbigl'=\nbigl'$.
We set
$\nbigl'':=
 \sum_{j=0}^{k-1}
 (y_{m_1}^{\ell}\Phi^{\ast})^j\nbigl'$.
Then, we obtain
$(y_{m_1}^{\ell}\Phi^{\ast})\nbigl''=\nbigl''$.
\hfill\qed

\subsubsection{Basic examples of pure isoclinic $\gminiq_m$-difference modules}
\label{subsection;18.12.21.2}

Let  $\omega\in\rnum$.
If $m\omega\in\seisuu$,
we obtain
$\LL_m(\omega)\in\Diff_m(\nbigk,\gminiq;\omega)$
by the $\nbigk_m$-vector space
$\nbigk_m\cdot\sle_{m,\omega}$
with the operator
$\Phi^{\ast}(\sle_{m,\omega})
=y_m^{-m\omega}\sle_{m,\omega}$.
For $\omega\in\rnum\setminus\frac{1}{m}\seisuu$,
we express  $\omega=\ell_0/m_0$
for $\ell_0\in\seisuu$ and $m_0\in\seisuu_{>0}$
with $\gcd(m_0,\ell_0)=1$.
Let $m_1$ be the least common multiple of
$m_0$ and $m$.
We obtain $\LL_{m_1}(\omega)\in\Diff_{m_1}(\nbigk,\gminiq)$.
We set
$\LL_m(\omega):=(\sfp_{m,m_1})_{\ast}\LL_{m_1}(\omega)
\in \Diff_{m}(\nbigk,\gminiq)$.

\begin{lem}
$\LL_m(\omega)\in\Diff_m(\nbigk,\gminiq;\omega)$.
\end{lem}
\pf
Set $b:=m_1/m$.
There exists a natural isomorphism
$(\sfp_{m,m_1})^{\ast}\LL_m(\omega)
\simeq
 \bigoplus_{\mu\in\Gal(m_1,m)}
 \mu^{\ast}\LL_{m_1}(\omega)
\simeq
 \bigoplus_{\mu\in\Gal(m_1,m)}
 \LL_{m_1}(\omega)
\otimes
 \VV_{m_1}(\mu^{b})$.
Then, the claim is clear.
\hfill\qed

\vspace{.1in}
Let $(\nbigv,\Phi^{\ast})\in\Diff_m(\nbigk,\gminiq;\omega)$.
If $\omega\in\frac{1}{m_1}\seisuu$
for $m_1\in m\seisuu_{>0}$,
then there exist
$(\nbigu^{\reg},\Phi^{\ast})
\in\Diff_{m_1}(\nbigk,\gminiq;0)$
and an isomorphism
$(\sfp_{m,m_1})^{\ast}(\nbigv,\Phi^{\ast})
 \simeq
 \LL_{m_1}(\omega)\otimes
 (\nbigu^{\reg},\Phi^{\ast})$.

\subsubsection{Slope decompositions}

\begin{df}
Let $(\nbigv,\Phi^{\ast})\in\Diff_m(\nbigk,\gminiq)$.
A decomposition
$(\nbigv,\Phi^{\ast})=
 \bigoplus_{\omega\in\rnum}(\nbigv_{\omega},\Phi^{\ast})$
in $\Diff_m(\nbigk,\gminiq_m)$
is called a slope decomposition
if 
$(\nbigv_{\omega},\Phi^{\ast})
\in\Diff_m(\nbigk,\gminiq)$.
\hfill\qed
\end{df}

We obtain the uniqueness of slope decompositions
from Lemma \ref{lem;18.8.7.12}.
\begin{lem}
If $\nbigv=\bigoplus \nbigv^{(i)}_{\omega}$ $(i=1,2)$
are slope decompositions of 
$(\nbigv,\Phi^{\ast})
 \in\Diff_m(\nbigk,\gminiq)$,
then $\nbigv^{(1)}_{\omega}=\nbigv^{(2)}_{\omega}$ hold
for any $\omega\in\rnum$.
\hfill\qed
\end{lem}
As a corollary, we obtain the following.
\begin{cor}
Let $(\nbigv,\Phi^{\ast})\in\Diff_m(\nbigk,\gminiq)$.
Let $\sfp_{m,m_1}^{\ast}(\nbigv,\Phi^{\ast})
=\bigoplus(\nbigv^{\langle m_1\rangle}_{\omega},\Phi^{\ast})$
be a slope decomposition of
$\sfp_{m,m_1}^{\ast}(\nbigv,\Phi^{\ast})$.
Then, the following holds.
\begin{itemize}
\item
$(\nbigv^{\langle m_1\rangle}_{\omega},\Phi^{\ast})$
is $\Gal(m_1,m)$-equivariant.
In particular, we obtain a decomposition
$(\nbigv,\Phi^{\ast})=
 \bigoplus_{\omega\in\rnum}
 (\nbigv_{\omega},\Phi^{\ast})$
as the descent.
\item
The decomposition
is a slope decomposition of
$(\nbigv,\Phi^{\ast})$.
\hfill\qed
\end{itemize}
\end{cor}

\begin{prop}
\label{prop;18.8.19.10}
Any $\gminiq_m$-difference $\nbigk_m$-module
has a slope decomposition.
\end{prop}

If $\gminiq_m$ is not a root of $1$,
Proposition \ref{prop;18.8.19.10}
is classically well known.
(See \cite{van-der-Put-Singer, Soibelman-Vologodsky, Sauloy2004}.)
We give an outline of the proof
only in the case $\gminiq_m$ is a root of $1$.

\begin{notation}
For any $(\nbigv,\Phi^{\ast})\in\Diff_m(\nbigk,\gminiq)$,
let $\Slope(\nbigv)$ denote the set of
$\omega\in\rnum$
such that $\nbigv_{\omega}\neq 0$.
\hfill\qed
\end{notation}

\subsubsection{Proof of Proposition \ref{prop;18.8.19.10}
 in the case where $\gminiq$ is a root of $1$}
\label{subsection;18.12.12.2}
\paragraph{Cyclic vectors}

Let $(\nbigv,\Phi^{\ast})\in\Diff_m(\nbigk,\gminiq)$.
For any $v\in\nbigv$,
we set 
\[
\langle\!\langle v\rangle\!\rangle:=
\sum_{j\in\seisuu} \nbigk_m\cdot(\Phi^{\ast})^j(v),
\quad\quad
\langle v\rangle:=
\sum_{j\geq 0} \nbigk_m\cdot (\Phi^{\ast})^j(v).
\]
Note that
$\langle v\rangle
=\langle\!\langle v\rangle\!\rangle$
holds.
Indeed, we clearly have 
$\Phi^{\ast}
 \bigl(\langle v\rangle\bigr)
\subset
 \langle v\rangle$.
Because 
$\dim_{\nbigk_p}\Phi^{\ast}\bigl(
 \langle v\rangle
 \bigr)
=\dim_{\nbigk_p}\langle v\rangle$,
we obtain 
$\Phi^{\ast}
 \bigl(\langle v\rangle\bigr)
=\langle v\rangle$.

\vspace{.1in}

An element $v\in\nbigv$ is called a cyclic vector
if $\langle\!\langle v\rangle\!\rangle=\nbigv$.
The following lemma is standard.

\begin{lem}
If $\nbigv$ has a cyclic vector $v$,
there exist $m_1\in m\seisuu_{>0}$,
$\ell\in\seisuu$
and a decomposition
$\sfp_{m,m_1}^{\ast}(\nbigv,\Phi^{\ast})
=(\nbigv_1,\Phi^{\ast})\oplus(\nbigv_2,\Phi^{\ast})$
such that 
(i) $(\nbigv_1,\Phi^{\ast})$
is pure isoclinic of slope $\ell/m_1$,
(ii) $(\nbigv_1,\Phi^{\ast})\neq 0$.
\end{lem}
\pf
We give only an indication.
Set $r:=\dim_{\nbigk_m}\nbigv$.
It is easy to see that
$v,\Phi^{\ast}(v),\ldots,
 (\Phi^{\ast})^{r-1}(v)$
induce a frame of $\nbigv$ over $\nbigk_m$.
There exists a relation
$(\Phi^{\ast})^r(v)=
\sum_{j=0}^{r-1} a_j\cdot (\Phi^{\ast})^j(v)$,
where $a_j\in\nbigk_m$.
Note that one of $a_j$ is not $0$.
We set
\[
 \ell/s:=\max\Bigl\{
 -\frac{\ord_{y_m}(a_j)}{r-j}\,\Big|\,
 j=0,\ldots,r-1
 \Bigr\},
\]
where $(\ell,s)\in\seisuu\times\seisuu_{>0}$
such that $\gcd(\ell,s)=1$.
Note that $\ord_{y_m}(0)=\infty$.
We set $m_1:=sm$.
Because
$(y_{m_1}^{\ell}\Phi^{\ast})^{j}(v)
=y_{m_1}^{\ell j}\gminiq_{m_1}^{\ell j(j-1)/2}
 (\Phi^{\ast})^j(v)$,
we obtain the following:
\[
 (y_{m_1}^{\ell}\Phi^{\ast})^r(v)
=\sum_{j=0}^{r-1}
 y_{m_1}^{\ell(r-j)}a_j
 \gminiq_{m_1}^{(\ell r(r-1)-\ell j(j-1))/2}
 (y_{m_1}^{\ell}\Phi^{\ast})^j(v).
\]
Note that
$b_j:=y_{m_1}^{\ell(r-j)}a_j\in\nbigr_{m_1}$,
and 
there exists $j_0$ such that
$b_{j_0}(0)\neq 0$.

Let $\nbigl\subset\nbigv\otimes\nbigk_{m_1}$ 
be the lattice generated by
$(y_{m_1}^{\ell}\Phi^{\ast})^j(v)$ $(j\in\seisuu)$.
Clearly,
$y_{m_1}^{\ell}\Phi^{\ast}(\nbigl)\subset\nbigl$ holds.
Moreover, 
the induced endomorphism 
$F:=\sigma(y_{m_1}^{\ell}\Phi^{\ast};\nbigl)$
of $\nbigl_{|0}$ is not nilpotent.
There exists the decomposition
$\nbigl_{|0}=L_1\oplus L_2$
such that 
(i) $F(L_i)\subset L_i$, 
(ii) $F_{|L_1}$ is invertible,
(iii) $F_{|L_2}$ is nilpotent,
(iv) $L_1\neq 0$.
There exists the decomposition
$\nbigl=\nbigl_1\oplus\nbigl_2$
such that
(i) $y_{m_1}^{\ell}\Phi^{\ast}(\nbigl_i)\subset\nbigl_i$,
(ii) $\nbigl_{i|0}=L_i$.
It induces a decomposition
$\nbigv=\nbigv_1\oplus\nbigv_2$
with the desired property.
\hfill\qed

\begin{rem}
If $\gminiq_m$ is not a root of $1$,
it is classically known that
any $(\nbigv,\Phi^{\ast})\in\Diff_m(\nbigk,\gminiq)$
has a cyclic vector.
(See {\rm\cite{DiVizio2002, Sauloy2000}}.)
If $\gminiq_m$ is a root of $1$,
a $\gminiq_m$-difference $\nbigk_m$-module
does not necessarily have a cyclic vector.
\hfill\qed
\end{rem}

\paragraph{Eigen decompositions}

As a preliminary to prove Proposition \ref{prop;18.8.19.10},
we recall a standard result in linear algebra.
Let $f$ be a $\nbigk_m$-automorphism of
a $\nbigk_m$-vector space $\nbigu$.
Set $r!m$.
Recall that the set of the eigenvalues
$\Sp(f)$ is contained in $\nbigk_{r!m}$.
We obtain the decomposition
\[
 \nbigu\otimes_{\nbigk_m}\nbigk_{r!m}
 =\bigoplus_{\gminia\in\Sp(f)}\nbigu^{\langle r!m\rangle}_{\gminia}.
\]
For any $\omega\in \frac{1}{r!m}\seisuu$,
we put
$\Sp(f,\omega):=
 \bigl\{
 \gminia\in\Sp(f)\,\big|\,
 \ord_{y_{r!m}}(\gminia)=r!m\omega
 \bigr\}$.
We set
\[
 \nbigu^{\langle r!m\rangle}_{\omega}:=
 \bigoplus_{\gminia\in\Sp(f,\omega)}
 \nbigu^{\langle r!m\rangle}_{\gminia}.
\]
Let $G$ denote the Galois group of $\nbigk_{r!m}$ over $\nbigk_m$.
There is a natural $G$-action on $\nbigv\otimes_{\nbigk_m}\nbigk_{r!m}$.
Because $\nbigu^{\langle r!m\rangle}_{\omega}$
is $G$-invariant,
we have a subspace $\nbigu_{\omega}\subset\nbigu$
such that 
$\nbigu^{\langle r!m\rangle}_{\omega}=
 \nbigu_{\omega}\otimes_{\nbigk_m}\nbigk_{r!m}$,
and we obtain the following decomposition:
\begin{equation}
\label{eq;18.12.12.1}
 \nbigu\otimes_{\nbigk_m}\nbigk_{r!m}
 =\bigoplus_{\omega\in\rnum}
 \nbigu^{\langle r!m\rangle}_{\omega}.
\end{equation}

\paragraph{Proof of Proposition \ref{prop;18.8.19.10}
in the case $\gminiq$ is a root of $1$}

Let $(\nbigv,\Phi^{\ast})$ be a $\gminiq_m$-difference 
$\nbigk_m$-module.
Let us prove that
$(\nbigv,\Phi^{\ast})$ has a slope decomposition
in the case 
where $\gminiq_m$ is an $s$-th root of $1$
for some $s\in\seisuu_{>0}$.
We  use an induction of $\dim_{\nbigk_m}\nbigv$.
We put $s_1:=r!s$. 
We set $\Psi^{\ast}:=(\Phi^{\ast})^{s_1}$.
Note that
$\Psi^{\ast}=\id$ on $\nbigk_{m'}$
for any $m'\in\{m,2m,\ldots,r!m\}$.
Hence, 
$\Psi^{\ast}$ on $\nbigv$
is $\nbigk_m$-linear,
and $\Psi^{\ast}$ on 
$\nbigv\otimes_{\nbigk_m}\nbigk_{m'}$
for $m'\in\{m,2m,\ldots,r!m\}$
are the induced $\nbigk_{m'}$-linear
automorphisms.
We obtain the decomposition
$\nbigv
=\bigoplus
 \nbigv_{\omega}$
as in (\ref{eq;18.12.12.1}).
Note that 
$(\nbigv_{\omega},\Psi^{\ast})$
has pure slope $\omega$.
By using the commutativity of
$\Phi^{\ast}$ and $\Psi^{\ast}$,
and by the construction of (\ref{eq;18.12.12.1}),
we obtain that 
$\Phi^{\ast}(\nbigv_{\omega})
=\nbigv_{\omega}$.
Let us prove that
$(\nbigv_{\omega},\Phi^{\ast})$
has pure slope $\omega/s_1$.

Suppose that
$\nbigv_{\omega}$ does not have a cyclic vector.
Take any $v\in \nbigv_{\omega}$.
Note that
$(\langle v\rangle,\Phi^{\ast})
 \subsetneq
 (\nbigv_{\omega},\Phi^{\ast})$.
Then, by the assumption of the induction,
we may assume that there exists
a decomposition
$\langle v\rangle
=\bigoplus_{\mu\in\rnum}
\langle v\rangle_{\mu}$,
where
$(\langle v\rangle_{\mu},\Phi^{\ast})$
has pure slope $\mu$.
Because 
$(\langle v\rangle_{\mu},\Psi^{\ast})$
has pure slope $\mu s_1$,
we obtain that 
$\langle v\rangle_{\mu}=0$
unless $s_1\mu=\omega$.
Hence, we obtain that 
$(\langle v\rangle,\Phi^{\ast})$ has pure slope
$\omega/s_1$.
By varying $v$, we obtain that 
$(\nbigv_{\omega},\Phi^{\ast})$
has pure slope $\omega/s_1$.

Suppose that
$\nbigv_{\omega}$ has a cyclic vector.
Then, there exist $m_1\in m\seisuu_{>0}$
and a decomposition
$\nbigv_{\omega}\otimes_{\nbigk_m}\nbigk_{m_1}
=\nbigv^{(1)}\oplus\nbigv^{(2)}$
such that
(i) $\Phi^{\ast}(\nbigv^{(i)})=\nbigv^{(i)}$,
(ii) $\nbigv^{(1)}\neq 0$,
(iii) $\nbigv^{(1)}$ has pure slope.
By using the hypothesis of the induction,
we may assume that
$\nbigv^{(2)}$ has a slope decomposition
with respect to $\Phi^{\ast}$.
Hence, $\nbigv\otimes_{\nbigk_m}\nbigk_{m_1}$
has a slope decomposition.
As in the previous paragraph,
we obtain that 
the slope of $(\nbigv^{(1)},\Phi^{\ast})$ is $\omega/s_1$,
and that 
$(\nbigv\otimes_{\nbigk_m}\nbigk_{m_1},\Phi^{\ast})$
has pure slope $\omega/s_1$.
Hence, we can conclude that
$(\nbigv_{\omega},\Phi^{\ast})$
has pure slope $\omega/s_1$.
\hfill\qed

\subsection{Filtered formal bundles}

We recall the notion of filtered bundles
on $\nbigk_m$-vector spaces.
Let $\nbigv$ be a finite dimensional vector space over 
$\nbigk_m$.
A filtered bundle over $\nbigv$
is an increasing sequence 
 $\nbigp_{\ast}\nbigv
=\bigl(\nbigp_a\nbigv\,|\,a\in\real\bigr)$
of $\nbigr_m$-lattices of
$\nbigv$ such that
(i) $\nbigp_a(\nbigv)=\bigcap_{a<b}\nbigp_b(\nbigv)$
 for any $a\in\real$,
(ii) $\nbigp_{a+n}(\nbigv)=y_m^{-n}\nbigp_a(\nbigv)$
 for any $a\in\real$ and $n\in\seisuu$.
We set
$\Gr^{\nbigp}_a(\nbigv):=
 \nbigp_{a}(\nbigv)\big/
 \nbigp_{<a}(\nbigv)$.
A morphism of filtered bundles
$F:\nbigp_{\ast}\nbigv_1\lrarr\nbigp_{\ast}\nbigv_2$
is defined to be a $\nbigk_m$-homomorphism $F$
satisfying
$F(\nbigp_a\nbigv_1)\subset\nbigp_a\nbigv_2$.
Let $\Mod(\nbigk_m)^{\Par}$ denote the category of
filtered bundles over finite dimensional $\nbigk_m$-vector spaces.

\subsubsection{Pull back}
\label{subsection;19.2.5.1}

Let $m_1\in m\seisuu_{>0}$.
Let $\nbigp_{\ast}\nbigv_1\in\Mod(\nbigk_m)^{\Par}$.
Recall that we obtain the induced filtered bundle
$\nbigp_{\ast}(\sfp_{m,m_1}^{\ast}\nbigv)$
over $\sfp_{m,m_1}^{\ast}\nbigv$
given as follows:
\begin{equation}
\label{eq;18.12.18.2}
 \nbigp_a\bigl(
 \sfp_{m,m_1}^{\ast}\nbigv
 \bigr)
=\sum_{(b,n)\in S(m_1,m)}
 y_{m_1}^{-n}\nbigp_b(\nbigv)
 \otimes_{\nbigr_m}\nbigr_{m_1},
\end{equation}
where
$S(m_1,m):=\bigl\{
 (b,n)\in\real\times\seisuu\,\big|\,
 \frac{m_1}{m}b+n\leq a
 \bigr\}$.
The filtered bundle
$\nbigp_{\ast}(\sfp_{m,m_1}^{\ast}\nbigv)$
is also denoted by
$\sfp_{m,m_1}^{\ast}\nbigp_{\ast}\nbigv$.
Thus, we obtain the pull back functor
\[
\sfp_{m,m_1}^{\ast}:
 \Mod(\nbigk_m)^{\Par}\lrarr
 \Mod(\nbigk_{m_1})^{\Par}.
\]

Let $\nbigp_{\ast}\nbigv\in\Mod(\nbigk_m)^{\Par}$.
We obtain the map
$\Gr^{\nbigp}_a(\nbigv)
\lrarr
\Gr^{\nbigp}_{(m_1a/m)+n}(
 \sfp_{m,m_1}^{\ast}\nbigv)$
as follows.
For $s\in\Gr^{\nbigp}_a(\nbigv)$,
take a lift $\stilde\in\nbigp_a(\nbigv)$ of $s$,
then we obtain the element 
in $\Gr^{\nbigp}_{(m_1a/m)+n}(\sfp_{m,m_1}^{\ast}\nbigv)$
induced by $y_{m_1}^{-n}\stilde$,
which is independent of a choice of $\stilde$.
This procedure induces an isomorphism
\begin{equation}
\label{eq;18.12.18.1}
\bigoplus_{(n,b)\in S_0(m_1,m,a)}
 \Gr^{\nbigp}_b(\nbigv)
\simeq
 \Gr^{\nbigp}_a\bigl(\sfp_{m,m_1}^{\ast}\nbigv\bigr)
\end{equation}
where
$S_0(m_1,m,a):=\bigl\{
 (n,b)\in\seisuu\times\real\,\big|\,
 0\leq n<\frac{m_1}{m},\,
 \frac{m_1}{m}b+n=a
 \bigr\}$.
Each $\nbigp_{a}(\sfp_{m,m_1}^{\ast}\nbigv)$
is preserved by the $\Gal(m_1,m)$-action,
and hence we obtain the $\Gal(m_1,m)$-action on
$\Gr^{\nbigp}_a(\sfp_{m,m_1}^{\ast}\nbigv)$.
The decomposition (\ref{eq;18.12.18.1})
is identified with the canonical decomposition
with respect to the $\Gal(m_1,m)$-action,
and 
$\Gr^{\nbigp}_b(\nbigv)$ 
is identified with the $\Gal(m_1,m)$-invariant part.

\subsubsection{Push-forward and descent}

Let $m$ and $m_1$ be positive integers such that
$m_1\in m\seisuu_{>0}$.
Let $\nbigp_{\ast}\nbigv\in\Mod(\nbigk_{m_1})^{\Par}$.
Recall that the filtered bundle
$\nbigp_{\ast}(\sfp_{m,m_1\ast}\nbigv)$
is induced as follows:
\[
 \nbigp_a(\sfp_{m,m_1\ast}\nbigv)
=\nbigp_{a(m_1/m)}\nbigv.
\]
The filtered bundle
$\nbigp_{\ast}(\sfp_{m,m_1\ast}\nbigv)$
is also denoted by
$\sfp_{m,m_1\ast}\nbigp_{\ast}\nbigv$.
Thus, we obtain
the push-forward
$\sfp_{m,m_1\ast}:
 \Mod(\nbigk_{m_1})^{\Par}
\lrarr
 \Mod(\nbigk_m)^{\Par}$.

Let $\nbigp_{\ast}(\nbigv)\in\Mod(\nbigk_{m_1})^{\Par}$.
There exists the natural isomorphism
\[
 \Gr^{\nbigp}_a(\sfp_{m,m_1\ast}\nbigv)
\simeq
 \Gr^{\nbigp}_{a(m_1/m)}(\nbigv).
\]

Let $\nbigv$ be a finite dimensional $\Gal(m_1,m)$-equivariant 
$\nbigk_{m_1}$-vector space.
We say that a filtered bundle
 $\nbigp_{\ast}\nbigv$ over $\nbigv$ is 
$\Gal(m_1,m)$-equivariant
if $\mu\nbigp_a(\nbigv)=\nbigp_a\nbigv$
for any $a\in\real$
and $\mu\in\Gal(m_1,m)$.
We obtain the $\nbigk_m$-vector space 
$\nbigv^{\Gal(m_1,m)}$
as the descent,
i..e., as the $\Gal(m_1,m)$-invariant part of $\nbigv$.
We have the induced filtered bundle
$\nbigp_{\ast}(\nbigv^{\Gal(m_1,m)})$
over $\nbigv^{\Gal(m_1,m)}$
as 
\[
 \nbigp_{a}(\nbigv^{\Gal(m_1,m)}):=
 \bigl(
 \nbigp_{a(m_1/m)}\nbigv
 \bigr)^{\Gal(m_1,m)}.
\]
The filtered bundle is denoted by
$\nbigp_{\ast}(\nbigv)^{\Gal(m_1,m)}$,
and called the decent of $\nbigp_{\ast}\nbigv$.

\begin{lem}
\mbox{{}}
\begin{itemize}
\item
$\nbigp_{\ast}\nbigv\in\Mod(\nbigk_m)^{\Par}$
is naturally isomorphic to
$(\sfp_{m,m_1}^{\ast}\nbigp_{\ast}\nbigv)^{\Gal(m_1,m)}$.
\item
Let $\nbigp_{\ast}\nbigv_1\in\Mod(\nbigk_{m_1})^{\Par}$.
Then,
$\sfp_{m_1,m}^{\ast}
 \sfp_{m,m_1\ast}
 \nbigp_{\ast}\nbigv_1$
is naturally isomorphic to 
$\bigoplus_{\mu\in\Gal(m_1,m)}
 \mu^{\ast}\nbigp_{\ast}(\nbigv_1)$.

\hfill\qed
\end{itemize}
\end{lem}

\subsubsection{Reduction}

For any $\nbigp_{\ast}\nbigv\in\Mod(\nbigk_m)^{\Par}$,
we set
\[
 \ttG(\nbigp_{\ast}\nbigv):=\bigoplus_{a\in\real}
 \Gr^{\nbigp}_a(\nbigv).
\]
The multiplication of $y_m$
induces $\cnum$-linear isomorphisms
$\Gr^{\nbigp}_a(\nbigv)
\lrarr
 \Gr^{\nbigp}_{a-1}(\nbigv)$ for any $a\in\real$.
Hence, we may naturally regard
$\ttG(\nbigp_{\ast}\nbigv)$ as 
a free $\cnum[y_m,y_m^{-1}]$-module
with an $\real$-grading.
It is also $\real$-graded.
For any $a\in\real$,
we set
\[
 \nbigp_a\ttG(\nbigp_{\ast}\nbigv):=
 \bigoplus_{b\leq a}
 \Gr^{\nbigp}_b(\nbigv).
\]
It is a $\cnum[y_m]$-lattice 
of $\ttG(\nbigp_{\ast}\nbigv)$.
By the construction,
there exists a natural isomorphism
\[
 \Gr^{\nbigp}_a\ttG(\nbigp_{\ast}\nbigv):=
 \nbigp_a\ttG(\nbigp_{\ast}\nbigv)\big/
 \nbigp_{<a}\ttG(\nbigp_{\ast}\nbigv)
\simeq
  \Gr^{\nbigp}_a(\nbigv).
\]

We set
$\ttGhat(\nbigp_{\ast}\nbigv):=
 \ttG(\nbigp\nbigv)\otimes_{\cnum[y_m,y_m^{-1}]} \nbigk_m$.
We also set
$\nbigp_a\ttGhat(\nbigp_{\ast}\nbigv)
:=\nbigp_a\ttG(\nbigp_{\ast}\nbigv)
 \otimes_{\cnum[y_m]}\nbigr_m$
for any $a\in\real$.
They give a filtered bundle
over $\ttGhat(\nbigp_{\ast}\nbigv)$.
In this way,
we also regard $\ttGhat(\nbigp_{\ast}\nbigv)$
as a filtered bundle.
For any $a\in\real$,
there exist the natural isomorphisms:
\begin{equation}
\label{eq;19.1.15.1}
 \Gr^{\nbigp}_a\ttGhat(\nbigp_{\ast}\nbigv)
\simeq
 \Gr^{\nbigp}_a\ttG(\nbigp_{\ast}\nbigv)
\simeq
 \Gr^{\nbigp}_a(\nbigv).
\end{equation}
\begin{rem}
There exist a (non-unique) isomorphism
of filtered bundles
$\ttGhat(\nbigp_{\ast}\nbigv)
\simeq
 \nbigp_{\ast}\nbigv$
which induces the isomorphisms
{\rm(\ref{eq;19.1.15.1})}.
\hfill\qed
\end{rem}

\subsection{Graded $\gminiq$-difference $\cnum[y,y^{-1}]$-modules}

A $\gminiq_m$-difference
free $\cnum[y_m,y_m^{-1}]$-module $(M,\Phi^{\ast})$
is a free $\cnum[y_m,y_m^{-1}]$-module $M$ of finite rank
equipped with a $\cnum$-linear automorphism $\Phi^{\ast}$
such that
$\Phi^{\ast}(y_ms)=\gminiq_my_m\Phi^{\ast}(s)$
for any $s\in M$.
A $(\rnum,\real)$-grading is
a decomposition
\[
 M=\bigoplus_{(\omega,a)\in\real\times\rnum}
 M_{\omega,a}
\]
such that the following holds:
\begin{itemize}
\item
 $y_m M_{\omega,a}=M_{\omega,a-1}$
 for any $(\omega,a)\in\rnum\times\real$.
\item
 $\Phi^{\ast}M_{\omega,a}=M_{\omega,a+m\omega}$.
\end{itemize}
A morphism of 
$\gminiq_m$-difference
free $\cnum[y_m,y_m^{-1}]$-modules
with $(\rnum,\real)$-grading
$(M^{(1)}_{\bullet,\bullet},\Phi^{\ast})
\lrarr
 (M^{(2)}_{\bullet,\bullet},\Phi^{\ast})$
is defined to be a morphism 
of $\gminiq_m$-difference $\cnum[y_m,y_m^{-1}]$-modules
preserving the gradings.
Let $\Diff_m(\cnum[y,y^{-1}],\gminiq)_{(\rnum,\real)}$
denote the category of 
$\gminiq_m$-difference
free $\cnum[y_m,y_m^{-1}]$-modules
equipped with $(\rnum,\real)$-grading.

\vspace{.1in}

For each $\omega$,
we have the expression
$m\omega=\ell(m\omega)/k(m\omega)$,
where 
$k(m\omega)$ and $\ell(m\omega)$ are uniquely determined by
the conditions
$k(m\omega)\in\seisuu_{>0}$,
$\ell(m\omega)\in \seisuu$
and $\gcd(k(m\omega),\ell(m\omega))=1$.
Let $\Lambda_{\omega}:=\frac{1}{k(m\omega)}\seisuu$.
Note that $\Lambda_{\omega}$
is the image of
the map
$\seisuu^2\lrarr \real$
defined by
$(n_1,n_2)\longmapsto
 n_1m\omega-n_2$.

We obtain the automorphism
$F_{\omega,a}:=y_m^{\ell(m\omega)}(\Phi^{\ast})^{k(m\omega)}$
on $M_{\omega,a}$
for any $(a,\omega)\in\real\times\rnum$.
We obtain the generalized eigen decomposition
\[
 (M_{\omega,a},F_{\omega,a})
=\bigoplus_{\alpha\in\cnum^{\ast}} 
 (M_{\omega,a,\alpha},F_{\omega,a,\alpha})
\]
where 
$F_{\omega,a,\alpha}$
has a unique eigenvalue $\alpha$.
It is easy to see
\[
 y_m \cdot M_{\omega,a,\alpha}
=M_{\omega,a-1,\alpha \gminiq_m^{k(m\omega)}},
\quad
 \Phi^{\ast}\cdot M_{\omega,a,\alpha}
=M_{\omega,a+m\omega,\alpha \gminiq_m^{-\ell(m\omega)}}.
\]
For 
$\omega\in\rnum$,
$-1/k(m\omega)<a\leq 0$
and $\alpha\in\cnum^{\ast}$,
we set
\[
 M(\omega,a,\alpha):=
 \bigoplus_{b\in\Lambda_1}
 M_{\omega,a+b,\alpha\gminiq^{-k(m\omega)b}}.
\]
Then, we obtain a decomposition
of $(\rnum,\real)$-graded
$\gminiq_m$-difference
$\cnum[y_m,y_m^{-1}]$-modules:
\[
 M=
 \bigoplus_{\omega\in\rnum}
 \bigoplus_{-k(m\omega)^{-1}<a\leq 0}
 \bigoplus_{\alpha\in\cnum^{\ast}}
 M(\omega,a,\alpha).
\]

\subsubsection{The induced nilpotent endomorphism 
and the weight filtration}

Each $M_{\omega,a,\alpha}$
is equipped with the nilpotent endomorphism $N_{\omega,a,\alpha}$
obtained as the logarithm of the unipotent part of
$F_{\omega,a,\alpha}$.
It induces the weight filtration
$W(M_{\omega,a,\alpha})$.
We obtain the nilpotent endomorphism 
$N=\bigoplus N_{\omega,a,\alpha}$ of $M$.
It commutes with $y_m$ and $\Phi^{\ast}$.
Hence, $N$ is a nilpotent endomorphism
of $M$ in $\Diff_m(\cnum[y,y^{-1}])_{(\rnum,\real)}$.
The weight filtration $W$ is a filtration of $M$
in $\Diff_m(\cnum[y,y^{-1}])_{(\rnum,\real)}$.

\subsubsection{Classification}

Let $\nbigc$ denote the category of
finite dimensional vector spaces $V$
equipped with a grading
\[
 V=
 \bigoplus_{\omega\in\rnum}
\bigoplus_{-k(m\omega)^{-1}<a\leq 0}
 \bigoplus_{\alpha\in\cnum^{\ast}}
 V_{\omega,a,\alpha}
\]
and a graded unipotent automorphism
$u=
 \bigoplus_{\omega\in\rnum}
 \bigoplus_{-k(m\omega)^{-1}<a\leq 0}
 \bigoplus_{\alpha\in\cnum^{\ast}}
 u_{\omega,a,\alpha}$.
A morphism
\[
 F:(V^{(1)}_{\bullet},u^{(1)}_{\bullet})
\lrarr (V^{(2)}_{\bullet},u^{(2)}_{\bullet})
\]
in $\nbigc$
is a $\cnum$-linear map
$F:V^{(1)}\lrarr V^{(2)}$
such that
(i) $F$ preserves the gradings,
(ii) $F\circ u_{\bullet}^{(1)}=u_{\bullet}^{(2)}\circ F$.

For any $M\in \Diff_m(\cnum[y,y^{-1}],\gminiq)_{(\rnum,\real)}$,
we obtain
the finite dimensional graded vector space
\[
 \bigoplus_{\omega\in\rnum}
 \bigoplus_{-k(m\omega)^{-1}<a\leq 0}
 \bigoplus_{\alpha\in\cnum^{\ast}}
 M_{\omega,a,\alpha}.
\]
Let $u_{\omega,a,\alpha}$ denote the unipotent part of
$F_{\omega,a,\alpha}$.
We obtain an object
$\bigoplus_{\omega,a,\alpha}
 (M_{\omega,a,\alpha},u_{\omega,a,\alpha})$
in $\nbigc$.
Thus, we obtain a functor
$\Diff_m(\cnum[y,y^{-1}],\gminiq)_{(\rnum,\real)}
\lrarr\nbigc$.
The following is easy to see.
\begin{lem}
The functor
$\Diff_m(\cnum[y,y^{-1}],\gminiq)_{(\rnum,\real)}
\lrarr \nbigc$
is an equivalence.
\hfill\qed
\end{lem}

\begin{rem}
Let $V_{\bullet}\in\nbigc$.
For each $(\omega,a,\alpha)$,
we obtain the nilpotent endomorphism
$N_{\omega,a,\alpha}$ of $V_{\omega,a,\alpha}$
as the logarithm of the unipotent automorphism
$u_{\omega,a,\alpha}$.
We obtain the weight filtration $W(V_{\omega,a,\alpha})$
with respect to $N_{\omega,a,\alpha}$.
Note that the conjugacy classes of $u_{\omega,a,\alpha}$
are determined by 
the filtrations $W(V_{\omega,a,\alpha})$.
\hfill\qed
\end{rem}

\subsubsection{Tensor product}

Let $M^{(i)}\in \Diff_m(\cnum[y,y^{-1}],\gminiq)_{(\rnum,\real)}$.
We obtain the $\gminiq_m$-difference
free $\cnum[y_m,y_m^{-1}]$-module
$M^{(1)}\otimes M^{(2)}$
by the tensor product over $\cnum[y_m,y_m^{-1}]$.
Let
$(M^{(1)}\otimes M^{(2)})_{\omega,a}$
be the image of the injective map
\begin{equation}
 \bigoplus_{\omega_1+\omega_2=\omega}
 \bigoplus_{\substack{a_1+a_2=a \\-1<a_1\leq 0}}
 M^{(1)}_{\omega_1,a_1}
\otimes_{\cnum}
 M^{(2)}_{\omega_2,a_2}
\lrarr
 M^{(1)}\otimes_{\cnum[y_m,y_m^{-1}]} M^{(2)}.
\end{equation}
Then, we obtain the grading
$M^{(1)}\otimes M^{(2)}
=\bigoplus (M^{(1)}\otimes M^{(2)})_{\omega,a}$.
We have the automorphisms
$F^{(i)}_{\omega,a}$ of $M^{(i)}_{\omega,a}$.
We also have the automorphism
$F_{\omega,a}$ of $(M^{(1)}\otimes M^{(2)})_{\omega,a}$.

\begin{lem}
Suppose that $M^{(1)}_{\omega,a}=0$ unless $\omega=0$.
Under the identification
\[
 (M^{(1)}\otimes M^{(2)})_{\omega,a}
=
 \bigoplus_{\substack{a_1+a_2=a \\-1<a_1\leq 0}}
 M^{(1)}_{0,a_1}
\otimes_{\cnum}
 M^{(2)}_{\omega,a_2},
\]
we have
$F_{\omega,a}
=\bigoplus
 (F^{(1)}_{0,a_1})^{k(m\omega)}
\otimes
 F^{(2)}_{\omega,a_2}$.
The nilpotent endomorphism
$N$ of $M^{(1)}\otimes M^{(2)}$
is equal to
$k(m\omega)N^{(1)}\otimes \id+\id\otimes N^{(2)}$,
where $N^{(i)}$ are the nilpotent endomorphism
of $M^{(i)}$.
The filtration
$W\bigl(
 (M^{(1)}\otimes M^{(2)})
 \bigr)$
is equal to the filtration induced by $W(M^{(i)})$
$(i=1,2)$.
\hfill\qed
\end{lem}

\subsubsection{Pull back}

Let $m_1\in m\seisuu_{>0}$.
Let $M\in\Diff_m(\cnum[y,y^{-1}],\gminiq)_{(\rnum,\real)}$.
We set
$\sfp_{m,m_1}^{\ast}(M):=
 M\otimes_{\cnum[y_m,y_m^{-1}]} \cnum[y_{m_1},y_{m_1}^{-1}]$
which is naturally a $\gminiq_{m_1}$-difference
$\cnum[y_{m_1},y_{m_1}^{-1}]$-module.
Set 
$S_0(m_1,m,a):=\bigl\{
 (n,b)\in\seisuu\times\real\,\big|\,
 0\leq n<\frac{m_1}{m},\,
 \frac{m_1}{m}b+n=a
 \bigr\}$
as in \S\ref{subsection;19.2.5.1}.
We define
$\sfp_{m,m_1}^{\ast}(M)_{\omega,a}$
as the image of the injection:
\begin{equation}
\label{eq;19.1.17.1}
\bigoplus_{(b,i)\in S_0(m,m_1,a)}
 y_{m_1}^{-i}
 M_{\omega,b}
\lrarr
 \sfp_{m,m_1}^{\ast}(M).
\end{equation}
Then, we obtain the grading
$\sfp_{m,m_1}^{\ast}(M)
=\bigoplus \sfp_{m,m_1}^{\ast}(M)_{\omega,a}$.
Thus, we obtain
\[
 \sfp_{m,m_1}^{\ast}:
 \Diff_m(\cnum[y,y^{-1}],\gminiq)_{(\rnum,\real)}
\lrarr
 \Diff_{m_1}(\cnum[y,y^{-1}],\gminiq)_{(\rnum,\real)}.
\]

Let $F^{(1)}_{\omega,a}$ be the automorphism
of $\sfp_{m,m_1}^{\ast}(M)_{\omega,a}$
induced by
$y_{m_1}^{\ell(m_1\omega)}(\Phi^{\ast})^{k(m_1\omega)}$.
Set
$d:=k(m\omega)/k(m_1\omega)\in\seisuu_{>0}$.
\begin{lem}
Under the identification of 
$\sfp_{m,m_1}^{\ast}(M)_{\omega,a}
\simeq
\bigoplus_{(b,i)\in S_0(m,m_1,a)}
 y_{m_1}^{-i}
 M_{\omega,b}$,
we have
\[
 \bigl(
 F^{(1)}_{\omega,a}
 \bigr)^{d}
=\bigoplus_{(b,i)\in S_0(m,m_1,a)}
 \gminiq_{m_1}^{\frac{1}{2}\ell(m_1\omega)k(m_1\omega)
 d(d-1)-ik(m_1\omega)d}
 F_{\omega,b}.
\]
Hence,
$dN^{(1)}=\sfp_{m_1,m}^{\ast}N$ holds,
where $N^{(1)}$ and $N$ are the nilpotent endomorphisms
of $\sfp_{m,m_1}^{\ast}(M)$ and $M$,
respectively.
We also obtain
$W(\sfp_{m,m_1}^{\ast}M)
=\sfp_{m,m_1}^{\ast}W(M)$.
\hfill\qed
\end{lem}

\subsubsection{Push-forward}

Let $M\in \Diff_{m_1}(\cnum[y,y^{-1}],\gminiq)_{(\rnum,\real)}$.
It naturally induces a $\gminiq_m$-difference
$\cnum[y_m,y_m^{-1}]$-module
$\sfp_{m,m_1\ast}(M)$.
We set 
$\sfp_{m,m_1\ast}(M)_{\omega,a}:= M_{\omega,am/m_1}$.
Thus, we obtain
\[
 \sfp_{m,m_1\ast}:
 \Diff_{m_1}(\cnum[y,y^{-1}],\gminiq)_{(\rnum,\real)}
\lrarr
 \Diff_{m}(\cnum[y,y^{-1}],\gminiq)_{(\rnum,\real)}.
\]

For any $(\omega,b)\in\rnum\times\real$,
let $F^{(1)}_{\omega,b}$ be the automorphism
of $\sfp_{m,m_1\ast}(M)_{\omega,b}$
induced by
$y_{m}^{\ell(m\omega)}(\Phi^{\ast})^{k(m\omega)}$.
Set
$d:=k(m\omega)/k(m_1\omega)\in\seisuu_{>0}$.
\begin{lem}
We have
\[
(F_{\omega,a})^d
=\gminiq_{m_1}^{\frac{1}{2}\ell(m_1\omega)k(m_1\omega)d(d-1)}
F^{(1)}_{\omega,am/m_1}.
\]
As a result,
$\sfp_{m,m_1\ast}(dN)=N^{(1)}$
and 
$\sfp_{m,m_1\ast}W(M)
=W(\sfp_{m,m_1\ast}M)$
hold.
\hfill\qed
\end{lem}

\subsubsection{Examples}

Let $\omega\in\rnum$
and 
 $-k(m\omega)^{-1}<a\leq 0$.
Let $\LL^{\ttG}_m(\omega,a)\in
\Diff_m(\cnum[y,y^{-1}],\gminiq)_{(\real,\rnum)}$ 
be the object corresponding to
\[
 \bigl(M_{\omega',a'},F_{\omega',a'}\bigr)
=\left\{
 \begin{array}{ll}
 \bigl(\cnum,
 \gminiq_{m}^{-\frac{1}{2}\ell(m\omega)k(m\omega)(k(m\omega)-1)}
 \bigr) & 
 \bigl((\omega',a')=(\omega,a)
 \bigr)\\
\mbox{{}}\\
 (0,1) & (\mbox{\rm otherwise}).
 \end{array}
 \right.
\]
There exists a natural isomorphism
$\sfp_{m,mk(m\omega)\ast}
 \LL^{\ttG}_{mk(\omega)}(\omega,ak(m\omega))
\simeq
 \LL^{\ttG}_{m}(\omega,a)$.

For a finite dimensional vector space $V$
with an automorphism $F$,
let $\VV^{\ttG}_m(V,F)\in \Diff_m(\cnum[y,y^{-1}],\gminiq)_{(\real,\rnum)}$
be the object corresponding to
\[
  \bigl(M_{\omega',a'},F_{\omega',a'}\bigr)
=\left\{
 \begin{array}{ll}
 \bigl(V,F\bigr) & 
 \bigl((\omega',a')=(0,0)
 \bigr)\\
 (0,1) & (\mbox{\rm otherwise}).
 \end{array}
 \right.
\]

Any $M\in \Diff_m(\cnum[y,y^{-1}],\gminiq)_{(\rnum,\real)}$
is isomorphic to
the object of the following form:
\[
 \bigoplus_{i=1}^{N}
 \LL^{\ttG}_{m}(\omega_i,a_i)
 \otimes
 \VV^{\ttG}_m(V_i,F_i).
\]

\subsection{Good filtered formal $\gminiq$-difference modules}

\subsubsection{Good filtered bundles}
\label{subsection;18.8.28.1}

Let $(\nbigv,\Phi^{\ast})\in\Diff_m(\nbigk,\gminiq)$.

\begin{df}
\label{df;18.8.28.4}
A filtered bundle $\nbigp_{\ast}\nbigv$ over $\nbigv$
is called good
if the following holds.
\begin{itemize}
\item
The filtration $\nbigp_{\ast}\nbigv$
is compatible with the slope decomposition
$\nbigv=\bigoplus_{\omega\in\Slope(\nbigv)}\nbigv_{\omega}$,
i.e.,
\[
 \nbigp_{\ast}\nbigv
=\bigoplus_{\omega\in\Slope(\nbigv)}
 \nbigp_{\ast}(\nbigv_{\omega}).
\]
\item
Take $m_1\in m\seisuu_{>0}$ such that
$m_1\omega\in\seisuu$
for any $\Slope(\nbigv)$.
Then, 
the following holds
for any $\omega\in\Slope(\nbigv)$
and for any $a\in\real$:
\[
 y_{m_1}^{m_1\omega}\Phi^{\ast}\Bigl(
 \nbigp_a\bigl(
 \sfp_{m,m_1}^{\ast}
 \nbigv_{\omega}
 \bigr)
 \Bigr)
=\nbigp_a\bigl(
 \sfp_{m,m_1}^{\ast}
 \nbigv_{\omega}
 \bigr).
\]
\end{itemize}
Such $(\nbigp_{\ast}\nbigv,\Phi^{\ast})$
is called a good filtered 
$\gminiq_m$-difference $\nbigk_m$-module.
\hfill\qed
\end{df}

\begin{rem}
As a special case,
$(\nbigp_{\ast}\nbigv,\Phi^{\ast})$
is called unramifiedly good (resp. regular)
if $\Slope(\nbigv)\subset\seisuu$
(resp. $\Slope(\nbigv)=\{0\}$).
\end{rem}

A morphism of good filtered $\gminiq_m$-difference
$\nbigk_m$-modules
$F:(\nbigp_{\ast}\nbigv_1,\Phi)
\lrarr
 (\nbigp_{\ast}\nbigv_2,\Phi)$
is 
defined to be a morphism of
$\gminiq_m$-difference
$\nbigk_m$-modules $F$
such that
$F(\nbigp_a\nbigv_1)\subset\nbigp_a\nbigv_2$
for any $a\in\real$.
Let $\Diff_m(\nbigk,\gminiq)^{\Par}$
denote the category of
good filtered $\gminiq_m$-difference
$\nbigk_m$-modules.
Let $\Diff_m(\nbigk,\gminiq;\omega)^{\Par}$
denote the full subcategory of
good filtered $\gminiq_m$-difference
$\nbigk_m$-modules $(\nbigp_{\ast}\nbigv,\Phi)$
such that $(\nbigv,\Phi)\in\Diff_m(\nbigk,\gminiq;\omega)$.

\begin{lem}
Let $(\nbigv,\Phi^{\ast})\in\Diff_m(\nbigk,\gminiq;\omega)$.
Let $\nbigp_{\ast}\nbigv$ be a filtered bundle
over $\nbigv$.
Then, $(\nbigp_{\ast}\nbigv,\Phi^{\ast})$ is good
if and only if
$\Phi^{\ast}(\nbigp_a\nbigv)=\nbigp_{a+m\omega}\nbigv$
for any $a\in\real$.
\end{lem}
\pf
It is clear if $\omega\in\seisuu$.
In general, we take $m_1\in m\seisuu_{>0}$
such that $m_1\omega\in\seisuu$.
By definition,
$\nbigp_{\ast}(\nbigv)$ is good 
if and only if
$\Phi^{\ast}\nbigp_a(\sfp_{m,m_1}^{\ast}(\nbigv))
=\nbigp_{a+m_1\omega}(\sfp_{m,m_1}^{\ast}\nbigv)$.
Because
$\nbigp_b(\nbigv)
=\nbigp_{b(m_1/m)}(\sfp_{m,m_1}^{\ast}\nbigv)^{\Gal}$,
we obtain the claim of the lemma.
\hfill\qed

\subsubsection{Reduction to $(\rnum,\real)$-graded $\gminiq$-difference modules}

Let
$(\nbigp_{\ast}\nbigv,\Phi^{\ast})\in
\Diff_{m}(\nbigk,\gminiq)^{\Par}$.
There exists the slope decomposition
$\nbigp_{\ast}\nbigv=
 \bigoplus_{\omega\in\rnum}
 \nbigp_{\ast}\nbigv_{\omega}$.
We have the induced isomorphisms:
\[
 \Phi^{\ast}:
 \Gr^{\nbigp}_{a}(\nbigv_{\omega})
\simeq
 \Gr^{\nbigp}_{a+m\omega}(\nbigv_{\omega})
\]
for any $a\in\real$.
Thus, we obtain a $\cnum$-linear automorphism
$\Phi^{\ast}$ on $\ttG(\nbigp_{\ast}\nbigv_{\omega})$.
It is easy to check
$\Phi^{\ast}(y_m s)=\gminiq_m y_m\Phi^{\ast}(s)$
for any $s\in\ttG(\nbigv_{\omega})$.
Thus, 
$(\ttG(\nbigp_{\ast}\nbigv),\Phi^{\ast})
=\bigoplus(\ttG(\nbigp_{\ast}\nbigv_{\omega}),\Phi^{\ast})
 \in\Diff_m(\cnum[y,y^{-1}],\gminiq)_{(\rnum,\real)}$.
Thus,
we obtain a functor
$\ttG:\Diff_m(\nbigk,\gminiq)^{\Par}
\lrarr
 \Diff_m(\cnum[y,y^{-1}],\gminiq)_{(\rnum,\real)}$.
We have natural equivalences
$\sfp_{m,m_1\ast}\circ\ttG
\simeq
 \ttG\circ\sfp_{m,m_1\ast}$
and 
$\sfp_{m,m_1}^{\ast}\circ\ttG
\simeq
 \ttG\circ\sfp_{m,m_1}^{\ast}$.

For any $(\nbigp_{\ast}\nbigv,\Phi^{\ast})\in
 \Diff_m(\nbigk,\gminiq)^{\Par}$,
by taking the formal completion
of $\ttG(\nbigp_{\ast},\Phi^{\ast})$,
we obtain 
$(\ttGhat(\nbigv),\Phi^{\ast})
 \in\Diff_m(\nbigk,\gminiq)$.
Moreover,
together with the induced filtered bundle
$\nbigp_{\ast}\ttGhat(\nbigv)$ over $\ttGhat(\nbigv)$,
we obtain
$(\nbigp_{\ast}\ttGhat(\nbigv),\Phi^{\ast})
\in\Diff_m(\nbigk,\gminiq)^{\Par}$.
Clearly,
$\ttG(\nbigp_{\ast}\ttGhat(\nbigv))
\simeq
 \ttG(\nbigp_{\ast}\nbigv)$.

\subsubsection{The generalized eigen decomposition and the weight filtration}

Let $(\nbigp_{\ast}\nbigv,\Phi^{\ast})
 \in\Diff_m(\nbigk,\gminiq)^{\Par}$.
For each $(\omega,a)\in\rnum\times\real$,
$\Gr^{\nbigp}_a(\nbigv_{\omega})$
is equipped with the automorphism 
$F_{a,\omega}$
induced
$y_m^{\ell(m\omega)}(\Phi^{\ast})^{k(m\omega)}$.
We obtain 
the generalized eigen decomposition
$\Gr^{\nbigp}_a(\nbigv_{\omega})
=\bigoplus_{\alpha\in\cnum^{\ast}}
 \EE_{\alpha}\Gr^{\nbigp}_a(\nbigv_{\omega})$.
Let $N_{a,\omega}$ denote the nilpotent endomorphism
of the unipotent part of $F_{a,\omega}$.
We obtain the weight filtration $W$
on $\Gr^{\nbigp}_a(\nbigv_{\omega})$
with respect to $N_{a,\omega}$.
It is compatible with the generalized eigen decomposition.

\subsubsection{Basic examples}
\label{subsection;18.12.21.3}

For $\alpha\in\cnum^{\ast}$,
let $\VV_m(\alpha)=\nbigk_m\,e$ be a Fuchsian $\gminiq_m$-difference
$\nbigk_m$-module
defined by $\Phi^{\ast}(e)=\alpha e$,
as in Example \ref{example;18.12.20.1}.
For $a\in\real$,
we define the filtered bundle $\nbigp_{\ast}^{(a)}\VV_m(\alpha)$
over $\VV_m(\alpha)$
as follows.
\[
 \nbigp^{(a)}_b\VV_m(\alpha)=
 \nbigr_m y^{-[b-a]}e.
\]
Here, we set $[c]:=\max\{n\in\seisuu\,|\,n\leq c\}$
for $c\in\real$.
Thus, we obtain
$(\nbigp^{(a)}_{\ast}\VV_m(\alpha),\Phi^{\ast})
\in \Diff_m(\nbigk,\gminiq;0)^{\Par}$.
\begin{lem}
$\ttG(\nbigp_{\ast}^{(a)}\VV_m(\alpha),\Phi^{\ast})
\simeq
\LL_m^{\ttG}(0,a)\otimes
\VV^{\ttG}_m(\alpha)$.
\hfill\qed
\end{lem}

Let $V$ be a finite dimensional $\cnum$-vector space
with a unipotent automorphism $u$.
For any $a\in\real$,
we define a filtered bundle $\nbigp_{\ast}^{(a)}\VV_m(V,u)$
over $\VV_m(V,u)$ by 
\[
 \nbigp^{(a)}_b\VV_m(V,u)
=\nbigr_m y^{-[b-a]}\otimes_{\cnum}V.
\]
Thus, we obtain 
$(\nbigp^{(a)}_{\ast}\VV_m(V,u),\Phi^{\ast})
\in\Diff_m(\nbigk,\gminiq;0)^{\Par}$.

\begin{lem}
$\ttG(\nbigp^{(a)}_{\ast}\VV_m(V,u),\Phi^{\ast})
\simeq
 \LL^{\ttG}_m(0,a)\otimes
 \VV^{\ttG}_m(V,u)$.
\hfill\qed
\end{lem}

Take $\omega\in\frac{1}{m}\seisuu$.
For any $a\in\real$,
we define the filtered bundle
$\nbigp_{\ast}^{(a)}\LL_m(\omega)$
as follows:
\[
 \nbigp_{b}^{(a)}\LL_m(\omega)
=\nbigr_m y^{-[b-a]}\,\sle_{m,\omega}.
\]
Thus, we obtain
$(\nbigp^{(a)}_{\ast}\LL_m(\omega),\Phi^{\ast})
\in\Diff_m(\nbigk,\gminiq)^{\Par}$.
More generally, for any $\omega\in \rnum$.
Set $m_1:=k(m\omega)\cdot m$.
We define
\[
\nbigp^{(a)}_{\ast}\LL_m(\omega):=
 \sfp_{m,m_1\ast}\bigl(
 \nbigp^{(am_1/m)}_{\ast}\LL_{m_1}(\omega)
 \bigr).
\]
The following is easy to see.
\begin{lem}
$\ttG\bigl(
 \nbigp^{(a)}_{\ast}\LL_m(\omega)
 \bigr)
\simeq
 \LL^{\ttG}_m(\omega,a)$.
\hfill\qed
\end{lem}

For $i=1,\ldots,N$,
we take
$\omega_i\in\rnum$,
$a_i\in\real$,
$\alpha_i\in\cnum^{\ast}$,
and finite dimensional vector spaces $V_i$
with a unipotent automorphism $u_i$.
Let us consider
\[
 (\nbigp_{\ast}\nbigv,\Phi^{\ast})
=\bigoplus_i
 \nbigp^{(a_i)}_{\ast}\LL_m(\omega_i)
 \otimes
 \nbigp^{(0)}_{\ast}\VV_m(\alpha_i)
\otimes
 \nbigp^{(0)}_{\ast}\VV_m(V_i,u_i).
\]
Then, we have
\[
 \ttG(\nbigp_{\ast}\nbigv,\Phi^{\ast})
\simeq
\bigoplus_i
 \LL^{\ttG}_m(\omega_i,a_i)
\otimes
 \VV_m^{\ttG}(\alpha_i)
\otimes
 \VV_m^{\ttG}(V_i,u_i).
\]

\section{Mini-complex manifolds}

\subsection{A twistor family of mini-complex manifolds}
\label{subsection;18.8.19.20}

\subsubsection{A hyperk\"ahler manifold $X$ equipped with 
$\real\times\seisuu^2$-action}
\label{subsection;18.8.5.2}

Take $\mu_1,\mu_2\in\cnum$
which are linearly independent over $\real$.
We assume that
$\Image(\mu_2/\mu_1)>0$.
Let $\Gamma$ denote the lattice of $\cnum$
generated by $\mu_1$ and $\mu_2$.
Let  $\Vol(\Gamma)$ denote 
the volume of $\cnum/\Gamma$ 
with respect to the volume form $\frac{\sqrt{-1}}{2}dz\,d\zbar$,
where $z$ is the standard coordinate of $\cnum$.
The following holds:
\[
 \Vol(\Gamma)=\frac{1}{2\sqrt{-1}}(\mu_2\mubar_1-\mubar_2\mu_1).
\]

We set $X:=\cnum_z\times\cnum_w$
with the Euclidean metric $dz\,d\zbar+dw\,d\wbar$.
It is a hyperk\"ahler manifold.
Let us consider the action of the group 
$\seisuu \tte_1\oplus\seisuu \tte_2$ on $X$
given by
\[
 n_i\tte_i(z,w)=(z,w)+n_i(\mu_i,0).
\]
We also consider the action of $\real\,\tte_0$ on $X$
given by
\[
 s\tte_0(z,w)=(z,w+s).
\]
Thus, we obtain an action of
$\real\tte_0\oplus\seisuu\tte_1\oplus\seisuu\tte_2$
on $X$.

\subsubsection{Complex manifold $X^{\lambda}$}

For each $\lambda\in\cnum$,
there exists the complex structure of $X$
given by the coordinate system
\[
 (\xi,\eta)=(z+\lambda\wbar,w-\lambda\zbar).
\]
The complex manifold is denoted by $X^{\lambda}$.
The action of 
$\real \tte_0\oplus
 \seisuu\tte_1\oplus \seisuu\tte_2$
is described as follows with respect to
the coordinate system $(\xi,\eta)$:
\[
 s\tte_0(\xi,\eta)
=(\xi,\eta)+(\lambda s,s),
\quad\quad
 n_i\tte_i(\xi,\eta)
=(\xi,\eta)+n_i(\mu_i,-\lambda \overline{\mu}_i)
\,\,\,\,\,(i=1,2).
\]

\subsubsection{Some calculations}

To introduce a more convenient complex coordinate system
of $X^{\lambda}$,
we make some calculations.

\begin{lem}
\label{lem;18.3.16.2}
There exist
$\tts_1\in\real$ and $\ttg_1\in\cnum$ with $|\ttg_1|=1$
such that
\begin{equation}
\label{eq;18.8.5.1}
 -\lambda\mubar_1+\tts_1=
 \ttg_1(\mu_1+\lambda \tts_1)\neq 0.
\end{equation}
\begin{itemize}
\item
If $|\lambda|\neq 1$,
there are two choices of $(\tts_1,\ttg_1)$.
One is contained in
$\real_{>0}\times S^1$,
and the other is contained in
$\real_{<0}\times S^1$.
Moreover, $1-\ttg_1\lambda\neq 0$ holds.
\item
If $|\lambda|=1$ and $\lambda\neq \pm\sqrt{-1}\mu_1|\mu_1|^{-1}$,
such $(\tts_1,\ttg_1)$ is 
uniquely determined as 
$(\tts_1,\ttg_1)=(0,-\lambda\mubar_1\mu_1^{-1})$.
Moreover, $1-\ttg_1\lambda\neq 0$ holds.
\item
If $\lambda=\pm\sqrt{-1}\mu_1|\mu_1|^{-1}$,
the set of such $(\tts_1,\ttg_1)$
is 
$\{(s,\lambda^{-1})\,|\,s\in\real\}$.
\end{itemize}
\end{lem}
\pf
Let us consider the condition 
$\bigl|
 -\lambda \mubar_1+\tts_1
 \bigr|
=\bigl|
 \mu_1+\lambda \tts_1
 \bigr|$
for $\tts_1\in\real$.
It is equivalent to the following:
\begin{equation}
\label{eq;18.3.16.1}
 (1-|\lambda|^2)\tts_1^2
-2(\lambda\mubar_1+\lambdabar\mu_1)\tts_1
-(1-|\lambda|^2)|\mu_1|^2=0.
\end{equation}

If $|\lambda|\neq 1$,
there exist two distinct solutions:
\[
 \tts_1=
 \frac{\lambda \mubar_1+\lambdabar\mu_1}{1-|\lambda|^2}
\pm
 \left(
  \frac{(\lambda \mubar_1+\lambdabar\mu_1)^2}{(1-|\lambda|^2)^2}
+|\mu_1|^2
 \right)^{1/2}
=
 \frac{\lambda \mubar_1+\lambdabar\mu_1}{1-|\lambda|^2}
\pm
 \frac{|\mu_1+\lambda^2\mubar_1|}{\bigl|1-|\lambda|^2\bigr|}.
\]
Hence, we obtain
\[
 \mu_1+\lambda\tts_1
=
 \frac{\mu_1+\lambda^2\mubar_1}{1-|\lambda|^2}
\pm
 \lambda\frac{|\mu_1+\lambda^2\mubar_1|}{\bigl|1-|\lambda|^2\bigr|}.
\]
Because $|\lambda|\neq 1$,
we obtain $\mu_1+\lambda\tts_1\neq 0$.
Once we choose $\tts_1$,
we obtain a unique complex number $\ttg_1$
satisfying $|\ttg_1|=1$
determined by the condition
(\ref{eq;18.8.5.1}).
Because $|\lambda|\neq 1$ and $|\ttg_1|=1$,
we obtain $1-\ttg_1\lambda\neq 0$.

If $|\lambda|=1$
and $\lambda\neq \pm\sqrt{-1}\mu_1|\mu_1|^{-1}$,
we obtain
$\lambda\mubar_1+\lambdabar\mu_1\neq 0$,
and hence the equation (\ref{eq;18.3.16.1})
has a unique solution $\tts_1=0$.
In this case, $\ttg_1$ is determined by
$-\lambda\mubar_1=\ttg_1\mu_1$,
i.e.,
$\ttg_1=-\lambda\mubar_1 \mu_1^{-1}$.
The following holds:
\[
1-\lambda\ttg_1=1+\lambda^2\mubar_1/\mu_1
=\lambda\mu_1^{-1}(\lambdabar\mu_1+\lambda\mubar_1)\neq 0.
\]

If $\lambda=\pm\sqrt{-1}\mu_1|\mu_1|^{-1}$,
we can check the claim by a direct computation.
\hfill\qed

\begin{lem}
The following holds:
\begin{equation}
\label{eq;18.8.8.2}
 \Image\Bigl(
 \frac{\ttg_1\mu_2+\lambda\mubar_2}{1-\ttg_1\lambda}
 \Bigr)
=\frac{\Vol(\Gamma)}{\Re(\ttg_1\mu_1)}
\neq 0.
\end{equation}
In particular, $\Re(\ttg_1\mu_1)\neq 0$.
If $|\lambda|\neq 1$,
we have
$\Re(\ttg_1\mu_1)
\cdot
 (1-|\lambda|^2)\tts_1>0$.
\end{lem}
\pf
Let us consider the case $|\lambda|\neq 1$.
Because
$(1-\ttg_1\lambda)^{-1}
=\tts_1(\ttg_1\mu_1+\lambda\mubar_1)^{-1}$,
the following holds:
\[
\Image\Bigl(
\frac{\ttg_1\mu_2
+\lambda\mubar_2}{1-\ttg_1\lambda}
\Bigr)
=
 \Image\Bigl(
 \frac{
 (\ttg_1\mu_2
+\lambda\mubar_2)\tts_1}
{\ttg_1\mu_1
+\lambda\mubar_1}
 \Bigr)
=\Image\Bigl(
 \frac{(1-|\lambda|^2)(\mu_2\mubar_1-\mu_1\mubar_2)\tts_1}
 {2|\ttg_1\mu_1+\lambda\mubar_1|^2}
 \Bigr)
=\frac{(1-|\lambda|^2)\Vol(\Gamma) \tts_1}
 {|\ttg_1\mu_1+\lambda\mubar_1|^2}.
\]
By using $(1-\ttg_1\lambda)^{-1}
=\tts_1(\ttg_1\mu_1+\lambda\mubar_1)^{-1}$ again,
we obtain
\begin{equation}
\label{eq;18.8.8.1}
 \frac{|\ttg_1\mu_1+\lambda\mubar_1|^2}{\tts_1(1-|\lambda|^2)}
=\frac{1}{1-|\lambda|^2}
 \Bigl(
 \ttgbar_1\mubar_1+\lambdabar\mu_1
-\lambda\mubar_1-|\lambda|^2\ttg_1\mu_1
 \Bigr).
\end{equation}
Because the left hand side of (\ref{eq;18.8.8.1}) is real,
it is equal to the following:
\[
 \frac{1}{2(1-|\lambda|^2)}
 \Bigl(
 \ttgbar_1\mubar_1+\ttg_1\mu_1
-|\lambda|^2\ttg_1\mu_1-|\lambda|^2\ttgbar_1\mubar_1
 \Bigr)
=\Re(\ttg_1\mu_1).
\]
Thus, we obtain (\ref{eq;18.8.8.2}).

Suppose $|\lambda|=1$.
Because $\ttg_1=-\lambda\mubar_1\mu_1^{-1}$,
the following holds:
\[
 \Image\Bigl(
 \frac{\ttg_1\mu_2+\lambda\mubar_2}{1-\ttg_1\lambda}
 \Bigr)
=\Image\Bigl(
 \frac{\mubar_2\mu_1-\mu_2\mubar_1}{\lambdabar\mu_1+\lambda\mubar_1}
 \Bigr)
=\frac{-2\Vol(\Gamma)}{(\lambdabar\mu_1+\lambda\mubar_1)}
=\frac{-2\Vol(\Gamma)}{-(\ttgbar_1\mubar_1+\ttg_1\mu_1)}
=\frac{\Vol(\Gamma)}{\Re(\ttg_1\mu_1)}.
\]
Thus, we are done.
\hfill\qed

\begin{lem}
$(\ttg_1-\lambdabar)(\mu_1+\lambda\tts_1)
=(1+|\lambda|^2)\Re(\ttg_1\mu_1)$ holds.
\end{lem}
\pf
We have
$(\ttg_1-\lambdabar)(\mu_1+\lambda\tts_1)
=-\lambdabar(\mu_1+\lambda\tts_1)
-\lambda\mubar_1+\tts_1
=(1-|\lambda|^2)\tts_1-\lambdabar\mu_1-\lambda\mubar_1$.
In particular, it is a real number.
We have the following:
\[
 (\ttg_1-\lambdabar)(\mu_1+\lambda\tts_1)
=\ttg_1(\mu_1+\lambda\tts_1)
-\lambdabar\ttgbar_1(-\lambda\mubar_1+s_1)
=\ttg_1\mu_1+|\lambda|^2\mubar_1\ttgbar_1
+\tts_1\bigl(\ttg_1\lambda-\lambdabar\ttgbar_1\bigr).
\]
Because it is a real number,
it is equal to
$\Re(\ttg_1\mu_1+|\lambda|^2\mubar_1\ttgbar_1)
=(1+|\lambda|^2)\Re(\ttg_1\mu_1)$.
\hfill\qed

\begin{lem}
\label{lem;18.12.14.2}
Suppose $|\lambda|\neq 1$.
Let $(\tts_1,\ttg_1)$ and $(\tts_1',\ttg_1')$
be two solutions of
the equation {\rm(\ref{eq;18.8.5.1})}.
Then, the following holds:
\[
 \Re(\mu_1\ttg_1)+\Re(\mu_1\ttg_1')=0.
\]
\end{lem}
\pf
The following holds:
\[
 \Re(\mu_1\ttg_1)
=\Re\Bigl(
 \frac{(-\lambda\mubar_1+s_1)\mu_1}{\mu_1+\lambda \tts_1}
 \Bigr)
=\frac{|\mu_1+\mubar_1\lambda^2|^2}{(1-|\lambda|^2)}
 \cdot
 \frac{\tts_1}{|\mu_1|^2+(\lambdabar\mu_1
 +\lambda\mubar_1)\tts_1+|\lambda|^2\tts_1^2}.
\]
The following holds:
\begin{multline}
\label{eq;18.12.14.1}
 \tts_1\bigl(
 |\mu_1|^2+(\lambdabar\mu_1+\lambda\mubar_1)\tts_1'+
 |\lambda|^2(\tts_1')^2
 \bigr)
+
 \tts_1'\bigl(
 |\mu_1|^2+(\lambdabar\mu_1+\lambda\mubar_1)\tts_1+
 |\lambda|^2\tts_1^2
 \bigr)
\\
=(\tts_1+\tts_1')(|\mu_1|^2+\tts_1\tts_1'|\lambda|^2)
+2\tts_1\tts_1'(\lambdabar\mu_1+\lambda\mubar_1).
\end{multline}
By using
$\tts_1+\tts_1'=2(\lambdabar\mu_1+\lambda\mubar_1)(1-|\lambda|^2)^{-1}$
and 
$\tts_1\tts_1'=-|\mu_1|^2$,
we obtain that (\ref{eq;18.12.14.1}) is $0$.
Then, we obtain the claim of the lemma by a direct calculation.
\hfill\qed

\begin{lem}
\label{lem;19.2.5.20}
We have
$(1+|\lambda|^2)\bigl|\Re(\mu_1\ttg_1)\bigr|
=|\mu_1+\lambda^2\mubar_1|$.
If $|\lambda|\neq 1$,
we have the following more precise formula:
\begin{equation}
\label{eq;19.2.5.101}
 (1+|\lambda|^2)\Re(\mu_1\ttg_1)
=\sgn(1-|\lambda|^2)\cdot\sgn(\tts_1)\cdot
\bigl|\mu_1+\lambda^2\mubar_1\bigr|.
\end{equation}
\end{lem}
\pf
If $|\lambda|=1$,
we have
$\Re(\mu_1\ttg_1)=\Re(-\lambda\mubar_1)$.
Because $|\lambda|=1$,
we also have
$\bigl|\mu_1+\lambda^2\mubar_1\bigr|
=\bigl|\lambdabar\mu_1+\lambda\mubar_1\bigr|
=2\bigl|\Re(\lambda\mubar_1)\bigr|$.
Hence, the claim of the lemma is clear.

Suppose $|\lambda|\neq 1$.
We have the following:
\begin{multline}
 |\mu_1|^2+(\lambdabar\mu_1+\lambda\mubar_1)\tts_1
+|\lambda|^2\tts_1^2
=|\mu_1|^2(1+|\lambda|^2)
+(\lambdabar\mu_1+\lambda\mubar_1)
 \tts_1
 \Bigl(
 1+\frac{2|\lambda|^2}{1-|\lambda|^2}
 \Bigr)
 \\
=(1+|\lambda|^2)
\Bigl(
 |\mu_1|^2+\frac{\lambdabar\mu_1+\lambda\mubar_1}{1-|\lambda|^2}
 \tts_1
\Bigr).
\end{multline}
We also have the following:
\[
 \frac{|\mu_1|^2}{\tts_1}
=
-\frac{\lambda\mubar_1+\lambdabar\mu_1}{1-|\lambda|^2}
\pm
\frac{|\mu_1+\lambda^2\mubar_1|}{\bigl|1-|\lambda|^2\bigr|}.
\]
Here, $\pm$ is equal to $\sgn(\tts_1)$.
Because
\[
 \frac{\Re(\ttg_1\mu_1)}{\bigl|\mu_1+\lambda^2\mubar_1\bigr|}
=\frac{\bigl|\mu_1+\lambda^2\mubar_1\bigr|}{(1-|\lambda|^2)}
 \frac{\tts_1}{|\mu_1|^2+(\lambdabar\mu_1+\lambda\mubar_1)\tts_1
 +|\lambda_1|^2\tts_1^2}
=\sgn(\tts_1)\frac{\bigl|\mu_1+\lambda^2\mubar_1\bigr|}{(1-|\lambda|^2)}
 \frac{\bigl|1-|\lambda|^2\bigr|}{(1+|\lambda|^2)|\mu_1+\lambda^2\mubar_1|},
\]
we obtain the claim of the lemma.
\hfill\qed

\subsubsection{Coordinate system $(\ttu,\ttv)$}

We introduce a more convenient complex coordinate system
of $X^{\lambda}$.

\begin{assumption}
In the following,
we suppose 
$\lambda\neq \pm\sqrt{-1}\mu_1|\mu_1|^{-1}$.
\hfill\qed
\end{assumption}

We take $\tts_1$ and $\ttg_1$ as in Lemma \ref{lem;18.3.16.2}.
We consider the $\cnum$-linear coordinate change
$\cnum_{\ttu}\times\cnum_{\ttv}
\simeq
 \cnum_{\xi}\times\cnum_{\eta}$
given by
\[
 (\xi,\eta)=(\ttu+\lambda\ttv,\ttg_1\ttu+\ttv),
\quad
 (\ttu,\ttv)=
 \frac{1}{1-\ttg_1\lambda}
 (\xi-\lambda\eta,-\ttg_1\xi+\eta).
\]
The action
of $\real \tte_0\oplus\seisuu \tte_1\oplus\seisuu \tte_2$
on $X^{\lambda}$
is described as follows
in terms of $(\ttu,\ttv)$:
\[
 s\tte_0(\ttu,\ttv)=(\ttu,\ttv)+(0,s),
\quad\quad
 n_i\tte_i(\ttu,\ttv)
=(\ttu,\ttv)
+\frac{n_i}{1-\ttg_1\lambda}
(\mu_i+\lambda^2\mubar_i,
 -\ttg_1\mu_i-\lambda\mubar_i)
\,\,\,(i=1,2).
\]

\begin{lem}
The following holds:
\begin{equation}
\label{eq;18.8.5.3}
 (\tte_1+\tts_1\tte_0)(\ttu,\ttv)
=(\ttu,\ttv)+\Bigl(
 \frac{\mu_1+\lambda^2\mubar_1}{1-\ttg_1\lambda},0
 \Bigr)
=(\ttu,\ttv)+(\mu_1+\lambda \tts_1,0).
\end{equation}
\end{lem}
\pf
Note that the following holds
by our choice of $\tts_1$ and $\ttg_1$:
\[
 \tts_1=\frac{\ttg_1\mu_1+\lambda\mubar_1}{1-\ttg_1\lambda}.
\]
Hence, we obtain the first equality in (\ref{eq;18.8.5.3}).
Note that 
\[
 1-\ttg_1\lambda=
 1-\frac{-\lambda^2\mubar_1+\tts_1\lambda}{\mu_1+\lambda\tts_1}
=\frac{\mu_1+\lambda^2\mubar_1}{\mu_1+\lambda\tts_1}.
\]
Hence,
we obtain
$(\mu_1+\lambda^2\mubar_1)(1-\ttg_1\lambda)^{-1}
=\mu_1+\lambda \tts_1$,
and the second equality in (\ref{eq;18.8.5.3}).
\hfill\qed

\begin{rem}
Let $(E,\nabla,h)$ is an instanton on $X^{\lambda}$.
Let $F(\nabla)=
F_{\xi\xibar}d\xi\,d\xibar+
F_{\xi\etabar}d\xi\,d\etabar+
F_{\eta\xibar}d\eta\,d\xibar+
F_{\eta\etabar}d\eta\,d\etabar$
denote the curvature.
For $\alpha\in\cnum$,
set $H_{\alpha}:=\{(\ttu,\ttv)\,|\,\ttv=\alpha\}
\subset X^{\lambda}$.
Because
$d\xi\,d\xibar=d\ttu\,d\ttubar$
and 
$d\eta\,d\etabar=d\ttu\,d\ttubar$
on $H_{\alpha}$,
the restriction of $F(\nabla)$
to $H_{\alpha}$
is equal to the restriction of
$F_{\xi\etabar}d\xi\,d\etabar+
F_{\eta\xibar}d\eta\,d\xibar$.

In the study of doubly periodic monopoles,
it is appropriate to assume the boundedness of
$F(\nabla)$.
In general, $F_{\xi\xibar}$ and $F_{\eta\etabar}$
are only bounded,
but $F_{\xi\etabar}$ and $F_{\eta\xibar}$
decay more rapidly.
We consider the above coordinate
$(\ttu,\ttv)$ to obtain an appropriate curvature decay
along $H_{\alpha}$.
\hfill\qed
\end{rem}

\subsubsection{Partial quotient $Y_p^{\lambda}$
and its partial compactification}
\label{subsection;18.11.21.2}

Take $p\in\seisuu_{>0}$.
Let $Y_p^{\lambda}$ denote the quotient space of
$X^{\lambda}$ by the action of $\seisuu\cdot p(\tte_1+\tts_1\tte_0)$.
There exists the following induced holomorphic function
on $Y_p^{\lambda}$:
\[
\ttU_p:=\exp\Bigl(
 2\pi\sqrt{-1}\frac{1-\ttg_1\lambda}{p(\mu_1+\lambda^2\mubar_1)}\ttu
 \Bigr)
=\exp\Bigl(
 2\pi\sqrt{-1}\frac{1}{p(\mu_1+\lambda \tts_1)}\ttu
 \Bigr).
\]
We obtain the holomorphic isomorphism
$Y_p^{\lambda}\simeq
 \cnum^{\ast}\times\cnum$ 
induced by $(\ttU_p,\ttv)$,
with which
we identify $Y_p^{\lambda}$ and $\cnum^{\ast}\times\cnum$.
We set $\Ybar_p^{\lambda}:=\proj^1\times\cnum$,
which is a partial compactification of $Y_p^{\lambda}$.
We set
\[
 \gminiq_p^{\lambda}:=
 \exp\Bigl(
 2\pi\sqrt{-1}
 \frac{\mu_2+\lambda^2\mubar_2}{p(\mu_1+\lambda^2\mubar_1)}
 \Bigr),
\quad\quad
 \gminit^{\lambda}:=
-\frac{\Vol(\Gamma)}{\Re(\ttg_1\mu_1)}.
\]
Then, the action of 
$\real \tte_0\oplus\seisuu \tte_2$ on $Y^{\lambda}$
is described as follows:
\[
 s\tte_0(\ttU_p,\ttv)=(\ttU_p,\ttv+s),
\quad\quad
  n\tte_2(\ttU_p,\ttv)
=\left(
 (\gminiq_p^{\lambda})^n\ttU_p,\,\,
 \ttv+\sqrt{-1}n \gminit^{\lambda}
\right).
\]

We set $\nbigg_p:=
(\seisuu/p\seisuu)\cdot (\tte_1+\tts_1\tte_0)$.
There exists the induced $\nbigg_p$-action on $Y_p$
described as follows:
\[
 (\tte_1+\tts_1\tte_0)(\ttU_p,\ttv)
=\Bigl(
 \exp(2\pi\sqrt{-1}/p)\cdot \ttU_p,\ttv
 \Bigr).
\]
The action naturally extends to an action on $\Ybar_p^{\lambda}$.

\begin{rem}
The following holds:
\begin{equation}
\label{eq;19.2.5.21}
 \Image\Bigl(
 \frac{\mu_2+\lambda^2\mubar_2}{\mu_1+\lambda^2\mubar_1}
 \Bigr)
=\frac{(1+|\lambda|^2)(1-|\lambda|^2)\Vol(\Gamma)}
 {|\mu_1+\lambda^2\mubar_1|^2}.
\end{equation}
In particular,
we obtain;
$|\gminiq_{p}^{\lambda}|<1$
in the case  $|\lambda|<1$;
$|\gminiq_{p}^{\lambda}|>1$
in the case $|\lambda|>1$;
$|\gminiq_{p}^{\lambda}|=1$
in the case $|\lambda|=1$.
\hfill\qed
\end{rem}

\subsubsection{Mini-complex manifolds
$\nbigm_p^{\lambda\cov}$ and
$\nbigmbar_p^{\lambda\cov}$}

Let $\nbigm_p^{\lambda\,\cov}$ be the quotient space of
$Y_p^{\lambda}$ by the action of $\real \tte_0$.
By setting $\ttt:=\Image(\ttv)$, 
we obtain the mini-complex coordinate system
$(\ttU_p,\ttt)$ 
of $\nbigm_p^{\lambda\,\cov}$.
The coordinate system induces
the identification
$\nbigm_p^{\lambda\,\cov}\simeq\cnum^{\ast}\times\real$.
The induced action of $\seisuu \tte_2$
is described as follows:
\[
 \tte_2(\ttU_p,\ttt)
=(\gminiq_p^{\lambda}\ttU_p,\ttt+\gminit^{\lambda}).
\]
Note that $\nbigm_p^{\lambda\,\cov}$
is naturally identified with
the quotient space of
$X^{\lambda}$
by the action of $\real\tte_0\oplus\seisuu\cdot  p\tte_1$.

Similarly, let $\nbigmbar_p^{\lambda\,\cov}$
denote the quotient space of
$\Ybar_p^{\lambda}$ by $\real \tte_0$.
It is naturally a mini-complex manifold
and naturally identified with
$\proj^1\times\real$.

We set
$H_{0,p}^{\lambda\cov}:=
 \{0\}\times\real$
and
$H_{\infty,p}^{\lambda\cov}:=
 \{\infty\}\times\real$
in $\proj^1\times\real$.
We set
$H_p^{\lambda\,\cov}:=
 H_{0,p}^{\lambda\cov}
\cup
 H_{\infty,p}^{\lambda\cov}$.

The $\nbigg_p$-action on $Y_p$
induces $\nbigg_p$-actions on $\nbigm^{\lambda\cov}_p$
and $\nbigmbar^{\lambda\cov}_p$.
We identify $\nbigg_p$ 
and $(\seisuu/p\seisuu)\cdot \tte_1$
by $\tte_1+\tts_1\tte_0\longmapsto \tte_1$.
Then, the action of $\nbigg_p$ on $\nbigmbar^{\lambda\cov}_p$
is identified with 
$\tte_1(\ttU_p,\ttt)=(e^{2\pi\sqrt{-1}/p}\ttU_p,\ttt)$.

\subsubsection{Mini-complex manifolds
$\nbigm_p^{\lambda}$ and $\nbigmbar_p^{\lambda}$}

Let $\nbigmlambda_p$ be the mini-complex manifold
obtained as the quotient space of
$\nbigm_p^{\lambda\cov}$ by the action of $\seisuu\tte_2$.
Similarly,
let $\nbigmbarlambda_p$ be the mini-complex manifold
obtained as the quotient space of
$\nbigmbar_p^{\lambda\cov}$ by the action of $\seisuu\tte_2$.
Let $H^{\lambda}_{\nu,p}$ $(\nu=0,\infty)$
denote the quotient of $H^{\lambda\cov}_{\nu,p}$
by $\seisuu\tte_2$.
We set
$H_p^{\lambda}:=H^{\lambda}_{0,p}\cup H^{\lambda}_{\infty,p}$,
which is the quotient of $H^{\lambda\cov}_p$ by $\seisuu\tte_2$.
We have 
$\nbigmbarlambda_p=\nbigmlambda_p
\cup H^{\lambda}_p$,
and it is compact.
There exist the naturally induced
$\nbigg_p$-actions
on $\nbigm^{\lambda}_p$
and $\nbigmbar^{\lambda}_p$.

Let 
$\ttP^{\lambda}_p:\nbigmbar_p^{\lambda\cov}
\lrarr\nbigmbar_p^{\lambda}$
denote the projections.
Let $\sfp_{p_1,p_2}^{\cov}:
 \nbigmbar_{p_2}^{\lambda\cov}
\lrarr
 \nbigmbar_{p_1}^{\lambda\cov}$
and 
$\sfp_{p_1,p_2}:
 \nbigmbar_{p_2}^{\lambda}
\lrarr
 \nbigmbar_{p_1}^{\lambda}$
denote the naturally induced morphisms.

\subsubsection{Neighbourhoods of $H^{\lambda}_{\nu,p}$}
\label{subsection;18.12.17.1}

Let 
$\Psi:\nbigm_p^0\lrarr \real_{y_0}$ be the proper map induced by
$(z,w)\longmapsto \Image(w)$.
Set $\nbigh_R:=\{y_0<R\}$ and 
$\nbigu^0_{p,R}:=\Psi^{-1}(\nbigh_R)\subset\nbigm^0$.
The corresponding open subset
in $\nbigm_p^{\lambda}$ is denoted by
$\nbigu^{\lambda}_{p,R}$.

\begin{lem}
The following holds:
\begin{equation}
 \label{eq;19.2.5.50}
\Psi(\ttU_p,\ttt)
=-\frac{p}{2\pi}\Re(\ttg_1\mu_1)\log|\ttU_p|
+\frac{1-|\lambda|^2}{1+|\lambda|^2}\ttt.
\end{equation}
\end{lem}
\pf
We have the following description of $\Image(w)$
in terms of $(\ttu,\ttv)$:
\[
 \Image(w)=\frac{1}{1+|\lambda|^2}
 \Bigl(
 \Image\bigl((\ttg_1-\lambdabar)\ttu\bigr)
+(1-|\lambda|^2)\Image(\ttv)
 \Bigr).
\]
The following holds:
\[
 \log|\ttU_p|
=\Re\Bigl(
2\pi\sqrt{-1}\frac{\ttu}{p(\mu_1+\lambda\tts_1)}
 \Bigr)
=-\frac{2\pi}{p}
 \Image\Bigl(
 \frac{\ttu}{(\mu_1+\lambda\tts_1)}
 \Bigr).
\]
Because
\[
-\frac{p}{2\pi}
 \frac{(\mu_1+\lambda\tts_1)}{\ttu}
\times
 (\ttg_1-\lambdabar)\ttu
=-\frac{p}{2\pi}(1+|\lambda|^2)\Re(\ttg_1\mu_1),
\]
the claim follows.
\hfill\qed

\vspace{.1in}

\begin{cor}
If $\Re(\ttg_1\mu_1)<0$,
$\nbigu^{\lambda}_{p,R}\cup H^{\lambda}_{0,p}$
is a neighbourhood of 
$H^{\lambda}_{0,p}$.
If $\Re(\ttg_1\mu_1)>0$,
$\nbigu^{\lambda}_{p,R}\cup H^{\lambda}_{\infty,p}$
is a neighbourhood of 
$H^{\lambda}_{\infty,p}$.
\hfill\qed
\end{cor}

\subsubsection{Complement in the case $|\lambda|\neq 1$}

Suppose $|\lambda|\neq 1$.
For simplicity we assume $p=1$.
We use the notation
$\ttU$, $\gminiq$, etc.,
instead of
$\ttU_p$, $\gminiq^{\lambda}_p$, etc.
According to Lemma {\rm\ref{lem;19.2.5.20}},
the following holds:
\[
 \gminit^{\lambda}=
 -\Vol(\Gamma)\frac{1+|\lambda|^2}{|\mu_1+\lambda^2\mubar_1|}
 \sgn(\tts_1)\sgn(1-|\lambda|^2).
\]
By {\rm(\ref{eq;19.2.5.21})},
the following holds:
\[
 \log\bigl|\gminiq^{\lambda}\bigr|
=-2\pi\frac{(1+|\lambda|^2)(|\lambda|^2-1)\Vol(\Gamma)}
 {|\mu_1+\lambda^2\mubar_1|^2}.
\]
We obtain the following:
\begin{equation}
\label{eq;19.2.5.100}
 \frac{\log|\gminiq^{\lambda}|}{\gminit^{\lambda}}
=\frac{2\pi\bigl|1-|\lambda|^2\bigr|}{|\mu_1+\lambda^2\mubar_1|}
 \sgn(\tts_1),
\quad\quad
 \frac{\log|\gminiq^{\lambda}|}
 {(\gminit^{\lambda})^2}
=\frac{2\pi(|\lambda|^2-1)}{\Vol(\Gamma)(1+|\lambda|^2)}.
\end{equation}
In particular, we obtain the following.
\begin{lem}
\label{lem;19.2.8.201}
If $|\lambda|\neq 1$,
$(\gminit^{\lambda})^{-2}\log|\gminiq^{\lambda}|$
is independent of the choice of
$(\tte_1,\tts_1)$.
\hfill\qed
\end{lem}

Let us rewrite (\ref{eq;19.2.5.50}).
For simplicity, we assume $p=1$.

\begin{lem}
\label{lem;19.2.8.200}
If $|\lambda|\neq 1$,
the following holds.
\[
\Psi(\ttU,\ttt)
=\frac{1-|\lambda|^2}{1+|\lambda|^2}
 \Bigl(
 \ttt-\gminit^{\lambda}\frac{\log|\ttU|}{\log|\gminiq^{\lambda}|}
 \Bigr).
\]
In particular,
$\ttt-\gminit^{\lambda}\frac{\log|\ttU|}{\log|\gminiq^{\lambda}|}$
is independent of the choice of $(\tte_1,\tts_1)$.
\end{lem}
\pf
By (\ref{eq;19.2.5.101}) and (\ref{eq;19.2.5.100}),
we obtain
\[
 \frac{\gminit^{\lambda}}{\log|\gminiq^{\lambda}|}
=\frac{1}{2\pi}\frac{1+|\lambda|^2}{1-|\lambda|^2}
 \Re(\ttg_1\mu_1).
\]
Together with (\ref{eq;19.2.5.50}),
we obtain the claim of the lemma.
\hfill\qed

\subsubsection{Two compactifications in the case $|\lambda|\neq 1$}

If $|\lambda|\neq 1$,
there are two solutions 
$(\tts_1,\ttg_1)$ and $(\tts_1',\ttg_1')$
of (\ref{eq;18.8.5.1}).
We obtain two mini-complex coordinate systems
$(\ttU_p,\ttt)$ and $(\ttU_p',\ttt')$
on $\nbigm_p^{\lambda\cov}$.
We obtain another partial compactification
$\nbigm_p^{\prime\lambda\cov}$
from $(\ttU_p',\ttt')$.
Let $\nbigm_p^{\prime\lambda}$
denote the quotient of
$\nbigm_p^{\prime\lambda\cov}$
by the action of $\seisuu \tte_2$.

By the construction,
we have $\ttU_p=\ttU_p'$.
We have the relation:
\[
 \ttt'=\ttt-2\gminit^{\lambda}\frac{\log|\ttU|}{\log|\gminiq^{\lambda}|}
\]
The identity on $\nbigm_p^{\lambda\cov}$
is not extended to an isomorphism
$\nbigmbar_p^{\lambda\cov}$
and 
$\nbigmbar_p^{\prime\lambda\cov}$.

\vspace{.1in}
We consider the automorphism $F$ of
$\nbigm_p^{\lambda\cov}$
defined by
$F^{\ast}(\ttU_p)=\ttU_p$
and 
$F^{\ast}(\ttv')=-\ttv$.
Then, $F$ is equivariant with respect to 
the $\seisuu \tte_2$-action
by Lemma \ref{lem;18.12.14.2}.
Moreover $F$ is extended to an isomorphism
$\nbigmbar^{\lambda\cov}_p
\simeq
\nbigmbar^{\prime\lambda\cov}_p$.
Hence, $F$ induces an isomorphism
$\nbigmbar^{\lambda}_p\simeq
 \nbigmbar^{\prime\lambda}_p$.

\subsection{Curvature of mini-holomorphic bundles
with Hermitian metric on $\nbigmlambda$}

\subsubsection{Mini-complex manifold $A^{\lambda}$}

We set $A:=X/\real \tte_0$.
For each $\lambda$,
it is equipped with the mini-complex structure
induced by the complex structure of $X^{\lambda}$.
(See \cite[\S2.6]{Mochizuki-difference-modules}.)
The mini-complex manifold is denoted by $A^{\lambda}$.
There exists the naturally induced action of
$\seisuu\tte_1\oplus\seisuu\tte_2$
on $A^{\lambda}$.
The quotient space of $A^{\lambda}$
by $p\seisuu \tte_1$ is naturally isomorphic to
$\nbigm_p^{\lambda\cov}$,
and 
the quotient space of $A^{\lambda}$
by $p\seisuu \tte_1\oplus\seisuu\tte_2$ 
is naturally isomorphic to
$\nbigm_p^{\lambda}$.

\subsubsection{Coordinate system
$(\alpha,\tau)$ on $A^{\lambda}$}
\label{subsection;18.12.17.20}

We have the complex coordinate system
$(\alpha,\beta)$ on $X^{\lambda}$
determined by the following relation:
\[
 (\xi,\eta)=\alpha(1,-\lambdabar)+\beta(\lambda,1)
=(\alpha+\beta\lambda,-\lambdabar \alpha+\beta),
\quad\quad
 (\alpha,\beta)
=\frac{1}{1+|\lambda|^2}
 \bigl(
 \xi-\lambda\eta,\,
 \eta+\lambdabar\xi
 \bigr).
\]
We can check the following
by direct computations.
\begin{lem}
We have
$d\alpha\,d\alphabar+d\beta\,d\betabar=
(1+|\lambda|^2)^{-1}(d\xi\, d\xibar+d\eta\,d\etabar)
=dz\,d\zbar+dw\,d\wbar$.
\hfill\qed
\end{lem}

The actions of $\real \tte_0\oplus\seisuu \tte_1\oplus \seisuu \tte_2$
are described as follows
with respect to $(\alpha,\beta)$:
\[
 s\tte_0(\alpha,\beta)=(\alpha,\beta)+(0,s),
\quad
 \tte_i(\alpha,\beta)=(\alpha,\beta)+\frac{1}{1+|\lambda|^2}
 \bigl(
 \mu_i+\lambda^2\mubar_i,
 -\lambda\mubar_i+\lambdabar\mu_i
 \bigr)
\,\,\,(i=1,2).
\]

Setting $\tau:=\Image(\beta)$,
we obtain a mini-complex coordinate
$(\alpha,\tau)$ on $A^{\lambda}$.
We have the complex vector fields
$\del_{\alpha}$, $\del_{\alphabar}$
and $\del_{\tau}$ on $A^{\lambda}$.
The induced complex vector fields on
$\nbigm_p^{\lambda}$
are also denoted by the same notation.

We have the following relation:
\[
 \alpha=\frac{1-\lambda\ttg_1}{1+|\lambda|^2}\ttu,
\quad\quad
 \tau=\frac{\Image((\ttg_1+\lambdabar)\ttu)}{1+|\lambda|^2}+\ttt.
\]
Hence, we have the following relation between 
the complex vector fields:
\[
 \del_{\ttubar}=\frac{1-\lambdabar\ttgbar_1}{1+|\lambda|^2}\del_{\alphabar}
-\frac{1}{2\sqrt{-1}}\frac{(\ttgbar_1+\lambda)}{1+|\lambda|^2}\del_{\tau},
\quad\quad
 \del_{\ttu}=\frac{1-\lambda\ttg_1}{1+|\lambda|^2}\del_{\alpha}
+\frac{1}{2\sqrt{-1}}\frac{(\ttg_1+\lambdabar)}{1+|\lambda|^2}\del_{\tau},
\quad\quad
 \del_{\ttt}=\del_{\tau}.
\]

\subsubsection{Monopoles and mini-holomorphic bundles}
\label{subsection;19.2.6.5}

Let $(E,h,\nabla,\phi)$ be a monopole on 
an open subset $\nbigu$ of $\nbigm^{\lambda}_p$,
i.e.,
$E$ is a vector bundle on $\nbigu$
with a Hermitian metric $h$,
a unitary connection $\nabla$,
and an anti-self-adjoint endomorphism $\phi$ of $E$
satisfying the Bogomolny equation 
\begin{equation}
\label{eq;19.2.6.1}
 F(\nabla)=\ast\nabla\phi.
\end{equation}
Here, $F(\nabla)$ denotes the curvature of $\nabla$,
and $\ast$ denotes the Hodge star operator
with respect to the Riemannian metric 
$d\alpha\,d\alphabar+d\tau\,d\tau$.
We have the expression 
$F(\nabla)=
 F(\nabla)_{\alpha\alphabar}\,d\alpha\,d\alphabar
+F(\nabla)_{\alpha\,\tau}\,d\alpha\,d\tau
+F(\nabla)_{\alphabar\,\tau}\,d\alphabar\,d\tau$.
Then, the Bogomolny equation is equivalent to
the pair of the following equations:
\begin{equation}
\label{eq;19.2.6.2}
 \bigl[\nabla_{\alphabar},\nabla_{\tau}\bigr]=0,
\end{equation}
\begin{equation}
\label{eq;19.2.6.3}
 F(\nabla)_{\alpha\alphabar}=
 \frac{\sqrt{-1}}{2}\nabla_{\tau}\phi.
\end{equation}
The equation (\ref{eq;19.2.6.2}) implies that
$\nabla_{\alphabar}$ and $\nabla_{\tau}$
determine a mini-holomorphic structure on $E$.
(See \cite[\S2.2]{Mochizuki-difference-modules}
for mini-holomorphic bundles).

Conversely,
Let $(E,\delbar_E)$ be a mini-holomorphic bundle on 
an $\nbigu$ of $\nbigm_p^{\lambda}$.
We have the differential operators
$\del_{E,\alphabar}$ and $\del_{E,\tau}$.
Let $h$ be a Hermitian metric of $E$.
Recall that we obtain the Chern connection $\nabla_h$
and the Higgs field $\phi_h$.
(See \cite[\S2.3]{Mochizuki-difference-modules}.)
Let $F(h)$ denote the curvature of $\nabla_h$.
We have the expression
$F(h)=F(h)_{\alpha\alphabar}\,d\alpha\,d\alphabar
+F(h)_{\alpha,\tau}\,d\alpha\,d\tau
+F(h)_{\alphabar,\tau}d\alphabar\,d\tau$.
Then, 
$(E,h,\nabla_h,\phi_h)$ is a monopole
if and only if
\begin{equation}
\label{eq;19.2.6.4}
 F(h)_{\alpha\alphabar}
=\frac{\sqrt{-1}}{2}\nabla_{h,\tau}\phi_h.
\end{equation}
If $(E,h,\nabla_h,\phi_h)$ is a monopole,
$(E,\delbar_E,h)$ is also called a monopole.

\subsubsection{Contraction of curvature and the analytic degree}

Let $(E,\delbar_E)$ be a mini-holomorphic bundle
with a Hermitian metric $h$
on an open subset $\nbigu\subset\nbigm^{\lambda}_p$.
We obtain $(E,h,\nabla_h,\phi_h)$ as in \S\ref{subsection;19.2.6.5}.
We set
\begin{equation}
 G(h):=F(h)_{\alpha\alphabar}
-\frac{\sqrt{-1}}{2}\nabla_{h,\tau}\phi_h.
\end{equation}
Note that the Bogomolny equation for 
$(E,h,\nabla_h,\phi_h)$ is equivalent to $G(h)=0$.

\begin{df}
Suppose that  $\Tr G(h)$ is expressed as
$g_1+g_2$,
where $g_1$ is an $L^1$-function on $U$,
and $g_2$ is non-positive everywhere.
Then, we set
$\deg(E,\delbar_E,h):=
 \int_U\Tr G(h)\,\dvol_U
\in \real\cup\{-\infty\}$,
which is called the analytic degree of $(E,\delbar_E,h)$.
\hfill\qed
\end{df}

Let us recall some formulas for $G(h)$.
See \cite[\S2.8]{Mochizuki-difference-modules}
for more detail.
\begin{lem}
Let $V$ be a mini-holomorphic bundle of $E$.
Let $h_V$ be the induced metric of $V$.
Let $p_V$ denote the orthogonal projection of $E$
onto $V$.
Then, the following holds:
\[
 \Tr G(h_V)=\Tr\bigl(G(h)p_V\bigr)
-\bigl|\del_{E,\alphabar}p_V\bigr|_h^2
-\frac{1}{4}\bigl|
 \del_{E,\tau}p_V
 \bigr|_h^2.
\]
In particular,
if $|G(h)|_h$ is integrable,
then $\deg(V,h_V)$ is well defined for any 
mini-holomorphic subbundles $V$ of $E$.
\hfill\qed
\end{lem}

\begin{lem}
Let $h_1$ be another Hermitian metric of $E$.
Let $s$ be the automorphism of $E$
determined by $h_1=h\cdot s$.
Then, the following holds.
\[
 G(h_1)=G(h)-
\del_{E,\alphabar}\bigl(
 s^{-1}\del_{E,h,\alpha}s
 \bigr)
-\frac{1}{4}
 \Bigl[
 \nabla_{h,\tau}-\sqrt{-1}\phi_h,
 \bigl[\nabla_{h,\tau}+\sqrt{-1}\phi_h,s\bigr]
 \Bigr].
\]
As a consequence, we obtain the following equality:
\[
 -\Bigl(
 \del_{\alpha}\del_{\alphabar}+\frac{1}{4}\del_{\tau}^2
 \Bigr)
 \Tr(s)
=\Tr\Bigl(
 s\bigl(G(h_1)-G(h)\bigr)
 \Bigr)
-\bigl|
 s^{-1/2}\del_{E,h,\alpha}s
 \bigr|_h^2
-\frac{1}{4}
 \bigl|
 s^{-1/2}
 \del'_{E,h,\tau}s
 \bigr|_h^2.
\]
The following equality also holds:
\[
  -\Bigl(
 \del_{\alpha}\del_{\alphabar}+\frac{1}{4}\del_{\tau}^2
 \Bigr)
 \log\bigl(
 \Tr(s)
\bigr)
\leq
 \Bigl|
 G(h_1)
 \Bigr|_{h_1}
+\Bigl|
 G(h)
 \Bigr|_h.
\]
If $\rank(E)=1$,
then
$G(h_1)-G(h)
=4^{-1}\Delta \log(s)$ holds on $U$.
\hfill\qed
\end{lem}

\subsubsection{Another expression of $G(h)$}

We introduce the following real vector fields on
$A^{\lambda}$:
\[
\gminiv:=
 (\ttg_1\lambda+\ttgbar_1\lambdabar)\del_{\tau}
+\sqrt{-1}(\ttgbar_1-\lambda^2\ttg_1)\del_{\alpha}
-\sqrt{-1}(\ttg_1-\lambdabar^2\ttgbar_1)\del_{\alphabar}.
\]
The induced vector fields
on $\nbigm_p^{\lambda}$ are also denoted by
$\gminiv$.

Let $(E,\delbar_E)$ be a mini-holomorphic bundle
on an open subset $U\subset \nbigmlambda_p$
with a Hermitian metric $h$.

\begin{prop}
\label{prop;18.9.2.1}
We have the following equality:
\begin{equation}
 G(h)=
 |1-\ttg_1\lambda|^{-2}
  (1+|\lambda|^2)^2
 \bigl[
 \del_{E,h,\ttu},
 \del_{E,\ttubar}
 \bigr]
\\
+|1-\ttg_1\lambda|^{-2}
 \nabla_{h,\gminiv}\phi_h.
\end{equation}
\end{prop}
\pf
We have the following formula 
for complex vector fields:
\[
 \del_{\alphabar}=
 \frac{1+|\lambda|^2}{1-\lambdabar\ttgbar_1}\del_{\ttubar}
+\frac{1}{2\sqrt{-1}}
 \frac{\ttgbar_1+\lambda}{1-\lambdabar\ttgbar_1}\del_{\ttt},
\quad\quad
 \del_{\alpha}=
 \frac{1+|\lambda|^2}{1-\lambda\ttg_1}\del_{\ttu}
-\frac{1}{2\sqrt{-1}}
 \frac{\ttg_1+\lambdabar}{1-\lambda\ttg_1}\del_{\ttt},
\quad\quad
 \del_{\tau}=\del_{\ttt}.
\]
Hence, we have the following formulas:
\[
 \frac{1+|\lambda|^2}{1-\lambdabar\ttgbar_1}
 \del_{E,\ttubar}
=\nabla_{h,\alphabar}
-\frac{1}{2\sqrt{-1}}
 \frac{\ttgbar_1+\lambda}{1-\lambdabar\ttgbar_1}
 (\nabla_{h,\tau}-\sqrt{-1}\phi),
\]
\[
  \frac{1+|\lambda|^2}{1-\lambda\ttg_1}
 \del_{E,h,\ttu}
=\nabla_{h,\alpha}
+\frac{1}{2\sqrt{-1}}
 \frac{\ttg_1+\lambdabar}{1-\lambda\ttg_1}
 (\nabla_{h,\tau}+\sqrt{-1}\phi).
\]
We recall the formulas
$[\nabla_{h,\alphabar},\nabla_{h,\tau}]
=\sqrt{-1}\nabla_{h,\alphabar}\phi_h$
and 
$[\nabla_{h,\alpha},\nabla_{h,\tau}]
=-\sqrt{-1}\nabla_{h,\alpha}\phi_h$.
(See \cite[\S2.8.2]{Mochizuki-difference-modules}.)
Then, we obtain the following:
\begin{multline}
\frac{(1+|\lambda|^2)^2}{|1-\lambda\ttg_1|^2}
 \bigl[
 \del_{E,h,\ttu},\del_{E,\ttubar}
 \bigr]
=\bigl[
 \nabla_{h,\alpha},
 \nabla_{h,\alphabar}
 \bigr]
-\frac{\ttg_1+\lambdabar}{1-\ttg_1\lambda}
 \nabla_{h,\alphabar}\phi
+\frac{\ttgbar_1+\lambda}{1-\ttgbar_1\lambdabar}
 \nabla_{h,\alpha}\phi
-\frac{\sqrt{-1}}{2}\frac{|\ttgbar_1+\lambda|^2}{|1-\lambda\ttg_1|^2}
 \nabla_{\tau}\phi\\
=G(h)
-\frac{\ttg_1+\lambdabar}{1-\ttg_1\lambda}
 \nabla_{h,\alphabar}\phi
+\frac{\ttgbar_1+\lambda}{1-\ttgbar_1\lambdabar}
 \nabla_{h,\alpha}\phi
-\sqrt{-1}
 \frac{\lambda\ttg_1+\lambdabar\ttgbar_1}
 {|1-\lambda\ttg_1|^2}
 \nabla_{\tau}\phi.
\end{multline}
Then, we obtain the desired formula.
\hfill\qed

\vspace{.1in}

Recall that
$(z,w)$ is the complex coordinate system of $X^0$.
By setting $y:=\Image(w)$,
we obtain a mini-complex coordinate system
$(z,y)$ of $A^{0}$.
We obtain the induced complex vector fields
$\del_z$, $\del_{\zbar}$ and $\del_y$
on $\nbigm_p^{0}$.

\begin{lem}
\label{lem;19.2.7.10}
$\gminiv=(1+|\lambda|^2)
 \bigl(
 \sqrt{-1}\ttgbar_1\del_z
-\sqrt{-1}\ttg_1\del_{\zbar}
 \bigr)$
holds.
\end{lem}
\pf
We obtain the following relations between complex vector fields:
\[
(1+|\lambda|^2)
 \del_{\alpha}=
\del_z+\lambdabar^2\del_{\zbar}
+\sqrt{-1}\lambdabar\del_y,
\]
\[
 (1+|\lambda|^2)\del_{\alphabar}
=\del_{\zbar}+\lambda^2\del_z
-\sqrt{-1}\lambda\del_y,
\]
\[
(1+|\lambda|^2)
 \del_{\tau}=2\sqrt{-1}\lambda\del_z
-2\sqrt{-1}\lambdabar\del_{\zbar}
+(1-|\lambda|^2)\del_y.
\]
Then, we obtain the claim of the lemma.
\hfill\qed

\vspace{.1in}
Let us give a consequence.
Suppose $U=\nbigm^{\lambda}_p\setminus Z$,
where $Z$ is a finite set.
We set $S^1_{\lambda}:=\real/\gminit^{\lambda}\seisuu$.
Let $\pi_p^{\cov}:\nbigm_p^{\lambda\cov}\lrarr\real$
be the map defined by
$\pi_p^{\cov}(\ttU_p,\ttt)=\ttt$.
It induces a map $\pi_p:\nbigm_p^{\lambda}\lrarr S^1_{\lambda}$.

\begin{prop}
\label{prop;19.2.8.40}
Suppose that 
$\Tr G(h)$
and 
$\Tr\bigl(
 [\del_{E,h,\ttu},\del_{E,\ttubar}]
 \bigr)$
are integrable on $\nbigm^{\lambda}_p\setminus Z$.
Then,
the following equality holds:
\begin{equation}
\label{eq;18.8.5.12}
 \int_{\nbigm^{\lambda}_p}
 \Tr G(h)\dvol
=
\int_{S^1_{\lambda}}d\ttt
\int_{\pi_{p}^{-1}(\ttt)}
\Tr\bigl(
 \bigl[
 \del_{E,h,\ttu},
 \del_{E,\ttubar}
 \bigr]
\bigr)
 \frac{\sqrt{-1}}{2}d\ttu\,d\ttubar.
\end{equation}
\end{prop}
\pf
By the assumption,
the following holds:
\begin{equation}
\label{eq;18.8.9.20}
 \int_{\nbigm^{\lambda}}
 \Tr G(h)\dvol
=
 \int_{\nbigm^{\lambda}}
 |1-\ttg_1\lambda|^{-2}(1+|\lambda|^2)^2
 \Tr\Bigl(
 \bigl[
 \del_{E,h,\ttu},\del_{E,\ttubar}
 \bigr]
 \Bigr)
\dvol
+
 \int_{\nbigm^{\lambda}}
 \Tr\Bigl(
 |1-\ttg_1\lambda|^{-2}
 \nabla_{h,\gminiv}\phi_h
 \Bigr)
 \dvol.
\end{equation}
Because
\[
 \dvol=\frac{\sqrt{-1}}{2}d\alpha\,d\alphabar\,d\tau
=\frac{|1-\ttg_1\lambda|^2}{(1+|\lambda|^2)^2}
 \frac{\sqrt{-1}}{2}\,d\ttu\,d\ttubar\,d\ttt,
\]
the first term of the right hand side of (\ref{eq;18.8.9.20})
is equal to
the right hand side of (\ref{eq;18.8.5.12}).
Let $T_p$ denote the quotient of
$\cnum$ by $p\seisuu\mu_1+\seisuu\mu_2$.
Because
$\dvol=\frac{\sqrt{-1}}{2}dz\,d\zbar\,dy$,
the following holds:
\[
 \int_{\nbigm^{\lambda}}
 \Tr\Bigl(
 |1-\ttg_1\lambda|^{-2}
 \nabla_{h,\gminiv}\phi_h
 \Bigr)
 \dvol
=\lim_{C\to\infty}
 \int_{-C}^Cdy
 \int_{T_p\times\{y\}}
  \Tr\Bigl(
 |1-\ttg_1\lambda|^{-2}
 \nabla_{h,\gminiv}\phi_h
 \Bigr)
 \frac{\sqrt{-1}}{2}dz\,d\zbar.
\]
Note that
$\int_{T_p\times\{y\}}
 \Tr\bigl(
 \nabla_{h,\gminiv}\phi_h
 \bigr)dz\,d\zbar=0$.
Hence, we obtain (\ref{eq;18.8.5.12}).
\hfill\qed

\section{Good filtered bundles with Dirac type singularity on
$(\nbigmbar^{\lambda};H^{\lambda},Z)$}
\label{section;19.1.24.1}

\subsection{Good filtered bundles on 
$(\Hhat^{\lambda}_{\nu,p},H^{\lambda}_{\nu,p})$}

\subsubsection{$\nbigo_{\Hhat_{\nu,p}^{\lambda}}
 (\ast H^{\lambda}_{\nu,p})$-modules}
\label{subsection;19.1.19.1}

For $\nu=0,\infty$,
let $\Hhat^{\lambda}_{\nu,p}$
denote the formal completion of $\nbigmbar_p^{\lambda}$
along $H^{\lambda}_{\nu,p}$.
Similarly,
let $\Hhat^{\lambda\cov}_{\nu,p}$
denote the formal completion of $\nbigmbar_p^{\lambda\cov}$
along $H^{\lambda\cov}_{\nu,p}$.
We have the natural $\seisuu \tte_2$-action on
$\Hhat^{\lambda\cov}_{\nu,p}$,
and 
$\Hhat^{\lambda}_{\nu,p}$
is naturally isomorphic to 
the quotient of $\Hhat^{\lambda\cov}_{\nu,p}$.
Hence,
$\nbigo_{\Hhat^{\lambda}_{\nu,p}}
 (\ast H^{\lambda}_{\nu,p})$-modules
are equivalent to
$\seisuu \tte_2$-equivariant
$\nbigo_{\Hhat^{\lambda\cov}_{\nu,p}}
 (\ast H^{\lambda\cov}_{\nu,p})$-modules.
Let $\LFM(\Hhat^{\lambda}_{\nu,p},H^{\lambda}_{\nu,p})$
(resp. $\LFM(\Hhat^{\lambda\cov}_{\nu,p},H^{\lambda\cov}_{\nu,p})$)
denote the category of locally free
$\nbigo_{\Hhat^{\lambda}_{\nu,p}}(\ast H^{\lambda}_{\nu,p})$-modules
(resp. $\nbigo_{\Hhat^{\lambda\cov}_{\nu,p}}
 (\ast H^{\lambda\cov}_{\nu,p})$-modules).

For $\nu=0,\infty$,
let $\nuhat_p$ denote the formal completion of
$\proj^1_{\ttU_p}$ at $\ttU_p=\nu$.
We have the natural isomorphism
$\Hhat^{\lambda\cov}_{\nu,p}\simeq
 \nuhat_{p}\times\real$.
Set $\ttU_{0,p}:=\ttU_p$
and $\ttU_{\infty,p}:=\ttU_p^{-1}$.
We also set
$\gminiq^{\lambda}_{0,p}:=\gminiq^{\lambda}_p$,
and 
$\gminiq^{\lambda}_{\infty,p}:=(\gminiq^{\lambda}_p)^{-1}$.
The $\seisuu\tte_2$-action on $\Hhat^{\lambda\cov}_{\nu,p}$
is described as
$\tte_2(\ttU_{\nu,p},\ttt)=
 (\gminiq_{\nu,p}^{\lambda}\ttU_{\nu,p},\ttt+\gminit^{\lambda})$.
The $\nbigg_p$-action on 
$\Hhat^{\lambda\cov}_{\nu,p}$
is described as
$\tte_1(\ttU_{\nu,p},\ttt)=(e^{\pm 2\pi\sqrt{-1}/p}\ttU_{\nu,p},\ttt)$,
where the signature is $+$ if $\nu=0$,
and $-$ if $\nu=\infty$.

Let 
$\pi_{\nu,p}^{\cov}:
 \Hhat^{\lambda\cov}_{\nu,p}
\lrarr
 \real$
denote the projection.
We have the natural identification
$(\pi_{\nu,p}^{\cov})^{-1}(\ttt)\simeq\nuhat_p$.
We set $S^1_{\lambda}:=\real/\gminit^{\lambda}\seisuu$.
We obtain the induced map
$\pi_{\nu,p}:
 \Hhat^{\lambda}_{\nu,p}
\lrarr S^1_{\lambda}$.
For each $\ttt\in S^1_{\lambda}$,
once we fix its lift to $\real$,
we obtain an isomorphism
$\pi_{\nu,p}^{-1}(\ttt)\simeq\nuhat_{p}$.

Set $\nbigk_{\nu,p}:=\cnum(\!(\ttU_{\nu,p})\!)$.
Let us observe that locally free
$\nbigo_{\Hhat^{\lambda}_{\nu,p}}
 (\ast H^{\lambda}_{\nu,p})$-modules
are equivalent to
$\gminiq_{\nu,p}^{\lambda}$-difference $\nbigk_{\nu,p}$-modules.
Let $\slq_{\nu,p}:\Hhat^{\lambda\cov}_{\nu,p}\lrarr \nuhat_p$
denote the projection.
Let $(\nbigv,\Phi^{\ast})$ be a $\gminiq_{\nu,p}^{\lambda}$-difference
$\nbigk_{\nu,p}$-module.
We obtain 
the $\nbigo_{\Hhat^{\lambda\cov}_{\nu,p}}(\ast H^{\lambda\cov}_{\nu,p})$-module
$\slq^{\ast}_{\nu,p}\nbigv$.
By the action of $\Phi^{\ast}$,
$\slq^{\ast}_{\nu,p}\nbigv$ is naturally $\seisuu\tte_2$-equivariant.
Hence, we obtain 
an $\nbigo_{\Hhat^{\lambda}_{\nu,p}}(\ast H^{\lambda}_{\nu,p})$-module
as the descent of $\slq^{\ast}_{\nu,p}\nbigv$,
which we denote by 
$\Upsilon^{\lambda}_{\nu,p}(\nbigv)$.
The following is easy to see.
\begin{lem}
$\Upsilon^{\lambda}_{\nu,p}$
induces an equivalence 
$\Diff_{p}(\nbigk_{\nu},\gminiq_{\nu}^{\lambda})
\simeq
 \LFM(\Hhat^{\lambda}_{\nu,p},H^{\lambda}_{\nu,p})$.
The quasi inverse is induced by the restriction
$\gbigv\longmapsto
(\Upsilon^{\lambda}_{\nu,p})^{-1}(\gbigv):=
\gbigv^{\cov}_{|\pi_{\nu,p}^{-1}(0)}$,
where $\gbigv^{\cov}$ is the pull back of $\gbigv$
by $\Hhat^{\lambda\cov}_{\nu,p}\lrarr \Hhat^{\lambda}_{\nu,p}$.
\hfill\qed
\end{lem}

\begin{df}
We say that 
a locally free $\nbigo_{\Hhat^{\lambda}_{\nu,p}}(\ast H^{\lambda}_{\nu,p})$-module
is pure isoclinic of slope $\omega$
if the corresponding 
$\gminiq_{\nu,p}^{\lambda}$-difference $\nbigk_{\nu,p}$-module
is pure isoclinic of slope $\omega$.
Let $\LFM_p(\Hhat^{\lambda}_{\nu,p},H^{\lambda}_{\nu,p};\omega)$
denote the full subcategory of
pure isoclinic modules of slope $\omega$.
A pure isoclinic modules of slope $0$ is also called Fuchsian
ore regular.
\hfill\qed
\end{df}

The following is a consequence of Proposition \ref{prop;18.8.19.10}.
\begin{prop}
Any  $\gbigv\in \LFM(\Hhat^{\lambda}_{\nu,p},H^{\lambda}_{\nu,p})$
has a decomposition
$\gbigv=
 \bigoplus_{\omega\in \rnum}
 \gbigv_{\omega}$
such that 
$\gbigv_{\omega}$
are pure isoclinic of slope $\omega$.
\hfill\qed
\end{prop}

For $p_2\in p_1\seisuu_{>0}$,
we may regard
$\Hhat^{\lambda}_{\nu,p_1}$
as the quotient of
$\Hhat^{\lambda}_{\nu,p_2}$
by the action of the subgroup
$(p_1\seisuu/p_2\seisuu)\tte_1
\subset
 (\seisuu/p_2\seisuu)\tte_1$.
We have the naturally induced morphisms
$\sfp_{p_1,p_2}:
 \Hhat^{\lambda}_{\nu,p_2}\lrarr 
 \Hhat^{\lambda}_{\nu,p_1}$.
We have the pull back
and the push-forward:
\[
 \sfp_{p_1,p_2}^{\ast}:
\LFM(\Hhat^{\lambda}_{\nu,p_1},H^{\lambda}_{\nu,p_1}) 
\lrarr
\LFM(\Hhat^{\lambda}_{\nu,p_2},H^{\lambda}_{\nu,p_2}),
\quad
 \sfp_{p_1,p_2\ast}:
\LFM(\Hhat^{\lambda}_{\nu,p_2},H^{\lambda}_{\nu,p_2})
\lrarr
\LFM(\Hhat^{\lambda}_{\nu,p_1},H^{\lambda}_{\nu,p_1}).
\]
They are compatible with
the pull back and push-forwards
for $\Diff_{p_1}(\nbigk_{\nu},\gminiq_{\nu}^{\lambda})$
between
$\Diff_{p_2}(\nbigk_{\nu},\gminiq_{\nu}^{\lambda})$.
We also have the descent of
$(p_1\seisuu/p_2\seisuu)\tte_1$-equivariant
locally free objects
in $\LFM_{p_2}(\Hhat^{\lambda}_{\nu,p_2},H^{\lambda}_{\nu,p_2})$.

\subsubsection{Filtered bundles on 
 $(\Hhat^{\lambda}_{\nu,p},H^{\lambda}_{\nu,p})$}

\begin{df}
For
any $\gbigv\in \LFM_p(\Hhat^{\lambda}_{\nu,p},H^{\lambda}_{\nu,p})$,
a filtered bundle over $\gbigv$ is defined to be
a family of filtered bundles
$\nbigp_{\ast}(\gbigv_{|\pi_{\nu,p}^{-1}(\ttt)})$ 
$(\ttt\in S^1_{\lambda})$.
Similarly,
for any $\gbigv^{\cov}\in 
 \LFM_p(\Hhat^{\lambda\cov}_{\nu,p},H^{\lambda\cov}_{\nu,p})$,
a filtered bundle over $\gbigv^{\cov}$ is defined to be
a family of filtered bundles
$\nbigp_{\ast}(\gbigv^{\cov}_{|(\pi^{\cov}_{\nu,p})^{-1}(\ttt)})$ 
$(\ttt\in \real)$.
Such families
are often denoted by $\nbigp_{\ast}\gbigv$
and $\nbigp_{\ast}\gbigv^{\cov}$.
\hfill\qed
\end{df}

Let $p_2\in p_1\seisuu_{>0}$.
For any filtered bundle
$\nbigp_{\ast}\bigl(\gbigv\bigr)$
over 
$\gbigv\in\LFM(\Hhat^{\lambda}_{\nu,p_1},H^{\lambda}_{\nu,p_1})$,
we obtain the induced filtered bundle
$\nbigp_{\ast}\bigl(
 \sfp_{p_1,p_2}^{\ast}(\gbigv)
 \bigr)$
over $\sfp_{p_1,p_2}^{\ast}(\gbigv)$.
For any filtered bundle
$\nbigp_{\ast}\bigl(\gbigv\bigr)$
over 
$\gbigv\in\LFM(\Hhat^{\lambda}_{\nu,p_2},H^{\lambda}_{\nu,p_2})$,
we obtain the induced filtered bundle
$\nbigp_{\ast}\bigl(
 \sfp_{p_1,p_2\ast}(\gbigv)
 \bigr)$
over $\sfp_{p_1,p_2\ast}(\gbigv)$.
For any $(p_1\seisuu/p_2\seisuu)$-equivariant
locally free filtered bundle
$\nbigp_{\ast}(\gbigv)$
over a $(p_1\seisuu/p_2\seisuu)$-equivariant
$\gbigv\in\LFM(\Hhat^{\lambda}_{\nu,p_2},H^{\lambda}_{\nu,p_2})$,
we obtain 
$\gbigv_1\in\LFM(\Hhat^{\lambda}_{\nu,p_1},H^{\lambda}_{\nu,p_1})$
as the descent of $\gbigv$,
and we obtain a filtered bundle
$\nbigp_{\ast}(\gbigv_1)$ over $\gbigv_1$
as the decent of $\nbigp_{\ast}(\gbigv)$.

\subsubsection{Good filtered bundles on
$(\Hhat^{\lambda}_{\nu,p},H^{\lambda}_{\nu,p})$}

Let $\gbigv$ be a locally free
$\nbigo_{\Hhat^{\lambda}_{\nu,p}}(\ast H^{\lambda}_{\nu,p})$-module.
\begin{df}
A filtered bundle $\nbigp_{\ast}(\gbigv)$ 
over $\gbigv$ is pure isoclinic of slope $\omega$
if the following holds.
\begin{itemize}
\item
Let $\gbigv^{\cov}\in
 \LFM(\Hhat^{\lambda\cov}_{\nu,p},H^{\lambda\cov}_{\nu,p})$
be the pull back of $\gbigv$.
Take $\ttt_1,\ttt_2\in \real\simeq H^{\lambda\cov}_{\nu,p}$.
Then, under the isomorphism
$\gbigv^{\cov}_{|\pi_{p,\nu}^{-1}(\ttt_1)}\simeq
\gbigv^{\cov}_{|\pi_{p,\nu}^{-1}(\ttt_2)}$
induced by the parallel transport
along the path,
\begin{equation}
\label{eq;19.2.6.10}
 \nbigp_a(\gbigv^{\cov}_{|\pi_{p,\nu}^{-1}(\ttt_1)})
=\nbigp_{a+p\omega(\ttt_2-\ttt_1)/\gminit^{\lambda}}
 (\gbigv^{\cov}_{|\pi_{p,\nu}^{-1}(\ttt_2)})
\end{equation}
holds for any $a\in\real$.
Note that the underlying $\gbigv$ is pure isoclinic of slope $\omega$.
Note also that
$\nbigp_{\ast}(\gbigv^{\cov}_{|\pi_{p,\nu}^{-1}(\ttt)})$
are uniquely determined by
$\nbigp_{\ast}(\gbigv^{\cov}_{|\pi_{p,\nu}^{-1}(0)})$.
\end{itemize}
Let $\LFM(\Hhat^{\lambda}_{\nu,p},H^{\lambda}_{\nu,p};\omega)^{\Par}$
denote the category of filtered flat bundles
over $(\Hhat^{\lambda}_{\nu,p},H^{\lambda}_{\nu,p})$
which are pure isoclinic of slope $\omega$.
\hfill\qed
\end{df}

\begin{rem}
If $\nbigp_{\ast}\gbigv$ has pure slope $0$,
it is also called a regular filtered bundle.
\hfill\qed
\end{rem}

\begin{df}
A filtered bundle $\nbigp_{\ast}\gbigv$
over $(\Hhat^{\lambda}_{\nu,p},H^{\lambda}_{\nu,p})$
is called good
if $\nbigp_{\ast}\gbigv=\bigoplus
 \nbigp_{\ast}\gbigv_{\omega}$,
where 
$\nbigp_{\ast}\gbigv_{\omega}
\in \LFM(\Hhat^{\lambda}_{\nu,p},H^{\lambda}_{\nu,p};\omega)^{\Par}$.
Let $\LFM(\Hhat^{\lambda}_{\nu,p},H^{\lambda}_{\nu,p})^{\Par}$
denote the category of 
good filtered bundles over
$(\Hhat^{\lambda}_{\nu,p},H^{\lambda}_{\nu,p})$.

\hfill\qed
\end{df}

For any
$\nbigp_{\ast}\gbigv\in
 \LFM(\Hhat^{\lambda}_{\nu,p},H^{\lambda}_{\nu,p})^{\Par}$,
the filtered bundle
$\nbigp_{\ast}\bigl(
 (\Upsilon^{\lambda}_{\nu,p})^{-1}(\gbigv)\bigr)$
is defined to be
$\nbigp_{\ast}(\gbigv^{\cov}_{|(\pi_{p,\nu}^{\cov})^{-1}(0)})$.
Conversely,
for any $(\nbigp_{\ast}\nbigv,\Phi)
=\bigoplus_{\omega}(\nbigp_{\ast}\nbigv_{\omega},\Phi)
 \in
 \Diff_{p}(\nbigk_{\nu},\gminiq_{\nu}^{\lambda})$,
the filtered bundle over
$\Upsilon^{\lambda}_{\nu,p}(\nbigv)
=\bigoplus \Upsilon^{\lambda}_{\nu,p}(\nbigv_{\omega})$
is defined by (\ref{eq;19.2.6.10})
and 
$\nbigp_{\ast}\bigl(
 \Upsilon^{\lambda}_{\nu,p}(\nbigv)^{\cov}
 _{|\pi_{p,\nu}^{-1}(0)}
 \bigr)
=\nbigp_{\ast}\nbigv$.
The following is clear.
\begin{lem}
$\LFM(\Hhat^{\lambda},H^{\lambda}_{\nu,p})^{\Par}$
and 
$\Diff_p(\nbigk_{\nu},\gminiq^{\lambda}_{\nu})^{\Par}$
are equivalent by
$\Upsilon^{\lambda}_{\nu,p}$ and 
$(\Upsilon^{\lambda}_{\nu,p})^{-1}$.
They also induce
equivalences between
$\LFM(\Hhat^{\lambda},H^{\lambda}_{\nu,p};\omega)^{\Par}$
and 
$\Diff_p(\nbigk_{\nu},\gminiq^{\lambda}_{\nu};\omega)^{\Par}$.
\hfill\qed
\end{lem}

For any $\nbigp_{\ast}\gbigv\in
 \LFM(\Hhat^{\lambda}_{\nu,p},H^{\lambda}_{\nu,p})^{\Par}$,
we define
$\ttG(\nbigp_{\ast}\gbigv):=
 \ttG(\Upsilon^{-1}(\nbigp_{\ast}\gbigv))
\in
 \Diff_m(\cnum[y,y^{-1}],\gminiq)_{(\rnum,\real)}$.

\subsubsection{Basic examples}
\label{subsection;18.12.21.10}

For any finite dimensional $\cnum$-vector space $V$
with an automorphism $f$,
we set
$\vecVV_{\nu,p}(V,f):=
 \Upsilon^{\lambda}_{\nu,p}\bigl(
 \VV_p(V,f)\bigr)$.
(See Example \ref{example;18.12.20.1}
for $\VV_p(V,f)$.)
Recall that we have constructed 
filtered bundles
$\nbigp^{(a)}_{\ast}\VV_p(V,f)$ over $\VV_p(V,f)$
in \S\ref{subsection;18.12.21.3}.
The $\nbigr_p$-lattices
$\nbigp^{(a)}_b\VV_p(V,f)$ naturally define
$\nbigo_{\Hhat^{\lambda}_{\nu,p}}$-lattices
$\nbigp^{(a)}_b\vecVV_{\nu,p}(V,f)$
of $\vecVV_{\nu,p}(V,f)$.
They induce 
a filtered bundle
$\nbigp^{(a)}_{\ast}\vecVV_{\nu,p}(V,f)$
over $\vecVV_{\nu,p}(V,f)$.

For any $A\in \GL_r(A)$
we set
$\vecVV_{\nu,p}(A):=
\Upsilon^{\lambda}_{\nu,p}\bigl(
 \VV_p(A)\bigr)$.
(See Example \ref{example;18.12.20.1}
for $\VV_p(A)$.)
For any $a\in\real$,
we obtain a filtered bundle
$\nbigp^{(a)}_{\ast}\vecVV_{\nu,p}(A)$
over $\vecVV_{\nu,p}(A)$ similarly.

\begin{lem}
$\nbigp^{(a)}_{\ast}\vecVV_{\nu,p}(V,f)$
and 
$\nbigp^{(a)}_{\ast}\vecVV_{\nu,p}(A)$
are objects in
$\LFM(\Hhat^{\lambda}_{\nu,p},H^{\lambda}_{\nu,p};0)^{\Par}$.
We have the natural isomorphisms
$\ttG(\nbigp^{(a)}_{\ast}\vecVV_{\nu,p}(V,f))
\simeq
 \LL^{\ttG}(0,a)\otimes
 \VV^{\ttG}(V,f)$
and
$\ttG(\nbigp^{(a)}_{\ast}\vecVV_{\nu,p}(A))
\simeq
 \LL^{\ttG}(0,a)\otimes
 \VV^{\ttG}(A)$.
\hfill\qed
\end{lem}

For any $\omega\in\rnum$,
we set 
$\vecLL_{\nu,p}(\omega):=\Upsilon^{\lambda}_{\nu,p}(\LL_p(\omega))$.
(See \S\ref{subsection;18.12.21.2}
for $\LL_p(\omega)$.)
Set $\nbigr_{\nu,p}:=\cnum[\![\ttU_{\nu,p}]\!]$.
If $p\omega\in\seisuu$,
the filtered bundle
$\nbigp^{(0)}_{\ast}\bigl(
 \slq_{\nu,p}^{\ast}\LL_p(\omega)
 \bigr)=
 \bigl(
 \nbigp^{(0)}_{\ast}\bigl(
 \slq_{\nu,p}^{\ast}\LL_p(\omega)_{|(\pi^{\cov}_{\nu,p})^{-1}(\ttt)}
 \bigr)\,\big|\,\ttt\in\real
 \bigr)$ over $\slq_{\nu,p}^{\ast}\LL_{p}(\omega)$
is given as follows:
\[
 \nbigp^{(a)}_{b}\Bigl(
 \slq_{\nu,p}^{\ast}\LL_p(\omega)_{|(\pi^{\cov}_{\nu,p})^{-1}(\ttt)}
 \Bigr)
=
 \ttU_{\nu,p}^{-[b-a-p\omega \ttt/\gminit^{\lambda}]}
 \nbigr_{\nu,p}\cdot \slq_{\nu,p}^{-1}(\sle_{p,\omega}).
\]
Here, we set $[c]:=\max\{n\in\seisuu\,|\,n\leq c\}$
for any $c\in\real$.
Because it is naturally $\seisuu \tte_2$-equivariant,
we obtain an induced filtered bundles
$\vecLL_{\nu,p}(\omega)$
denoted by 
$\nbigp^{(a)}_{\ast}(\vecLL_{\nu,p}(\omega))=\bigl(
 \nbigp^{(a)}_{\ast}\bigl(\vecLL_{\nu,p}(\omega)_{|\pi_{\nu,p}^{-1}(\ttt)}\bigr)
 \,\big|\,
\ttt\in S^1_{\lambda}
 \bigr)$.

For general $\omega\in\rnum$,
we take set $k_1:=k(p\omega)$, $\ell_1:=\ell(p\omega)$
and $p_1:=p\cdot k_1$.
A filtered bundle
$\nbigp^{(a)}_{\ast}\vecLL_{\nu,p}(\omega)$
over $\vecLL_{\nu,p}(\omega)$
is obtained as the push-forward of
$\nbigp^{(k_1a)}_{\ast}\vecLL_{\nu,p_1}(\omega)$.

\begin{lem}
$\nbigp^{(a)}_{\ast}\vecLL_{\nu,p}(\omega)$
is an object in
$\LFM(\Hhat^{\lambda}_{\nu,p},H^{\lambda}_{\nu,p})$.
We have the natural isomorphism
$\ttG(\nbigp^{(a)}_{\ast}\vecLL_{\nu,p}(\omega))
\simeq
 \LL^{\ttG}(\omega,a)$.
\hfill\qed
\end{lem}

Let $\nbigp_{\ast}\gbigv\in
 \LFM(\Hhat^{\lambda}_{\nu,p},H^{\lambda}_{\nu,p})^{\Par}$.
There exists the slope decomposition
$\gbigv=\bigoplus_{\omega\in\Slope(\gbigv)}\gbigv_{\omega}$,
where each $\gbigv_{\omega}$ has pure slope $\omega$.
We take $p_1\in p\seisuu_{>0}$
such that $p_1\omega\in\seisuu$
for any $\omega\in \Slope(\gbigv)$.
There exists an isomorphism
\begin{equation}
\label{eq;18.8.28.3}
 \sfp_{p,p_1}^{\ast}\gbigv
\simeq
 \bigoplus_{\omega\in\Slope(\gbigv)}
 \vecLL_{\nu,p_1}(\omega)\otimes
\gbigu^{\langle p_1\rangle}_{\omega},
\end{equation}
where 
$\gbigu^{\langle p_1\rangle}_{\omega}$
are Fuchsian.
Then, we have
\begin{equation}
 \label{eq;18.8.28.2}
 \nbigp_{\ast}\bigl(
 \sfp_{p,p_1}^{\ast}\gbigv_{|\pi_{p,\nu}^{-1}(\ttt)}
 \bigr)
\simeq
 \bigoplus_{\omega\in\Slope(\gbigv)}
 \nbigp^{(0)}_{\ast}\bigl(
 \vecLL_{\nu,p_1}(\omega)_{|\pi_{p,\nu}^{-1}(\ttt)}
 \bigr)
\otimes
 \nbigp_{\ast}\bigl(
 \gbigu^{\langle p_1\rangle}_{\omega|\pi_{p,\nu}^{-1}(\ttt)}
 \bigr),
\end{equation}
where
$\nbigp_{\ast}\bigl(
 \gbigu^{\langle p_1\rangle}_{\omega}
 \bigr)$
are isoclinic of pure slope $0$.

\subsubsection{Decomposition and weight filtration
on the associated graded vector spaces}
\label{subsection;18.12.21.20}

Let $\nbigp_{\ast}\gbigv$
be a good filtered bundle on
$(\Hhat^{\lambda}_{p,\nu},H^{\lambda}_{p,\nu})$
with the slope decomposition
$\nbigp_{\ast}\gbigv
=\bigoplus_{\omega\in\rnum} \nbigp_{\ast}\gbigv_{\omega}$.
Let $\gbigv^{\cov}=\bigoplus \gbigv^{\cov}_{\omega}$ denote 
the locally free 
$\nbigo_{\Hhat^{\lambda\cov}_{p,\nu}}
 (\ast H^{\lambda\cov}_{p,\nu})$-module
obtained as the pull back of $\gbigv$.
Let $\nbigp_{\ast}(\gbigv^{\cov})=
 \bigoplus\nbigp_{\ast}\gbigv^{\cov}_{\omega}$
denote the induced filtered bundle over
$\gbigv^{\cov}$.

By the parallel transport along the path connecting
$\ttt_1,\ttt_2\in\real$,
we obtain the isomorphism
\begin{equation}
\label{eq;18.12.21.11}
 \Gr^{\nbigp}_{a}
\bigl(
 \gbigv^{\cov}_{\omega|(\pi^{\cov}_{p,\nu})^{-1}(\ttt_1)}
\bigr)
\simeq
  \Gr^{\nbigp}_{a+p\omega(\ttt_2-\ttt_1)/\gminit^{\lambda}}
\bigl(
 \gbigv^{\cov}_{\omega|(\pi^{\cov}_{p,\nu})^{-1}(\ttt_2)}
\bigr).
\end{equation}

Recall that
$\ttG(\nbigp_{\ast}\gbigv)$
is $(\rnum,\real)$-graded
$\ttG(\nbigp_{\ast}\gbigv)
=\bigoplus_{\omega,a}
 \ttG(\nbigp_{\ast}\gbigv)_{\omega,a}$.
Each $\ttG(\nbigp_{\ast}\gbigv)_{\omega,a}$
is equipped with 
the automorphism $F_{\omega,a}$
and a generalized eigen decomposition
$\ttG(\nbigp_{\ast}\gbigv)_{\omega,a}
=\bigoplus_{\alpha\in\cnum^{\ast}}
\ttG(\nbigp_{\ast}\gbigv)_{\omega,a,\alpha}$.
Moreover, it is equipped with
the nilpotent endomorphism $N_{a,\omega}$
and the weight filtration $W$.
By the construction,
$\ttG(\nbigp_{\ast}\gbigv)_{\omega,a}$
is naturally identified with
$\Gr^{\nbigp}_{a}\bigl(
 \gbigv^{\cov}_{\omega|(\pi^{\cov}_{p,\nu})^{-1}(0)}
 \bigr)$.
Hence, 
each
$\Gr^{\nbigp}_{a}\bigl(
 \gbigv^{\cov}_{\omega|(\pi^{\cov}_{p,\nu})^{-1}(0)}
 \bigr)$
is equipped with the automorphism $F_{\omega,a}$
and the generalized eigen decomposition
$\Gr^{\nbigp}_{a}\bigl(
 \gbigv^{\cov}_{\omega|(\pi^{\cov}_{p,\nu})^{-1}(0)}
 \bigr)
=\bigoplus
 \EE_{\alpha}
\Gr^{\nbigp}_{a}\bigl(
 \gbigv^{\cov}_{\omega|(\pi^{\cov}_{p,\nu})^{-1}(0)}
 \bigr)$.
Moreover, it is equipped with
the nilpotent endomorphism
$N_{\omega,a,\alpha}$
and the weight filtration $W$.

By the isomorphisms (\ref{eq;18.12.21.11}),
each
$\Gr^{\nbigp}_{a+p\omega\ttt/\gminit^{\lambda}}\bigl(
 \gbigv^{\cov}_{\omega|(\pi^{\cov}_{p,\nu})^{-1}(\ttt)}
 \bigr)$
is equipped with the automorphism $F_{\omega,a}$
and the generalized eigen decomposition
$\Gr^{\nbigp}_{a+p\omega\ttt/\gminit^{\lambda}}\bigl(
 \gbigv^{\cov}_{\omega|(\pi^{\cov}_{p,\nu})^{-1}(\ttt)}
 \bigr)
=\bigoplus
 \EE_{\alpha}
\Gr^{\nbigp}_{a+p\omega\ttt/\gminit^{\lambda}}\bigl(
 \gbigv^{\cov}_{\omega|(\pi^{\cov}_{p,\nu})^{-1}(\ttt)}
 \bigr)$.
Moreover, it is equipped with
the nilpotent endomorphism
$N_{\omega,a,\alpha}$
and the weight filtration $W$.

\subsubsection{The associated local systems}

By using the isomorphisms (\ref{eq;18.12.21.11}),
we obtain a local system
$\ttL^{\cov}_{\omega,a}(\nbigp_{\ast}\gbigv)$ on
$H^{\lambda\cov}_{\nu,p}$
by setting
\[
 \ttL^{\cov}_{\omega,a}(\nbigp_{\ast}\gbigv)_{\ttt}:=
 \Gr^{\nbigp}_{a+p\omega\ttt/\gminit^{\lambda}}
 (\gbigv^{\cov}_{\omega|(\pi_{p,\nu}^{\cov})^{-1}(\ttt)}).
\]
We obtain the automorphism $F_{a,\omega}$,
the decomposition
$\ttL^{\cov}_{\omega,a}(\nbigp_{\ast}\gbigv)=\bigoplus
 \EE_{\alpha}\ttL^{\cov}_{\omega,a}(\nbigp_{\ast}\gbigv)$,
the nilpotent endomorphism
$N_{\omega,a}=\bigoplus N_{\omega,a,\alpha}$
and the weight filtration $W$.

The multiplication of $\ttU_{\nu,p}$
induces isomorphisms
$\ttL^{\cov}_{\omega,a}(\nbigp_{\ast}\gbigv)
\simeq
\ttL^{\cov}_{\omega,a-1}(\nbigp_{\ast}\gbigv)$.
We also have the isomorphisms
\[
 \tte_{2}^{\ast}
 \ttL^{\cov}_{\omega,a}(\nbigp_{\ast}\gbigv)
\simeq
 \ttL^{\cov}_{\omega,a+p\omega}(\nbigp_{\ast}\gbigv).
\]
Therefore, 
the multiplication of
$\ttU_{\nu,p}^{\ell(p\omega)}$
induces an isomorphism
\[
 (\tte_2^{\ast})^{k(p\omega)}
 \ttL^{\cov}_{\omega,a}(\nbigp_{\ast}\gbigv)
\simeq
 \ttL^{\cov}_{\omega,a}(\nbigp_{\ast}\gbigv).
\]
Hence, we obtain systems
$\ttL_{\omega,a}(\nbigp_{\ast}\gbigv)$
on 
$S^1_{\lambda,\omega}:=
 H^{\lambda\cov}_{\nu,p}/k(p\omega)\seisuu\tte_2$.

We obtain the monodromy $F_{\omega,a}$ on $\ttL_{\omega,a}$.
We obtain the generalized eigen decomposition
$\ttL_{\omega,a}(\nbigp_{\ast}\gbigv)
=\bigoplus_{\alpha\in\cnum^{\ast}}
 \EE_{\alpha}\ttL_{\omega,a}(\nbigp_{\ast}\gbigv)$
with respect to $F_{\omega,a}$.
Let $N_{\omega,a}=\bigoplus N_{\omega,a,\alpha}$
be the nilpotent endomorphism obtained
as the logarithm of the unipotent part of $F_{\omega,a}$.
Let $W$ be the weight filtration of $N_{\omega,a}$.

Their pull back to $\ttL^{\cov}_{\omega,a}(\nbigp_{\ast}\gbigv)$
are equal to 
the automorphism
the decomposition,
the nilpotent endomorphism
and the weight filtration 
on $\ttL^{\cov}_{\omega,a}(\nbigp_{\ast}\gbigv)$.

\subsubsection{Local filtrations by lattices}
\label{subsection;18.12.17.100}

Let $\nbigp_{\ast}\gbigv$ be a good filtered bundle
on $(\Hhat^{\lambda}_{\nu,p},H^{\lambda}_{\nu,p})$.
Take $\ttt_0\in \real$.
Take $a\in\real$.
Take a small $\epsilon>0$.
Set $I(\nu_p,\ttt_0,\epsilon):=\{\ttt\,|\,|\ttt-\ttt_0|<\epsilon\}
\subset H^{\lambda\cov}_{\nu,p}$
and 
$\Ihat(\nu_p,\ttt_0,\epsilon):=\nuhat_{p}\times I(\nu_p,\ttt_0,\epsilon)$.
We obtain
$\vecP^{(\ttt_0)}_a\gbigv^{\cov}
\subset
 \gbigv^{\cov}_{|\Ihat(\nu_p,\ttt_0,\epsilon)}$
determined by the following
for $\ttt\in I(\nu_p,\ttt_0,\epsilon)$:
\[
 \vecP^{(\ttt_0)}_a\gbigv^{\cov}_{|(\pi^{\cov}_{p,\nu})^{-1}(\ttt)}
=\bigoplus_{\omega\in\Slope(\gbigv)}
 \nbigp_{a+p\omega(\ttt-\ttt_0)/\gminit^{\lambda}}
 \gbigv^{\cov}_{\omega|(\pi^{\cov}_{p,\nu})^{-1}(\ttt)}.
\]
We may naturally regard $I(\nu_p,\ttt_0,\epsilon)\subset H^{\lambda}_{\nu,p}$.
We obtain a filtration
$\vecP^{(\ttt_0)}_{\ast}$
of $\gbigv_{|\Ihat(\nu,\ttt_0,\epsilon)}$.
We obtain the local systems
$\Gr^{\vecP^{(\ttt_0)}}_{a}(\gbigv_{|\Ihat(\nu_p,\ttt_0,\epsilon)})$
on $I(\nu_p,\ttt_0,\epsilon)$.
We have the weight filtration $W$
on $\Gr^{\vecP^{(\ttt_0)}}_a(\gbigv_{|\Ihat(\nu_p,\ttt_0,\epsilon)})$.

\subsection{Good filtered bundles with Dirac type singularity
on $(\nbigmbar^{\lambda}_{p};H^{\lambda}_p,Z)$}

Let 
$\pi^{\cov}_p:\nbigmbar_p^{\lambda\cov}\lrarr \real$
denote the projection
$\pi^{\cov}_p(\ttU_p,\ttt)=\ttt$.
It induces
$\pi_p:\nbigmbar_p^{\lambda}\lrarr 
S^1_{\lambda}$.
The fibers $(\pi^{\cov}_p)^{-1}(\ttt)\subset \nbigmbar_p^{\lambda\cov}$ 
$(\ttt\in\real)$
are identified with $\proj^1$.
For each $\ttt\in S^1_{\lambda}$,
by fixing its lift to $\real$,
we obtain the isomorphism
$\pi_p^{-1}(\ttt)\simeq\proj^1$.

Let $Z\subset\nbigm^{\lambda}_p$ be a finite subset.
Let $\gbigv$ be a locally free 
$\nbigo_{\nbigmbarlambda_p\setminus Z}(\ast H^{\lambda}_p)$-module.
A filtered bundle over $\gbigv$
is a family of filtered bundles
$\nbigp_{\ast}(\gbigv)=
 \bigl(
 \nbigp_{\ast}(\gbigv_{|\pi_p^{-1}(\ttt)})\,\big|\,
 \ttt\in S^1_{\lambda}
 \bigr)$
over $\gbigv_{|\pi_p^{-1}(\ttt)}$.
It induces filtered bundles
$\nbigp_{\ast}(\gbigv_{|\Hhat^{\lambda}_{\nu,p}})$
$(\nu=0,\infty)$
over 
$\gbigv_{|\Hhat^{\lambda}_{\nu,p}}$.
\begin{df}
$\nbigp_{\ast}(\gbigv)$ is called good
if the induced filtered bundles
$\nbigp_{\ast}(\gbigv_{|\Hhat^{\lambda}_{\nu,p}})$
are good.
If moreover each point of $Z$ is Dirac type singularity of $\gbigv$,
we say that 
$\nbigp_{\ast}\gbigv$ is a good filtered bundle
with Dirac type singularity 
over $(\nbigmbar^{\lambda}_p;H^{\lambda}_p,Z)$.
\hfill\qed
\end{df}

\subsubsection{Degree and stability condition}

Let $\nbigp_{\ast}(\gbigv)$ be a good filtered bundle
with Dirac type singularity on
$(\nbigmbar^{\lambda}_p;H^{\lambda}_p,Z)$.
We define the degree of
$\nbigp_{\ast}(\gbigv)$ as follows:
\[
 \deg(\nbigp_{\ast}\gbigv):=
\int_{S^1_{\lambda}}
 \deg\bigl(
 \nbigp_{\ast}(\gbigv_{|\pi^{-1}_{p}(\ttt)})
 \bigr)\,d\ttt.
\]

Let $\gbigv_1\subset\gbigv$ be
an $\nbigo_{\nbigmbar^{\lambda}_p}(\ast H^{\lambda}_p)$-submodule.
Then, it is also locally free,
and each point of $Z$ is with Dirac type singularity.
The induced filtered bundle
$\nbigp_{\ast}(\gbigv_1)$ is good.
We say that 
$\nbigp_{\ast}(\gbigv)$ is stable
if 
\[
 \deg(\nbigp_{\ast}\gbigv_1)/\rank(\gbigv_1)
<
 \deg(\nbigp_{\ast}\gbigv)/\rank(\gbigv)
\]
for any saturated submodules $\gbigv_1$ of $\gbigv$
such that
$\gbigv_1\neq 0,\gbigv$.
We say that 
$\nbigp_{\ast}(\gbigv)$ is semistable
if 
\[
 \deg(\nbigp_{\ast}\gbigv_1)/\rank(\gbigv_1)
\leq
 \deg(\nbigp_{\ast}\gbigv)/\rank(\gbigv)
\]
for any non-trivial submodules $\gbigv_1$ of $\gbigv$.
We say that 
$\nbigp_{\ast}(\gbigv)$ is polystable
if it is semistable and a direct sum of stable ones.

\subsection{Good filtered bundles on neighbourhoods
of $H^{\lambda}_{\nu,p}$}
\label{subsection;18.9.1.1}

For $\nu=0,\infty$,
let $\nbigubar^{\lambda}_{\nu,p}$
be a neighbourhood of
$H^{\lambda}_{\nu,p}$.
We set
$\nbigu^{\lambda}_{\nu,p}:=
 \nbigubar^{\lambda}_{\nu,p}
\setminus H^{\lambda}_{\nu,p}$.
The induced map
$\nbigubar^{\lambda}_{\nu,p}\lrarr
 S^1_{\lambda}$
is denoted by $\pi_p$.
Let $\gbigv$ be a locally free
$\nbigo_{\nbigubar^{\lambda}_{\nu,p}}(\ast
H^{\lambda}_{\nu,p})$-module.
A filtered bundle over $\gbigv$
be a family of filtered bundles
$\nbigp_{\ast}(\gbigv_{|\pi_p^{-1}(\ttt)})$
over $\gbigv_{|\pi_p^{-1}(\ttt)}$
$(\ttt\in S^1_{\lambda})$.
The tuple
$(\nbigp_{\ast}(\gbigv_{|\pi_p^{-1}(\ttt)})\,|\,\ttt\in S^1_{\lambda})$
is denoted by
$\nbigp_{\ast}(\gbigv)$.
A filtered bundle $\nbigp_{\ast}\gbigv$ over $\gbigv$
is called good if the induced filtered bundle
over $\gbigv_{|\Hhat_{\nu,p}^{\lambda}}$ is good.

\subsubsection{Filtrations by local lattices}

For $\ttt_0\in S^1_{\lambda}$,
we set
$I(\ttt_0,\epsilon):=\{\ttt\,|\,|\ttt-\ttt_0|<\epsilon\}$.
For $a\in\real$,
we obtain the lattice
$\vecP^{(\ttt_0)}_a(\gbigv_{|\pi_p^{-1}(I(\ttt_0,\epsilon))})
\subset
 \gbigv_{|\pi_p^{-1}(I(\ttt_0,\epsilon))}$
from 
$\vecP^{(\ttt_0)}_a(\gbigv_{|\Ihat(\nu_p,\ttt_0,\epsilon)})$.
Thus, we have the filtration
$\vecP^{(\ttt_0)}_{\ast}(\gbigv_{|\pi_p^{-1}(I(\ttt_0,\epsilon))})$.
The induced local system
$\Gr^{\vecP^{(\ttt_0)}}_a(\gbigv_{|\pi_p^{-1}(I(\ttt_0,\epsilon))})$
on $I(\ttt_0,\epsilon)$
is equipped with the weight filtration $W$.
We also have the decomposition
\[
\Gr^{\vecP^{(\ttt_0)}}\bigl(
 \gbigv_{|\pi_p^{-1}(I(\ttt_0,\epsilon))}\bigr)
=\bigoplus
\Gr^{\vecP^{(\ttt_0)}}\bigl(
 \gbigv_{\omega|\pi_p^{-1}(I(\ttt_0,\epsilon))}\bigr)
\]
induced by
$\gbigv_{|\Hhat^{\lambda}_{\nu,p}}
=\bigoplus_{\omega}
 \gbigv_{\omega}$.
The decomposition and the filtration $W$ are compatible.

\subsubsection{Compatible frame}

We continue to use the notation in \S\ref{subsection;18.9.1.1}.
Set $r:=\rank(\gbigv)$.

\begin{df}
Let $\vecv=(v_i\,|\,i=1,\ldots,r)$ be a frame of 
$\vecP^{(\ttt_0)}_a\gbigv$ 
on a neighbourhood of $\pi_p^{-1}(\ttt_0)$.
We say that $\vecv$ is compatible with 
the filtration $\vecP^{(\ttt_0)}_{\ast}\gbigv$
and the slope decomposition
if there exists a decomposition
 $\{1,\ldots,r\}=
 \coprod_{\omega\in\Slope(\gbigv)}
 \coprod_{a-1<b\leq a} I_{\omega,b}$
 such that 
 $(v_i\,|\,i\in I_{\omega,b})$ induces a frame of
 $\Gr^{\vecP^{(\ttt_0)}}_{b}(\gbigv_{\omega})$
for $a-1<b\leq a$.

We say that $\vecv$ is compatible with
the slope decomposition,
the filtration $\vecP^{(\ttt_0)}_{\ast}\gbigv$ 
and the filtration $W$
if there exists a decomposition
 $\{1,\ldots,r\}=
 \coprod_{\omega\in\Slope(\gbigv)}
 \coprod_{a-1<b\leq a}
 \coprod_{k\in\seisuu} I_{\omega,b,k}$
 such that 
 $(v_i\,|\,i\in I_{\omega,b,k})$ induces a frame of
 $\Gr^W_k\Gr^{\vecP^{(\ttt_0)}}_b(\gbigv_{\omega})$.
\hfill\qed
\end{df}

Take a local frame $\vecv$ of 
$\vecP_a^{(\ttt_0)}\gbigv$
compatible with 
the slope decomposition
and the filtration $\vecP^{(\ttt_0)}_{\ast}$.
We set $b(v_i):=b$
and $\omega(v_i):=\omega$
if $i\in I_{\omega,b}$.
If moreover $\vecv$ is compatible with $W$,
we also set $k(v_i):=k$
if $i\in I_{\omega,b,k}$.

\subsubsection{Adaptedness and norm estimate}

Let $\nbigp_{\ast}\gbigv$ be a good filtered bundle 
over $\gbigv$.
Let $V$ be the mini-holomorphic bundle on $\nbigu^{\lambda}_{\nu,p}$
obtained as the restriction
$\gbigv_{|\nbigu^{\lambda}_{\nu,p}}$.
Let $P$ be a point of $H^{\lambda}_{\nu,p}$.
Let $U_P$ be a neighbourhood of $P$
in $\nbigubar^{\lambda}_{\nu,p}$.
Let $\vecv$ be a frame of  $\vecP^{(\ttt_0)}_a\gbigv$
on $U_P$
compatible with the slope decomposition
and the filtration $\vecP^{(\ttt_0)}_{\ast}\gbigv$.
Let $h_{P,\vecv}$ be the Hermitian metric 
of $V_{|U_P\setminus H^{\lambda}_{\nu,p}}$
determined by
\[
h_{P,\vecv}(v_i,v_j):=
 \left\{
\begin{array}{ll}
 |\ttU_{\nu,p}|^{-2b(v_i)-2p\omega(v_i)(\ttt-\ttt_0)/\gminit^{\lambda}}
 & (i=j)\\
 0 & (i\neq j).
\end{array}
 \right.
\]
If moreover $\vecv$ is compatible with the filtration $W$,
then let $\htilde_{P,\vecv}$ be the Hermitian metric
of $V_{|U_P\setminus H^{\lambda}_{\nu,p}}$
determined by 
\[
\htilde_{P,\vecv}(v_i,v_j):=
 \left\{
 \begin{array}{ll}
 |\ttU_{\nu,p}|^{-2b(v_i)-2p\omega(v_i)(\ttt-\ttt_0)/\gminit^{\lambda}}
 (-\log|\ttU_{\nu,p}|)^{k(v_i)}
 & (i=j)\\
 0 & (i\neq j). 
 \end{array}
 \right.
\]

The following is easy to see.
\begin{lem}
Let $\vecv$ and $\vecv'$ be frames of $\vecP^{(\ttt_0)}_a(\gbigv)$
on $U_P$
compatible with
the slope decomposition
and the filtration $\vecP^{(\ttt_0)}_{\ast}(\gbigv)$.
Take a relative compact neighbourhood $U_P'$ of $P$ in $U_P$.
Then,
$h_{P,\vecv}$ and $h_{P,\vecv'}$ are mutually bounded
on $U_P'\setminus H^{\lambda}_{\nu,p}$.
If moreover both $\vecv$ and $\vecv'$ are compatible with $W$,
then $\htilde_{P,\vecv}$ and $\htilde_{P,\vecv'}$
are mutually bounded 
on $U_P'\setminus H^{\lambda}_{\nu,p}$.
\hfill\qed
\end{lem}

\begin{df}
A Hermitian metric $h$ of $V$
is called 
adapted to $\nbigp_{\ast}\gbigv$ around $P$
if the following holds.
\begin{itemize}
\item
Let $\vecv$ be a frame of $\vecP^{(\ttt_0)}_{a}\gbigv$
on a neighbourhood $U_P$ of $P$
compatible with the slope decomposition
and the filtration $\vecP^{(\ttt_0)}_{\ast}\gbigv$.
Then, for any smaller neighbourhood $U_P'\subset U_P$
and for any $\epsilon$,
there exists $C_{\epsilon}>1$ such that 
\[
 C_{\epsilon}^{-1}|\ttU_{\nu,p}|^{\epsilon}h_{P,\vecv}
\leq
 h
\leq
 C_{\epsilon}|\ttU_{\nu,p}|^{-\epsilon}h_{P,\vecv}
\]
on $U_P'\setminus H^{\lambda}_{\nu,p}$.
\end{itemize}

We say that $\nbigp_{\ast}\gbigv$ is adapted to $h$
if it is adapted to $h$ around any point of $H^{\lambda}_{\nu,p}$.
\hfill\qed
\end{df}

\begin{df}
Let $h$ be a Hermitian metric of $V$.
We say that 
the norm estimate holds for 
$\nbigp_{\ast}\gbigv$ and $h$
around $P$,
if the following holds.
\begin{itemize}
\item
Let $\vecv$ be a frame of $\vecP^{(\ttt_0)}_{a}\gbigv$
on a neighbourhood $U_P$ of $P$
compatible with the slope decomposition,
the filtration $\vecP^{(\ttt_0)}_{\ast}\gbigv$
and $W$.
Then, for any smaller neighbourhood $U_P'\subset U_P$
there exists $C>1$ such that 
\[
 C^{-1}
 \htilde_{P,\vecv}
\leq
 h
\leq
 C\htilde_{P,\vecv}
\]
on $U_P'\setminus H^{\lambda}_{\nu,p}$.
\end{itemize}
We say that 
the norm estimate holds for 
$\nbigp_{\ast}\gbigv$ and $h$
if the norm estimate holds 
around any point of $H^{\lambda}_{\nu,p}$.
\hfill\qed
\end{df}

\subsection{Approximation}

We use the notation in \S\ref{subsection;18.9.1.1}.
Let $\nbigc^{\infty}_{\nbigubar^{\lambda}_{\nu,p}}$
denote the sheaf of $C^{\infty}$-functions 
on $\nbigc^{\infty}_{\nbigubar^{\lambda}_{\nu,p}}$.
For good filtered bundles 
$\nbigp_{\ast}\gbigv^{(i)}$ $(i=1,2)$
over
$(\nbigubar^{\lambda}_{\nu,p},H^{\lambda}_{\nu,p})$,
a $C^{\infty}$-isomorphism of
$f:\nbigp_{\ast}\gbigv^{(1)}\simeq
 \nbigp_{\ast}\gbigv^{(2)}$
means an isomorphism
$\gbigv^{(1)}\otimes\nbigc^{\infty}_{\nbigubar^{\lambda}_{\nu,p}}
\simeq
 \gbigv^{(2)}\otimes\nbigc^{\infty}_{\nbigubar^{\lambda}_{\nu,p}}$
such that 
the restriction to
$\pi^{-1}_{p}(\ttt)$ 
preserve the induced filtrations.

The following lemma is clear.
\begin{lem}
Let $\nbigp_{\ast}\gbigv^{(i)}$ $(i=1,2)$
be good filtered bundles 
$(\nbigubar^{\lambda}_{\nu,p},H^{\lambda}_{\nu,p})$.
If there exists an isomorphism
$\fhat:
 \nbigp_{\ast}\gbigv^{(1)}_{|\Hhat^{\lambda}_{\nu,p}}
\simeq
 \nbigp_{\ast}\gbigv^{(2)}_{|\Hhat^{\lambda}_{\nu,p}}$,
then there exists an isomorphism
\[
 f_{C^{\infty}}:
\nbigp_{\ast}\gbigv^{(1)}\otimes
 \nbigc^{\infty}_{\nbigubar^{\lambda}_{\nu,p}}
\simeq
\nbigp_{\ast}\gbigv^{(2)}\otimes
 \nbigc^{\infty}_{\nbigubar^{\lambda}_{\nu,p}}
\]
whose restriction to 
$\Hhat^{\lambda}_{\nu,p}$ is equal to $\Fhat$.
\hfill\qed
\end{lem}

\begin{lem}
Let $\nbigp_{\ast}\gbigv^{(i)}$ $(i=1,2)$
be good filtered bundles 
$(\nbigubar^{\lambda}_{\nu,p},H^{\lambda}_{\nu,p})$.
If there exists an isomorphism
$f^{\ttG}:\ttG(\nbigp_{\ast}\gbigv^{(1)})\simeq
 \ttG(\nbigp_{\ast}\gbigv^{(2)})$,
there exists an isomorphism
\[
 f:
\nbigp_{\ast}\gbigv^{(1)}\otimes
 \nbigc^{\infty}_{\nbigubar^{\lambda}_{\nu,p}}
\simeq
\nbigp_{\ast}\gbigv^{(2)}\otimes
 \nbigc^{\infty}_{\nbigubar^{\lambda}_{\nu,p}}
\]
such that the following holds.
\begin{itemize}
\item
 For each $\ttt\in S^1_{\lambda}$,
 the restriction of $f$ to $\pi_p^{-1}(\ttt)$
 is holomorphic and preserves 
 the filtrations.
\item
 The induced morphism
 $\Gr^{\nbigp}_a(f_{|\pi^{-1}_p(\ttt)})$
 preserves the decomposition
$\Gr^{\nbigp}_a(\gbigv^{(i)})
=\bigoplus_{\omega}\Gr^{\nbigp}_a(\gbigvhat^{(i)}_{\omega})$
induced by the slope decomposition
$\gbigv^{(i)}_{|\Hhat^{\lambda}_{\nu,p}}
=\bigoplus
 \gbigvhat^{(i)}_{\omega}$.
As a result,
we obtain the decomposition
 $\Gr^{\nbigp}_a(f_{|\pi^{-1}_p(\ttt)})
=\bigoplus_{\omega}
\Gr^{\nbigp}_a(f_{|\pi^{-1}_p(\ttt)})_{\omega}$.
\item
If $\ttt_1-\ttt_2$ is small,
$\Gr^{\nbigp}_{a+p\omega(\ttt_2-\ttt_1)/\gminit^{\lambda}}
 (f_{|\pi^{-1}_p(\ttt_2)})_{\omega}$
and 
$\Gr^{\nbigp}_{a}
 (f_{|\pi^{-1}_p(\ttt_1)})_{\omega}$
are equal
under the natural isomorphism
\begin{equation}
 \label{eq;19.1.18.1}
 \Gr^{\nbigp}_{a+p\omega(\ttt_2-\ttt_1)/\gminit^{\lambda}}
 (\gbigvhat_{\omega|\pi_p^{-1}(\ttt_2)}^{(i)})
\simeq
 \Gr^{\nbigp}_a(\gbigvhat_{\omega|\pi^{-1}_p(\ttt_1)}^{(i)}).
\end{equation}
\end{itemize}
\end{lem}
\pf
The isomorphism $f^{\ttG}$
induces an isomorphism
$f^{\ttG}_{\omega,a,\ttt}:
 \Gr^{\nbigp}_{a}(\gbigv^{(i)}_{\omega|\pi_p^{-1}(\ttt)})$
for any $a\in\real$, $\omega\in\rnum$
and $\ttt\in S^1_{\lambda}$
satisfying
$f^{\ttG}_{\omega,a,\ttt_1}
=f^{\ttG}_{\omega,a+p\omega(\ttt_2-\ttt_1),\ttt_2}$
under (\ref{eq;19.1.18.1}).
For any $\ttt_0\in S^1_{\lambda}$,
we take a small neighbourhood $I(\ttt_0)$
in $S^1_{\lambda}$.
We can take a holomorphic isomorphism
$f_{I(\ttt_0)}:
 \gbigv^{(1)}_{|\pi_{p}^{-1}(I(\ttt_0))}
 \simeq
 \gbigv^{(2)}_{|\pi_p^{-1}(I(\ttt_0))}$
such that the following holds:
\begin{itemize}
\item
For each $\ttt\in I(\ttt_0)$,
the restriction to
$\pi_p^{-1}(\ttt)$
preserves the filtrations.
\item
The induced isomorphism
$\Gr^{\nbigp}_a(\gbigvhat^{(1)}_{\omega|\pi_p^{-1}(\ttt)})
\simeq
 \Gr^{\nbigp}_a(\gbigvhat^{(2)}_{\omega|\pi_p^{-1}(\ttt)})$
is equal to $f^{\ttG}_{\omega,a,\ttt}$.
\end{itemize}
We take a finite covering 
$S^1_{\lambda}=\bigcup_{i=1}^N I(\ttt_0^{(i)})$
and a partition of unity
$\{\chi_i\}$ subordinate to the covering.
We construct 
a $C^{\infty}$-isomorphism $f$
as $f=\sum_{i=1}^N \chi_if_{I(\ttt_0^{(i)})}$.
Then, $f$ satisfies the conditions.
\hfill\qed

\section{Basic examples of doubly periodic monopoles}

\subsection{Examples (1)}
\label{subsection;18.12.14.1}

\subsubsection{Construction}
\label{subsection;18.11.24.30}

On $A^{0}$, we have the mini-complex coordinate system $(z,y)$,
where $y:=\Image(w)$.
Let $\underline{\cnum}\cdot \gminie$ denote the product
line bundle on $A^{0}$ with a global frame $\gminie$.
Let $h$ be the metric given by $h(\gminie,\gminie)=1$.
We consider the $\seisuu \tte_1$-action
on $\underline{\cnum}\,\gminie$
given by
$\tte_1^{\ast}(\gminie)=\gminie$.
It induces an action of $\seisuu (m\tte_1)$
for any $m\in\seisuu_{>0}$ as the restriction.

Take a positive integer $p$ and a rational number 
$\omega\in \frac{1}{p}\seisuu$.
We have the expression
$\omega=\ell(\omega)/k(\omega)$,
where $k(\omega)\in\seisuu_{>0}$,
$\ell(\omega)\in\seisuu$ and $\gcd(k(\omega),\ell(\omega))=1$.
We set
\[
 \alpha(\omega):=
 \frac{2\pi\omega}{\Vol(\Gamma)}.
\]
We define the $\seisuu \tte_2$-action
on $\underline{\cnum}\,\gminie$
by
\[
 \tte_2^{\ast}(\gminie)
\longmapsto
\gminie\cdot
\exp\Bigl(
-\sqrt{-1}\Vol(\Gamma)
 \alpha(\omega)|\mu_1|^{-2}
 \Re(\mubar_1z)
 \Bigr)
=
\gminie\cdot
\exp\Bigl(
-2\pi\sqrt{-1}\omega|\mu_1|^{-2}
 \Re(\mubar_1z)
 \Bigr).
\]
\begin{lem}
The actions of $\seisuu (k(\omega)\tte_1)$
and $\seisuu\tte_2$ are commutative,
i.e.,
the action of $\seisuu (k(\omega)\tte_1)\oplus\seisuu\tte_2$
on $\underline{\cnum}\,\gminie$ is well defined.
\end{lem}
\pf
It follows from
$\exp\Bigl(
 -2\pi\sqrt{-1}\omega|\mu_1|^{-2}\Re(\mubar_1\cdot k(\omega)\mu_1)
 \Bigr)
=\exp\Bigl(
 -2\pi\sqrt{-1}k(\omega)\omega
 \Bigr)=1$.
\hfill\qed

\vspace{.1in}

Let $\phi_{p,\omega}$ be the Higgs field
given as $\phi_{p,\omega}=\sqrt{-1}\alpha(\omega) y$.
We define the connection $\nabla_{p,\omega}$ by
\[
 \nabla_{p,\omega}\gminie
=\gminie
 \Bigl(-\frac{\alpha(\omega)}{4}\Bigr)
 |\mu_1|^{-2}
 \bigl(\mubar_1z-\mu_1\zbar\bigr)
 (\mubar_1\,dz+\mu_1\,d\zbar).
\]
\begin{lem}
The Bogomolny equation 
$F(\nabla_{p,\omega})=
\ast\nabla_{p,\omega}\phi_{p,\omega}$
is satisfied.
\end{lem}
\pf
We have
$\nabla_{p,\omega}\phi_{p,\omega}=\sqrt{-1}\alpha(\omega) dy$,
and hence
$\ast\nabla_{p,\omega}\phi_{p,\omega}
=-\frac{1}{2}\alpha(\omega)\,dz\,d\zbar$.
We also have
\[
 F(\nabla_{p,\omega})
=-\frac{\alpha(\omega)}{4}
 |\mu_1|^{-2}
 \bigl(
 \mubar_1\,dz-\mu_1\,d\zbar
 \bigr)
 (\mu_1\,d\zbar+\mubar_1\,dz)
=-\frac{\alpha(\omega)}{2}\,dz\,d\zbar.
\]
Hence, the Bogomolny equation is satisfied.
\hfill\qed

\begin{lem}
$(k(\omega)\tte_1)^{\ast}\nabla_{p,\omega}=\nabla_{p,\omega}$
and 
$\tte_2^{\ast}\nabla_{p,\omega}=\nabla_{p,\omega}$.
\end{lem}
\pf
The claim $(k(\omega)\tte_1)^{\ast}\nabla_{p,\omega}=\nabla_{p,\omega}$
is clear.
Because
\begin{multline}
 \tte_2^{\ast}(\nabla_{p,\omega})
 \tte_2^{\ast}(\gminie)
=\tte_2^{\ast}(\gminie)
 \cdot
 \Bigl(
 -\frac{\alpha(\omega)}{4}
 \Bigr)
|\mu_1|^{-2}
 \Bigl(
 (\mubar_1z-\mu_1\zbar)
+(\mubar_1\mu_2-\mu_1\mubar_2)
 \Bigr)
(\mubar_1\,dz+\mu_1\,d\zbar)
 \\
=\tte_2^{\ast}(\gminie)
 \cdot
 \Bigl(
 -\frac{\alpha(\omega)}{4}
 \Bigr)
|\mu_1|^{-2}
 \Bigl(
 (\mubar_1z-\mu_1\zbar)
+2\sqrt{-1}\Vol(\Gamma)
 \Bigr)
(\mubar_1\,dz+\mu_1\,d\zbar)
\end{multline}
we obtain $\tte_2^{\ast}\nabla_{p,\omega}=\nabla_{p,\omega}$.
\hfill\qed

\vspace{.1in}

The monopole
$(\underline{\cnum}\gminie,h,\nabla_{p,\omega},\phi_{p,\omega})$
on $A^0$
is denoted by $\sfL_p(\omega)$.
Because it is equivariant with respect to 
$\seisuu k(\omega)\tte_1\oplus\seisuu \tte_2$,
we obtain 
a monopole
$\vecsfL^{\cov}_p(\omega)$ on $\nbigm_p^{0\cov}$,
and a monopole $\vecsfL_p(\omega)$ on $\nbigm_p^0$.
Moreover,
the monopoles are equivariant with respect to the 
$\bigl(k(\omega)\seisuu/p\seisuu\bigr)\tte_1$-action.

Let $\nbigl_p^{\lambda\cov}(\omega)$
be the mini-holomorphic bundle on $\nbigm_p^{\lambda\cov}$
underlying $\vecsfL_p^{\cov}(\omega)$,
which is naturally $\seisuu \tte_2$-equivariant.
Let $\nbigl^{\lambda}_p(\omega)$
be the mini-holomorphic bundle on $\nbigm^{\lambda}_p$
underlying $\vecsfL_p(\omega)$,
which is obtained as the descent of 
$\nbigl^{\lambda\cov}_p(\omega)$.
The mini-holomorphic bundles are equivariant
with respect to 
the $\bigl( k(\omega)\seisuu/p\seisuu\bigr)\tte_1$-action.

\subsubsection{Corresponding instantons on $X$}

Let $\sfLtilde_p(\omega)=\bigl(
 \underline{\cnum}\cdot\gminietilde,
 \htilde,\nablatilde_{p,\omega} \bigr)$
denote the $\real \tte_0\oplus 
 \seisuu k(\omega)\tte_1\oplus \seisuu \tte_2$-equivariant
instanton on $X$ corresponding to $\sfL_p(\omega)$.
We obtain $\htilde(\gminietilde,\gminietilde)=1$
and 
\[
 \nablatilde_{p,\omega}\gminietilde
=\gminietilde
\Bigl(
 -\frac{\alpha(\omega)}{4}
\Bigr)
\Bigl(
|\mu_1|^{-2}
 (\mubar_1z-\mu_1\zbar)
 (\mubar_1\,dz+\mu_1d\zbar)
-(w-\wbar)(dw+d\wbar)
 \Bigr).
\]

Let $(\nbigltilde^{\lambda}_p(\omega),\delbar^{\lambda})$
be the underlying holomorphic vector bundle
on $X^{\lambda}$,
which is equivariant with respect to
the $\real \tte_0\oplus\seisuu k(\omega)\tte_1
 \oplus\seisuu \tte_2$-action.

\begin{lem}
The following holds:
\begin{multline}
 \delbar^{\lambda}\gminietilde=\gminietilde
 \frac{\alpha(\omega)}{4}
 \frac{1}{(1+|\lambda|^2)^2}
 \Bigl(
 -(1+|\lambda|^2)\xi\,d\xibar
+(1+|\lambda|^2)\eta\,d\etabar
+(\mubar_1^2|\mu_1|^{-2}\lambda -\lambdabar)\xi\,d\etabar
 \\
+(-\mu_1^2|\mu_1|^{-2}\lambdabar+\lambda)
 \eta d\xibar
+(\mu_1^2|\mu_1|^{-2}+\lambda^2)\xibar d\xibar
-(\lambda^2\mubar_1^2|\mu_1|^{-2}+1)\etabar\,d\etabar
 \Bigr).
\end{multline}
\end{lem}
\pf
In the proof,
$\alpha(\omega)$ is denoted by $\alpha$.
Because
\[
 z=\frac{1}{1+|\lambda|^2}(\xi-\lambda\etabar),
\quad
 w=\frac{1}{1+|\lambda|^2}(\eta+\lambda\xibar),
\]
the following holds:
\[
 \mubar_1z-\mu_1\zbar
=\frac{1}{1+|\lambda|^2}
 \bigl(\mubar_1\xi-\lambda\mubar_1\etabar
-\mu_1\xibar+\mu_1\lambdabar\eta\bigr),
\]
\[
 \mubar_1\,dz+\mu_1d\zbar
=\frac{1}{1+|\lambda|^2}
 \bigl(
 \mubar_1d\xi
-\mubar_1\lambda d\etabar
+\mu_1\,d\xibar
-\lambdabar\mu_1d\eta
 \bigr).
\]
Hence,
we obtain
\begin{multline}
 (\mubar_1z-\mu_1\zbar)(\mubar_1\,dz+\mu_1d\zbar)
=\frac{1}{(1+|\lambda|^2)^2}
 \Bigl(
 \mubar_1^2\xi d\xi
+|\mu_1|^2\xi d\xibar
-\mubar_1^2\lambda\xi d\etabar
-|\mu_1|^2\lambdabar \xi d\eta
 \\
-|\mu_1|^2\xibar\,d\xi
-\mu_1^2\xibar d\xibar 
+|\mu_1|^2\lambda \xibar d\etabar
+\lambdabar\mu_1^2\xibar d\eta
 \\
-\lambda\mubar_1^2\etabar d\xi
-\lambda|\mu_1|^2\etabar d\xibar
+\lambda^2\mubar_1^2\etabar d\etabar
+|\lambda|^2|\mu_1|^2\etabar d\eta
 \\
+|\mu_1|^2\lambdabar \eta d\xi
+\mu_1^2\lambdabar \eta d\xibar
-|\mu_1|^2|\lambda|^2 \eta d\etabar
-\mu_1^2\lambdabar^2\eta d\eta
 \Bigr).
\end{multline}

Note that the following also holds:
\begin{multline}
 (w-\wbar)(dw+d\wbar)
=\frac{1}{(1+|\lambda|^2)^2}
 \Bigl(
 \eta d\eta
+\lambda \eta d\xibar
+\eta d\etabar
+\eta\lambdabar d\xi\\
+\lambda \xibar d\eta+\lambda^2\xibar d\xibar
+\lambda \xibar d\etabar +|\lambda|^2 \xibar d\xi
-\etabar d\eta-\lambda\etabar d\xibar
-\etabar d\etabar -\lambdabar \etabar d\xi\\
-\lambdabar \xi d\eta
-|\lambda|^2\xi d\xibar
-\lambdabar \xi d\etabar
-\lambdabar^2\xi d\xi
 \Bigr).
\end{multline}
Hence, we obtain
\begin{multline}
 \Bigl(
 -\frac{\alpha}{4}|\mu_1|^{-2}
 \bigl(\mubar_1 z-\mu_1\zbar\bigr)
 (\mubar_1\,dz+\mu_1d\zbar)
+\frac{\alpha}{4}(w-\wbar)(dw+d\wbar)
 \Bigr)^{0,1}
=\\
 \frac{\alpha}{4}
 \frac{1}{(1+|\lambda|^2)^2}
 \Bigl(
 -(1+|\lambda|^2)\xi\,d\xibar
+(1+|\lambda|^2)\eta\,d\etabar
+(\mubar_1^2|\mu_1|^{-2}\lambda -\lambdabar)\xi\,d\etabar
 \\
+(-\mu_1^2|\mu_1|^{-2}\lambdabar+\lambda)
 \eta d\xibar
+(\mu_1^2|\mu_1|^{-2}+\lambda^2)\xibar d\xibar
-(\lambda^2\mubar_1^2|\mu_1|^{-2}+1)\etabar\,d\etabar
 \Bigr).
\end{multline}
Thus, we obtain the claim of the lemma.
\hfill\qed

\subsubsection{Holomorphic frame of $\nbigltilde^{\lambda}_p(\omega)$}

We consider the following holomorphic frame
of $\nbigltilde^{\lambda}_p(\omega)$ on $X^{\lambda}$:
\begin{multline}
 \sfvtilde_{p,\omega}^{\lambda}:=
 \gminietilde\exp\Bigl(
 \frac{\alpha(\omega)}{4}
 \frac{1}{(1+|\lambda|^2)^2}
 \Bigl(
 (1+|\lambda|^2)\xi\xibar
-(1+|\lambda|^2)\eta\etabar
 \\
-(\mubar_1^2|\mu_1|^{-2}\lambda-\lambdabar)\xi\etabar
-(-\mu_1^2|\mu_1|^{-2}\lambdabar+\lambda)\eta\xibar
-(\mu_1^2|\mu_1|^{-2}+\lambda^2)\frac{1}{2}\xibar^2
+(\lambda^2\mubar_1^2|\mu_1|^{-2}+1)\frac{1}{2}\etabar^2
 \\
-(-\mubar_1^2|\mu_1|^{-2}+\lambdabar^2)\frac{1}{2}\xi^2
+(-\lambdabar^2\mu_1^2|\mu_1|^{-2}+1)\frac{1}{2}\eta^2\\
-(\xi-\lambda\eta)2\lambdabar\eta
-\frac{1}{\mu_1+\lambda^2\mubar_1}
 \bigl(2|\lambda|^2\mubar_1-\lambdabar^2\mu_1+\mubar_1\bigr)
 (\xi-\lambda\eta)^2
 \Bigr)
 \Bigr).
\end{multline}

\begin{lem}
\label{lem;18.8.20.10}
We have
$(\tte_0)^{\ast}\sfvtilde_{p,\omega}^{\lambda}
=\sfvtilde^{\lambda}_{p,\omega}$
and
$(k(\omega)\tte_1)^{\ast}\sfvtilde_{p,\omega}^{\lambda}
=\sfvtilde_{p,\omega}^{\lambda}$.
We also have
\begin{equation}
\label{eq;18.8.20.1}
 \tte_2^{\ast}\sfvtilde_{p,\omega}^{\lambda}
=\sfvtilde_{p,\omega}^{\lambda}
 \cdot
 \ttU_p^{-p\omega}
 \exp\Bigl(
 \frac{\alpha(\omega)}{4}
 (\mu_1+\lambda^2\mubar_1)^{-1}
 (\mu_1-\lambda^2\mubar_1)
 \cdot|\mu_1|^{-2}
 \bigl(
 |\mu_1|^2|\mu_2|^2
-\mu_2^2\mubar_1^2/2
-\mubar_2^2\mu_1^2/2
 \bigr)
 \Bigr).
\end{equation}
\end{lem}
\pf
We can check
$\tte_0^{\ast}\sfvtilde^{\lambda}_{p,\omega}
=\sfvtilde^{\lambda}_{p,\omega}$
and 
$(k(\omega)\tte_1)^{\ast}\sfvtilde^{\lambda}_{p,\omega}
=\sfvtilde^{\lambda}_{p,\omega}$
by direct computations.
We give an indication
to check the formula (\ref{eq;18.8.20.1}).
We have
\[
 \tte_2^{\ast}\sfvtilde^{\lambda}_{p,\omega}
=\sfvtilde^{\lambda}_{p,\omega}
 \exp\Bigl(
 \frac{-\alpha(\omega)}{4}|\mu_1|^{-2}
 (\mubar_1\mu_2-\mu_1\mubar_2)
 (\mubar_1z+\mu_1\zbar)
 \Bigr)
\exp\Bigl(
 \frac{\alpha(\omega)}{4(1+|\lambda|^2)^2}G
 \Bigr),
\]
where
\begin{multline}
G=(1+|\lambda|^2)
\bigl(
(\xi+\mu_2)(\xibar+\mubar_2)
-\xi\xibar
\bigr)
-(1+|\lambda|^2)
\Bigl(
 \bigl(\eta-\lambda\mubar_2\bigr)(\etabar-\lambdabar\mu_2)
-\eta\etabar
\Bigr)
\\
-\Bigl(
 \frac{\mubar_1^2\lambda}{|\mu_1|^2}
-\lambdabar
 \Bigr)
\Bigl(
 (\xi+\mu_2)(\etabar-\lambdabar\mu_2)
-\xi\etabar
\Bigr)
-\Bigl(
 -\frac{\mu_1^2\lambdabar}{|\mu_1|^2}+\lambda
 \Bigr)
 \Bigl(
 (\eta-\lambda\mubar_2)(\xibar+\mubar_2)
-\eta\xibar
 \Bigr)
\\
-\Bigl(
  \frac{\mu_1^2}{|\mu_1|^2}+\lambda^2
 \Bigr)
 \frac{1}{2}
\Bigl(
 (\xibar+\mubar_2)^2
-\xibar^2
\Bigr)
+\Bigl(
 \lambda^2\frac{\mubar_1^2}{|\mu_1|^2}+1
 \Bigr)
 \frac{1}{2}
\Bigl(
 (\etabar-\lambdabar\mu_2)^2
-\etabar^2
\Bigr)
 \\
-\Bigl(
 -\frac{\mubar_1^2}{|\mu_1|^2}+\lambdabar^2
 \Bigr)
 \frac{1}{2}
\Bigl(
 (\xi+\mu_2)^2
-\xi^2
\Bigr)
+\Bigl(
 -\lambdabar^2\frac{\mu_1^2}{|\mu_1|^2}+1
 \Bigr)\frac{1}{2}
 \Bigl(
 (\eta-\lambda\mubar_2)^2
-\eta^2
 \Bigr)
 \\
-(\xi-\lambda\eta+\mu_2+\lambda^2\mubar_2)
 2\lambdabar(\eta-\lambda\mubar_2)
+(\xi-\lambda\eta)2\lambdabar\eta
 \\
-(\mu_1+\lambda^2\mubar_1)^{-1}
 \bigl(2|\lambda|^2\mubar_1-\lambdabar^2\mu_1+\mubar_1\bigr)
\Bigl(
 (\xi-\lambda\eta+\mu_2+\lambda^2\mubar_2)^2
-(\xi-\lambda\eta)^2
\Bigr).
\end{multline}
We set
\begin{multline}
 F:=-|\mu_1|^{-2}(\mubar_1\mu_2-\mu_1\mubar_2)
  (1+|\lambda|^2)
 \bigl(\mubar_1(\xi-\lambda\etabar)+\mu_1(\xibar-\lambdabar\eta)\bigr)
 \\
+(1+|\lambda|^2)\bigl(
 \xi\mubar_2+\mu_2\xibar+|\mu_2|^2
+\lambdabar\mu_2\eta+\lambda\mubar_2\etabar
-|\lambda|^2|\mu_2|^2
 \bigr)
-\Bigl(
 \frac{\mubar_1^2\lambda}{|\mu_1|^2}-\lambdabar
 \Bigr)
 (-\xi\lambdabar\mu_2+\mu_2\etabar-\lambdabar\mu_2^2)
 \\
-\Bigl(
 -\frac{\mu_1^2\lambdabar}{|\mu_1|^2}+\lambda
 \Bigr)(\mubar_2\eta-\lambda\mubar_2\xibar-\lambda\mubar_2^2)
-\Bigl(
 \frac{\mu_1^2}{|\mu_1|^2}+\lambda^2
 \Bigr)(\xibar\mubar_2+\mubar_2^2/2)
 \\
+\Bigl(
 \lambda^2\frac{\mubar_1^2}{|\mu_1|^2}+1
\Bigr)(-\etabar\lambdabar\mu_2+\lambdabar^2\mu_2^2/2)
-\Bigl(
 -\frac{\mubar_1^2}{|\mu_1|^2}+\lambdabar^2
 \Bigr)
 (\xi\mu_2+\mu_2^2/2)
+\Bigl(
 -\lambdabar^2\frac{\mu_1^2}{|\mu_1|^2}+1
 \Bigr)(-\eta\lambda\mubar_2+\lambda^2\mubar_2^2/2)
 \\
-(\xi-\lambda\eta)2\lambdabar(-\lambda\mubar_2)
-(\mu_2+\lambda^2\mubar_2)
 2\lambdabar\eta
+(\mu_2+\lambda^2\mubar_2)2|\lambda|^2\mubar_2
 \\
-(\mu_1+\lambda^2\mubar_1)^{-1}
 \bigl(
 2|\lambda|^2\mubar_1-\lambdabar^2\mu_1+\mubar_1
 \bigr)
 \bigl(
 2(\xi-\lambda\eta)(\mu_2+\lambda^2\mubar_2)
+(\mu_2+\lambda^2\mubar_2)^2
 \bigr).
\end{multline}
Then, we have
\[
 \tte_2^{\ast}\sfvtilde^{\lambda}_{p,\omega}
=\sfvtilde^{\lambda}_{p,\omega}\exp
 \Bigl(
 \frac{\alpha(\omega)}{4(1+|\lambda|^2)^2}F
 \Bigr).
\]
We have the expression
$F=A_1\xibar+A_2\etabar+A_3\xi+A_4\eta+A_5$
for some constants $A_i$.
Because 
$\sfvtilde^{\lambda}$ and
$\tte_2^{\ast}\sfvtilde^{\lambda}$
are holomorphic and $\tte_0$-invariant,
we have $A_1=A_2=0$
and $A_4=-\lambda A_3$.
By a direct computation,
we obtain that
\[
 A_3=2(1+|\lambda|^2)^2
 (\mu_1\mubar_2-\mu_2\mubar_1)(\mu_1+\lambda^2\mubar_1)^{-1}
=-4(1+|\lambda|^2)^2\sqrt{-1}\Vol(\Gamma)
 (\mu_1+\lambda^2\mubar_1)^{-1}.
\]
We can also obtain the following by a direct computation:
\[
 A_5=(1+|\lambda|^2)^2
 (\mu_1+\lambda^2\mubar_1)^{-1}
 (\mu_1-\lambda^2\mubar_1)
 |\mu_1|^{-2}
 \bigl(
 |\mu_2|^2|\mu_1|^2
-\mu_2^2\mubar_1^2/2
-\mubar_2^2\mu_1^2/2
 \bigr).
\]
Then, we obtain the desired formula.
\hfill\qed

\vspace{.1in}

Let us study the growth order of $|\sfvtilde_{p,\omega}^{\lambda}|$
as $\ttU_p\to 0$ or $\ttU_p\to\infty$.
Recall
$\ttU_p=\exp\bigl(
2\pi\sqrt{-1}p^{-1}(\mu_1+\lambda \tts_1)^{-1}\ttu
 \bigr)$.
We describe
\[
 \ttu=p(\mu_1+\lambda \tts_1)\frac{c+\sqrt{-1}\sigma}{\sqrt{-1}}
\]
for real numbers $c$ and $\sigma$.
\begin{lem}
\label{lem;18.8.20.22}
We have
\[
 |\sfvtilde_{p,\omega}^{\lambda}|
\sim
 \exp\Bigl(
 \alpha(\omega)\Image(\ttv)
 \Re(\ttg_1\mu_1)pc
 \Bigr)
=
 \exp\Bigl(
 p\omega \Image(\ttv)
 \frac{\Re(\ttg_1\mu_1)}{\Vol(\Gamma)}
 2\pi c
 \Bigr)
=
 \exp\Bigl(
 -p\omega \Image(\ttv)
  (\gminit^{\lambda})^{-1}
 2\pi c
 \Bigr).
\]

\end{lem}
\pf
We have
\begin{multline}
 \bigl|
 \sfvtilde_{p,\omega}^{\lambda}
 \bigr|
=\exp\Bigl(
 \frac{\alpha(\omega)}{4}\frac{1}{(1+|\lambda|^2)^2}
 \Re\Bigl(
 (1+|\lambda|^2)(\ttu+\lambda\ttv)(\ttubar+\lambdabar \ttvbar)
-(1+|\lambda|^2)(\ttg_1\ttu+\ttv)(\ttgbar_1\ttubar+\ttvbar)
 \\
-\lambdabar^2(\ttu+\lambda \ttv)^2
+(\ttg_1\ttu+\ttv)^2
-\ttu(1-\ttg_1\lambda)2\lambdabar(\ttg_1\ttu+\ttv)
 \\
-(\mu_1+\lambda^2\mubar_1)^{-1}
 (2|\lambda|^2\mubar_1-\lambdabar^2\mu_1+\mubar_1)
 \ttu^2(1-\ttg_1\lambda)^2
 \Bigr)
 \Bigr).
\end{multline}

Let us look at the quadratic term with respect to $\ttu$.
\begin{equation}
\label{eq;18.3.16.20}
 -\lambdabar^2\ttu^2+\ttg_1^2\ttu^2
-2\lambdabar\ttg_1\ttu^2
+2\ttg_1^2|\lambda|^2\ttu^2
-(\mu_1+\lambda^2\mubar_1)^{-1}
 \bigl(
 2|\lambda|^2\mubar_1-\lambdabar^2\mu_1+\mubar_1
 \bigr)
 (1-\ttg_1\lambda)^2\ttu^2.
\end{equation}
We have
\[
 \ttg_1\ttu=(-\lambda\mubar_1+\tts_1)
 \frac{p(c+\sqrt{-1}\sigma)}{\sqrt{-1}}.
\]
Hence,
we can rewrite (\ref{eq;18.3.16.20})
as follows:
\begin{multline}
\Bigl(
 \lambdabar^2(\mu_1+\lambda \tts_1)^2
-(-\lambda\mubar_1+\tts_1)^2
+2\lambdabar(\mu_1+\lambda \tts_1)
 (-\lambda\mubar_1+\tts_1)
-2|\lambda|^2(-\lambda\mubar_1+\tts_1)^2
+(\mu_1+\lambda^2\mubar_1)
 (2|\lambda|^2\mubar_1-\lambdabar^2\mu_1+\mubar_1)
 \Bigr) \\
\times p^2(c+\sqrt{-1}\sigma)^2.
\end{multline}
It is equal to
$p^2(c+\sqrt{-1}\sigma)^2(1+|\lambda|^2)
 \Bigl(
 (|\lambda|^2-1)\tts_1^2
+2(\lambda\mubar_1+\lambdabar\mu_1)\tts_1
+(1-|\lambda|^2)|\mu_1|^2
 \Bigr)$.
By our choice of $\tts_1$,
it is $0$.

Let us study the linear term
with respect to $\ttu$ and $\ttubar$.
\begin{multline}
 \Re\Bigl(
 (1+|\lambda|^2)
 \bigl(\lambdabar \ttu\ttvbar
 +\lambda \ttubar \ttv
 -\ttg_1\ttu\ttvbar-\ttgbar_1\ttv\ttubar
 \bigr)
-2\lambdabar|\lambda|^2\ttu\ttv
+2\ttg_1\ttu\ttv
-(1-\ttg_1\lambda)2\lambdabar \ttu\ttv
 \Bigr) \\
=-2(1+|\lambda|^2)
 \Re\Bigl(
 \ttu(\lambdabar-\ttg_1)(\ttv-\ttvbar)
 \Bigr).
\end{multline}
Because
$\ttu=p(\mu_1+\lambda \tts_1)(c+\sqrt{-1}\sigma)/\sqrt{-1}$,
it is rewritten as follows:
\begin{equation}
\label{eq;18.3.16.30}
 -4(1+|\lambda|^2)\Image(\ttv)
 \Re\bigl(
 p(c+\sqrt{-1}\sigma)
 (\mu_1+\lambda \tts_1)(\lambdabar-\ttg_1)
 \bigr).
\end{equation}
We have the following:
\[
 (\mu_1+\lambda \tts_1)\lambdabar
-\ttg_1(\mu_1+\lambda \tts_1)
=\lambdabar(\mu_1+\lambda \tts_1)
-(-\lambda\mubar_1+\tts_1)
=(|\lambda|^2-1)\tts_1+\lambdabar\mu_1+\lambda\mubar_1
\in\real.
\]
We also have
\begin{equation}
\label{eq;18.8.8.10}
  (\mu_1+\lambda \tts_1)\lambdabar
-\ttg_1(\mu_1+\lambda \tts_1)
=\ttgbar_1(-\lambda\mubar_1+\tts_1)\lambdabar
-\ttg_1(\mu_1+\lambda \tts_1)
=-|\lambda|^2\ttgbar_1\mubar_1
-\ttg_1\mu_1
+\tts_1(\lambdabar\ttgbar_1-\lambda\ttg_1).
\end{equation}
Because (\ref{eq;18.8.8.10}) is real,
it is equal to
\[
 \frac{1}{2}\Bigl(
-|\lambda|^2\ttgbar_1\mubar_1
-|\lambda|^2\ttg_1\mu_1
-\ttg_1\mu_1
-\ttgbar_1\mubar_1
 \Bigr)
=-(1+|\lambda|^2)\Re(\ttg_1\mu_1).
\]
Hence, (\ref{eq;18.3.16.30}) is rewritten as 
$4(1+|\lambda|^2)^2\Re(\ttg_1\mu_1)\Image(\ttv)pc$.
Thus, we obtain the claim of the lemma.
\hfill\qed

\subsubsection{Mini-holomorphic frames of $\nbigl^{\lambda\cov}_{p}(\omega)$}

Because 
$\tte_0^{\ast}\sfvtilde^{\lambda}_{p,\omega}=\sfvtilde^{\lambda}_{p,\omega}$
and 
$(k(\omega)\tte_1)^{\ast}
 \sfvtilde^{\lambda}_{p,\omega}=\sfvtilde^{\lambda}_{p,\omega}$,
we obtain a mini-holomorphic frame
$\sfv^{\lambda}_{p,\omega}$
of $\nbigl^{\lambda\cov}_p(\omega)$
on $\nbigm_p^{\lambda\cov}$.
By Lemma \ref{lem;18.8.20.10},
we have
\begin{equation}
\label{eq;18.8.20.21}
  \tte_2^{\ast}\sfv_{p,\omega}^{\lambda}
=\sfv_{p,\omega}^{\lambda}
 \cdot
 \ttU_p^{-p\omega}
 \exp\Bigl(
 \frac{\alpha(\omega)}{4}
 (\mu_1+\lambda^2\mubar_1)^{-1}
 (\mu_1-\lambda^2\mubar_1)
 \cdot|\mu_1|^{-2}
 \bigl(
 |\mu_1|^2|\mu_2|^2
-\mu_2^2\mubar_1^2/2
-\mubar_2^2\mu_1^2/2
 \bigr)
 \Bigr).
\end{equation}
By Lemma \ref{lem;18.8.20.22},
we have
\[
 |\sfv_{p,\omega}^{\lambda}|
\sim
 \exp\Bigl(
 \alpha(\omega)\ttt
 \Re(\ttg_1\mu_1)pc
 \Bigr)
=
 \exp\Bigl(
 p\omega\ttt
 \frac{\Re(\ttg_1\mu_1)}{\Vol(\Gamma)}
 2\pi c
 \Bigr)
=|\ttU_p|^{-p\omega \ttt/\gminit^{\lambda}}
=|\ttU_p|^{p\omega \ttt \Re(\ttg_1\mu_1)/\Vol(\Gamma)}.
\]

\subsubsection{Associated filtered bundles}
\label{subsection;19.2.8.30}

By using the frame $\sfv^{\lambda}_{p,\omega}$,
we extend $\nbigl^{\lambda\cov}_p(\omega)$
to a locally free
$\nbigo_{\nbigmbar_p^{\lambda\cov}}(\ast H_p^{\lambda\cov})$-module
$\nbigp\nbigl^{\lambda\cov}_p(\omega)$.
Because $\nbigl^{\lambda\cov}_p(\omega)$ is
$\seisuu k(\omega)\tte_1\oplus\seisuu \tte_2$-equivariant,
we obtain the induced locally free 
$\nbigo_{\nbigmbar_p^{\lambda}}(\ast H_p^{\lambda})$-module
$\nbigp\nbigl^{\lambda}_p(\omega)$,
which is $(k(\omega)\seisuu/p\seisuu)\tte_1$-equivariant.

We define a filtered bundle
$\nbigp_{\ast}(\nbigl^{\lambda\cov}_p(\omega)_{|\pi_p^{-1}(\ttt)})$
over 
$\nbigl^{\lambda}_p(\omega)_{|\pi_p^{-1}(\ttt)}$
as follows:
for $\veca=(a_0,a_{\infty})\in\real^2$,
\[
 \nbigp_{\veca}\bigl(
\nbigl^{\lambda\cov}_p(\omega)_{|\pi_p^{-1}(\ttt)}
 \bigr)
=
 \nbigo_{\proj^1}\Bigl(
 \bigl[
 a_0-p\omega\ttt(\gminit^{\lambda})^{-1}
 \bigr]\cdot\{0\}
+
 \bigl[
 a_{\infty}+p\omega\ttt(\gminit^{\lambda})^{-1}
 \bigr]\cdot\{\infty\}
 \Bigr)
 \sfv^{\lambda}_{p,\omega}.
\]
We obtain a filtered bundle
$\nbigp_{\ast}\nbigl^{\lambda}_p(\omega)$
over $\nbigp\nbigl^{\lambda}_p(\omega)$
as the descent,
which is $(k(\omega)\seisuu/p\seisuu)\tte_1$-equivariant.

\begin{lem}
$\nbigp_{\ast}\nbigl^{\lambda}_p(\omega)_{|\Hhat^{\lambda}_{p,0}}$
is isomorphic to
$\nbigp_{\ast}^{(0)}
 \vecLL_{p,0}(\omega)
 \otimes
 \nbigp_{\ast}^{(0)}\vecVV_{p,0}\Bigl(
 \beta(\omega)
 \Bigr)$,
where
\[
 \beta(\omega):=
\exp\Bigl(
 \frac{\alpha(\omega)}{4}
 (\mu_1+\lambda^2\mubar_1)^{-1}
 (\mu_1-\lambda^2\mubar_1)
 \cdot|\mu_1|^{-2}
 \bigl(
 |\mu_1|^2|\mu_2|^2
-\mu_2^2\mubar_1^2/2
-\mubar_2^2\mu_1^2/2
 \bigr)
 \Bigr).
\]
Similarly,
$\nbigp_{\ast}\nbigl^{\lambda}_p(\omega)_{|\Hhat^{\lambda}_{p,\infty}}$
is isomorphic to
$\nbigp_{\ast}^{(0)}
 \vecLL_{p,\infty}(-\omega)
 \otimes
 \nbigp_{\ast}^{(0)}\vecVV_{p,\infty}\Bigl(
 \beta(\omega)
 \Bigr)$.
\hfill\qed
\end{lem}

\subsection{Examples (2)}
\label{subsection;18.12.21.22}

\subsubsection{Preliminary}

We define the action of
$\real\,\tte_3$ on $\cnum^2$
by
\[
 \tte_3(z,w)=(z,w+\sqrt{-1}).
\]
It is described as follows
in terms of $(\xi,\eta)$:
\[
 \tte_3(\xi,\eta)=(\xi,\eta)
+(-\lambda\sqrt{-1},\sqrt{-1}).
\]

Let $(\ttx,\tty)$ be 
the complex coordinate system
determined by
$(\xi,\eta)=
 \ttx(-\lambda,1)+\tty(1,\lambdabar)$.
Note that
$d\xi\,d\xibar+d\eta\,d\etabar
=(1+|\lambda|^2)
 \bigl(
d\ttx\,d\ttxbar
+d\tty\,d\ttybar
\bigr)$.
The following holds:
\[
 \left\{
 \begin{array}{l}
 \xi=-\lambda \ttx+\tty\\
 \eta=\ttx+\lambdabar \tty,
 \end{array}
 \right.
\quad\quad
 \left\{
 \begin{array}{l}
 \ttx=(1+|\lambda|^2)^{-1}
 \bigl(\eta-\lambdabar\xi\bigr)\\
 \tty=(1+|\lambda|^2)^{-1}
 \bigl(\xi+\lambda\eta\bigr).
 \end{array}
 \right.
\]
We have the following formulas:
\[
 \tte_3(\ttx,\tty)=(\ttx,\tty)+(\sqrt{-1},0),
\]
\[
 \tte_0(\ttx,\tty)=
 (\ttx,\tty)
+\frac{1}{1+|\lambda|^2}(1-|\lambda|^2,2\lambda),
\]
\[
 \tte_i(\ttx,\tty)=(\ttx,\tty)
+\frac{1}{1+|\lambda|^2}
 \bigl(
 -\lambda\mubar_i-\lambdabar\mu_i,\,
 \mu_i-\lambda^2\mubar_i
 \bigr)
\quad(i=1,2).
\]

We have the following relations:
\[
 \left\{
 \begin{array}{l}
 \ttu=(1-\ttg_1\lambda)^{-1}
 \bigl(
 -2\lambda\ttx+(1-|\lambda|^2)\tty
 \bigr)
\\
 \ttv=(1-\ttg_1\lambda)^{-1}
 \bigl(
 (1+\ttg_1\lambda)\ttx+
 (\lambdabar-\ttg_1)\tty
 \bigr),
 \end{array}
 \right.
\quad
 \left\{
 \begin{array}{l}
 \ttx=(1+|\lambda|^2)^{-1}
 \bigl((\ttg_1-\lambdabar)\ttu
+(1-|\lambda|^2)\ttv
\bigr)
  \\
 \tty=(1+|\lambda|^2)^{-1}
 \bigl(
(1+\lambda\ttg_1)\ttu
+2\lambda \ttv
 \bigr).
 \end{array}
 \right.
\]

\begin{lem}
There exists a unique solution $(\ttA,\ttB)\in\cnum^2$
of the equation
\begin{equation}
\label{eq;18.3.23.1}
 \tte_0^{\ast}\bigl(\ttybar+\ttA\ttx+\ttB\tty\bigr)
 -(\ttybar+\ttA\ttx+\ttB\tty)=0,
\quad
 \tte_1^{\ast}\bigl(\ttybar+\ttA\ttx+\ttB\tty\bigr)
-(\ttybar+\ttA\ttx+\ttB\tty)
 =0.
\end{equation}
Indeed, we have
\begin{equation}
\label{eq;18.8.17.1}
 \ttA=\frac{2(\lambda\mubar_1-\lambdabar\mu_1)}{\mu_1+\lambda^2\mubar_1},
\quad\quad
 \ttB=\frac{-(\mubar_1+\lambdabar^2\mu_1)}{\mu_1+\lambda^2\mubar_1}.
\end{equation}
For such $\ttA$ and $\ttB$,
the following holds:
\begin{equation}
\label{eq;18.8.17.2}
\ttC:=
 \tte_2^{\ast}\bigl(\ttybar+\ttA\ttx+\ttB\tty\bigr)
-(\ttybar+\ttA\ttx+\ttB\tty)
=
-2\sqrt{-1}(1+|\lambda|^2)
 \frac{\Vol(\Gamma)}{\mu_1+\lambda^2\mubar_1}
\neq 0.
\end{equation}
\end{lem}
\pf
The equation (\ref{eq;18.3.23.1})
is equivalent to the following equation:
\[
 \left\{
 \begin{array}{l}
  2\lambdabar+\ttA(1-|\lambda|^2)+2\ttB\lambda=0\\
 \mubar_1-\lambdabar^2\mu_1
+\ttA(-\lambda\mubar_1-\lambdabar\mu_1)
+\ttB(\mu_1-\lambda^2\mubar_1)=0.
 \end{array}
 \right.
\]
Because
$(1-|\lambda|^2)(\mu_1-\lambda^2\mubar_1)
-2\lambda(-\lambda\mubar_1-\lambdabar\mu_1)
=(1+|\lambda|^2)(\mu_1+\lambda^2\mubar_1)\neq 0$,
we have a unique solution $(\ttA,\ttB)$.
We obtain (\ref{eq;18.8.17.1})
and  (\ref{eq;18.8.17.2})
by direct computations.
\hfill\qed

\vspace{.1in}

Recall
$\ttU_p=\exp\Bigl(
 \frac{2\pi\sqrt{-1}}{p(\mu_1+\lambda \tts_1)}
 \ttu
 \Bigr)$.
We consider
\begin{equation}
\label{eq;18.8.20.30}
 \ttu=\frac{\mu_1+\lambda \tts_1}{\sqrt{-1}}p(c+\sqrt{-1}\sigma)
\sim
\frac{\mu_1+\lambda \tts_1}{\sqrt{-1}}pc.
\end{equation}
We have
\[
 \ttx\sim
 \frac{1}{1+|\lambda|^2}
 (\ttg_1-\lambdabar)
 \frac{\mu_1+\lambda \tts_1}{\sqrt{-1}}pc
=\frac{\Re(\ttg_1\mu_1)}{\sqrt{-1}}pc,
\quad\quad
 \tty\sim
 \frac{1}{1+|\lambda|^2}
 \bigl(\mu_1-\lambda^2\mubar_1+2\lambda\tts_1\bigr)
 \frac{pc}{\sqrt{-1}}.
\]

\subsubsection{Construction}

For $(\sfa,\sfb)\in\real\times\cnum$,
let $\Ltilde(\lambda,\sfa,\sfb)$
be the line bundle on $X^{\lambda}$
with a global frame $\vtilde_{0,(\sfa,\sfb)}$.
Let $\htilde$ be the metric determined by
$\htilde(\vtilde_{0,(\sfa,\sfb)},\vtilde_{0,(\sfa,\sfb)})=1$.
Let $\delbar_{\Ltilde(\lambda,\sfa,\sfb)}$
be the holomorphic structure determined by
\[
 \del_{\ttxbar}\vtilde_{0,(\sfa,\sfb)}=
 \vtilde_{0,(\sfa,\sfb)}\sqrt{-1}\sfa,
\quad
 \del_{\ttybar}\vtilde_{0,(\sfa,\sfb)}=
 \vtilde_{0,(\sfa,\sfb)}\sfb.
\]
The holomorphic bundle $\Ltilde(\lambda,\sfa,\sfb)$
with the metric is equivariant with respect to
the action of 
$\real \tte_0\oplus\seisuu \tte_1\oplus\seisuu \tte_2$
by $\tte_i^{\ast}(\vtilde_{0,(\sfa,\sfb)})=\vtilde_{0,(\sfa,\sfb)}$.
It induces a mini-holomorphic bundle
$L_p^{\cov}(\lambda,\sfa,\sfb)$ of rank $1$
with the induced metric $h^{\cov}$
on $\nbigm_p^{\lambda\cov}$,
which is a $(\seisuu/p\seisuu)\tte_1\times \seisuu \tte_2$-equivariant
monopole.
We also obtain a monopole
$(L_p(\lambda,\sfa,\sfb),h)$ on $\nbigmlambda_p$
as the descent,
which is $(\seisuu/p\seisuu)\tte_1$-equivariant.

\subsubsection{Underlying mini-holomorphic bundles}
\label{subsection;18.12.21.21}

We have the holomorphic section $\vtilde_{1,(\sfa,\sfb)}$ of 
$\Ltilde(\lambda,\sfa,\sfb)$
given as follows:
\[
 \vtilde_{1,(\sfa,\sfb)}:=
 \vtilde_{0,(\sfa,\sfb)}\cdot\exp\Bigl(
 (\ttx-\ttxbar)\sqrt{-1}\sfa
-(\ttybar+\ttA\ttx+\ttB\tty)\sfb
 \Bigr).
\]
We have
$\tte_0^{\ast}\vtilde_{1,(\sfa,\sfb)}=\vtilde_{1,(\sfa,\sfb)}$
and 
$\tte_1^{\ast}\vtilde_{1,(\sfa,\sfb)}=\vtilde_{1,(\sfa,\sfb)}$.
We also have 
\[
 \tte_2^{\ast}\vtilde_{1,(\sfa,\sfb)}
=\vtilde_{1,(\sfa,\sfb)}\cdot\exp\Bigl(-\ttC\sfb\Bigr).
\]
We obtain the induced mini-holomorphic frame
$v_{1,(\sfa,\sfb)}$ of 
$L_p^{\cov}(\lambda,\sfa,\sfb)$
on $\nbigm_p^{\lambda\cov}$
for which the following holds:
\[
  \tte_2^{\ast}v_{1,(\sfa,\sfb)}
=v_{1,(\sfa,\sfb)}\cdot\exp\Bigl(-\ttC\sfb\Bigr).
\]
Because
$|\vtilde_{1,(\sfa,\sfb)}|_{\htilde}
=\exp\Bigl(
 \Re\bigl((\ttx-\ttxbar)\sqrt{-1}\sfa\bigr)
-\Re\bigl(
 (\ttybar+\ttA\ttx+\ttB\tty)\sfb
 \bigr)
 \Bigr)$,
we have
\[
 |\vtilde_{1,(\sfa,\sfb)}|_{\htilde}
\sim
 \exp\Bigl(
 \Re\Bigl(
2\Re(\ttg_1\mu_1)\sfa
-\frac{2}{\sqrt{-1}}
 \frac{\sfb}{\mu_1+\lambda^2\mubar_1}
 \Bigl(
 (|\lambda|^2-1)|\mu_1|^2
-(\lambdabar\mu_1+\lambda\mubar_1)\tts_1
+(\lambda\mubar_1-\lambdabar\mu_1)
 \Re(\ttg_1\mu_1)
 \Bigr)pc
 \Bigr),
\]
where $c$ is introduced as in (\ref{eq;18.8.20.30}).

\subsubsection{Associated filtered bundles}
\label{subsection;19.2.8.31}

By using the frame $v_{1,(\sfa,\sfb)}$,
we extend $L_p^{\cov}(\lambda,\sfa,\sfb)$
to a locally free
$\nbigo_{\nbigmbar_p^{\lambda\cov}}(\ast H^{\lambda\cov})$-module
$\nbigp L^{\cov}(\lambda,\sfa,\sfb)$.
We set
\[
 \paramap(\lambda,\sfa,\sfb):=
  \frac{1}{2\pi}
 \Re\Bigl[
2\Re(\ttg_1\mu_1)\sfa
-\frac{2}{\sqrt{-1}}
 \frac{\sfb}{\mu_1+\lambda^2\mubar_1}
 \Bigl(
 (|\lambda|^2-1)|\mu_1|^2
-(\lambdabar\mu_1+\lambda\mubar_1)\tts_1
+(\lambda\mubar_1-\lambdabar\mu_1)
 \Re(\ttg_1\mu_1)
 \Bigr)
 \Bigr].
\]
We obtain the filtered bundle
$\nbigp_{\ast}L^{\cov}_p(\lambda,\sfa,\sfb)$
over $\nbigp L_p^{\cov}(\lambda,\sfa,\sfb)$
determined as follows;
for $\veca=(a_0,a_1)\in\real^2$,
\[
 \nbigp_{\veca}
\bigl(
 L_p^{\cov}(\lambda,\sfa,\sfb)
_{|\pi_p^{-1}(\ttt)}
\bigr)
:=\nbigo_{\proj^1}\Bigl(
 [a_0-p\paramap(\lambda,\sfa,\sfb)]\cdot\{0\}
+[a_{\infty}+p\paramap(\lambda,\sfa,\sfb)]\cdot\{\infty\}
 \Bigr)
 \cdot
 v_{1,(\sfa,\sfb)|\pi_p^{-1}(\ttt)}.
\]
Because 
$\nbigp(L_p^{\cov}(\lambda,\sfa,\sfb))$
and $\nbigp_{\ast}(L_p^{\cov}(\lambda,\sfa,\sfb))$
are equivariant with respect to
the $\seisuu \tte_2$-action,
we obtain 
a locally free
$\nbigo_{\nbigmbar_p^{\lambda}}(\ast H_p^{\lambda})$-module
$\nbigp(L_p(\lambda,\sfa,\sfb))$
and a filtered bundle
$\nbigp_{\ast}(L_p(\lambda,\sfa,\sfb))$
over 
$\nbigp(L_p(\lambda,\sfa,\sfb))$.

\begin{lem}
\label{lem;19.1.7.10}
$\nbigp_{\ast}L_p(\lambda,\sfa,\sfb)_{|\Hhat^{\lambda}_{p,0}}$
is isomorphic to
$\nbigp^{(p\paramap(\lambda,\sfa,\sfb))}_{\ast}
 \vecVV_{p,0}\bigl(e^{-\ttC\sfb}\bigr)$,
and
$\nbigp_{\ast}L_p(\lambda,\sfa,\sfb)_{|\Hhat^{\lambda}_{p,\infty}}$
is isomorphic to
$\nbigp^{(-p\paramap(\lambda,\sfa,\sfb))}_{\ast}
 \vecVV_{p,\infty}(e^{-\ttC\sfb})$.
\hfill\qed
\end{lem}

\subsubsection{Isomorphisms}
\label{subsection;19.1.7.1}

For any $\vecn=(n_1,n_2)\in\seisuu^2$,
we set 
\[
 \chi_{\vecn}(z):=
 \exp\Bigl(
 \frac{\pi}{\Vol(\Gamma)}
 \Bigl(
-n_1p(\mubar_1z-\mu_1\zbar)
+n_2(\mubar_2z-\mu_2\zbar)
 \Bigr)
 \Bigr).
\]
It induces a function $\chi_{\vecn}$
on $\nbigm_p^0=\nbigm_p^{\lambda}$.
We have the isomorphism of monopoles
$F_{\vecn}:
 L_p(\lambda,\sfa,\sfb)
\simeq
 L_p(\lambda,\sfa',\sfb')$
induced by
$F_{\vecn}(\chi_{\vecn}\vtilde_{0,(\sfa,\sfb)})
=\vtilde_{0,(\sfa',\sfb')}$,
where
\[
 (\sfa',\sfb')
:=(\sfa,\sfb)
+\frac{\pi}{(1+|\lambda|^2)\Vol(\Gamma)}
 \Bigl[
 pn_1
\Bigl(
 -\sqrt{-1}(\mubar_1\lambda-\mu_1\lambdabar),
 \mu_1+\lambda^2\mubar_1
 \Bigr)
-n_2
\Bigl(
 -\sqrt{-1}(\mubar_2\lambda-\mu_2\lambdabar),
 \mu_2+\lambda^2\mubar_2
 \Bigr)
\Bigr].
\]
We have 
$F(v_{1,(\sfa,\sfb)})
=\ttU_p^{pn_2}v_{1,(\sfa',\sfb')}$.

\begin{rem}
The numbers
$\exp\Bigl(-\ttC\sfb\Bigr)$
and 
$p\paramap(\lambda,\sfa,\sfb)$
determine $(\sfa,\sfb)$
up to the induced action of $\seisuu p\tte_1$.
\hfill\qed
\end{rem}

\subsubsection{Comparison with $\lambda=0$}

We define the bijection
$\ttF^{\lambda}:\real\times\cnum\simeq \real\times\cnum$
by
\[
 \ttF^{\lambda}(\sfa,\sfb):=
 \Bigl(
 \sfa+\frac{2\Image(\sfbbar\lambda)}{1+|\lambda|^2},
 \,\,
 \frac{\sfb+\lambda^2\sfbbar}{1+|\lambda|^2}
 \Bigr).
\]

\begin{lem}
$\vecsfL_p(0,\sfa,\sfb)
=\vecsfL_p\bigl(\lambda,\ttF^{\lambda}(\sfa,\sfb)\bigr)$
holds
on $\nbigm^0_p=\nbigm^{\lambda}_p$.
\end{lem}
\pf
It is enough to compare the corresponding instantons
on $X$.
Let $\vtilde^0_{0,(\sfa,\sfb)}$ be the global frame
of $\Ltilde(0,\sfa,\sfb)$.
The unitary connection is given as
\[
 \nablatilde\vtilde^0_{0,(\sfa,\sfb)}
=\vtilde^0_{0,(\sfa,\sfb)}
 \bigl(
 \sqrt{-1}\sfa\,d\wbar
-\sqrt{-1}\sfa\,dw
+\sfb\,d\zbar
-\sfbbar\,dz
 \bigr).
\]
We have the following relation:
\[
 z=\frac{1}{1+|\lambda|^2}
 (-\lambda \ttx-\lambda\ttxbar+\tty-\lambda^2\ttybar),
\quad
 w=\frac{1}{1+|\lambda|^2}
 \bigl(\ttx-|\lambda|^2\ttxbar+\lambdabar\tty+\lambda\ttybar\bigr).
\]
By a direct computation,
we obtain
\begin{multline}
 \sqrt{-1}\sfa\,d\wbar
-\sqrt{-1}\sfa\,dw
+\sfb\,d\zbar
-\sfbbar\,dz
=\\
\sqrt{-1}\sfa(d\ttxbar-d\ttx)
+\frac{1}{1+|\lambda|^2}
 \Bigl(
 2\sqrt{-1}\Image(\sfbbar\lambda)\,d\ttx
+2\sqrt{-1}\Image(\sfbbar\lambda)\,d\ttxbar
+(\sfb+\lambda^2\sfbbar)\,d\ttybar
-(\sfbbar+\lambdabar^2\sfb)\,d\tty
 \Bigr).
\end{multline}
Thus, we obtain the claim of the lemma.
\hfill\qed

\subsubsection{Twist}

Recall that we constructed a monopole
$\vecsfL_p(\omega)$ on $\nbigm_p^0$
for $\omega\in \frac{1}{p}\seisuu$
in \S\ref{subsection;18.12.14.1}.

\begin{lem}
\label{lem;19.1.5.2}
We set
$\sfb_0:=-\frac{\pi\omega\mu_1}{\Vol(\Gamma)}$.
Then,
$\tte_1^{\ast}\vecsfL_p(\omega)$
is isomorphic to
$\vecsfL_p(\omega)
\otimes
 L_p(0,0,\sfb_0)$.
The isomorphism is induced by
$\tte_1^{\ast}(\gminie)
\longmapsto
\gminie\otimes v_{1,(0,\sfb_0)}$.
\end{lem}
\pf
Note that  $-\ttC\sfb_0=-2\pi\sqrt{-1}\omega$.
Let $v_{1,(0,\sfb_0)}$ be the mini-holomorphic frame 
of $L_p^{\cov}(0,0,\sfb_0)$
as in \S\ref{subsection;18.12.21.21}.
Then, we have
$|v_{1,(0,\sfb_0)}|=1$
and 
$\tte_2^{\ast}v_{1,(0,\sfb_0)}=
 v_{1,(0,\sfb_0)}
 \exp(-2\pi\omega)$.
Because
\[
 \tte_2^{\ast}(\tte_1^{\ast}\gminie)
\longmapsto
 \tte_1^{\ast}(\gminie)
 \exp\Bigl(
 -2\pi\sqrt{-1}\omega|\mu_1|^{-2}\Re(\mubar_1z)
 \Bigr)
\cdot
 \exp(-2\pi\omega),
\]
we obtain the claim of the lemma.
\hfill\qed

\subsection{Examples (3)}
\label{subsection;19.1.7.21}

\subsubsection{Neighbourhoods}

We continue to use the notation in \S\ref{subsection;18.12.21.22}.
Let $R>0$.
We set
$\nbigutilde_{-,R}=\{(\ttx,\tty)\in X^{\lambda}\,|\,\Image(\ttx)<-R\}$
and 
$\nbigutilde_{+,R}=\{(\ttx,\tty)\in X^{\lambda}\,|\,\Image(\ttx)>R\}$.
Let $\nbigu^{\cov}_{p,\pm,R}$
denote the quotient of $\nbigutilde_{\pm,R}$
by the action of $\real\tte_0\oplus \seisuu p\tte_1$.
Let $\nbigu_{p,\pm,R}$ denote the quotient of
$\nbigu^{\cov}_{p,\pm,R}$
by the action of $\seisuu\tte_2$.

If $\Re(\ttg_1\mu_1)>0$,
we set
\[
 \nbigu^{\lambda}_{p,\infty,R}:=\nbigu_{p,-,R},
\quad
 \nbigu^{\lambda\cov}_{p,\infty,R}:=\nbigu^{\cov}_{p,-,R},
\quad
 \nbigu^{\lambda}_{-,0,R}:=\nbigu_{p,+,R},
\quad
 \nbigu^{\lambda\cov}_{p,0,R}:=\nbigu^{\cov}_{p,+,R}.
\]
If $\Re(\ttg_1\mu_1)<0$,
we set
\[
 \nbigu^{\lambda}_{p,\infty,R}:=\nbigu_{p,+,R},
\quad
 \nbigu^{\lambda\cov}_{p,\infty,R}:=\nbigu^{\cov}_{p,+,R},
\quad
 \nbigu^{\lambda}_{-,0,R}:=\nbigu_{p,-,R},
\quad
 \nbigu^{\lambda\cov}_{p,0,R}:=\nbigu^{\cov}_{p,-,R}.
\]
Then, 
$\nbigubar^{\lambda}_{p,\nu,R}:=
 \nbigu^{\lambda}_{p,\nu,R}
 \cup
 H^{\lambda}_{p,\nu}$
is a neighbourhood of $H^{\lambda}_{p,\nu}$.
Similarly,
$\nbigubar^{\lambda\cov}_{p,\nu,R}:=
 \nbigu^{\lambda\cov}_{p,\nu,R}
 \cup
 H^{\lambda\cov}_{p,\nu}$
is a neighbourhood of $H^{\lambda\cov}_{p,\nu}$.

\subsubsection{Examples of monopoles of rank $2$
with unipotent monodromy}

Let $\Vtilde_{\pm}(\lambda,2)$ be the holomorphic vector bundle
on $\nbigutilde_{\pm,R}$
with a global frame $\vecvtilde=(\vtilde_1,\vtilde_2)$
with the holomorphic structure
determined by
\[
 \del_{\ttxbar}\vecvtilde=0,
\quad
 \del_{\ttybar}\vecvtilde=\vecvtilde
 N_2,
\quad
\mbox{\rm where }
N_2:=
 \left(
 \begin{array}{cc}
 0 & 0 \\
 1 & 0
 \end{array}
 \right).
\]
Let $\htilde$ be the metric of $\Vtilde_{\pm}(\lambda,2)$
determined by $\htilde(\vtilde_1,\vtilde_2)=0$,
$\htilde(\vtilde_1,\vtilde_1)=\bigl|\Image(\ttx)\bigr|$
and 
$\htilde(\vtilde_2,\vtilde_2)=\bigl|\Image \ttx\bigr|^{-1}$.
The holomorphic bundles
with a Hermitian metric are instantons
on $\nbigutilde_{\pm,R}$.

We have the holomorphic frame
$\vecutilde:=\vecvtilde\cdot
 \exp\Bigl(
 -(\ttybar+\ttA\ttx+\ttB\tty)N_2
 \Bigr)$
of $\Vtilde_{\pm}(\lambda,2)$.
We have
$\tte_0^{\ast}\vecutilde=\vecutilde$
and $\tte_1^{\ast}\vecutilde=\vecutilde$.
We also have
\[
 \tte_2^{\ast}\vecutilde=\vecutilde\exp(-\ttC N_2).
\]
\begin{lem}
Let $\htilde_0$ be the metric of $V_{\pm}(\lambda,2)$
determined by
\[
\htilde_0(\utilde_1,\utilde_1)=
 \bigl|\Image(\ttx)\bigr|,
\quad
\htilde_0(\utilde_2,\utilde_2)=
\bigl|\Image(\ttx)\bigr|^{-1},
\quad
\htilde_0(\utilde_1,\utilde_2)=0.
\]
Then, we have
$\htilde_0$ and $\htilde$ are mutually bounded.
\hfill\qed
\end{lem}

We define the action of
$\real \tte_0\oplus\seisuu \tte_1\oplus \seisuu \tte_2$
on $\Vtilde_{\pm}(\lambda,2)$
by $\tte_i^{\ast}(\vecvtilde)=\vecvtilde$,
and the holomorphic structure and the metric
are preserved by the action.
We obtain the corresponding mini-holomorphic bundles
$V^{\cov}_{p,\nu}(\lambda,2)$
on $\nbigu^{\lambda\cov}_{p,\nu,R}$
and 
$V_{p,\nu}(\lambda,2)$ 
on $\nbigu^{\lambda}_{p,\nu,R}$
for $\nu=0,\infty$.
They are equipped with the induced metrics
$h^{\cov}$ and $h$, respectively.
With the metrics,
they are monopoles.

We obtain the induced mini-holomorphic frame
$\vecu$ of $V_{p,\nu}^{\cov}(\lambda,2)$,
with which 
$V_{p,\nu}^{\cov}(\lambda,2)$
extends to a mini-holomorphic bundle
$\nbigp_0V^{\cov}_{p,\nu}(\lambda,2)$
on $\nbigubar^{\lambda\cov}_{p,\nu,R}$.
It induces a filtered bundle
$\nbigp_{\ast}V^{\cov}_{p,\nu}(\lambda,2)$
over $(\nbigubar_{p,\nu,R}^{\lambda\cov},H^{\lambda\cov}_{p,\nu})$
such that
$\Gr^{\nbigp}_aV^{\cov}_{p,\nu}(\lambda,2)=0$
unless $a\in\seisuu$.
We obtain the induced frame 
$[\vecu]$ of $\Gr^{\nbigp}_0(V^{\cov}_{p,\nu}(\lambda,2))$,
for which
$\tte_2^{\ast}[\vecu]
=[\vecu]\exp(-\ttC N_2)$ holds.

Because $\nbigp_{\ast}V^{\cov}_{p,\nu}(\lambda,2)$
is $\seisuu\tte_2$-equivariant,
we obtain an induced filtered bundle
$\nbigp_{\ast}V^{\cov}_{p,\nu}(\lambda,2)$
on $(\nbigubar^{\lambda}_{p,\nu,R},H^{\lambda}_{p,\nu})$,
which is an extension of
$V_{p,\nu}(\lambda,2)$.
The conjugacy class of the monodromy of
$\Gr^{\nbigp}_0(V_{p,\nu}(\lambda,2))$
is $\exp(-\ttC N_2)$.

\subsubsection{Examples with any monodromy at infinity}

Take $(\sfa_i,\sfb_i)\in\real\times\cnum$ $(i=1,\ldots,m)$
and $\ell_i\in\seisuu_{\geq 0}$ $(i=1,\ldots,m)$.
We obtain the following monopole:
\[
 E=
 \bigoplus_{i=1}^m
 L(\lambda,\sfa_i,\sfb_i)
\otimes
 \Sym^{\ell_i} V(\lambda,2).
\]
We have
\[
 \Gr^{\nbigp}_a(E)
=\bigoplus_i
 \Gr^{\nbigp}_a\bigl(L(\lambda,\sfa_i,\sfb_i)\bigr)
 \otimes
 \Gr^{\nbigp}_0(V(\lambda,2)).
\]
The conjugacy class 
of the monodromy on 
$\Gr^{\nbigp}_{\paramap(\lambda,\sfa_i,\sfb_i)}
 \bigl(L(\lambda,\sfa_i,\sfb_i)\bigr)
 \otimes
 \Gr^{\nbigp}_0(V(\lambda,2))$
 is
\[
 \exp(-\ttC\sfb_i)\cdot
 \exp(-\ttC N_{\ell_i+1}).
\]
Here, $N_{\ell_i+1}$ is a $(\ell_i+1)$-square matrix
such that
$(N_{\ell_i+1})_{j+1,j}=1$ $(j=1,\ldots,\ell_i)$
and $(N_{\ell_i+1})_{i,j}=0$ $(i\neq j+1)$.

\subsubsection{Another expression}

Suppose that $\vecA=(A_1,A_2,A_3)\in\su(n)$ $(i=1,2,3)$ 
satisfy $[A_i,A_j]+A_k=0$
for any cyclic permutation $(i,j,k)$ of $(1,2,3)$.
Let $\nbigvtilde_{\pm}$ be
a product bundle 
$\nbigutilde_{\pm,R}\times\cnum^n$
on $\nbigutilde_{\pm,R}$
with a global frame $\vece=(e_1,\ldots,e_n)$.

Let $h_{\nbigvtilde_{\pm}}$ be the Hermitian metric 
of $\nbigvtilde_{\pm}$ for which
the frame $\vece$ is orthonormal.
We define operators
$\del_{\nbigvtilde_{\pm},\ttxbar}$
and 
$\del_{\nbigvtilde_{\pm},\ttybar}$
on $\nbigvtilde_{\pm}$
by 
\[
 \del_{\nbigvtilde_{\pm},\ttxbar}\vece=
 \vece\cdot
 \frac{1}{2\Image(\ttx)}A_3,
\quad
 \del_{\nbigvtilde_{\pm},\ttybar}\vece
=\vece\cdot
 \frac{1}{2\Image(\ttx)}\bigl(
 A_1+\sqrt{-1}A_2
 \bigr).
\]
Then, the operators give a holomorphic structure
$\delbar_{\nbigvtilde_{\pm}}$ of $\nbigvtilde_{\pm}$,
and 
$(\nbigvtilde_{\pm},\delbar_{\nbigvtilde_{\pm}},h_{\nbigvtilde_{\pm}})$
are instantons on $\nbigutilde_{\pm,R}$.
It is naturally equivariant
with respect to the action of
$\real\tte_0\oplus\seisuu\tte_1\oplus\seisuu\tte_2$
determined by
$\tte_i^{\ast}\vece=\vece$.
Hence, we obtain monopoles
$\nbigv^{\cov}_{p,\pm}(\lambda,\vecA)$
on $\nbigu^{\cov}_{p,\pm,R}$,
and 
$\nbigv_{p,\pm}(\lambda,\vecA)$
on $\nbigu_{p,\pm,R}$.

The following is easy to check.
\begin{lem}
If $(k_1,\ldots,k_m)$ be the weight decomposition of
the $\su(2)$-representation determined by $\vecA$,
then 
$\nbigv(\lambda,\vecA)$
is naturally isomorphic to
$\bigoplus \Sym^{k_i}V(\lambda,2)$.
\hfill\qed
\end{lem}
The following is easy to see.
\begin{lem}
$\nbigv(\lambda,\vecA)$
is isomorphic to 
$\nbigv(0,\vecA)$.
\hfill\qed
\end{lem}

\subsection{Example (4)}
\label{subsection;19.1.5.10}

We set $\Gamma^{\lor}:=\bigl\{
 \sfb\in\cnum\,\big|\,
 \Image(\chi \sfbbar)\in\pi\seisuu
 \,\,(\forall \chi\in\Gamma)
 \bigr\}$.
We set $\mu_i^{\lor}:=\pi\mu_i/\Vol(\Gamma)$.
Then,
$\Gamma^{\lor}:=\seisuu \mu_1^{\lor}\oplus\seisuu\mu_2^{\lor}$.
We set
$\Gamma_p:=\seisuu (p\mu_1)\oplus\seisuu \mu_2$.
We have
$(\Gamma_p)^{\lor}=
 \seisuu \mu_1^{\lor}
\oplus
 \seisuu\cdot (\mu_2^{\lor}/p)$.

Let $\omega\in\rnum$.
We set $k(\omega):=\min\{p\in\seisuu_{>0}\,|\,p\omega\in\seisuu\}$.
We have the action of
$(\seisuu/k(\omega)\seisuu)\cdot
 (\omega\mu_1^{\lor})$
on 
$\cnum/(\Gamma_{k(\omega)})^{\lor}$
induced by the addition.
It naturally induces
an action of 
$(\seisuu/k(\omega)\seisuu)
\cdot (\omega\mu_1^{\lor})$
on $\real\times
\cnum/(\Gamma_{k(\omega)})^{\lor}$.

Let $I\subset \rnum$
be a finite subset.
For each $\omega\in I$,
let $\nbigs_{\omega}\subset
 \real\times\cnum/(\Gamma_{k(\omega)})^{\lor}$
which is preserved by 
the action of 
$(\seisuu/k(\omega)\seisuu)\cdot (\omega\mu_1^{\lor})$.
For each $(\sfa,\sfb)\in \nbigs_{\omega}$,
let $n(\omega,\sfa,\sfb)\in\seisuu_{\geq 0}$,
and let 
\[
 \vecA_{\omega,\sfa,\sfb}
=(A_{1,\omega,\sfa,\sfb},A_{2,\omega,\sfa,\sfb},A_{3,\omega,\sfa,\sfb})
\in\su(n(\omega,\sfa,\sfb))^3
\]
such that 
$[A_{i,\omega,\sfa,\sfb},A_{j,\omega,\sfa,\sfb}]
+A_{k,\omega,\sfa,\sfb}=0$
for any cyclic permutation $(i,j,k)$ of $(1,2,3)$.
We assume 
\[
 \vecA_{\omega,\sfa,\sfb}
=\vecA_{\omega,\sfa,\sfb+\omega\mu_1^{\lor}}.
\]
Let $\nbigstilde_{\omega}\subset\real\times\cnum$
be a lift of $\nbigs_{\omega}$,
i.e., the projection
$\real\times\cnum\lrarr
\real\times(\cnum/\Gamma_{k(\omega)}^{\lor})$
induces a bijection
$\pi:\nbigstilde_{\omega}\simeq\nbigs_{\omega}$.
For each $(\sfatilde,\sfbtilde)\in\nbigstilde_{\omega}$,
we set
$\vecA_{\omega,\sfatilde,\sfbtilde}:=
 \vecA_{\omega,\pi(\sfatilde,\sfbtilde)}$.

We obtain the following monopole on 
$\nbigu_{k(\omega),\nu,R}$ $(\nu=0,\infty)$:
\begin{equation}
\label{eq;19.1.7.2}
\ttM_{k(\omega),\nu}(\omega,\nbigs_{\omega},\{\vecA_{\omega,\sfa,\sfb}\})
:=
 \bigoplus_{(\sfatilde,\sfbtilde)\in\nbigstilde_{\omega}}
 \vecsfL_{k(\omega)}(\omega)
\otimes
 L_{0}(0,\sfatilde,\sfbtilde)
\otimes
 \nbigv(0,\vecA_{\omega,\sfatilde,\sfbtilde}).
\end{equation}
Recall that we have the isomorphism
$L_0(0,\sfatilde,\sfbtilde+n_1\mu_1^{\lor}+n_2\mu_2^{\lor})
\simeq
L_0(0,\sfatilde,\sfbtilde)$
as explained in \S\ref{subsection;19.1.7.1}.
We also have the isomorphism
$\tte_1^{\ast}\vecsfL_{k(\omega)}
\simeq
 \vecsfL_{k(\omega)}(\omega)
\otimes
 L_{k(\omega)}(0,0,-\mu_1^{\lor})$
as in Lemma \ref{lem;19.1.5.2}.
By the isomorphisms,
the monopole
$\ttM_{k(\omega),\nu}(\omega,\nbigs_{\omega},\{\vecA_{\omega,\sfa,\sfb}\})$
is naturally equivariant with respect to 
the action of $(\seisuu/k(\omega)\seisuu)\tte_1$.
We obtain monopoles
\[
 \ttM_{\nu}(\omega,\nbigs_{\omega},\{\vecA_{\omega,\sfatilde,\sfbtilde}\})
\]
on $\nbigu_{1,\nu,R}$
as the descent of
$\ttM_{k(\omega),\nu}(\nbigs_{\omega},\{\vecA_{\omega,\sfa,\sfb}\})$.
By taking the direct sum,
we obtain a monopole
\[
 \ttM_{\nu}(I,\{\nbigs_{\omega}\},\{\vecA_{\omega,\sfatilde,\sfbtilde}\})
:=
 \bigoplus_{\omega\in I}
\ttM_{\nu}(\omega,\nbigs_{\omega},\{\vecA_{\omega,\sfatilde,\sfbtilde}\})
\]
on $\nbigu_{1,\nu,R}$.

\section{Asymptotic behaviour of doubly periodic monopoles}
\label{subsection;18.12.17.50}

\subsection{Statements}

Let $(y_0,y_1,y_2)$ be the standard coordinate
of $\real^3$.
We consider the Euclidean metric
$\sum_{i=0,1,2} dy_i\,dy_i$.
Let $\Gamma\subset\{0\}\times\real^2$
be a lattice.
The volume of $\real^2/\Gamma$
is denoted by $\Vol(\Gamma)$.
We may assume that
$\Gamma=\seisuu\cdot(0,a,0)\oplus
 \seisuu\cdot(0,b,c)$,
where $a$ and $c$ are positive numbers.
We consider the action of $\seisuu \tte_1\oplus\seisuu \tte_2$
on $\real^3$ by
$\tte_1(y_0,y_1,y_2)=(y_0,y_1+a,y_2)$
and 
$\tte_2(y_0,y_1,y_2)=(y_0,y_1+b,y_2+c)$.

For any $R\in\real$,
we set $\nbigutilde_{R}:=\{(y_0,y_1,y_2)\in\real^3\,|\,y_0<-R\}$.
Let $\nbigu_{R}$ denote the quotient space of
$\nbigutilde_{R}$
by the action of $\seisuu \tte_1\oplus\seisuu\tte_2$.

\vspace{.1in}
Let $(E,h,\nabla,\phi)$ be a monopole
on $\nbigu_{R_0}$ for some $R_0>0$.
By the pull back,
we obtain the $\seisuu\tte_1\oplus\seisuu\tte_2$-equivariant
monopole
$(\Etilde,\htilde,\nablatilde,\phitilde)$
on $\nbigutilde_{R_0}$.
\begin{assumption}
We assume that the curvature 
$F(\nabla)$ is bounded.
It particularly implies
$|\phi|_h=O(|y_0|)$.
\end{assumption}

\subsubsection{First reduction}
\label{subsection;18.9.1.11}

We shall prove the following proposition in
\S\ref{subsection;18.11.23.10}.
\begin{prop}
\label{prop;18.8.14.1}
There exists a finite subset $I(\phi)\subset\rnum$,
and positive numbers $R_1>0$ and $C_1>0$
such that the following holds
for $(y_0,y_1,y_2)\in\nbigu_{R_1}$:
\begin{itemize}
\item
For any eigenvalue $\alpha$
of $\phi_{|(y_0,y_1,y_2)}$,
there exists $\omega\in I(\phi)$
such that 
\begin{equation}
\label{eq;18.8.14.2}
 \Bigl|
 \alpha
-\frac{2\pi\sqrt{-1}\omega y_0}{\Vol(\Gamma)}
 \Bigr|<C_1.
\end{equation}
\end{itemize}
In particular,
if $R_1>0$ is sufficiently large,
we obtain the orthogonal decomposition
\begin{equation}
\label{eq;18.9.1.10}
 (E,h,\phi)_{|\nbigu_{R_1}}=
 \bigoplus_{\omega\in I(\phi)}
 (E^{\bullet}_{\omega},h^{\bullet}_{\omega},\phi^{\bullet}_{\omega})
\end{equation}
such that 
any eigenvalue of $\phi^{\bullet}_{\omega|(y_0,y_1,y_2)}$
satisfies {\rm(\ref{eq;18.8.14.2})}.
\end{prop}

We obtain a decomposition
$\nabla=\nabla^{\bullet}+\rho$,
where $\nabla^{\bullet}$ is a direct sum
of unitary connections
$\nabla^{\bullet}_{\omega}$ on $E^{\bullet}_{\omega}$,
and $\rho$ is a section of
$\bigoplus_{\omega_1\neq\omega_2}
 \Hom(E^{\bullet}_{\omega_1},E^{\bullet}_{\omega_2})
 \otimes\Omega^1$.
The inner product of $\rho$ and $\del_{y_i}$
are denoted by $\rho_i$.
Similarly, for any section $s$ of
$\End(E)\otimes\Omega^p$,
we obtain a decomposition
$s=s^{\bullet}+s^{\top}$,
where $s^{\bullet}$ is a section of
$\bigoplus \End(E^{\bullet}_{\omega})\otimes\Omega^p$,
and 
$s^{\top}$ is a section of
$\bigoplus_{\omega_1\neq\omega_2}
 \Hom(E^{\bullet}_{\omega_1},E^{\bullet}_{\omega_2})\otimes\Omega^p$.
Note that
$(\nabla\phi)^{\bullet}
=\nabla^{\bullet}\phi$
and 
$(\nabla\phi)^{\top}=[\rho,\phi]$.

We shall prove the following proposition
in \S\ref{subsection;18.11.23.11}.
\begin{thm}
\label{thm;18.8.15.1}
There exist positive constants $R_2$, $C_2$ and $\epsilon_2$
such that 
$|\rho|_h\leq 
 C_2\exp(-\epsilon_2y_0^2)$ on $\nbigu_{R_2}$.
Moreover,
for any positive integer $k$,
there exist positive constants $C_2(k)$
and $\epsilon_2(k)$
such that 
\[
 \bigl|
 \nabla^{\bullet}_{\kappa_1}\circ
 \cdots\circ
 \nabla^{\bullet}_{\kappa_{k}}
 \rho
 \bigr|_h
\leq
 C_2(k)\cdot
 \exp\bigl(
 -\epsilon_2(k)y_0^2
 \bigr) 
\]
on $\nbigu_{R_2}$
for any $(\kappa_1,\ldots,\kappa_k)\in\{0,1,2\}^k$.
\end{thm}

As a direct consequence,
we obtain the following corollary.
\begin{cor}
For any $k$, there exist positive constants
$C_3(k)$ and $\epsilon_3(k)$
such that
\[
 \bigl|
 \nabla^{\bullet}_{\kappa_1}\circ\cdots\circ\nabla^{\bullet}_{\kappa_k}
 (F(\nabla^{\bullet})-\nabla^{\bullet}\phi)
 \bigr|_h
\leq C_3(k)
 \exp(-\epsilon_3(k)y_0^2)
\]
on $\nbigu_{R_2}$
for any $(\kappa_1,\ldots,\kappa_k)\in\{0,1,2\}^k$.
Moreover,
\[
 \bigl|
  \nabla^{\bullet}_{\omega,\kappa_1}\circ
 \cdots\circ \nabla^{\bullet}_{\omega,\kappa_k}
 \nabla^{\bullet}_{\omega} \phi^{\bullet}_{\omega}
 \bigr|
+
 \bigl|
  \nabla^{\bullet}_{\omega,\kappa_1}\circ\cdots 
 \circ \nabla^{\bullet}_{\omega,\kappa_k}
 F(\nabla^{\bullet}_{\omega})
 \bigr|
\]
is bounded on $\nbigu_{R_2}$
for any $(\kappa_1,\ldots,\kappa_k)\in\{0,1,2\}^k$.
\hfill\qed
\end{cor}

For each $\omega\in I(\phi)$,
let $p$ be determined by
$\min\{p'\in\seisuu_{>0}\,|\,p'\omega\in\seisuu\}$.
For any $R>0$,
let $\nbigu_{p,R}$ denote the quotient of
$\nbigutilde_R$ by the action of 
$p\seisuu \tte_1\oplus \seisuu \tte_2$.
Let $\sfp_p:\nbigu_{p,R}\lrarr\nbigu_R$
denote the projection.
On $\nbigu_{p,R_1}$,
we set 
\begin{equation}
\label{eq;19.1.6.11}
 (E_{\omega},h_{\omega},\nabla_{\omega},\phi_{\omega}):=
 \sfp_p^{-1}(E^{\bullet}_{\omega},h^{\bullet}_{\omega},
 \nabla^{\bullet}_{\omega},\phi^{\bullet}_{\omega})
 \otimes
 \vecsfL_p(-\omega). 
\end{equation}

\begin{prop}
\label{prop;18.8.27.1}
For any $k\in\seisuu_{\geq 0}$
and for any
$(\kappa_1,\ldots,\kappa_k)\in \{0,1,2\}^k$,
we obtain
\[
 \Bigl|
\nabla_{\omega,\kappa_1}\circ\cdots
 \circ\nabla_{\omega,\kappa_k}
 \bigl(
F(\nabla_{\omega})
 \bigr)
\Bigr|_{h_{\omega}}
+
 \Bigl|
\nabla_{\omega,\kappa_1}\circ\cdots
 \circ\nabla_{\omega,\kappa_k}
 \bigl(
 \nabla_{\omega}\phi_{\omega}
 \bigr)
 \Bigr|_{h_{\omega}}
\lrarr 0
\]
as $|y_0|\lrarr\infty$.
\end{prop}

\subsubsection{Second reduction}
\label{subsection;18.9.1.12}

For any $R>0$,
we set $\nbigh_R:=\{y_0\in\real\,|\,y_0<-R\}\subset\real$.
Let $\Psi:\nbigu_{p,R}\lrarr\nbigh_R$
denote the projection.
Let $\nbiga$ be the ring of 
the non-commutative polynomials
of four variables.
We obtain the following proposition
from Proposition \ref{prop;18.8.27.21},
Proposition \ref{prop;18.8.27.22},
and Proposition \ref{prop;18.8.27.23}
below.

\begin{prop}
\label{prop;19.1.5.1}
There exist a finite subset $S_{\omega}\subset\real^3$,
a graded vector bundle
$V_{\omega}=\bigoplus_{\veca\in S_{\omega}}
 V_{\omega,\veca}$ on $\nbigh_{R}$,
a graded Hermitian metrics
$h_{V_{\omega}}=\bigoplus_{\veca\in S_{\omega}}
 h_{V_{\omega,\veca}}$,
a graded unitary connection
$\nabla_{V_{\omega}}=
 \bigoplus_{\veca\in S_{\omega}}\nabla_{V_{\omega,\veca}}$,
graded anti-Hermitian endomorphisms
$\phi_{i,\omega}=
 \bigoplus_{\veca\in S_{\omega}} \phi_{i,\omega,\veca}$
$(i=1,2,3)$,
and an isomorphism
$E_{\omega}\simeq
 \Psi^{-1}(V_{\omega})$
such that the following holds:
\begin{itemize}
\item
 Let $b_{\omega}$ be the automorphism of $E_{\omega}$
 determined by 
 $h_{\omega}=\Psi^{-1}(h_{V_{\omega}})b_{\omega}$.
 Then, for any $P\in\nbiga$,
there exists $\epsilon(P)>0$ such that
\[
 \bigl|
 P(\nabla_{\omega,y_0},\nabla_{\omega,y_1},\nabla_{\omega,y_2},\phi_{\omega})
 (b_{\omega}-\id)
 \bigr|
=O\Bigl(
 e^{\epsilon(P)y_0}
 \Bigr).
\]
\item
For any $P\in\nbiga$,
there exists $\epsilon(P)>0$ such that 
\[
 \Bigl|
  P(\nabla_{\omega,y_0},\nabla_{\omega,y_1},\nabla_{\omega,y_2},\phi_{\omega}) 
 (\phi_{\omega}-\Psi^{-1}(\phi_{3,\omega}))
 \Bigr|
=O\Bigl(
 e^{\epsilon(P)y_0}
 \Bigr).
\]
\item
We set
$\ttR_{\omega,i}:=\nabla_{\omega,y_i}-
 (\del_{y_i}+\Psi^{-1}(\phi_{\omega,i}))$
$(i=1,2)$,
where $\del_{y_i}$ are the naturally induced
operators of $\Psi^{-1}(V_{\omega})$.
Then,  for any $P\in \nbiga$,
there exists $\epsilon(P)>0$ such that
\[
 \Bigl|
  P(\nabla_{\omega,y_0},\nabla_{\omega,y_1},\nabla_{\omega,y_2},\phi_{\omega}) 
 \ttR_{\omega,i}
 \Bigr|
=O\Bigl(
 e^{\epsilon(P)y_0}
 \Bigr).
\]
\item
There exist anti-Hermitian endomorphisms
$A_{i,\omega,\veca}$ $(i=1,2,3)$
of $V_{\omega,\veca}$
such that
$\nabla_{V_{\veca,\omega}}A_{i,\omega,\veca}=0$
and 
\[
 \phi_{i,\omega,\veca}=
 \sqrt{-1}a_i\id_{V_{\omega,\veca}}
 +y_0^{-1}A_{i,\omega,\veca}
+O(y_0^{-2}).
\]
Moreover, 
$\bigl[
 A_{i,\omega,\veca},
 A_{j,\omega,\veca}
 \bigr]
+A_{k,\omega,\veca}=0$ holds
for any cyclic permutation $(i,j,k)$ of $(1,2,3)$.
\end{itemize}
\end{prop}

Set $e_1:=(a,b)$ and $e_2:=(0,c)$.
Let $e_i^{\lor}\in\real^2$ 
be determined by
$(e_i^{\lor},e_i)=\pi$
and $(e_i^{\lor},e_j)=0$ $(i\neq j)$.
More explicitly,
\[
 e_1^{\lor}=(a^{-1}\pi,-\Vol(\Gamma)^{-1}b\pi),
\quad
e_2^{\lor}=(0,a\Vol(\Gamma)^{-1}\pi).
\]
Let $\Gamma^{\lor}_p:=
 \seisuu p^{-1}e_1^{\lor}\oplus \seisuu e_2^{\lor}$.
There exists the action of
$(\seisuu/p\seisuu ) \omega e_2^{\lor}$
on $\real^2/\Gamma_p^{\lor}$
induced by 
$\omega e_2^{\lor}\bullet (a_1,a_2)=
(a_1,a_2)+\omega e_2^{\lor}$.
The following will be clear by the choice of $S_{\omega}$.
\begin{prop}
Let $[S_{\omega}]\subset (\real^2/\Gamma^{\lor}_p)\times\real$
denote the image of $S_{\omega}$ by the projection
$\real^3\lrarr (\real^2/\Gamma^{\lor}_p)\times\real$.
Then, $[S_{\omega}]$ is well defined 
for $(E,h,\nabla,\phi)$,
and 
$[S_{\omega}]$ is naturally preserved 
by the above action of 
$(\seisuu/p\seisuu) \omega e_2^{\lor}$.
Moreover,
if $\veca\equiv \omega e_2^{\lor}\bullet\veca'$
in $(\real^2/\Gamma_p^{\lor})\times\real$,
then
$\vecA_{\omega,\veca}=\vecA_{\omega,\veca'}$
holds.
\end{prop}

\subsubsection{A consequence}

We obtain the following consequence.

\begin{cor}
\label{cor;18.12.17.3}
We obtain 
$\bigl|
 \nabla_{y_1}\phi
 \bigr|
+\bigl|
 \nabla_{y_2}\phi
 \bigr|_h
=O\bigl(y_0^{-2}\bigr)$.
Equivalently,
we obtain 
$\bigl|
 F(\nabla)_{y_0,y_i}
 \bigr|_h=O\bigl(y_0^{-2}\bigr)$
 $(i=1,2)$.
\hfill\qed
\end{cor}

\begin{rem}
Note that $\nabla_{y_0}\phi$ is not necessarily
$O(y_0^{-2})$.
Equivalently,
$|F(\nabla)_{y_1,y_2}|_h$ is not necessarily
$O(y_0^{-2})$.
See the examples in 
{\rm\S\ref{subsection;18.11.24.30}.}
\hfill\qed
\end{rem}

\subsection{Vector bundles with a connection on $S^1$}

\subsubsection{Statement}

Let $r$ be a positive integer.
Let $C_0>0$ be a constant.
Let $A_0$ be an $r$-square Hermitian matrix.
Set $S^1:=\real/\seisuu$.
Let $A_1:S^1\lrarr M_r(\cnum)$ be a continuous function
such that $|A_1|\leq C_0$.
Let $V$ be a $C^{\infty}$-vector bundle of rank $r$
on $S^1$
with a frame $\vecv$.
We have the connection $\nabla$
determined by
\[
 \nabla\vecv=\vecv\cdot (A_0+A_1)\,dt,
\]
where $t$ is the standard coordinate of $\real$,
and $dt$ is the induced $1$-form on $S^1$.
We have the monodromy 
$M(A_0+A_1):V_{|0}\lrarr V_{|1}=V_{|0}$
of the connection $\nabla$,
and let $\Sp(M(A_0+A_1))$ denote 
the set of eigenvalues.
We shall prove the following proposition
in \S\ref{subsection;18.8.27.2}--\ref{subsection;18.8.27.3}.
\begin{prop}
\label{prop;18.8.11.2}
There exists $R>0$ depending only on $C_0$
such that the following holds.
\begin{itemize}
\item
 For any $\alpha\in \Sp(M(A_0+A_1))$,
 there exists
 $\beta\in \Sp(M(A_0))$
 such that 
$|\alpha\beta^{-1}|\leq R$
and 
$|\alpha^{-1}\beta|\leq R$
Conversely, for any $\alpha\in \Sp(M(A_0))$,
 there exists
 $\beta\in \Sp(M(A_0+A_1))$
 such that 
$|\alpha\beta^{-1}|\leq R$
and 
$|\alpha^{-1}\beta|\leq R$.
\end{itemize}
\end{prop}

\subsubsection{Decomposition of a finite tuple of real numbers}
\label{subsection;18.8.27.2}

We consider a finite tuple
$(a_1,\ldots,a_N)$ of real numbers.
We assume $a_i\leq a_j$ for $i<j$.
We fix a positive number $c_0>0$.
We take any $c_1>10N$.

\begin{lem}
\label{lem;18.8.11.1}
There exist $k\geq 0$ and a decomposition
$\{1,\ldots,N\}=\coprod_{\ell=1}^m\gbigi_{\ell}$
such that the following holds.
\begin{itemize}
\item
 If $i,j\in \gbigi_{\ell}$,
 then
 $|a_i-a_j|\leq 3Nc_1^{k}c_0$.
\item
If $i\in \gbigi_{\ell_1}$ and $j\in \gbigi_{\ell_2}$
with $\ell_1\neq\ell_2$,
then
 $|a_i-a_j|\geq\frac{1}{2}c_1^{k+1}c_0$.
\end{itemize}
\end{lem}
\pf
We set $m(0):=N$.
We shall construct a finite decreasing sequence
$m(0)>m(1)>\ldots>m(k)$,
order preserving injective maps
$G_n:\{1,\ldots,m(n)\}\lrarr\{1,\ldots,N\}$
$(n=0,\ldots,k)$,
and order preserving surjective maps
$F_n:\{1,\ldots,m(n)\}
\lrarr \{1,\ldots,m(n+1)\}$
$(n=0,\ldots,k-1)$
by an inductive procedure.
Suppose that we have already constructed
$m(n)$,
$G_n:\{1,\ldots,m(n)\}\lrarr \{1,\ldots,N\}$.
We set
$J^{(n)}:=\{i\,|\,a_{G_n(i+1)}-a_{G_n(i)}>c_1^{n+1}c_0\}
 \cup\{m(n)\}$.
If $J^{(n)}:=\{1,\ldots,m(n)\}$,
we stop the procedure.
If $J^{(n)}\neq\{1,\ldots,m(n)\}$,
we set
$m(n+1):=|J^{(n)}|$.
We have the natural order preserving bijection
$\varphi_{n+1}:\{1,\ldots,m(n+1)\}\simeq J^{(n)}$.
Because $J^{(n)}\subset \{1,\ldots,m(n)\}$,
we obtain an injection
$G_{n+1}:\{1,\ldots,m(n+1)\}
\lrarr\{1,\ldots,N\}$
from $\varphi_{n+1}$ and $G_{n}$.
For $i\in\{1,\ldots,m(n)\}$,
there exists $j\in\{1,\ldots,m(n+1)\}$
such that  $\varphi_{n+1}(j-1)<i\leq\varphi_{n+1}(j)$,
where we formally set $\varphi_{n+1}(0)=0$.
We define $F_{n+1}(i)=j$ for such $i$ and $j$.
Thus, we obtain the order preserving surjection
$F_{n+1}:\{1,\ldots,m(n)\}\lrarr\{1,\ldots,m(n+1)\}$.
The procedure will stop after finite steps.

By the construction,
$|a_{G_k(i)}-a_{G_k(j)}|>c_1^{k+1}c_0$ holds
for $i,j\in\{1,\ldots,m(k)\}$
with $i\neq j$.
Let $F:\{1,\ldots,N\}\lrarr\{1,\ldots,m(k)\}$
be the map obtained as the composite of
$F_0,\ldots,F_{k-1}$.
For $\ell\in F^{-1}(i)$,
the following holds:
\[
  \bigl|
 a_{\ell}-a_{G_k(i)}
 \bigr|
\leq
 N(c_1+\cdots c_1^{k})c_0
=N(c_1^{k+1}-c_1)(c_1-1)^{-1}c_0.
\]
Hence, if $\ell_1,\ell_2\in F^{-1}(i)$,
then
\[
 \bigl|
 a_{\ell_1}-a_{\ell_2}
 \bigr|
\leq
 2N(c_1^{k+1}-c_1)(c_1-1)^{-1}c_0
\leq
 3Nc_1^kc_0.
\]
For $\ell_1\in F^{-1}(j)$
and $\ell_2\in F^{-1}(i)$ with $i\neq j$,
the following holds:
\[
  \bigl|
 a_{\ell_1}-a_{\ell_2}
 \bigr|
\geq
 c_1^{k+1}c_0-
 2N(c_1^{k+1}-c_1)(c_1-1)^{-1}c_0
\geq
 \frac{1}{2}c_1^{k+1}c_0.
\]
Thus, we are done.
\hfill\qed

\subsubsection{An estimate}

Let $a$ be a non-zero real number.
For any $C^0$-function $g$ on $S^1$,
we have a unique $C^1$-function $f$
such that $(\del_t+a)f=g$.

\begin{lem}
\label{lem;18.12.14.10}
We have $\sup|f|\leq 2|a|^{-1}\sup|g|$.
\end{lem}
\pf
It  is enough to consider the case $a>0$.
Let $f=\sum f_{n}e^{2\pi\sqrt{-1}n\theta}$
and $g=\sum g_ne^{2\pi\sqrt{-1}n\theta}$
be the Fourier expansions.
Because $(2\pi\sqrt{-1}n+a)f_n=g_n$,
we obtain 
$\int_{0}^1|f|^2dt=\sum|f_n|^2
\leq a^{-2}\sum|g_n|^2=a^{-2}\int_0^1|g|^2\,dt$.
Hence, there exists $t_0\in S^1$
such that 
$|f(t_0)|\leq a^{-1}\sup|g|$.
We may assume that $t_0=0$
by a coordinate change.
Because 
$\del_t(e^{at}f)=e^{at}g$,
we have
\[
 \Bigl|
 e^{at}f(t)-f(0)
 \Bigr|
\leq
 \int_0^1e^{as}|g(s)|\,ds
\leq
 \sup|g|\cdot a^{-1}e^{at}.
\]
Hence, we obtain the claim of the lemma.
\hfill\qed

\subsubsection{Solving a non-linear equation}

Let $m$ be a positive integer.
Let $D_0$ be an $m$-square Hermitian matrix.
Let $C_{10}$ be a positive constant.
Let $B_0(t)$ be a $C^0$-map $S^1\lrarr \cnum^m$
such that $|B_0(t)|\leq C_{10}/3$.
Let $B_1(t)$ be a $C^0$-map $S^1\lrarr M_m(\cnum)$
such that $|B_1(t)|\leq C_{10}/3$.
Let $B_2(t,x)$ be a $C^0$-map
$S^1\times \cnum^m\lrarr \cnum^m$
such that the following holds.
\begin{itemize}
\item $|B_2(t,x)|=o(|x|)$ as $|x|\to 0$.
\item
For any $\epsilon>0$,
there exists $\delta>0$ such that
$|B_2(t,x)-B_2(t,y)|\leq\epsilon|x-y|$
if $\max\{|x|,|y|\}\leq \delta$.
\end{itemize}

We take $T>1$ such that the following holds.
\begin{itemize}
\item
If $|x|<T^{-1}$, then $|B_2(t,x)|\leq C_{10}/3$.
\item
If $\max\{|x|,|y|\}<T^{-1}$,
then
$|B_2(t,x)-B_2(t,y)|\leq C_{10}|x-y|/3$.
\end{itemize}

\begin{lem}
\label{lem;18.8.27.10}
Assume that any eigenvalues $a$ of $D_0$
satisfies
$|a|\geq 10mTC_{10}$.
Then, there exists
$f:S^1\lrarr \cnum^m$
such that
(i) $(\del_t+D_0)f(t)+B_0(t)+B_1(t)\cdot f(t)+B_2(t,f(t))=0$,
(ii) $|f|\leq T^{-1}$.
Such a function $f$ is unique.
\end{lem}
\pf
We take any $C^0$-function $f_0:S^1\lrarr\cnum^m$
such that $|f_0|\leq T^{-1}$.
Inductively, we define $f_i$
as a unique solution of
$(\del_t+D_0)f_i(t)+B_0(t)+B_1(t)f_{i-1}(t)+B_2(t,f_{i-1}(t))=0$.
Because
$|B_1(t)f_{i-1}(t)|\leq C_{10}/3$
and 
$|B_2(t,f_{i-1}(t))|\leq C_{10}/3$,
we obtain $|f_i|\leq (C_{10}T)^{-1}C_{10}\leq T^{-1}$
by Lemma \ref{lem;18.12.14.10}.
Note that
\[
 (\del_t+D_0)(f_{i+1}(t)-f_i(t))
+B_1(t)(f_i(t)-f_{i-1}(t))+
 B_2(t,f_i(t))-B_2(t,f_{i-1}(t))=0.
\]
Because
$|B_1(t)(f_i(t)-f_{i-1}(t))+B_2(t,f_i(t))-B_2(t,f_{i-1}(t))|
\leq C_{10}|f_i(t)-f_{i-1}(t)|$,
we obtain
$\sup|f_{i+1}-f_{i}|\leq T^{-1}\sup|f_{i}-f_{i-1}|$
by Lemma \ref{lem;18.12.14.10}.
Hence, the sequence $f_i$ is convergent,
and the limit $f_{\infty}=\lim f_i$
satisfies the desired conditions.
We also obtain the uniqueness.
\hfill\qed

\subsubsection{Proof of Proposition \ref{prop;18.8.11.2}}
\label{subsection;18.8.27.3}

We may assume that $A_0$ is diagonal.
Let $a_i$ denote the $(i,i)$-th entries.
We may assume that $a_i\leq a_j$ for $i\leq j$.
Take a sufficiently large constant $C_1$.
We have $k\geq 0$ and a decomposition
$\{1,\ldots,r\}=\coprod\gbigi_{\ell}$
as in Lemma \ref{lem;18.8.11.1}.
We choose $i(\ell)\in\gbigi_{\ell}$,
and set $\alpha_{\ell}:=a_{i(\ell)}$.
We put $r(\ell):=|\gbigi_{\ell}|$.
We set
$\Atilde_{0}:=\bigoplus \alpha_{\ell}I_{r(\ell)}$
and 
$\Atilde_{1}:=A_0-\Atilde_0+A_1$.
We have
$|\Atilde_{1}|\leq 4rC_1^kC_0$.

According to the decomposition
$\{1,\ldots,r\}=\coprod \gbigi_{\ell}$,
we have the decomposition
$\cnum^r=\bigoplus \cnum^{r(\ell)}$.
It induces
$\End(\cnum^r)
=\bigoplus_{\ell}
 \End(\cnum^{r(\ell)})
\oplus
 \bigoplus_{\ell_1\neq\ell_2}
 \Hom(\cnum^{r(\ell_1)},\cnum^{r(\ell_2)})$.
For any matrix $D\in\End(\cnum^r)$,
we have the decomposition
$D=D^{\circ}+D^{\bot}$.

We consider the following equation
for $G:S^1\lrarr \GL(r,\cnum)$
and $U:S^1\lrarr \bigoplus \End(\cnum^{r(\ell)})$:
\begin{equation}
\label{eq;18.12.14.20}
 G^{-1}\circ
 \bigl(
 \del_t+\Atilde_0+\Atilde_1
 \bigr)\circ G
=\del_t+\Atilde_0+\Atilde_1^{\circ}+U.
\end{equation}
We impose that $G^{\circ}=I_{r}$,
and we regard (\ref{eq;18.12.14.20})
as an equation for $G^{\bot}$ and $U$.
It is equivalent to the following equations:
\[
 (\Atilde_1^{\bot}G^{\bot})^{\circ}=U,
\quad\quad
 \del_tG^{\bot}+[\Atilde_0,G^{\bot}]
+[\Atilde_1^{\circ},G]
+(\Atilde_1^{\bot}G^{\bot})^{\bot}
+G^{\bot}U=0.
\]
By eliminating $U$,
we obtain the following equation for $G^{\bot}$:
\[
 \del_tG^{\bot}
+[\Atilde_0,G^{\bot}]
+[\Atilde_1^{\circ},G^{\bot}]
+(\Atilde_1^{\bot}G^{\bot})^{\bot}
+\Atilde_1^{\bot}+G^{\bot}(\Atilde_1^{\bot}G^{\bot})^{\circ}=0.
\]
For a large $C_1$,
we set
$C_{10}:=400r^3C_1^kC_0$
and $T:=(1000r^3)^{-1}C_1$.
By using Lemma \ref{lem;18.8.27.10},
if $C_1$ is sufficiently large,
we have a solution 
$G^{\bot}$ with $|G^{\bot}|\leq T^{-1}$.
We also obtain $U$
such that $|U|\leq C_{10}T^{-1}$.

By considering the eigenvalues of the monodromy of
$\del_t+\Atilde_0+\Atilde_1^{\circ}+U$,
we obtain the claim of the proposition.

\hfill\qed

\subsection{First reduction}

\subsubsection{Proof of Proposition \ref{prop;18.8.14.1}}
\label{subsection;18.11.23.10}

We take the mini-holomorphic structure
determined by the decomposition
$\real^3=\real\cdot(0,a,0)\times
 \bigl(\real\cdot(0,a,0)\bigr)^{\bot}$.
We take
$\real^3\simeq \cnum\oplus\real$
given by
\[
 (y_0,y_1,y_2)\longmapsto
 \Bigl(
 2\pi c^{-1}(y_0+\sqrt{-1}y_2),y_1
 \Bigr).
\]
The action of $\seisuu \tte_1\oplus\seisuu \tte_2$
on $\cnum\times \real$
are described as
\[
 \tte_1(\zetatilde,y_1)=(\zetatilde,y_1+a),
\quad
 \tte_2(\zetatilde,y_1)
=(\zetatilde+2\pi\sqrt{-1},y_1+b).
\]

For any $R$,
we set
$\Utilde_{R}:=\bigl\{
 \zetatilde\in\cnum\,\big|\,
 \frac{c}{2\pi}\Re(\zetatilde)<-R
 \bigr\}$.
We have
$\nbigutilde_{R}= \Utilde_{R}\times\real$
under the above identification
$\real^3\simeq\cnum\times\real$.

We have the associated mini-holomorphic bundle
$(\Etilde,\delbar_{\Etilde})$ on $\nbigutilde_{R_0}$
with respect to the above mini-complex structure.
By considering the flat sections along $\{\zetatilde\}\times\real$
for each $\zetatilde$,
we obtain a holomorphic vector bundle $(\Vtilde,\delbar_{\Vtilde})$
on $\Utilde_{R_0}$.
The action of $\seisuu \tte_1$ induces a holomorphic automorphism
$\Ftilde$ of $\Vtilde$.
The action of $\seisuu \tte_2$ induces
an isomorphism
$\tte_2^{\ast}\Vtilde\simeq \Vtilde$,
where $\tte_2:\cnum\lrarr\cnum$ is given by
$\tte_2(\zetatilde)=\zetatilde+2\pi\sqrt{-1}$.

We identify
$\cnum/(2\pi\sqrt{-1}\seisuu)\simeq\cnum^{\ast}$
by $\zetatilde\longmapsto \zeta=e^{\zetatilde}$.
For any $R$, we set
$U_{R}:=\bigl\{\zeta\in\cnum^{\ast}\,\big|\,
 c\log|\zeta|<-R
 \bigr\}$.
We obtain the induced holomorphic bundle 
$(V,\delbar_V)$ on $U_{R_0}$.
Because the actions of $\tte_1$ and $\tte_2$ are commutative,
we obtain the induced automorphism $F$
of $(V,\delbar_V)$.
We obtain the spectral curve 
$\Sp(F)$ of $F$
contained in
$U_{R_0}\times\cnum^{\ast}$.
We set $\Ubar_{R_0}:=U_{R_0}\cup\{0\}$.

\begin{lem}
The closure $\overline{\Sp(F)}$ of $\Sp(F)$
in $\Ubar_{R_0}\times\proj^1$
is complex analytic.
\end{lem}
\pf
Let $\htilde$ denote the Hermitian metric of
$\Etilde$ induced by $h$.
Let $s$ be a flat section of 
$\Etilde_{|\{\zetatilde\}\times\real}$
with respect to
$\nabla_{y_1}-\sqrt{-1}\phi$.
Then, we have
\[
\del_{y_1}\htilde(s,s)
=\htilde\bigl(s, (\nabla_{y_1}+\sqrt{-1}\phi)s
 \bigr)
=\htilde\bigl(s,2\sqrt{-1}\phi s\bigr).
\]
Hence, there exists $C>0$, which is independent of $\zetatilde$,
such that
$\bigl|\del_{y_1}\htilde(s,s)\bigr|
\leq
 C\bigl|\Re(\zetatilde)\bigr|\cdot \htilde(s,s)$.
It implies that
$\Bigl|
 \log|\Ftilde|_{\htilde}\Bigr|=
O\Bigl(
 \bigl|\Re(\zetatilde)\bigr|
 \Bigr)$.
Then, we obtain the claim of the lemma.
\hfill\qed

\vspace{.1in}
By replacing $R_{0}$ with a larger number,
we may assume to have the decomposition
\begin{equation}
\label{eq;18.12.14.30}
 \overline{\Sp(F)}
=\coprod_{\omega\in\rnum}
 \overline{\Sp(F)}_{\omega},
\end{equation}
where $\overline{\Sp(F)}_{\omega,\alpha}\lrarr \Ubar_{R_0}$
are the union of graphs of ramified meromorphic functions $g$
such that $|\zeta|^{\omega} g$ are bounded.

The group $\seisuu\tte_1$ acts on
$\{\zetatilde\}\times\real$.
Let $S^1_{\zetatilde,a}$ denote the quotient space.
For $\zetatilde\in \Utilde_{R_0}$,
there exists a naturally induced injection
$S^1_{\zetatilde,a}\lrarr \nbigu_{R_0}$.
We obtain the induced vector bundle
$E^{\zetatilde}$ on $S^1_{\zetatilde}$
with the metric $h^{\zetatilde}$,
the unitary connection
$\nabla^{\zetatilde}$,
and the anti-Hermitian endomorphism
$\phi^{\zetatilde}$.
There exists an orthonormal frame 
$\vecu=(u_1,\ldots,u_r)$ of $E^{\zetatilde}$
such that the following holds.
\begin{itemize}
\item
 There exists a constant anti-Hermitian matrix $A$
 such that
 $\nabla^{\zetatilde}\vecu=\vecu\cdot A$.
 Moreover, the eigenvalues of $A$
 are contained in 
 $\{\sqrt{-1}\rho\,|\,0\leq \rho a\leq 2\pi\}$.
\end{itemize}
Because $\nabla^{\zetatilde}\phi^{\zetatilde}$
is bounded independently from $\zetatilde$,
there exists a constant $C_{10}$,
which is independent of $\zetatilde$,
and a decomposition
$\phi^{\zetatilde}=\psi_0+\psi_1$,
such that the following holds.
\begin{itemize}
\item
 There exists a constant anti-Hermitian matrix $\Psi_0$
 such that
 $\psi_0\vecu=\vecu\Psi_0$.
\item
 $|\psi_1|_{h^{\zetatilde}}\leq C_{10}$.
\end{itemize}
Let $\nbigs(\zetatilde)$ be the set of the eigenvalues of
$\Psi_0$.
Then, there exists $C_{11}>0$,
which is independent of $\zetatilde$,
such that 
the following holds for any $y_1$.
\begin{itemize}
\item
 For any eigenvalue $\alpha$ of
 $\phitilde_{|(\zetatilde,y_1)}$,
 there exists $\beta\in\nbigs(\zetatilde)$
 such that
 $|\alpha-\beta|<C_{11}$.
 Conversely, for any $\beta\in\nbigs(\zetatilde)$,
 there exists an eigenvalue $\alpha$
 of  $\phitilde_{|(\zetatilde,y_1)}$
 such that 
 $|\alpha-\beta|<C_{11}$.
\end{itemize}
By Proposition \ref{prop;18.8.11.2},
there exists $C_{12}>0$,
which is independent of $\zetatilde$,
such that the following holds.
\begin{itemize}
\item
 For any eigenvalue $\gamma$ of $F_{|\zetatilde}$,
 there exists $\beta\in\nbigs(\zetatilde)$
 such that
 $\bigl|
 \log|\gamma|+a\sqrt{-1}\beta
 \bigr|<C_{12}$.
\end{itemize}
Then, the claim of Proposition \ref{prop;18.8.14.1}
follows from the decomposition (\ref{eq;18.12.14.30}).
\hfill\qed

\subsubsection{Proof of Theorem \ref{thm;18.8.15.1}}
\label{subsection;18.11.23.11}

\begin{lem}
For any $k\in\seisuu_{\geq 0}$
and  any $(\kappa_1,\ldots,\kappa_k)\in\{0,1,2\}^k$,
$\bigl|
 \nabla_{\kappa_1}\circ\cdots
 \circ\nabla_{\kappa_k}(\nabla\phi)
 \bigr|_h$ is bounded on $\nbigu_{2R_0}$.
\end{lem}
\pf
Take a positive number $\epsilon_0>0$,
and we take $\epsilon_1>0$ such that
$2C_0\epsilon_1^2<\epsilon_0$.
For any $(y_0,y_1,y_2)\in \nbigu_{2R_0}$,
let 
$S_{y_0,y_1,y_2}=\big\{
 (z_0,z_1,z_2)\,|\,|z_0-y_0|<\epsilon_1\}$.
We have
$G_{y_0,y_1,y_2}:
 \{(x_0,x_1,x_2)\in|\,|x|<1\}\lrarr
 S_{y_0,y_1,y_2}$
by
$(x_0,x_1,x_2)\longmapsto
 (y_0,y_1,y_2)+\epsilon_1(x_0,x_1,x_2)$.
We have 
$\bigl|
 G_{y_0,y_1,y_2}^{-1}F(h)
\bigr|\leq \epsilon_0$
and 
$\bigl|
 \epsilon_1G_{y_0,y_1,y_2}^{-1}\nabla\phi
\bigr|\leq \epsilon_0$.
Set $\nabla':=G_{y_0,y_1,y_2}^{-1}(\nabla)$.
For any $k$
and $(\kappa_1,\ldots,\kappa_k)\in\{0,1,2\}^k$,
there exists $B_1(k)$
which is independent of $(y_0,y_1,y_2)$,
such that 
\[
 \Bigl|
 \nabla'_{x_{\kappa_1}}\circ
 \cdots
 \nabla'_{x_{\kappa_k}}
 G_{y_0,y_1,y_2}^{-1}(\nabla\phi)
 \Bigr|
\leq
 B_1(k).
\]
Then, we obtain the desired estimate for
the derivatives of $\nabla\phi$.
\hfill\qed

\vspace{.1in}
We obtain the following lemma
as in the case of 
\cite[Lemma 6.15]{Mochizuki-difference-modules}.
\begin{lem}
We have 
$|\rho_{\kappa}|=
 O\bigl(
 \bigl|
 (\nabla_{\kappa}\phi)^{\top}
 \bigr|_h
 \bigr)$
for $\kappa=0,1,2$.
We also have 
\[
 \bigl|
 \nabla^{\bullet}_{\kappa_1}
 \rho_{\kappa_2}
 \bigr|_h
= O\Bigl(
 \bigl|(\nabla_{\kappa_2}\phi)^{\top}\bigr|
+\bigl|
 \nabla^{\bullet}_{\kappa_1}
 (\nabla_{\kappa_2}\phi)^{\top}
\bigr|
 \Bigr)
\]
for any $\kappa_1,\kappa_2\in\{0,1,2\}$.
\hfill\qed
\end{lem}

By the argument in \cite[\S6.3.3]{Mochizuki-difference-modules},
we obtain the following estimates:
\begin{equation}
\label{eq;18.12.15.1}
 h\Bigl(
 \nabla^2_{\kappa_1}(\nabla_{\kappa_2}\phi)^{\bullet},
 (\nabla_{\kappa_2}\phi)^{\top}
 \Bigr)
=O\Bigl(
 \Bigl(
 \bigl|(\nabla_{\kappa_1}\phi)^{\top}\bigr|_h
+
 \bigl|\nabla^{\bullet}_{\kappa_1}(\nabla_{\kappa_1}\phi)^{\top}
 \bigr|
 \Bigr)
 \cdot
 \bigl|
 (\nabla_{\kappa_2}\phi)^{\top}
 \bigr|
 \Bigr),
\end{equation}
\begin{equation}
\label{eq;18.12.15.2}
 \sum_{\kappa_1=0,1,2}
 h\bigl(
 \nabla_{\kappa_1}^2
 (\nabla_{\kappa_2}\phi),
 (\nabla_{\kappa_2}\phi)^{\top}
 \bigr)
=
 \Bigl|
 \bigl[\phi,(\nabla_{\kappa_2}\phi)^{\top}\bigr]
 \Bigr|^2_h
+
O\Bigl(
 \bigl|
 (\nabla\phi)^{\top}
 \bigr|_h
\cdot
 \bigl|(\nabla_{\kappa_2}\phi)^{\top}\bigr|
 \Bigr).
\end{equation}
By using the estimates (\ref{eq;18.12.15.1}), (\ref{eq;18.12.15.2})
and the argument in \cite[\S6.3.3]{Mochizuki-difference-modules},
we obtain the following inequality
on $\nbigu_{R_{20}}$ for some 
$C_{20}>0$ and $R_{20}>R_0$:
\[
  -(\del_{y_0}^2+\del_{y_1}^2+\del_{y_2}^2)
 \bigl|(\nabla\phi)^{\top}\bigr|_h^2
\leq
 -C_{20}
 \bigl|
 (\nabla\phi)^{\top}
 \bigr|_h^2
 \cdot y_0^2.
\]

Let $\int_{T^2}|(\nabla\phi)^{\top}|^2$
denote the function 
on $\nbigh_{R_0}$
obtained as the fiber integral 
of $|(\nabla\phi)^{\top}|^2$
with respect to
$\nbigu_{R_0}\lrarr\nbigh_{R_0}$.
We obtain the following inequality
on $\nbigh_{R_{20}}$:
\[
-\del_{y_0}^2\int_{T^2}\bigl|(\nabla\phi)^{\top}\bigr|_h^2
\leq
 -C_{20}y_0^2\int_{T^2}\bigl|(\nabla\phi)^{\top}\bigr|_h^2.
\]

We take $B_1>0$ such that
$\int_{T^2}\bigl|(\nabla\phi)^{\top}\bigr|_h^2
\leq 
 B_1e^{-(C_{11}^{1/2}/2)y_0^2}$
at $y_0=-R_{20}$.
For any $\delta>0$, we set
$F_{\delta}:=B_1e^{-(C_{11}^{1/2}/2)y_0^2}
-\delta(y_0+R_{20})$.
The following holds on $\nbigh_{R_{20}}$:
\[
 -\del_{y_0}^2F_{\delta}
\geq 
-C_{11}F_{\delta}.
\]
We also have
$\int_{T^2}
 \bigl|(\nabla\phi)^{\top}\bigr|_h^2
\leq
 F_{\delta}$
at $y_0=-R_1$.
By an argument as in Ahlfors lemma \cite{a, Simpson90},
for any $\delta>0$ we obtain 
\[
 \int_{T^2}\bigl|(\nabla\phi)^{\top}\bigr|_h^2
\leq 
 F_{\delta}
\]
on $\nbigh_{R_1}$.
Then, by taking the limit $\delta\to 0$,
we obtain
\[
 \int_{T^2}\bigl|(\nabla\phi)^{\top}\bigr|_h^2
\leq 
 B_1e^{-(C_{11}^{1/2}/2)y_0^2}.
\]
Then, by the argument in \cite[\S6.3.4]{Mochizuki-difference-modules},
we obtain the Theorem \ref{thm;18.8.15.1}.
\hfill\qed

\vspace{.1in}

\begin{cor}
For any $k$,
there exist positive constants $C(k)$ and $\epsilon(k)$
such that
\begin{equation}
\label{eq;18.8.14.21}
\Bigl|
 \nabla^{\bullet}_{\kappa_1}\circ\cdots\circ
 \nabla^{\bullet}_{\kappa_k}
 \Bigl(
 (\nabla_i^{\bullet})^2\nabla^{\bullet}_a\phi
-4\bigl[\nabla_b^{\bullet}\phi,\nabla_c^{\bullet}\phi\bigr]
+\bigl[
 \phi,[\phi,\nabla^{\bullet}_a\phi]
 \bigr]
 \Bigr)
\Bigr|_h
\leq
 C(k)e^{-\epsilon(k)y_0^2}.
\end{equation}
Here, $(a,b,c)$ is a cyclic permutation of $(0,1,2)$.
\end{cor}
\pf
Recall the following equalities:
\begin{equation}
\label{eq;18.12.15.100}
 \bigl(
 \nabla_0^2+\nabla_1^2+\nabla_2^2
 \bigr)
 (\nabla_i\phi)
=4[\nabla_j\phi,\nabla_k\phi]
-\bigl[
 \phi,[\phi,\nabla_i\phi]
 \bigr],
\end{equation}
where
$(i,j,k)$ is a cyclic permutation of $(0,1,2)$.
(For example,
see \cite[Lemma 6.16]{Mochizuki-difference-modules}.)
Then, the corollary follows from
Theorem \ref{thm;18.8.15.1}
and (\ref{eq;18.12.15.100}).
\hfill\qed

\subsubsection{Proof of Proposition \ref{prop;18.8.27.1}}

Let $\gbiga_3$ denote the set of the permutations of
$(0,1,2)$.
The following holds:
\begin{equation}
\label{eq;18.8.14.20}
\sum_i
 \nabla_{\omega,i}\nabla_{\omega,i}\phi_{\omega}
+\sum_{\sigma\in\gbiga_3}
\nabla_{\omega,\sigma(0)}
 \bigl(
 F(\nabla_{\omega})_{\sigma(1)\sigma(2)}
-\nabla_{\omega,\sigma(0)}\phi_{\omega}
 \bigr)
=
\sum_{\sigma\in\gbiga_3}
\nabla_{\omega,\sigma(0)}
 F(\nabla_{\omega})_{\sigma(1)\sigma(2)}
=0.
\end{equation}
Because
$\nabla_{\omega,\sigma(0)}
 \bigl(
 F(\nabla_{\omega})_{\sigma(1)\sigma(2)}
-\nabla_{\omega,\sigma(0)}\phi_{\omega}
 \bigr)
=O\bigl(e^{-\epsilon_{30} y_0^2}\bigr)$
for some $\epsilon_{30}>0$,
we obtain the following estimate
for some $\epsilon_{31}>0$:
\[
-\bigl(
 \del_{y_0}^2
+\del_{y_1}^2
+\del_{y_2}^2
 \bigr)|\phi_{\omega}|^2
=
-2\bigl|\nabla_{\omega}\phi_{\omega}\bigr|^2
+O\bigl(e^{-\epsilon_{31} y_0^2}\bigr).
\]
We obtain the following for some $\epsilon_{32}>0$:
\[
 -\del_{y_0}^2\int_{T^2}|\phi_{\omega}|^2
+O\bigl(e^{-\epsilon_{32} y_0^2}\bigr)
=-2\int_{T^2}\bigl|\nabla_{\omega}\phi_{\omega}\bigr|^2
\leq 0.
\]
Because 
$|\phi_{\omega}|$ is bounded,
we obtain
$\del_{y_0}\int_{T^2}|\phi_{\omega}|^2\to 0$
as $y_0\to -\infty$.
We also obtain
\[
 \int_{-\infty}^{R_1}dy_0
 \int_{T^2}
 |\nabla_{\omega}\phi_{\omega}|^2<\infty.
\]
By using (\ref{eq;18.8.14.21}) with $k=0$,
we obtain that
$\bigl|\nabla_{\omega}\phi_{\omega}\bigr|\to 0$
as $y_0\to-\infty$.
By a standard bootstrapping argument,
we obtain that 
the norms of the higher derivatives of
$\nabla_{\omega}\phi_{\omega}$
also converge to $0$ as $y_0\to-\infty$.
We also obtain that
the norms of $F(\nabla_{\omega})$
and its higher derivatives converge to $0$
as $y_0\to-\infty$.
\hfill\qed

\subsection{Asymptotically spectral decomposition}
\label{subsection;18.8.27.30}

\subsubsection{Setting}

Let $E$ be a $C^{\infty}$-vector bundle
on $\nbigu_{R_{40}}$
with a Hermitian metric $h$,
a unitary connection $\nabla$,
and an anti-Hermitian endomorphism $\phi$
such that the following holds.
\begin{itemize}
\item
For any $k\geq 0$,
there exist $B(k)>0$ and $\epsilon(k)>0$
such that
\[
\Bigl|
 \nabla_{\kappa_1}\circ\cdots\circ
 \nabla_{\kappa_k}
 \bigl(
 F(\nabla)-\ast\nabla\phi
 \bigr)
\Bigr|
\leq 
 B(k)e^{-\epsilon(k)y_0^2}
\]
for any $(\kappa_1,\ldots,\kappa_k)\in\{0,1,2\}^k$.
\item
$|\phi|$ is bounded.
\item
For any $k\geq 0$,
$\bigl|
 \nabla_{\kappa_1}\circ\cdots\circ\nabla_{\kappa_k}
 (\nabla\phi)
 \bigr|\to 0$
as $y_0\to-\infty$.
\end{itemize}

\subsubsection{Modification to mini-holomorphic structures}

We set $\sfz:=y_1+\sqrt{-1}y_2$.
For any $k\geq 0$,
there exists $\epsilon_1(k)>0$
such that
\[
\Bigl|
 \nabla_{\kappa_1}\circ\cdots\circ
 \nabla_{\kappa_k}
 \bigl(
 \bigl[
 \nabla_{\sfzbar},
 \nabla_{y_0}-\sqrt{-1}\phi
 \bigr]
\bigr)
\Bigr|_h
=O\bigl(
 e^{-\epsilon_1(k)y_0^2}
 \bigr)
\]
for any $(\kappa_1,\ldots,\kappa_k)\in\{0,1,2\}^k$.

\begin{lem}
There exists $A\in\End(E)$ with the following property.
\begin{itemize}
\item
$\bigl[
 \nabla_{\sfzbar}+A,
 \nabla_{y_0}-\sqrt{-1}\phi
 \bigr]=0$.
\item
For any $k\in\seisuu_{\geq 0}$,
there exists $\epsilon_2(k)>0$
such that
$\Bigl|
 \nabla_{\kappa_1}\circ\cdots\circ
 \nabla_{\kappa_k}
 A
 \Bigr|
=O\bigl(e^{-\epsilon_2(k)y_0^2}\bigr)$
for any
$(\kappa_1,\ldots,\kappa_k)\in\{0,1,2\}^k$.
\end{itemize}
\end{lem}
\pf
It is enough to take the integral of 
$\bigl[
 \nabla_{\sfzbar},
 \nabla_{y_0}-\sqrt{-1}\phi
 \bigr]$
along $y_0$
by using the parallel transport with respect to 
$\nabla_{y_0}-\sqrt{-1}\phi$.
\hfill\qed

\vspace{.1in}

The bundle $E$ has the mini-holomorphic structure
$\delbar_{E}$
given by
$\del_{E,y_0}=\nabla_{y_0}-\sqrt{-1}\phi$
and $\del_{E,\sfzbar}:=\nabla_{\sfzbar}+A$.
By the construction,
the following holds:
\[
G(h)
=\bigl[
 \nabla_{\sfz}-A^{\dagger},\,
 \nabla_{\sfzbar}+A
 \bigr]
-\frac{\sqrt{-1}}{2}
 \nabla_{y_0}\phi.
\]
Hence, 
for any $k\geq 0$,
there exists $\epsilon_3(k)>0$
such that 
\begin{equation}
\label{eq;18.12.15.40}
\bigl|
 \nabla_{\kappa_1}\circ\cdots\circ\nabla_{\kappa_k}
 G(h)
 \bigr|_h=O\bigl(e^{-\epsilon_3(k) y_0^2}\bigr)
\end{equation}
for any $(\kappa_1,\ldots,\kappa_k)\in\{0,1,2\}^k$.

\subsubsection{Spectral decomposition}

We have the decomposition of the mini-holomorphic bundle
\[
 (E,\delbar_{E})
=\bigoplus_{\alpha\in (T^2)^{\lor}}
 (E_{\alpha},\delbar_{E_{\alpha}}),
\]
where 
$\Spec(E_{\alpha},\delbar_{E_{\alpha}})=\{\alpha\}$.
Let $\Psi:\nbigu_{R}\lrarr \nbigh_{R}$
denote the projection for any $R$.

\begin{lem}
\label{lem;18.12.15.10}
If $R_{41}$ is sufficiently large,
there exists a vector bundle
$V=\bigoplus V_{\alpha}$ on $\nbigh_{R_{41}}$
with a graded connection $\nabla_V=\bigoplus \nabla_{V_{\alpha}}$,
a graded endomorphism $f=\bigoplus f_{\alpha}$,
and a graded isomorphism
$\Psi^{-1}(V)=
 \bigoplus\Psi^{-1}(V_{\alpha})\simeq 
 E_{|\nbigu_{R_{41}}}=
 \bigoplus E_{\alpha|\nbigu_{R_{41}}}$
such that the following holds.
\begin{itemize}
\item
$f_{\alpha}$ has a unique eigenvalue $\alpha$.
\item
 $\nabla_{V_{\alpha}}(f_{\alpha})=0$.
\item
$\del_{E,\sfzbar}=\del_{\sfzbar}+\Psi^{-1}(f)$,
where $\del_{\sfzbar}$ is the naturally defined operator
on $\Psi^{-1}(V)$.
\item
$\nabla_{y_0}-\sqrt{-1}\phi=
\Psi^{-1}(\nabla_{V,y_0})$.
\end{itemize}
\end{lem}
\pf
If $R_{41}$ is sufficiently large,
$E_{|\Psi^{-1}(y_0)}$ is semistable of degree $0$
for any $y_0<-R_{41}$.
We may assume it from the beginning.

We set
$\nbigu_R^{\star}:=\nbigu_R\times\real_{y_3}$.
We introduce the complex coordinate system
$\sfz:=y_1+\sqrt{-1}y_2$
and $\sfw:=y_3+\sqrt{-1}y_0$.
We also set
$\nbigh_R^{\star}:=\nbigh_R\times\real_{y_3}$,
on which we have the complex coordinate
$\sfw:=y_3+\sqrt{-1}y_0$.
Let $\sfQ_1:\nbigu_R^{\star}\lrarr\nbigu_R$ 
and $\sfQ_0:\nbigh_R^{\star}\lrarr\nbigh_R$
denote the projections.
Let $\Psi^{\star}:\nbigu_R^{\star}\lrarr\nbigh_R^{\star}$
denote the projection.

We set
$E^{\star}:=\sfQ_1^{-1}(E)$,
which is naturally $\real_{y_3}$-equivariant.
Let $\del_{E^{\star},\sfzbar}$
denote the derivative on $E^{\star}$ with respect to $\del_{\sfzbar}$
induced by $\del_{E,\sfzbar}$.
Let $\del_{E^{\star},y_0}$ denote the derivative on
$E^{\star}$ induced by 
$\del_{E,y_0}=\nabla_{y_0}-\sqrt{-1}\phi$.
Let $\del_{E^{\star},y_3}$ denote the naturally induced derivative on
$E^{\star}$ with respect to $\del_{y_3}$.
Then,
we set
$\del_{E^{\star},\sfwbar}:=
 \frac{1}{2}\bigl(
 \del_{y_3}+\sqrt{-1}\del_{E^{\star},y_0}\bigr)$.
They determine a holomorphic structure
$\delbar_{E^{\star}}$ of $E^{\star}$,
which is $\real$-equivariant.
We have the spectral decomposition
$(E^{\star},\delbar_{E^{\star}})
 =\bigoplus(E^{\star}_{\alpha},\delbar_{E^{\star}_{\alpha}})$
corresponding to the spectral decomposition of
$E$.

According to \cite[\S2.1]{Mochizuki-doubly-periodic},
there exists an $\real_{y_3}$-equivariant
 graded holomorphic vector bundle
$(V^{\star},\delbar_{V^{\star}})
=\bigoplus
 (V^{\star}_{\alpha},\delbar_{V^{\star}_{\alpha}})$
on $\nbigh^{\star}_{R_{40}}$
with an $\real_{y_3}$-equivariant
 holomorphic graded endomorphism
$f^{\star}=\bigoplus f^{\star}_{\alpha}$
and an $\real_{y_3}$-equivariant graded isomorphism
\[
 (\Psi^{\star})^{-1}(V^{\star})
=\bigoplus
 (\Psi^{\star})^{-1}(V^{\star}_{\alpha})
\]
such that 
the following holds:
\begin{itemize}
\item
 $f^{\star}_{\alpha}$
 has a unique eigenvalue $\alpha$.
\item
 $\del_{E^{\star},\sfzbar}=\del_{\sfzbar}+(\Psi^{\star})^{-1}(f^{\star})$,
 where $\del_{\sfzbar}$ is the naturally induced
 derivative on $(\Psi^{\star})^{-1}(V^{\star})$.
\item
 $\del_{E^{\star},\sfwbar}$
 is equal to the operator induced by
 $\del_{V^{\star},\sfwbar}$.
\end{itemize}

By the $\real_{y_3}$-equivariance of $V^{\star}$,
we obtain a graded $C^{\infty}$-vector bundle 
$V=\bigoplus V_{\alpha}$ on $\nbigh_{R_{40}}$.
The $\real_{y_3}$-equivariant holomorphic structure induces 
a graded flat connection 
$\nabla_{V}=\bigoplus \nabla_{V_{\alpha}}$.
The $\real_{y_3}$-equivariant
holomorphic graded endomorphism $f^{\star}$
induces a flat graded endomorphism
$f=\bigoplus f_{\alpha}$.
The $\real_{y_3}$-equivariant graded isomorphism
induces
a graded isomorphism
$E\simeq \Psi^{-1}(V)$.
Then, it is easy to see that they have the desired property.
\hfill\qed

\vspace{.1in}

We obtain the Hermitian metric
$h_{\alpha}$ of $V_{\alpha}$
as follows:
\[
 h_{\alpha}
=\frac{1}{\vol(T^2)}
 \int_{T^2}
 h\bigl(
 \Psi^{-1}(u_1),\Psi^{-1}(u_2)
 \bigr)\,
 dy_1\,dy_2.
\]
We set $h^{\circ}:=\bigoplus_{\alpha}\Psi^{-1}(h_{\alpha})$
on $E$.
We obtain the automorphism $b$ 
which is self-adjoint with respect to both $h$ and $h^{\circ}$,
determined by $h=h^{\circ}\cdot b$.
The following estimate can be proved
by arguments in 
\cite{Mochizuki-doubly-periodic, Mochizuki-difference-modules}.
\begin{prop}
\label{prop;18.8.27.21}
For any $P\in\nbiga$,
there exist  $C(P)>0$ and $\epsilon(P)>0$ such that
\[
\bigl|
 P(\nabla_{y_0},\nabla_{y_1},\nabla_{y_2},\phi)
 (b-\id)
\bigr|
\leq C(P)
 e^{\epsilon(P)y_0}.
\]
\end{prop}
\pf
We give an outline of the proof.
We use the notation in the proof of Lemma \ref{lem;18.12.15.10}.
Let $h^{\star}$ be the metric of $E^{\star}$ induced by $h$.
We obtain Hermitian metrics $h^{\star}_{\alpha}$
of $V^{\star}_{\alpha}$ in a way similar to the construction
of $h_{\alpha}$.
We set
$h^{\star\circ}:=\bigoplus (\Psi^{\star})^{-1}(h_{\alpha}^{\star})$.
We obtain $b^{\star}$
by $h^{\star}=h^{\star\circ}b^{\star}$.
The metrics $h^{\star}$, $h_{\alpha}^{\star}$,
and $h^{\star\circ}$ are $\real_{y_3}$-equivariant,
and hence $b^{\star}$ is also $\real_{y_3}$-equivariant.

Let $F(h^{\star})$ denote the curvature
of the Chern connection $\nabla^{\star}$ of
$(E^{\star},\delbar_{E^{\star}},h^{\star})$.
We have the expression
$F(h^{\star})=
 F_{\sfz\sfzbar}\,d\sfz\,d\sfzbar
+F_{\sfz\sfwbar}\,d\sfz\,d\sfwbar
+F_{\sfw\sfzbar}\,d\sfw\,d\sfzbar
+F_{\sfw\sfwbar}\,d\sfw\,d\sfwbar$.

Let $U$ be any open subset of $\nbigh_{R_{41}}$.
Let $U^{\star}:=U\times\real_{y_3}$.
The fiber integral induces the map
\[
 C^{\infty}\bigl(
 (\Psi^{\star})^{-1}(U^{\star}),\End(E^{\star}_{\alpha})\bigr)
\lrarr
 C^{\infty}\bigl(U^{\star},\End(V_{\alpha}^{\star})\bigr).
\]
Let 
$C^{\infty}((\Psi^{\star})^{-1}(U^{\star}),\End(E^{\star}_{\alpha}))_0$
denote the kernel.
There exists the injection
\[
 C^{\infty}\bigl(U^{\star},\End(V_{\alpha}^{\star})\bigr)
\lrarr
 C^{\infty}\bigl(
 (\Psi^{\star})^{-1}(U^{\star}),\End(E^{\star}_{\alpha})\bigr)
\]
induced by the pull back.
Thus, we obtain the decomposition
\[
  C^{\infty}((\Psi^{\star})^{-1}(U^{\star}),\End(E^{\star}_{\alpha}))
=C^{\infty}(U^{\star},\End(V_{\alpha}^{\star}))
\oplus
  C^{\infty}((\Psi^{\star})^{-1}(U^{\star}),\End(E^{\star}_{\alpha}))_0.
\]
We set
\[
 C^{\infty}((\Psi^{\star})^{-1}(U^{\star}),\End(E^{\star}))^{\circ}
:=\bigoplus_{\alpha}
 C^{\infty}(U^{\star},\End(V^{\star})),
\]
\[
 C^{\infty}((\Psi^{\star})^{-1}(U^{\star}),\End(E^{\star}))^{\bot}
:=\bigoplus_{\alpha}
   C^{\infty}((\Psi^{\star})^{-1}(U^{\star}),\End(E^{\star}_{\alpha}))_0
\oplus
 \bigoplus_{\alpha\neq\beta}
 C^{\infty}\bigl(
 (\Psi^{\star})^{-1}(U^{\star}),
 \Hom(E^{\star}_{\alpha},E^{\star}_{\beta})
 \bigr).
\]
We obtain the decomposition
\[
 C^{\infty}((\Psi^{\star})^{-1}(U^{\star}),\End(E^{\star}))
=
 C^{\infty}((\Psi^{\star})^{-1}(U^{\star}),\End(E^{\star}))^{\circ}
\oplus
 C^{\infty}((\Psi^{\star})^{-1}(U^{\star}),\End(E^{\star}))^{\bot}.
\]
For any sections 
$s\in C^{\infty}((\Psi^{\star})^{-1}(U^{\star}),\End(E^{\star}))$,
we obtain the decomposition
$s=s^{\circ}+s^{\bot}$.
We also obtain 
a function $\|s\|^2$ on $U^{\star}$
by the fiber integral of
$|s|^2_{h^{\star}}$.
There exists $C>0$ such that
$\|\nabla^{\star}_{\sfzbar}s^{\bot}\|\geq C\|s^{\bot}\|$
and 
$\|\nabla^{\star}_{\sfz}s^{\bot}\|\geq C\|s^{\bot}\|$
for any 
$s^{\bot}\in C^{\infty}((\Psi^{\star})^{-1}(U^{\star}),\End(E^{\star}))^{\bot}$.

By using the argument in the proof of
\cite[\S5.5.2]{Mochizuki-doubly-periodic},
we obtain the following estimates
for some small $\epsilon_i>0$:
\begin{multline}
 -\del_{\sfw}\del_{\sfwbar}
 \bigl\|F_{\sfz\sfzbar}^{\bot}\bigr\|^2
\leq
 -\bigl\|\nabla^{\star}_{\sfzbar}F_{\sfz\sfzbar}^{\bot}\bigr\|^2
 -\bigl\|\nabla^{\star}_{\sfzbar}F_{\sfzbar\sfzbar}^{\bot}\bigr\|^2
 -\bigl\|\nabla^{\star}_{\sfwbar}F_{\sfz\sfzbar}^{\bot}\bigr\|^2
 -\bigl\|\nabla^{\star}_{\sfwbar}F_{\sfz\sfzbar}^{\bot}\bigr\|^2
\\
+
O\Bigl(
 \epsilon_1\|F_{\sfz\sfzbar}^{\bot}\|^2
+\epsilon_1\|F_{\sfz\sfzbar}^{\bot}\| \|F_{\sfw\sfzbar}^{\bot}\|
+\epsilon_1\|\nabla^{\star}_{\sfwbar}F_{\sfw\sfzbar}^{\bot}\|
 \|F_{\sfz\sfzbar}\|
+\epsilon_1\|\nabla^{\star}_{\sfzbar}F_{\sfz\sfzbar}^{\bot}\|
 \|F_{\sfz\sfzbar}^{\bot}\|
 \Bigr)
\\
+O\Bigl(
 \epsilon_1\|\nabla^{\star}_{\sfwbar}F_{\sfz\sfzbar}^{\bot}\| 
  \|F_{\sfz\sfzbar}^{\bot}\|
+\epsilon_1\|F_{\sfw\sfzbar}^{\bot}\|^2
+\epsilon_1\|F_{\sfw\sfzbar}^{\bot}\| \|\nabla^{\star}_{\sfwbar}F_{\sfz\sfzbar}^{\bot} \|
 \Bigr)
+O\Bigl(
 \exp\bigl(-\epsilon_2y_0^2\bigr)
 \Bigr).
\end{multline}
\begin{multline}
-\del_{\sfw}\del_{\sfwbar}\|F^{\bot}_{\sfz\sfwbar}\|^2
\leq
-\|\nabla^{\star}_{\sfzbar}F^{\bot}_{\sfz\sfwbar}\|^2
-\|\nabla^{\star}_{\sfz}F^{\bot}_{\sfz\sfwbar}\|^2
-\|\nabla^{\star}_{\sfwbar}F^{\bot}_{\sfz\sfwbar}\|^2
-\|\nabla^{\star}_{\sfw}F^{\bot}_{\sfz\sfwbar}\|^2
 \\
+O\Bigl( 
 \epsilon_1\|F_{\sfz\sfwbar}^{\bot}\|^2\|F_{\sfz\sfzbar}^{\bot}\|
+\epsilon_1\|\nabla^{\star}_{\sfwbar}F_{\sfw\sfzbar}^{\bot}\|
 \|F_{\sfz\sfwbar}^{\bot}\|
+\epsilon_1\|F_{\sfz\sfwbar}^{\bot}\|\|F_{\sfw\sfzbar}^{\bot}\|
+\epsilon_1\|\nabla^{\star}_{\sfwbar}F_{\sfz\sfzbar}^{\bot}\|
 \|F_{\sfz\sfwbar}^{\bot} \|
 \Bigr)
\\
+O\Bigl(
 \epsilon_1\|\nabla^{\star}_{\sfzbar}F_{\sfz\sfwbar}^{\bot}\|
 \|F_{\sfz\sfzbar}^{\bot}\|
+\epsilon_1 \|F_{\sfz\sfwbar}^{\bot}\|^2
 \Bigr)
+O\Bigl(
 \exp\bigl(-\epsilon_2y_0^2\bigr)
 \Bigr).
\end{multline}
\begin{multline}
 -\del_{\sfw}\del_{\sfwbar}\|F_{\sfw\sfzbar}^{\bot}\|^2
\leq
-\bigl\|
 \nabla^{\star}_{\sfz}F_{\sfw\sfzbar}^{\bot}
 \bigr\|^2
-\bigl\|
 \nabla^{\star}_{\sfzbar}F_{\sfw\sfzbar}^{\bot}
 \bigr\|^2
-\bigl\|
 \nabla^{\star}_{\sfw}F_{\sfw\sfzbar}^{\bot}
 \bigr\|^2
-\bigl\|
 \nabla^{\star}_{\sfwbar}F_{\sfw\sfzbar}^{\bot}
 \bigr\|^2
\\
+O\Bigl(
 \epsilon_1\|F_{\sfw\sfzbar}^{\bot}\|\|F_{\sfz\sfzbar}^{\bot}\|
+\epsilon_1\|\nabla^{\star}_{\sfwbar}F_{\sfw\sfzbar}^{\bot}\|
 \|F_{\sfw\sfzbar}^{\bot}\|
+\epsilon_1\|F_{\sfw\sfzbar}^{\bot}\|^2
+\epsilon_1\|\nabla^{\star}_{\sfwbar}F_{\sfz\sfzbar}^{\bot}\|
 \|F_{\sfw\sfzbar}^{\bot}\|
 \Bigr)
 \\
+O\Bigl(
 \epsilon_1\|\nabla^{\star}_{\sfzbar}F_{\sfw\sfzbar}^{\bot}\|
 \|F_{\sfz\sfzbar}\|
+\epsilon_1\|F_{\sfw\sfzbar}^{\bot}\| \|F_{\sfw\sfwbar}\|
\Bigr)
+O\Bigl(
 \exp\bigl(-\epsilon_2y_0^2\bigr)
 \Bigr).
\end{multline}
From the estimate for $G(h)$,
we obtain
\begin{equation}
 -\del_{\sfw}\del_{\sfwbar}\|F_{\sfw\sfwbar}^{\bot}\|^2
=-\del_{\sfw}\del_{\sfwbar}\|F_{\sfz\sfzbar}^{\bot}\|^2
+O\Bigl(
 \exp\bigl(-\epsilon_2y_0^2\bigr)
 \Bigr).
\end{equation}
We set
$g:= \|F^{\bot}_{\sfz\sfzbar}\|^2
+\|F^{\bot}_{\sfz\sfwbar}\|^2
+\|F^{\bot}_{\sfw\sfzbar}\|^2
+\|F^{\bot}_{\sfw\sfwbar}\|^2$.
From these estimates,
we obtain the following for some $C_i>0$:
\[
 -\del_{\sfw}\del_{\sfwbar}g
\leq
-C_1g+
 C_2\exp\bigl(-\epsilon_2y_0^2\bigr).
\]
Note that $g$ depends only on $y_0$
by the $\real_{y_3}$-equivariance.
Hence, we obtain the following:
\[
  -\del_{y_0}^2g
\leq
-C_1'g+
 C_2'\exp\bigl(-\epsilon_2y_0^2\bigr).
\]
By a standard argument of Ahlfors lemma \cite{a,Simpson90}
we obtain that 
$g=O\bigl(e^{\epsilon_3y_0}\bigr)$
for some $\epsilon_3>0$.

Set
$F(h^{\star})^{\bot}:=
 F_{\sfz\sfzbar}^{\bot}\,d\sfz\,d\sfzbar
+ F_{\sfz\sfwbar}^{\bot}\,d\sfz\,d\sfwbar
+F_{\sfw\sfzbar}^{\bot}\,d\sfw\,d\sfzbar
+F_{\sfw\sfwbar}^{\bot}\,d\sfw\,d\sfwbar$.
By using a standard bootstrapping argument
as in the proof of \cite[Proposition 5.8]{Mochizuki-doubly-periodic},
we obtain the following.
\begin{itemize}
\item
For any $P\in\nbiga$,
we have $C(P)>0$ and $\epsilon(P)>0$ such that
\[
 \bigl|
 P(\nabla^{\star}_{\sfz},\nabla^{\star}_{\sfzbar},
 \nabla^{\star}_{\sfw},\nabla^{\star}_{\sfwbar})
 F(h^{\star})^{\bot}
 \bigr|
\leq C(P)\exp\bigl(\epsilon(P)y_0\bigr).
\]
\end{itemize}
By \cite[Lemma 10.13]{Mochizuki-doubly-periodic},
we obtain the following.
\begin{itemize}
\item
For any $P\in\nbiga$,
we have $C(P)>0$ and $\epsilon(P)>0$ such that
\[
 \bigl|
 P(\nabla^{\star}_{\sfz},\nabla^{\star}_{\sfzbar},
 \nabla^{\star}_{\sfw},\nabla^{\star}_{\sfwbar})
 (b^{\star}-\id)
 \bigr|
\leq C(P)\exp\bigl(\epsilon(P)y_0\bigr).
\]
\end{itemize}
Note that for any $\real_{y_3}$-invariant section $s$
of $\End(E)$,
we have
\[
 \nabla^{\star}_{\sfwbar}(s)
=\frac{\sqrt{-1}}{2}\bigl(
 \nabla_{y_0}-\sqrt{-1}\phi\bigr)s,
\quad
  \nabla^{\star}_{\sfw}(s)
=\frac{\sqrt{-1}}{2}\bigl(
 \nabla_{y_0}+\sqrt{-1}\phi\bigr)s.
\]
We also have
$\nabla^{\star}_{\sfzbar}s
=(\nabla_{\sfzbar}+A)s$
and 
$\nabla^{\star}_{\sfz}s
=(\nabla_{\sfz}-A^{\dagger})s$.
Hence, we obtain the desired estimate for $b$.
\hfill\qed

\subsubsection{Anti-Hermitian endomorphisms}

We have $\phi_{3,\alpha}$
determined by
$\nabla_{V_{\alpha}}
=\nabla^u_{V_{\alpha}}
-\sqrt{-1}\phi_{3,\alpha}\,dy_0$,
where
$\nabla^u_{V_{\alpha}}$ 
is a unitary connection of
$(V_{\alpha},h_{\alpha})$,
and 
$\phi_{\alpha}$ is an anti-Hermitian endomorphism
of $(V_{\alpha},h_{\alpha})$.
Set $\phi_3:=\bigoplus \phi_{3,\alpha}$.
\begin{prop}
For any $P\in\nbiga$,
there exist $C(P)>0$ and $\epsilon(P)>0$
such that
\[
\Bigl|
 P(\nabla_{y_0},\nabla_{y_1},\nabla_{y_2},\phi)
 (\phi-\Psi^{-1}(\phi_{3}))
\Bigr|
\leq C(P)e^{\epsilon(P)y_0}.
\]
\end{prop}
\pf
It follows from Proposition \ref{prop;18.8.27.21}.
\hfill\qed

\vspace{.1in}

We define the anti-Hermitian endomorphisms
$\phi_i=\bigoplus \phi_{i,\alpha}$ $(i=1,2)$ of 
$(V,h_V)=\bigoplus (V_{\alpha},h_{V_{\alpha}})$
by $f=\frac{1}{2}(\phi_1+\sqrt{-1}\phi_2)$.
\begin{lem}
$\nabla_{V,y_0}\phi_1-[\phi_2,\phi_3]=0$
and
$\nabla_{V,y_0}\phi_2-[\phi_3,\phi_1]=0$
hold.
\end{lem}
\pf
It follows from the flatness
$[\nabla_{V,y_0}-\sqrt{-1}\phi_3,f]=0$.
\hfill\qed

\begin{prop}
For any $P\in\nbiga$,
there exist $C(P)>0$ and $\epsilon(P)>0$
such that
\[
 \Bigl|
 P(\nabla_{V,y_0},\phi_1,\phi_2,\phi_3)
 \bigl(
 \nabla_{V,y_0}\phi_3
-[\phi_1,\phi_2]
 \bigr)
 \Bigr|
\leq
 C(P)e^{\epsilon(P)y_0}.
\]
\end{prop}
\pf
It follows from Proposition \ref{prop;18.8.27.21}
and the estimate (\ref{eq;18.12.15.40}).
\hfill\qed

\begin{lem}
$\phi_i$ $(i=1,2)$ are bounded.
\end{lem}
\pf
Let $\delbar_{V^{\star}}+\del_{V^{\star}}$
be the Chern connection of
$(V^{\star},\delbar_{V^{\star}},h^{\star})$.
Set $\theta^{\star}:=f^{\star}d\sfw$,
which is a Higgs field of
$(V^{\star},\delbar_{V^{\star}})$.
Let $(\theta^{\star})^{\dagger}$
denote the adjoint of $\theta^{\dagger}$
with respect to $h^{\star}$.
We obtain
$[\delbar_{V^{\star}},\del_{V^{\star}}]
+[\theta^{\star},(\theta^{\star})^{\dagger}]
=O(e^{\epsilon y_0})$
for some $\epsilon>0$.
Note that the eigenvalues of $\ftilde$
are constant.
Hence, as a variant of Simpson's main estimate 
(\cite{Simpson90} and \cite[Proposition 2.10]{Mochizuki-KHI}),
we obtain that
$|\ftilde|_{h^{\star}}$ is bounded.
Then, the claim of the lemma follows.
\hfill\qed

\subsection{Approximate solutions of Nahm equations}
\label{subsection;18.8.27.31}

\subsubsection{Reduction}

\label{subsection;18.8.30.1}

Let $V$ be a $C^{\infty}$-vector bundle on 
$\nbigh_R$ with a Hermitian metric $h$,
a unitary connection $\nabla$,
and bounded anti-self-adjoint endomorphisms $\phi_i$ $(i=1,2,3)$.
We introduce a condition.
\begin{condition}
\label{condition;18.12.16.1}
For any $P\in\nbiga$,
there exist $\epsilon(P)>0$ and $B(P)>0$
such that 
\begin{equation}
\label{eq;18.8.15.2}
\Bigl|
 P(\nabla_{y_0},\phi_1,\phi_2,\phi_3)
 \bigl(
 \nabla_{y_0}\phi_i
-[\phi_j,\phi_k]
 \bigr)
\Bigr|
\leq
 B(P)e^{\epsilon(P)y_0},
\end{equation}
where $(i,j,k)$ denotes any cyclic permutation
of $(1,2,3)$.
Moreover,
$\nabla_{y_0}\phi_i\to 0$ as $y_0\to-\infty$.
\hfill\qed
\end{condition}

\begin{prop}
\label{prop;18.8.27.22}
There exist
a finite subset
$S\subset(\sqrt{-1}\real)^3$,
an orthogonal decomposition
$V=\bigoplus_{\vecb\in S}V_{\vecb}$,
a graded unitary connection 
$\nabla^{\sharp}=\bigoplus\nabla_{\vecb}^{\sharp}$,
and graded anti-self-adjoint endomorphisms
$\phi^{\sharp}_i=\bigoplus \phi^{\sharp}_{i,\vecb}$
such that the following holds:
\begin{itemize}
\item
The eigenvalues of $\phi^{\sharp}_{i,\vecb|y_0}$ converge to $b_i$
as $y_0\to-\infty$.
\item
Set $\rho^{\sharp}:=\nabla-\nabla^{\sharp}$.
For any $k$, there exist $B(k)>0$ and $\epsilon(k)>0$
such that
\[
\Bigl|
 \bigl(
 \nabla^{\sharp}_{y_0}
 \bigr)^k\rho^{\sharp}
\Bigr|_h
+\sum
 \Bigl|
 \bigl(
 \nabla^{\sharp}_{y_0}
 \bigr)^k(\phi_i-\phi_i^{\sharp})
 \Bigr|_h
\leq
 B(k)e^{\epsilon(k)y_0}.
\]
\end{itemize}
As a result,
$V_{\vecb}$, the induced metric $h_{\vecb}$,
the induced connection $\nabla^{\sharp}_{\vecb}$,
and the anti-Hermitian endomorphisms
$\phi^{\sharp}_{i,\vecb}$ $(i=1,2,3)$
satisfy Condition {\rm\ref{condition;18.12.16.1}}.
\end{prop}
\pf
We begin with a preliminary.

\begin{lem}
\label{lem;18.12.16.10}
For each $i$,
there exist a finite subset $S(\phi_i)\subset\sqrt{-1}\real$
such that the following holds.
\begin{itemize}
\item
Let $\Sp(\phi_{i|y_0})$ be the set of eigenvalues of
$\phi_{i|y_0}$.
For any $\alpha\in\cnum$ and $\delta>0$,
set $B_{\alpha}(\delta):=\{\beta\in\cnum\,|\,|\alpha-\beta|<\delta\}$.
Then, 
for any $\delta_1>0$,
there exists $R_2$ such that
the following holds for any $y_0<-R_2$:
\[
 \Sp(\phi_{i|y_0})
 \subset
 \bigcup_{\alpha\in S(\phi_i)}
 B_{\alpha}(\delta),
\quad
 S(\phi_i)
\subset
  \bigcup_{\alpha\in\Sp(\phi_{i|y_0})}
 B_{\alpha}(\delta).
\]
\end{itemize}
\end{lem}
\pf
We set $F:=\phi_2+\sqrt{-1}\phi_3$.
Then,
$(\nabla_{y_0}-\sqrt{-1}\phi_1)F=O(e^{\epsilon y_0})$.
There exist 
$\ttA\in\End(V)$ such that
$(\nabla_{y_0}-\sqrt{-1}\phi_1)\ttA
=(\nabla_{y_0}-\sqrt{-1}\phi_1)F$
and that
$\ttA=O(e^{\epsilon_1y_0})$ for some $\epsilon_1>0$.
We set $\Ftilde:=F-\ttA$.
Because
$(\nabla_{y_0}-\sqrt{-1}\phi_1)\Ftilde=0$,
the eigenvalues of $\Ftilde$ are constant with respect to $y_0$.
Then, we obtain the claim for $\phi_2$ and $\phi_3$.
Similarly, we obtain the claim for $\phi_1$.
\hfill\qed

\vspace{.1in}

Let 
$(V,\phi_3)=\bigoplus_{\alpha\in S(\phi_3)}
 (V_{\alpha},\phi_{3,\alpha})$
be the decomposition
satisfying the following condition.
\begin{itemize}
\item
For any $\delta_1>0$,
there exists $R_2$ such that
eigenvalues $\beta$ of $\phi_{3,\alpha|y_0}$  $(y_0<-R_2)$
satisfy
$|\alpha-\beta|\leq \delta_1$.
\end{itemize}
We obtain the decomposition
$\nabla=\nabla^{\bullet}+\rho$,
where $\nabla^{\bullet}=\bigoplus \nabla^{\bullet}_{\alpha}$ 
is the direct sum of unitary connections $\nabla^{\bullet}_{\alpha}$
on $V_{\alpha}$,
and $\rho$ is a section  of
$\bigoplus_{\alpha\neq\beta}
 \Hom(V_{\alpha},V_{\beta})\,dy_0$.
We also obtain the decomposition
$\phi_i=\phi_i^{\bullet}+\phi_i^{\top}$ $(i=1,2,3)$
according to the decomposition
$\End(V)=\bigoplus \End(V_{\alpha})
 \oplus \bigoplus_{\alpha\neq\beta}
 \Hom(V_{\alpha},V_{\beta})$.
Clearly, $\phi_3^{\bullet}=\phi_3$ holds.
We have the decomposition
$\phi^{\bullet}_i=\bigoplus \phi^{\bullet}_{i,\alpha}$.

We obtain the following estimate
by an argument similar to the proof of 
Theorem \ref{thm;18.8.15.1}.
\begin{lem}
For any $k$,
there exists $\epsilon(k)>0$ such that
$\bigl|
 (\nabla^{\bullet}_{y_0})^k
 \phi_{1}^{\top}
\bigr|_h
+
\bigl|
 (\nabla^{\bullet}_{y_0})^k
 \phi_{2}^{\top}
\bigr|_h
+
 \bigl|
 (\nabla^{\bullet}_{y_0})^k
 \rho
\bigr|_h
=O\bigl(
 e^{\epsilon(k)y_0}
 \bigr)$.
\end{lem}
\pf
We give only an outline.
By using an argument in the proof of 
$\nabla^k_{y_0}\phi_i$
are bounded for any $k$

We obtain a bundle
$\Vtilde:=\Psi^{-1}(V)$ on $\nbigu_{R_0}$,
with the metric $\htilde=\Psi^{-1}(h)$,
the unitary connection
$\nablatilde:=
\Psi^{-1}(\nabla)+\phi_1\,dy_1
+\phi_2\,dy_2$
and the anti-Hermitian metric $\phitilde:=\Psi^{-1}(\phi_3)$.
Let $F(\nablatilde)$ denote the curvature of $\nablatilde$.
We have the following:
\begin{itemize}
\item
For any $k\in\seisuu_{\geq 0}$,
there exists $C(k)>0$ and $\epsilon(k)>0$
such that
\[
 \Bigl|
 \nablatilde_{\kappa_1}\circ\cdots\circ
 \nablatilde_{\kappa_k}
 \bigl(
 F(\nablatilde)-\ast\nablatilde\phitilde
 \bigr)
\Bigr|\leq
C(k)e^{\epsilon(k) y_0}
\]
for any $(\kappa_1,\ldots,\kappa_k)\in\{0,1,2\}^k$.
\item
For any $(\kappa_1,\ldots,\kappa_k)\in\{0,1,2\}^k$,
$\bigl|\nablatilde_{\kappa_1}\circ\cdots\circ\nablatilde_{\kappa_k}
 (\nablatilde\phitilde)\bigr|\to 0$ 
as $y_0\to-\infty$.
\end{itemize}
Corresponding to the decomposition
$\End(\Vtilde)=\bigoplus \End(\Vtilde_{\alpha})
\oplus
 \bigoplus_{\alpha\neq\beta} \Hom(\Vtilde_{\alpha},\Vtilde_{\beta})$,
we obtain 
$\nablatilde=\nablatilde^{\bullet}+\rhotilde$.
Note that 
$\rhotilde=\Psi^{-1}(\rho)+
 \Psi^{-1}(\phi_1^{\top})dy_1
+\Psi^{-1}(\phi_2^{\top})dy_2$.
Any section $s$ of $\End(\Vtilde)$
is decomposed into
$s^{\bullet}+s^{\top}$.

By using the argument in the proof of
\cite[Lemma 6.16]{Mochizuki-difference-modules},
we obtain the following:
\[
\sum_{i=0,1,2}
 \nablatilde_{y_i}^2(\nablatilde_a\phitilde)
=4\bigl[
 \nablatilde_b\phitilde,\nablatilde_c\phitilde
 \bigr]
-\bigl[
 \phitilde,[\phitilde,\nablatilde_a\phitilde]
 \bigr]
+O(e^{\epsilon_0 y_0}),
\]
where $(a,b,c)$ is a cyclic permutation of $(0,1,2)$.
By the argument in the proof of
\cite[Lemma 6.17]{Mochizuki-difference-modules},
for any $\delta>0$ there exists $R_{10}$
such that the following holds 
on $\nbigu_{R_{10}}$:
\[
 \htilde\bigl(
 \nablatilde_{\kappa_1}^2
 (\nablatilde_{\kappa_2}\phitilde)^{\bullet},
 (\nablatilde_{\kappa_2}\phitilde)^{\top}
 \bigr)
=O\Bigl(
 \delta\cdot
 \Bigl(
 \bigl|
 (\nablatilde_{\kappa_1}\phitilde)^{\top}
 \bigr|_h
+ \bigl|
 \nablatilde^{\bullet}_{\kappa_1}
 (\nablatilde_{\kappa_1}\phitilde)^{\top}
 \bigr|_h
 \Bigr)
\cdot 
 \bigl|
 (\nablatilde_{\kappa_2}\phitilde)^{\top}
 \bigr|_h
 \Bigr)
+O(e^{\epsilon_0 y_0}).
\]
By the argument in the proof of
\cite[Lemma 6.18]{Mochizuki-difference-modules},
for any $\delta>0$ there exists $R_{10}$
such that the following holds 
on $\nbigu_{R_{10}}$:
\[
 \sum_{\kappa_1=0,1,2}
 \htilde\bigl(
 \nablatilde_{\kappa_1}^2
 (\nablatilde_{\kappa_2}\phitilde),
 (\nablatilde_{\kappa_2}\phitilde)^{\top}
 \bigr)
=
 \Bigl|
 \bigl[\phitilde,(\nablatilde_{\kappa_2}\phitilde)^{\top}\bigr]
 \Bigr|^2_h
+
O\Bigl(
\delta
 \bigl|
 (\nablatilde\phitilde)^{\top}
 \bigr|_h
\cdot
 \bigl|(\nablatilde_{\kappa_2}\phitilde)^{\top}\bigr|
 \Bigr)
+O(e^{\epsilon_0y_0})
\]
By using the argument in the proof of
\cite[Lemma 6.19]{Mochizuki-difference-modules},
we obtain
\begin{multline}
-\sum_{i=0,1,2}
 \del_i^2\bigl|
 (\nablatilde\phitilde)^{\top}
 \bigr|^2
=-\sum_{i=0,1,2}
 2\bigl|
 \nablatilde(\nablatilde_i\phi)^{\top}
 \bigr|^2
-2\bigl|
 \bigl[
 \phitilde,(\nablatilde\phitilde)^{\top}
 \bigr]
 \bigr|^2
+O\Bigl(
 \delta\cdot
 \bigl|(\nablatilde\phitilde)^{\top}\bigr|^2
 \Bigr)
\\
+O\Bigl(
 \delta
 \sum_i
 \bigl|
 \nablatilde_i(\nablatilde\phitilde)^{\top}
 \bigr|
 \cdot
 \bigl|(\nablatilde\phitilde)^{\top}\bigr|
 \Bigr)
+O(e^{\epsilon_0y_0})
\end{multline}
Note that
$\Bigl|
 \bigl[\phitilde,(\nabla_{\kappa_2}\phitilde)^{\top}\bigr]
 \Bigr|
\geq
 c\Bigl|
 (\nabla_{\kappa_2}\phitilde)^{\top}
 \Bigr|$
for some $c>0$,
which we may assume to be independent of $R_{10}$.
Hence, we obtain the following if $R_{10}$ is large enough:
\[
 -\sum_{i=0,1,2}
 \del_i^2\bigl|
 (\nablatilde\phitilde)^{\top}
 \bigr|^2
\leq
 -c_1 \bigl|
 (\nablatilde\phitilde)^{\top}
 \bigr|^2
+O(e^{\epsilon_0y_0}).
\]
We set
$g:=
 \bigl|
 (\nabla\phi_3)^{\top}
 \bigr|^2
+\bigl|
 [\phi_3,\phi_1^{\top}]
 \bigr|^2
+\bigl|
 [\phi_3,\phi_2^{\top}]
 \bigr|^2$.
Because
$(\Vtilde,\htilde,\nablatilde,\phitilde)$
is equivariant with respect to 
the natural action of $\real_{y_1}\oplus\real_{y_2}$,
we obtain the following:
\[
 -\del_{y_0}^2g
\leq -c_2g+O(e^{\epsilon_0y_0}).
\]
By a standard argument,
we obtain $g=O(e^{\epsilon_1y_0})$
for some $\epsilon_1>0$.
We obtain
$|\phi^{\top}_1|+|\phi^{\top}_2|+|\rho^{\top}|
=O(e^{\epsilon_1y_0})$.
By a bootstrapping argument,
we obtain the estimates for higher derivatives.
\hfill\qed

\vspace{.1in}

We obtain
$(V_{\alpha},h_{\alpha})$
with 
a unitary connection $\nabla^{\bullet}_{\alpha}$ and 
bounded anti-Hermitian endomorphisms
$\phi^{\bullet}_{i,\alpha}$ $(i=1,2,3)$
satisfying Condition \ref{condition;18.12.16.1}.
Moreover, the eigenvalues of $\phi_{3,\alpha}$ converges to $\alpha$
as $y_0\to-\infty$.
By an applying similar argument to
$\phi_{2,\alpha}$ and $\phi_{1,\alpha}$ inductively,
we obtain the claim of Proposition \ref{prop;18.8.27.22}.
\hfill\qed

\vspace{.1in}
We shall study the behaviour of 
$\phi^{\sharp}_{i,\vecb}-b_i\id_{V_{\vecb}}$
in the next subsection.

\subsubsection{Decay}
\label{subsection;18.12.16.11}

Let $(V,h,\nabla,\{\phi_i\}_{i=1,2,3})$
be as in \S\ref{subsection;18.8.30.1}
satisfying Condition \ref{condition;18.12.16.1}.
Moreover, we assume that
the eigenvalues of $\phi_i$ are convergent to $0$
as $y_0\to-\infty$.

\begin{prop}
\label{prop;18.8.27.23}
For any $k\in\seisuu_{\geq 0}$, 
$\bigl| y_0^{k+1}
 \nabla_{y_0}^k\phi_i\bigr|$
are bounded.
In particular,
we obtain the expression
$\phi_i=y_0^{-1}A_i+O(y_0^{-2})$
for endomorphisms $A_i$ such that
$\nabla A_i=0$,
and the tuple $(A_1,A_2,A_3)$ satisfies
$[A_i,A_j]=A_k$,
where $(i,j,k)$ are cyclic permutation of $(1,2,3)$.
\end{prop}
\pf
Let $F$, $\ttA$ and $\Ftilde$
be as in the proof of Lemma \ref{lem;18.12.16.10}.
Because the eigenvalues of $\phi_i$ converges to $0$,
we obtain that $\Ftilde$ is nilpotent.
By the construction,
we have
$\bigl[
 \nabla_{y_0}-\sqrt{-1}\phi_1,\Ftilde
 \bigr]=0$,
and
\[
 \bigl[
 \Ftilde^{\dagger},\Ftilde
 \bigr]
+2\sqrt{-1}\nabla_{y_0}\phi_1
=O(e^{\epsilon y_0}).
\]
We obtain the following:
\[
 \del_{y_0}^2\bigl|\Ftilde\bigr|_h^2
=\bigl|
 [\nabla_{y_0}+\sqrt{-1}\phi_1,\Ftilde]
\bigr|^2_h
+\bigl|
 [\Ftilde^{\dagger},\Ftilde]
 \bigr|^2_h
+O\bigl(e^{\epsilon y_0}|\Ftilde|_h^2\bigr).
\]
Hence, 
we obtain
\[
 -\del_{y_0}^2\log|\Ftilde|_h^2
\leq
 -\frac{\bigl|[\Ftilde^{\dagger},\Ftilde]\bigr|_h^2}{|\Ftilde|_h^2}
+O(e^{\epsilon_0 y_0}).
\]
Because $\Ftilde$ is nilpotent,
there exists a positive constant $c_1$
depending only on $\rank E$
such that 
$\bigl|
 [\Ftilde^{\dagger},\Ftilde]
 \bigr|_h^2
\geq
 c_1|\Ftilde|_h^2$.
Hence, we obtain the following for some $c_2>0$:
\[
 -\del_{y_0^2}
 \log\bigl|\Ftilde\bigr|_h^2
\leq
 -c_2|\Ftilde|^2_h
+O(e^{\epsilon y_0}).
\]
By a standard argument of  Ahlfors lemma \cite{a, Simpson90}
we obtain that 
$|\Ftilde|_h^2=O(y_0^{-2})$.
We obtain $|\phi_i|=O(y_0^{-1})$ $(i=2,3)$.
Similarly, we obtain $|\phi_1|=O(y_0^{-1})$.
Then, we obtain
$\nabla \phi_i=-[\phi_j,\phi_k]=O(y_0^{-2})$.
By an inductive argument,
we obtain the estimates for the higher derivatives of
$\phi_i$.
\hfill\qed

\subsubsection{Norm estimate and the conjugacy class of the nilpotent map}
\label{subsection;18.11.24.1}

Let $(V,h,\nabla)$ and $\phi_i$ $(i=1,2,3)$
be as in \S\ref{subsection;18.12.16.11}.
Let $\Ftilde$ be the endomorphism of $V$
as in the proof of Lemma \ref{lem;18.12.16.10}.
If is flat with respect to $\nabla_{y_0}-\sqrt{-1}\phi_1$.
In this case, $\Ftilde$ is nilpotent.
We obtain the weight filtration $W$ of $V$
with respect to $\Ftilde$,
which is preserved by $\nabla_{y_0}-\sqrt{-1}\phi_1$.

Let $\vece=(e_1,\ldots,e_r)$ be a frame of $V$
satisfying the following conditions.
\begin{itemize}
\item
 $(\nabla_{y_0}-\sqrt{-1}\phi_1)\vece=0$.
\item
 $\vece$ is compatible with $W$,
 i.e.,
 there is a decomposition
 $\vece=\bigcup_{k\in\seisuu} \vece_k$
 such that
 $\bigcup_{k\leq\ell} \vece_k$
 is a frame of $W_{\ell}$.
\end{itemize}
If $e_i\in\vece_k$, we set $k(i):=k$.
Let $h_0$ be the Hermitian metric of $V$
defined by
$h_0(e_i,e_i)=(-y_0)^{k(i)}$
and $h_0(e_i,e_j)=0$ $(i\neq j)$.

\begin{prop}
\label{prop;18.12.16.20}
$h$ and $h_0$ are mutually bounded.
\end{prop}
\pf
Set $\Delta^{\ast}_R:=\{\slw\in\cnum^{\ast}\,|\,\log|\slw|<-R\}$.
Let $\slQ:\Delta^{\ast}_R\lrarr \nbigh_R$
be the map defined by
$\slQ(\slw)=\log|\slw|$.
We set $\Vtilde:=\slQ^{-1}(V)$
and $\htilde:=\slQ^{-1}(h)$.
They are naturally $S^1$-equivariant.
We define the derivative 
$\del_{\Vtilde,\slwbar}$
on $\Vtilde$ with respect to $\slwbar$
by
$\slwbar\del_{\Vtilde,\slwbar}\slQ^{-1}(s)
=\slQ^{-1}\bigl((\nabla_{y_0}-\sqrt{-1}\phi_1)s\bigr)$.
It induces an $S^1$-equivariant holomorphic structure 
$\delbar_{\Vtilde}$ on $\Vtilde$.
Let $\ftilde$ be the holomorphic endomorphism
of $\Vtilde$
induced by $\Ftilde$.
We set $\thetatilde:=\ftilde\,d\slw/\slw$.
Let $\nablatilde$ denote the Chern connection of
$(\Vtilde,\delbar_{\Vtilde},\htilde)$,
and let $F(\nablatilde)$ be the curvature of $\nablatilde$.
Let $\thetatilde^{\dagger}$ denote the adjoint of $\thetatilde$
with respect to $\htilde$.
Then, we have
$F(\nablatilde)+[\thetatilde,\thetatilde^{\dagger}]
=O(|\slw|^{\epsilon-2})d\slw\,d\slwbar$
for some $\epsilon>0$.
We also have
$F(\nablatilde)=
 O\bigl(|\slw|^{-2}(\log|\slw|)^{-2}\bigr)
 d\slw\,d\slwbar$.

Set $\Delta_R:=\Delta_R^{\ast}\cup\{0\}$.
We have the associated filtered bundle
$\nbigp_{\ast}\Vtilde$ on $(\Delta_R,0)$.
Let us observe that
$\Gr^{\nbigp}_a(\Vtilde)=0$ unless $a\in\seisuu$.
Indeed, let $\vecetilde=(\etilde_i)$ denote
the $S^1$-equivariant holomorphic frame of $\Vtilde$
induced by $\vece$.
Let $H(\htilde,\vecetilde)$ be the Hermitian-matrix valued
function whose $(i,j)$-entries are $\htilde(\vtilde_i,\vtilde_j)$.
Then, it is easy to see that
$C^{-1}(-\log|\slw|)^{-N}
 <H(\htilde,\vecetilde)
 <C(-\log|\slw|)^N$ 
for some $C>1$
and $N>0$.
Thus, we obtain 
$\Gr^{\nbigp}_a(\Vtilde)=0$ unless $a\in\seisuu$.
Then, the claim of Proposition \ref{prop;18.12.16.20}
follows from the norm estimate in \cite{Simpson90}.
\hfill\qed

\vspace{.1in}
Let $C_0$ be the matrix 
determined by
$(C_0)_{i,i}=k(i)/2$
and $(C_0)_{i,j}=0$ $(i\neq j)$

\begin{prop}
\label{prop;18.12.16.30}
The conjugacy class of $-\sqrt{-1}A_1$ 
is represented by $C_0$.
\end{prop}
\pf
Let $\vecv$ be an orthonormal frame of $V$
such that $\nabla_{y_0}\vecv=0$.
We obtain the matrix valued function $\nbiga_1$
determined by
$\phi_1\vecv=\vecv\nbiga_1$.
There is a constant matrix $\nbiga_{1,0}$
such that $\nbiga_1-y_0^{-1}\nbiga_{1,0}=O(y_0^{-2})$.
We have
$(\nabla_{y_0}-\sqrt{-1}\phi_1)\vecv=
\vecv\cdot(-\sqrt{-1}\nbiga_1)$.
We may assume that $\nbiga_{1,0}$ is diagonal.

We set $e_i':=(-y_0)^{-k(i)/2}e_i$.
We obtain a frame $\vece'=(e_i')$.
Let $B$ be the $\GL(r)$-valued function 
determined by
$\vecv=\vece'\cdot B$.
By Proposition \ref{prop;18.12.16.20},
$B$ and $B^{-1}$ are bounded.
Because
$(\nabla_{y_0}-\sqrt{-1}\phi_1)\vece'
=\vece' C_0y_0^{-1}$,
we obtain the relation
\[
 y_0\del_{y_0}B+C_0B+\sqrt{-1}B\cdot \nbiga_{1,0}
+\sqrt{-1}B\cdot (y_0\nbiga_{1}-\nbiga_{1,0})=0.
\]
Note that the eigenvalues of $\sqrt{-1}\nbiga_{1,0}$
are contained in $\frac{1}{2}\seisuu$
because $(A_1,A_2,A_3)$ induces an $\su(2)$-representation.
It is easy to check the following lemma.
\begin{lem}
\label{lem;18.12.16.31}
Let $a\in \frac{1}{2}\seisuu$.
Let $g$ be a bounded $C^{\infty}$-function 
on $\nbigh_R$
satisfying
$y_0\del_{y_0}g+ag=O(|y_0|^{-1})$.
Then,
the following holds.
\begin{itemize}
\item
 If $a=1$, then
 $g=O(|y_0|^{-1}\log|y_0|)$.
\item
 If $a=1/2$, then
 $g=O(|y_0|^{-1/2})$.
\item
 If $a=0$, 
 there exists $g_0\in\cnum$
 such that
 $g-g_0=O(|y_0|^{-1})$.
\item
 Otherwise,
 $g=O(|y_0|^{-1})$.
\hfill\qed
\end{itemize}
\end{lem}

By Lemma \ref{lem;18.12.16.31},
we obtain the following.
\begin{itemize}
\item
There exists $B_{i,j,0}\in\cnum$
such that 
$B_{i,j}-B_{i,j,0}=O(|y_0|^{-1/2})$.
Moreover,
$B_{i,j,0}=0$
unless $-\sqrt{-1}(\nbiga_{1,0})_{j,j}= k(i)/2$.
\end{itemize}
Then, the claim of Proposition \ref{prop;18.12.16.30}
follows from the boundedness of $B$ and $B^{-1}$.
\hfill\qed

\begin{prop}
\label{prop;19.2.8.2}
The conjugacy class of $\Ftilde$
is equal to the conjugacy class of
$A_2-\sqrt{-1}A_3$.
\end{prop}
\pf
Let $\nbign_0$ be the matrix valued function determined by
$\Ftilde\vece'=\vece'\cdot \nbign_0$.
There exists a constant matrix $N_0$
such that
$\nbign_0-N_0y_0^{-1}=O(|y_0|^{-3/2})$.
It is easy to observe that 
the conjugacy class of $\Ftilde$ is represented by $N_0$.

Let $\nbign_1$ be the matrix valued function 
determined by 
$(\phi_2-\sqrt{-1}\phi_3)\vece'=\vece'\cdot \nbign_1$.
Because $\Ftilde-(\phi_2-\sqrt{-1}\phi_3)=O(e^{\epsilon y_0})$,
we obtain
$\nbign_1-N_0y_0^{-1}=O(|y_0|^{-3/2})$.

Let $\vecv$ and $B$ be as in the proof of 
Proposition \ref{prop;18.12.16.30}.
Let $\nbign_2$ be the matrix valued function
determined by
$(\phi_2-\sqrt{-1}\phi_3)\vecv=\vecv\cdot \nbign_2$.
We have the constant matrix $N_2$
such that
$\nbign_2-N_2=O(|y_0|^{-1})$.
The conjugacy class of $A_2-\sqrt{-1}A_3$
is represented by $N_2$.
We have the relation 
$\nbign_2=B^{-1}\nbign_1B$.
Then, we obtain that $N_2$ and $N_0$ are conjugate.
\hfill\qed

\section{Hermitian metrics and filtered prolongation}

\subsection{Prolongation of monopoles with bounded curvature}
\label{subsection;19.1.27.2}

\subsubsection{Prolongation of mini-holomorphic bundles
with Hermitian metric}

\label{subsection;18.8.27.40}

We use the notation in \S\ref{subsection;18.8.19.20}.
Let $\nu$ denote $0$ or $\infty$.
Let $p$ be any positive integer.
Let $\nbigubar^{\lambda}_{\nu,p}$ be a neighbourhood of
$H^{\lambda}_{\nu,p}$ in $\nbigmbar^{\lambda}_{\nu,p}$.
We set $\nbigu^{\lambda}_{\nu,p}:=
 \nbigubar^{\lambda}_{\nu,p}
 \setminus
 H^{\lambda}_{\nu,p}$.
For any $\ttt\in S^1_{\lambda}$,
we put
$\nbigubar^{\lambda}_{\nu,p}(\ttt):=
 \pi_p^{-1}(\ttt)\cap\nbigubar^{\lambda}_{\nu,p}$
and 
$\nbigu^{\lambda}_{\nu,p}(\ttt):=
 \pi_p^{-1}(\ttt)\cap\nbigu^{\lambda}_{\nu,p}$.
We also set
$\nbigubar^{\lambda\cov}_{\nu,p}:=
 \ttP^{-1}(\nbigubar^{\lambda}_{\nu,p})$
and 
$\nbigu^{\lambda\cov}_{\nu,p}:=
 \ttP^{-1}(\nbigu^{\lambda}_{\nu,p})$.
For any $\ttt\in\real$,
we put 
$\nbigubar^{\lambda\cov}_{\nu,p}(\ttt):=
 (\pi^{\cov}_p)^{-1}(\ttt)\cap\nbigubar^{\lambda\cov}_{\nu,p}$
and 
$\nbigu^{\lambda\cov}_{\nu,p}(\ttt):=
 (\pi^{\cov}_p)^{-1}(\ttt)\cap\nbigu^{\lambda\cov}_{\nu,p}$.

Let $(E,\delbar_E)$ be a mini-holomorphic bundle
on $\nbigu^{\lambda}_{\nu,p}$
with a Hermitian metric $h$.
We have the Chern connection $\nabla_h$
and the Higgs field $\phi_h$.
We set $(E^{\cov},\delbar_{E^{\cov}},h^{\cov}):=
\ttP^{-1}(E,\delbar_{E},h)$
on $\nbigu^{\lambda\cov}_{\nu,p}$.
Suppose the following.
\begin{condition}
\label{condition;18.11.24.10}
 $\Bigl|
 \bigl[
 \del_{E,\ttubar},\del_{E,h,\ttu}
 \bigr]\Bigr|_{h}=O(y_0^{-2})$
and $|\phi_{h}|_{h}=O(|y_0|)$
around any point of
$H^{\lambda}_{\nu,p}$.
\hfill\qed
\end{condition}

Note that 
 $\Bigl|
 \bigl[
 \del_{E,\ttubar},\del_{E,h,\ttu}
 \bigr]\Bigr|_h=O(y_0^{-2})$
and $|\phi_h|_h=O(|y_0|)$
implies the acceptability of
the holomorphic bundles
with a Hermitian metric
$(E,\delbar_E,h)_{|\nbigu^{\lambda}_{\nu,p}(\ttt)}$
and 
$(E^{\cov},\delbar_{E^{\cov}},h^{\cov})
_{|\nbigu^{\lambda}_{\nu,p}(\ttt)}$.
Hence,
for any $\ttt\in\real$,
$E^{\cov}_{|\nbigu^{\lambda\cov}_{\nu,p}(\ttt)}$
naturally extends to a filtered bundle
$\nbigp_{\ast}\bigl(
 E^{\cov}_{|\nbigu^{\lambda\cov}_{\nu,p}(\ttt)}
 \bigr)$
over a locally free
$\nbigo_{\nbigubar^{\lambda\cov}_{\nu,p}(\ttt)}
 (\ast \nu)$-module
$\nbigp\bigl(
 E^{\cov}_{|\nbigu^{\lambda\cov}_{\nu,p}(\ttt)}
 \bigr)$.
For any $\ttt\in S^1_{\lambda}$,
$E_{|\nbigu^{\lambda}_{\nu,p}(\ttt)}$
naturally extends to a filtered bundle
$\nbigp_{\ast}\bigl(
 E_{|\nbigu^{\lambda}_{\nu,p}(\ttt)}
 \bigr)$
over a locally free
$\nbigo_{\nbigubar^{\lambda}_{\nu,p}(\ttt)}(\ast \nu)$-module
$\nbigp\bigl(
 E_{|\nbigu^{\lambda}_{\nu,p}(\ttt)}
 \bigr)$.

\begin{lem}
$(E^{\cov},\delbar_{E^{\cov}})$ uniquely extends to
$\seisuu \tte_2$-equivariant
$\nbigo_{\nbigubar^{\lambda\cov}_{\nu,p}}
 (\ast H^{\lambda\cov}_{\nu,p})$-module
$\nbigp(E^{\cov})$
such that 
\[
\nbigp(E^{\cov})_{|\nbigubar^{\lambda\cov}_{\nu,p}(\ttt)}
=\nbigp\bigl(
 E^{\cov}_{|\nbigu^{\lambda\cov}_{\nu,p}(\ttt)}
 \bigr) 
\]
for any $\ttt\in\real$.
Similarly,
$(E,\delbar_E)$ uniquely extends to
a locally free 
$\nbigo_{\nbigubar^{\lambda}_{\nu,p}}(\ast H^{\lambda}_{\nu,p})$-module
$\nbigp(E)$
such that
$\nbigp(E)_{|\nbigubar^{\lambda}_{\nu,p}(\ttt)}
=\nbigp(E_{|\nbigu^{\lambda}_{\nu,p}(\ttt)})$
for any $\ttt\in S^1_{\lambda}$.
\end{lem}
\pf
The uniqueness is clear.
Because of $|\phi_h|_h=O(|y_0|)$,
the scattering map induces an isomorphism
$\nbigp\bigl(
 E^{\cov}_{|\nbigu^{\lambda}_{\nu,p}(\ttt_1)}
 \bigr)
\simeq
 \nbigp\bigl(
 E^{\cov}_{|\nbigu^{\lambda}_{\nu,p}(\ttt_2)}
 \bigr)$
for any $\ttt_1,\ttt_2\in\real$.
Hence, the claim is clear.
\hfill\qed

\vspace{.1in}
In all,
from $(E,\delbar_E,h)$ satisfying 
Condition \ref{condition;18.11.24.10},
we obtain a locally free 
$\nbigo_{\nbigubar^{\lambda}_{\nu,p}}(\ast H^{\lambda}_{\nu,p})$-module
$\nbigp(E)$
and a filtered bundle
$\nbigp_{\ast}(E)=\bigl(
 \nbigp_{\ast}(E_{|\nbigu^{\lambda}_{\nu,p}(\ttt)})\,\big|\,\ttt\in S^1_{\lambda}
 \bigr)$
over $\nbigp(E)$.
We also obtain a locally free 
$\nbigo_{\nbigubar^{\lambda\cov}_{\nu,p}}
 (\ast H^{\lambda\cov}_{\nu,p})$-module
$\nbigp(E^{\cov})$
and a filtered bundle
$\nbigp_{\ast}(E^{\cov})=\bigl(
 \nbigp_{\ast}(E^{\cov}_{|\nbigu^{\lambda\cov}_{\nu,p}(\ttt)})
 \,\big|\,\ttt\in S^1_{\lambda}
 \bigr)$.

\subsubsection{Statements}

Let $(E,h,\nabla,\phi)$ be a monopole
with bounded curvature on 
$\nbigu^{\lambda}_{\nu,1}$.
We obtain the mini-holomorphic bundle
$(E,\delbar_E)$ with the metric $h$
on $\nbigu^{\lambda}_{\nu,1}$.
According to 
Proposition \ref{prop;18.9.2.1},
Lemma \ref{lem;19.2.7.10}
and Corollary \ref{cor;18.12.17.3},
$(E,\delbar_E,h)$ satisfies Condition \ref{condition;18.11.24.10}.
Hence, we obtain 
the locally free 
$\nbigo_{\nbigubar^{\lambda}_{\nu,1}}(\ast H^{\lambda}_{\nu,1})$-module
$\nbigp E^{\lambda}$
and a filtered bundle
$\nbigp_{\ast}(E^{\lambda})$ over $\nbigp(E^{\lambda})$.
We shall prove the following theorem
in \S\ref{subsection;18.12.17.111}
after some preliminaries.
\begin{thm}
\label{thm;18.12.17.110}
The  filtered bundle $\nbigp_{\ast}E^{\lambda}$
is good.
Moreover,
the norm estimate holds for
$(\nbigp_{\ast}E^{\lambda},h)$.
\end{thm}

There exist
$I(\phi)\subset\rnum$
and a decomposition (\ref{eq;18.9.1.10})
as in Proposition \ref{prop;18.8.14.1}.
For each $\omega$,
there exist a finite subset
$S_{\omega}\in \real^3$
as in Proposition \ref{prop;19.1.5.1}.
Moreover, for each $\veca\in S_{\omega}$,
there exits 
the $\su(2)$-representation $H_{\omega,\veca}$
determined by
$A_{i,\omega,\veca}$ $(i=1,2,3)$
in Proposition \ref{prop;19.1.5.1}.
As in \S\ref{subsection;19.1.5.10},
we obtain a monopole
\[
(E_0,h_0,\nabla_0,\phi_0)
:=
 \bigoplus_{\omega\in I(\phi)}
 \ttM(\omega,S_{\omega},\{\vecA_{\omega,\veca}\})
\]
on $\nbigu^{\lambda}_{\nu,1}$.
We obtain a good filtered bundle
$\nbigp_{\ast}E_0^{\lambda}$.

\begin{thm}
\label{thm;19.1.6.10}
There exists an isomorphism
$\ttG(\nbigp_{\ast}E^{\lambda})
\simeq
\ttG(\nbigp_{\ast}E_0^{\lambda})$.
\end{thm}

\subsection{Prolongation of asymptotically mini-holomorphic bundles}

\label{subsection;18.8.27.50}

Let $E$ be a $C^{\infty}$-vector bundle on 
$\nbigu^{\lambda}_{\nu,p}$
with a Hermitian metric $h$,
a unitary connection $\nabla$
and an anti-Hermitian endomorphism $\phi$.
Let $(\alpha,\tau)$ denote
the local mini-complex coordinate system
on $\nbigu^{\lambda}_{\nu,p}$
as in \S\ref{subsection;18.12.17.20}.
We define
differential operators
$\del_{E,\ttubar}$,
$\del_{E,h,\ttu}$
and 
$\del_{E,\ttt}$
by the following formula:
\begin{equation}
 \del_{E,\ttubar}
=\frac{1-\lambdabar\ttg_1}{1+|\lambda|^2}
 \nabla_{\alphabar}
-\frac{1}{2\sqrt{-1}}
 \frac{\ttgbar_1+\lambda}{1+|\lambda|^2}
 (\nabla_{\tau}-\sqrt{-1}\phi),
\end{equation}
\begin{equation}
 \del_{E,h,\ttu}
=\frac{1-\lambda\ttgbar_1}{1+|\lambda|^2}
 \nabla_{\alpha}
+\frac{1}{2\sqrt{-1}}
 \frac{\ttg_1+\lambdabar}{1+|\lambda|^2}
 (\nabla_{\tau}+\sqrt{-1}\phi),
\end{equation}
\begin{equation}
\del_{E,\ttt}:=\nabla_{\tau}-\sqrt{-1}\phi.
\end{equation}

We assume that 
$(E,h,\nabla,\phi)$ satisfies 
the following condition
in \S\ref{subsection;19.2.7.30}--\ref{subsection;19.1.6.1}.
\begin{condition}
\label{condition;19.1.6.20}
$[\del_{E,\ttubar},\del_{E,h,\ttu}]=O(y_0^{-2})$,
and 
$|\phi|_h$ is bounded.
\hfill\qed
\end{condition}

\subsubsection{Case 1}
\label{subsection;19.2.7.30}

In this subsection,
we assume the following additional condition.
\begin{condition}
\label{condition;18.11.24.20}
 For any $k\geq 0$,
 there exists $\epsilon(k)>0$
 such that the following holds for 
 $(\kappa_1,\ldots,\kappa_k)\in\{0,1,2\}^k$:
\[
 \Bigl|
 \nabla_{\kappa_1}\circ\cdots\circ\nabla_{\kappa_k}
 \bigl(
 [\del_{E,\ttubar},\del_{E,\ttt}]
 \bigr)
 \Bigr|_h
 =O\bigl(e^{-\epsilon(k) y_0^2}\bigr).
\]
\hfill\qed
\end{condition}

By taking the pull back by $\ttP$,
we obtain 
$(E^{\cov},h^{\cov})$
with 
the differential operators
$\del_{E^{\cov},\ttubar}$
and $\del_{E^{\cov},\ttt}$.
The restrictions
$(E^{\cov},\del_{E^{\cov},\ttubar},h^{\cov})
 _{|\nbigu^{\lambda\cov}_{\nu,p}(\ttt)}$
and 
$(E,\del_{E,\ttubar},h)
 _{|\nbigu^{\lambda}_{\nu,p}(\ttt)}$
are holomorphic vector bundles
with a Hermitian metric.
By the assumption
$\bigl[
 \del_{E,\ttubar},\del_{E,h,\ttu}
 \bigr]=O\bigl(y_0^{-2}\bigr)$,
$E^{\cov}_{|\nbigu^{\lambda\cov}_{\nu,p}(\ttt)}$
extends to a filtered bundle
$\nbigp_{\ast}\bigl(
 E^{\cov}_{|\nbigu^{\lambda\cov}_{\nu,p}(\ttt)}
 \bigr)$
for any $\ttt\in\real$.
Similarly,
for any $\ttt\in S^1_{\lambda}$,
$E_{|\nbigu^{\lambda}_{\nu,p}(\ttt)}$
extends to a filtered bundle
$\nbigp_{\ast}\bigl(
 E_{|\nbigu^{\lambda}_{\nu,p}(\ttt)}
 \bigr)$.

\begin{lem}
\mbox{}
\begin{itemize}
\item
For each $a\in\real$,
$E^{\cov}$ uniquely extends to
a $C^{\infty}$-bundle
$\nbigp_a(E^{\cov})$ on $\nbigubar^{\lambda\cov}_{\nu,p}$
such that
$\nbigp_a(E^{\cov})_{|\nbigubar^{\lambda\cov}_{\nu,p}(\ttt)}
=\nbigp_a\bigl(E^{\cov}_{|\nbigu^{\lambda\cov}_{\nu,p}(\ttt)}\bigr)$.
\item
$\del_{E^{\cov},\ttt}$
and $\del_{E^{\cov},\ttubar}$ extend
to $C^{\infty}$-differential operators
on $\nbigp_a(E^{\cov})$.
\item
$\bigl[\del_{E^{\cov},\ttt},\del_{E^{\cov},\ttubar}\bigr]
 _{|\Hhat^{\lambda\cov}_{\nu,p}}
 =0$.
\end{itemize}
Similar claims hold for $(E,\delbar_E,h)$
on $\nbigubar^{\lambda}_{\nu,p}$.
\end{lem}
\pf
Take a holomorphic frame $\vecv$ of
$\nbigp_a(E^{\cov}_{|\nbigu^{\lambda\cov}_{\nu,p}(0)})$.
We obtain a $C^{\infty}$-frame
$\vecvtilde$ of $E^{\cov}$
such that
(i) $\del_{E^{\cov},\ttt}\vecvtilde=0$,
(ii) $\vecvtilde_{|\nbigu^{\lambda\cov}_{\nu,p}(0)}=\vecv$.
We have the matrix valued function $\nbiga$
on $\nbigu^{\lambda\cov}_{\nu,p}$
determined by
$\del_{E,\ttubar}\vecvtilde=\vecvtilde\nbiga$.
For each $(\ell_1,\ell_2,\ell_3)\in \seisuu_{\geq 0}^3$,
there exists
$\epsilon(\ell_1,\ell_2,\ell_3)>0$
such that
$\del_{\ttt}^{\ell_1}\del_{\ttubar}^{\ell_2}\del_u^{\ell3}
 \nbiga
=O\bigl(e^{-\epsilon(\ell_1,\ell_2,\ell_3) y_0^2}\bigr)$.
It implies that
for each $(\ell_1,\ell_2,\ell_3)\in \seisuu_{\geq 0}^3$,
there exists
$\epsilon_1(\ell_1,\ell_2,\ell_3)>0$
such that
$\del_{\ttt}^{\ell_1}
 \del_{\ttU_{\nu,p}}^{\ell_2}
 \del_{\ttUbar_{\nu,p}}^{\ell_3}\nbiga
=O\bigl(e^{-\epsilon_1(\ell_1,\ell_2,\ell_3)y_0^2}\bigr)$.
Hence, $\nbiga$ extends to
a $C^{\infty}$-function on
$\nbigubar^{\lambda\cov}_{\nu,p}$.
Moreover we have
$\nbiga_{|\Hhat^{\lambda\cov}_{\nu,p}}=0$.

We extend $E^{\cov}$
to $\nbigp_a(E^{\cov})$
by using the frame $\vecvtilde$.
The bundle $\nbigp_a(E^{\cov})$
is independent of the choice of $\vecv$.
The operator $\del_{E^{\cov},\ttt}$ naturally induces
a $C^{\infty}$-differential operator
on $\nbigp_a(E^{\cov})$.
Because $\nbiga$ extends to a $C^{\infty}$-function 
on $\nbigubar^{\lambda\cov}_{\nu,p}$,
$\del_{E^{\cov},\ttubar}$ also 
induces a $C^{\infty}$-differential operator
on $\nbigp_a(E^{\cov})$.
Because 
$\nbiga_{|\Hhat^{\lambda\cov}_{\nu,p}}=0$,
we obtain
$\bigl[\del_{E^{\cov},\ttt},\del_{E^{\cov},\ttubar}\bigr]
 _{|\Hhat^{\lambda\cov}_{\nu,p}}
 =0$.
It is easy to see that
$\nbigp_a(E^{\cov})_{|\nbigubar^{\lambda\cov}_{\nu,p}(\ttt)}
=\nbigp_a\bigl(E^{\cov}_{|\nbigu^{\lambda\cov}_{\nu,p}(\ttt)}\bigr)$
in a natural way
for any $\ttt\in\real$.
\hfill\qed

\begin{cor}
If Condition {\rm\ref{condition;18.11.24.20}}
is satisfied, we obtain 
a locally free $\nbigo_{\Hhat^{\lambda}_{\nu,p}}$-module
$\nbigp_a(E)_{|\Hhat^{\lambda}_{\nu,p}}$
for each $a\in\real$,
and hence a regular filtered bundle
$\nbigp_{\ast}(E)_{|\Hhat^{\lambda}_{\nu,p}}$
over 
$(\Hhat^{\lambda}_{\nu,p},H^{\lambda}_{\nu,p})$.
\hfill\qed
\end{cor}

\subsubsection{Case 2}
\label{subsection;19.1.6.1}

In this subsection,
we assume the following additional condition
which is weaker than Condition \ref{condition;18.11.24.20}.

\begin{condition}
\label{condition;18.11.24.50}
 For any $k\geq 0$,
 there exists $\epsilon(k)>0$
 such that the following holds for 
 $(\kappa_1,\ldots,\kappa_k)\in\{0,1,2\}^k$:
\begin{equation}
\label{eq;18.12.17.40}
 \Bigl|
 \nabla_{\kappa_1}\circ\cdots\circ\nabla_{\kappa_k}
 \bigl(
 [\del_{E,\ttubar},\del_{E,\ttt}]
 \bigr)
 \Bigr|_h
 =O\bigl(e^{-\epsilon(k) |y_0|}\bigr).
\end{equation}
\hfill\qed
\end{condition}

We obtain the filtered bundles
$\nbigp_{\ast}(E_{|\nbigu^{\lambda}_{\nu,p}(\ttt)})$
$(\ttt\in S^1_{\lambda})$,
and the induced vector spaces
for any $\ttt\in S^1_{\lambda}$ and $a\in\real$:
\[
 \Gr^{\nbigp}_a(E,\ttt):=
\nbigp_{a}(E_{|\nbigu^{\lambda}_{\nu,p}(\ttt)})
 \big/
\nbigp_{<a}(E_{|\nbigu^{\lambda}_{\nu,p}(\ttt)}).
\]
Similarly,
we obtain the vector spaces
$\Gr^{\nbigp}_a(E^{\cov},\ttt):=
 \nbigp_a(E^{\cov}_{|\nbigu^{\lambda\cov}_{\nu,p}(\ttt)})
 \big/
  \nbigp_{<a}(E^{\cov}_{|\nbigu^{\lambda\cov}_{\nu,p}(\ttt)})$
for any $a\in\real$
and $\ttt\in \real$.

\begin{lem}
For any $a\in\real$,
and for any $\ttt_1,\ttt_2\in\real$,
we have natural isomorphisms
\begin{equation}
\label{eq;18.8.16.10}
 \Gr^{\nbigp}_a(E^{\cov},\ttt_1)
\simeq
 \Gr^{\nbigp}_a(E^{\cov},\ttt_2).
\end{equation}
\end{lem}
\pf
We take a section $s^{\ttt_1}$ of
$\nbigp_a(E^{\cov}_{|\nbigu^{\lambda\cov}_{\nu,p}(\ttt_1)})$.
By definition,
we have
$|s^{\ttt_1}|_h=O(|\ttU_{\nu,p}|^{-a-\epsilon})$
for any $\epsilon>0$.
By the parallel transport with respect to $\del_{E,\ttt}$,
we obtain an induced $C^{\infty}$-section $s^{\ttt_2}$ of
$E^{\cov}_{|\nbigu^{\lambda\cov}_{\nu,p}(\ttt_2)}$.
By (\ref{eq;18.12.17.40}),
we have
$\del_{E^{\cov},\ttubar}(s^{\ttt_2})=O(|\ttU_{\nu,p}|^{-a+\delta})$
for some $\delta>0$,
which implies
$\del_{E^{\cov},\ttUbar_{\nu,p}}(s^{\ttt_2})
=O\bigl(|\ttU_{\nu,p}|^{-a-1+\delta}\bigr)$.
There exists a $C^{\infty}$-section $b^{\ttt_2}$ of 
$E^{\cov}_{|\nbigu^{\lambda\cov}_{\nu,p}(\ttt_2)}$
such that
$\del_{E^{\cov},\ttubar}(s^{\ttt_2}+b^{\ttt_2})=0$
and $|b^{\ttt_2}|=O\bigl(|\ttU_{\nu,p}|^{-a+\delta_1}\bigr)$
for some $\delta_1>0$.
Then, $\stilde^{\ttt_2}=s^{\ttt_2}+b^{\ttt_2}$ is a section of
$\nbigp_a(E^{\cov}_{|\nbigu^{\lambda\cov}_{\nu,p}(\ttt_2)})$,
which induces an element of
$\Gr^{\nbigp}_a(E^{\cov},\ttt_2)$.
It induces a well defined isomorphism
(\ref{eq;18.8.16.10}).
\hfill\qed

\vspace{.1in}

Thus, for any $a\in\real$,
we obtain a local system
$\Gr^{\nbigp}_a(E^{\cov})$ on $\real$,
which is naturally $\seisuu \tte_2$-equivariant.
Thus, we obtain a local system
$\Gr^{\nbigp}_a(E)$ on $S_{\lambda}^1$
for any $a\in\real$.

We obtain the filtration $W$
on $\Gr^{\nbigp}_a(E)$
as the weight filtration of 
the nilpotent endomorphism
obtained as the logarithm
of the unipotent part of the monodromy.

Let $\ttt\in S^1_{\lambda}$.
Let $\vecv$ be a holomorphic frame of
$\nbigp_a(E_{|\nbigu^{\lambda}_{\nu,p}(\ttt)})$
compatible with the filtrations $\nbigp$
and $W$.
We obtain the numbers
$b(v_i):=\deg^{\nbigp}(v_i)$
and $k(v_i):=\deg^{W}(v_i)$.
Let $h_0$ be the metric of
$E_{|\nbigu^{\lambda}_{\nu,p}(\ttt)}$
determined by
$h_0(v_i,v_j)=0$ $(i\neq j)$
and
$h_0(v_i,v_i)=|\ttU_{\nu,p}|^{-2b(v_i)}
 \bigl|\log|\ttU_{\nu,p}|\bigr|^{k(v_i)}$.
We say that the norm estimate holds for
$(\nbigp_{\ast}(E_{|\nbigu^{\lambda}_{\nu,p}(\ttt)}),h)$
if $h_0$ and $h_{|\nbigu^{\lambda}_{\nu,p}}$
are mutually bounded.
The following lemma is easy to see.
\begin{lem}
If  the norm estimate holds for
$(\nbigp_{\ast}(E_{|\nbigu^{\lambda}_{\nu,p}(\ttt_0)}),h)$
at some $\ttt_0$,
then the norm estimate holds for
$(\nbigp_{\ast}(E_{|\nbigu^{\lambda}_{\nu,p}(\ttt)}),h)$
for any $\ttt\in S^1_{\lambda}$.
\hfill\qed
\end{lem}

\subsubsection{Comparison}

Let $(E,h,\nabla,\phi)$ be as in
\S\ref{subsection;19.1.6.1}.
Let $E^{\circledcirc}$ be a $C^{\infty}$-vector bundle
on $\nbigu^{\lambda}_{\nu,p}$
with a Hermitian metric $h^{\circledcirc}$,
a unitary connection $\nabla^{\circledcirc}$
and an anti-Hermitian endomorphism $\phi^{\circledcirc}$.
Let $F:E\simeq E^{\circledcirc}$ be a $C^{\infty}$-isomorphism.
Let $b^{\circ}$ be the endomorphism of $E$
determined by
$h=F^{\ast}(h^{\circledcirc}) b^{\circledcirc}$.
Assume the following condition on $F$.
\begin{condition}
\label{condition;19.2.8.1}
For any $k\in\seisuu_{\geq 0}$,
there exists $\epsilon(k)>0$
such that the following holds
for any $(\kappa_1,\ldots,\kappa_k)\in\{0,1,2\}^k$:
\[
\bigl|
 \nabla_{\kappa_1}\circ\cdots\circ
 \nabla_{\kappa_k}(b^{\circledcirc}-\id)
\bigr|_h
=O\bigl(e^{-\epsilon(k)|y_0|}\bigr),
\]
\[
 \bigl|
 \nabla_{\kappa_1}\circ\cdots\circ
 \nabla_{\kappa_k}(\nabla-F^{\ast}\nabla^{\circledcirc})
\bigr|_h
=O\bigl(e^{-\epsilon(k)|y_0|}\bigr),
\]
\[
  \bigl|
 \nabla_{\kappa_1}\circ\cdots\circ
 \nabla_{\kappa_k}(\phi-F^{\ast}\phi^{\circledcirc})
\bigr|_h
=O\bigl(e^{-\epsilon(k)|y_0|}\bigr).
\]
\end{condition}
Note that
$(E^{\circledcirc},h^{\circledcirc},
 \nabla^{\circledcirc},\phi^{\circledcirc})$
also satisfies Condition \ref{condition;18.11.24.50}.

\begin{lem}
\label{lem;19.1.6.21}
For any $a\in\real$,
there exists a naturally induced isomorphism
of the local systems
$\Gr^{\nbigp}_a(E)\simeq
 \Gr^{\nbigp}_a(E^{\circledcirc})$.
Moreover, if 
the norm estimate holds for
$(\nbigp_{\ast}E^{\circledcirc}
 _{|\nbigu^{\lambda}_{\nu,p}(\ttt)},
 h^{\circledcirc})$,
then the norm estimate also holds for 
$(\nbigp_{\ast}E
 _{|\nbigu^{\lambda}_{\nu,p}(\ttt)},
 h)$.
\end{lem}
\pf
Let $s$ be a holomorphic section of
$\nbigp_a(E_{|\nbigu^{\lambda}_{\nu,p}(\ttt)})$.
Let $[s]$ be the induced element of
$\Gr^{\nbigp}_a(E_{|\nbigu^{\lambda}_{\nu,p}(\ttt)})$.
There exists a $C^{\infty}$-section $c$
of $E^{\circledcirc}_{|\nbigu^{\lambda}_{\nu,p}(\ttt)}$
such that 
(i) $\stilde:=F(s)-c$ is a holomorphic section of
$\nbigp_a(E^{\circledcirc})$,
(ii) $|c|=O(|\ttU_{\nu,p}|^{-a+\epsilon})$ for some $\epsilon>0$.
Let $[\stilde]$ denote the induced element of
$\Gr^{\nbigp}_a(E^{\circledcirc}_{|\nbigu^{\lambda}_{\nu,p}(\ttt)})$.
Then, $[\stilde]$ depends only on $[s]$.
Thus, we obtain
$\Gr_a^{\nbigp}(E_{|\nbigu^{\lambda}_{\nu,p}(\ttt)})
\lrarr
 \Gr^{\nbigp}_a(E^{\circledcirc}_{|\nbigu^{\lambda}_{\nu,p}(\ttt)})$.
This procedure induces the desired isomorphism.
The claim for the norm estimate is easy to check.
\hfill\qed

\subsection{Prolongation to good filtered bundles}
\label{subsection;18.9.1.20}

Let $(E,\delbar_E)$ be a mini-holomorphic bundle
on $\nbigu^{\lambda}_{\nu,p}$
with a Hermitian metric $h$.
The Chern connection $\nabla$
and the Higgs field $\phi$
are associated to $(E,\delbar_E,h)$.
Let $(y_0,y_1,y_2)$ be the local coordinate system of $\nbigm^0$ 
induced by $z=y_1+\sqrt{-1}y_2$ and $\Image(w)=y_0$,
as in \S\ref{subsection;18.12.17.50}.
Let $\nabla_{h,i}$ denote $\nabla_{h,y_i}$.
Suppose that the following condition is satisfied.
\begin{condition}
\label{condition;18.11.24.11}
Condition {\rm\ref{condition;18.11.24.10}} is satisfied.
Moreover, there exists
an orthogonal decomposition
\[
 (E,h,\phi)=
 \bigoplus_{\omega\in\frac{1}{p}\seisuu}
 \bigl(E^{\bullet}_{\omega},h^{\bullet}_{\omega},
 \phi^{\bullet}_{\omega}\bigr)
\]
such that
the following holds.
\begin{itemize}
\item
 $\phi^{\bullet}_{\omega}-
 \bigl(2\pi\sqrt{-1}\omega/\Vol(\Gamma)\bigr)
 y_0
 \id_{E^{\bullet}_{\omega}}$
are bounded.
\item
We have the decomposition
 $\nabla=\nabla^{\bullet}+\rho$,
 where $\nabla^{\bullet}$ is the direct sum 
 of connections $\nabla^{\bullet}_{\omega}$
 of $E^{\bullet}_{\omega}$,
and $\rho$ is a section of
 $\bigoplus_{\omega_1\neq\omega_2}
 \Hom(E^{\bullet}_{\omega_1},E^{\bullet}_{\omega_2})\otimes\Omega^1$.
 Then,
 for any $k\in\seisuu_{\geq 0}$,
 we have $\epsilon(k)>0$ such that
the following holds for any $(\kappa_1,\ldots,\kappa_k)
 \in\{0,1,2\}^k$:
\begin{equation}
\label{eq;18.12.17.20}
 \bigl|
 \nabla^{\bullet}_{\kappa_1}\circ
 \cdots\circ
 \nabla^{\bullet}_{\kappa_k}
 \rho
 \bigr|=O(e^{-\epsilon(k)y_0^2}).
\end{equation}
\hfill\qed
\end{itemize}
\end{condition}

\begin{prop}
\label{prop;18.8.27.51}
$\nbigp_{\ast}(E)$
is a good filtered bundle
over $\nbigp(E)$.
\end{prop}
\pf
Let $\tti_{\omega}:E^{\bullet}_{\omega}\lrarr E$
denote the inclusion,
and let $\ttp_{\omega}:E\lrarr E^{\bullet}_{\omega}$
denote the orthogonal projection.
We set
$\del_{E^{\bullet}_{\omega},\ttubar}:=
  \ttp_{\omega}\circ \del_{E,\ttubar}\circ\tti_{\omega}$
and 
$\del_{E^{\bullet}_{\omega},\ttt}:=
  \ttp_{\omega}\circ \del_{E,\ttt}\circ\tti_{\omega}$.
Similarly, we obtain connection 
$\nabla_{E^{\bullet}_{\omega}}$
on $E^{\bullet}_{\omega}$.

We set
$E_{\omega}:=\vecsfL_p(-\omega)\otimes
 E^{\bullet}_{\omega}$.
(See \S\ref{subsection;18.11.24.30}
for the monopole $\vecsfL_{p}(-\omega)$
and the 
 $\nbigo_{\nbigubar^{\lambda}_{\nu,p}}
 (\ast H^{\lambda}_{\nu,p})$-module
$\nbigp\nbigl^{\lambda}_p(-\omega)$.)
Let $h_{\omega}$ be the induced metric
on $E_{\omega}$.
We obtain the differential operators
$\del_{E_{\omega},\ttubar}$
and $\del_{E_{\omega},\ttt}$
from 
the mini-holomorphic structure of 
$\nbigl^{\lambda}_p(-\omega)$,
and the operators
$\del_{E^{\bullet}_{\omega},\ttubar}$
and $\del_{E^{\bullet}_{\omega},\ttt}$.
We obtain the connection 
$\nabla_{E_{\omega}}$ 
of $E_{\omega}$
from 
$\nabla_{E^{\bullet}_{\omega}}$ 
and the connection of
$\vecsfL_p(-\omega)$.
Similarly, we obtain the anti-Hermitian
endomorphism $\phi_{\omega}$
from $\phi_{\omega}^{\bullet}$
and the anti-Hermitian endomorphism of
$\vecsfL_p(-\omega)$.
Then, 
$\bigl(E_{\omega},h_{\omega},
\nabla_{\omega},\phi_{\omega}
\bigr)$ 
satisfies Condition \ref{condition;18.11.24.20},
and the operators 
$\del_{E_{\omega},\ttubar}$
and $\del_{E_{\omega},\ttt}$
are induced by $\nabla_{\omega}$
and $\phi_{\omega}$
as in \S\ref{subsection;18.8.27.50}.
We obtain $C^{\infty}$-bundles
$\nbigp^{C^{\infty}}_a(E_{\omega})$
for each $a\in\real$.
We may regard them as 
$\nbigc^{\infty}_{\nbigubar^{\lambda}_{\nu,p}}$-modules.

Because
$E^{\bullet}_{\omega}
=\vecsfL_p(\omega)\otimes E_{\omega}$,
we have the following natural $C^{\infty}$-identification
on $\nbigu^{\lambda}_{\nu,p}$:
\begin{equation}
\label{eq;18.11.24.40}
 E\simeq
 \bigoplus_{\omega}
 \vecsfL_p(\omega)\otimes
 E_{\omega}.
\end{equation}
\begin{lem}
\label{lem;18.12.17.30}
The isomorphism {\rm(\ref{eq;18.11.24.40})}
extends to an isomorphism
of $\nbigc^{\infty}_{\nbigubar^{\lambda}_{\nu,p}}$-modules:
\[
F:
 \nbigp(E)\otimes_{\nbigo_{\nbigubar^{\lambda}_{\nu,p}}}
 \nbigc^{\infty}_{\nbigubar^{\lambda}_{\nu,p}}
\simeq
 \bigoplus
 \nbigp\bigl(\nbigl^{\lambda}_p(\omega)\bigr)
\otimes_{\nbigo_{\nbigubar^{\lambda}_{\nu,p}}}
 \nbigp^{C^{\infty}}_0(E_{\omega}).
\]
Moreover, 
$F_{|\Hhat^{\lambda}_{\nu,p}}$ 
is mini-holomorphic.
\end{lem}
\pf
We take $\ttt_0\in S^{1}_{\lambda}$
and a neighbourhood $I$ of $\ttt_0$ in $S^1_{\lambda}$.
We set
$\nbigubar^{\lambda}_{\nu,p}(I):=
 \pi_p^{-1}(I)\cap
\nbigubar^{\lambda}_{\nu,p}$,
and 
$\nbigu^{\lambda}_{\nu,p}(I):=
\nbigubar^{\lambda}_{\nu,p}(I)
\setminus H^{\lambda}_{\nu,p}$.
We also put
$\Hhat^{\lambda}_{\nu,p}(I):=
 \Hhat^{\lambda}_{\nu,p}
\cap
 \pi_p^{-1}(I)$.
We take a $C^{\infty}$-frame
$\vecv_{\omega}$ of 
$\nbigp^{C^{\infty}}_0(E_{\omega})_{|\nbigubar^{\lambda}_{\nu,p}(I)}$
such that 
$\vecv_{\omega|\Hhat^{\lambda}_{\nu,p}(I)}$
is mini-holomorphic.
Fixing a lift of $\ttt_0$ to $\real$,
we obtain the mini-holomorphic frame
$\sfv^{\lambda}_{p,\omega}$
of $\nbigp\nbigl^{\lambda}_p(\omega)
 _{|\nbigubar^{\lambda}_{\nu,p}(I)}$.
We obtain a $C^{\infty}$-frame
$\sfv^{\lambda}_{p,\omega}\otimes
 \vecv_{\omega}$
of $\nbigp\nbigl^{\lambda}_p(\omega)\otimes
 \nbigp^{C^{\infty}}_0E_{\omega}$
on $\nbigubar^{\lambda}_{\nu,p}(I)$.
They induce a frame $\vecu$ of
$\bigoplus
\nbigp(\nbigl^{\lambda}_p(\omega))
 \otimes
 \nbigp^{C^{\infty}}_0(E_{\omega})$
on $\nbigubar^{\lambda}_{\nu,p}(I)$.
By the frame $\vecu$,
we also obtain a $C^{\infty}$-vector bundle
$V$ on $\nbigubar^{\lambda}_{\nu,p}(I)$
with an isomorphism
$V_{|\nbigu^{\lambda}_{\nu,p}(I)}
\simeq E$.
We may naturally regard $V$
as a $\nbigc^{\infty}_{\nbigubar^{\lambda}_{\nu,p}}$-submodule
of 
$\bigoplus
\nbigp\nbigl^{\lambda}_p(\omega)
 \otimes
 \nbigp^{C^{\infty}}_0E_{\omega}$.

Let $A_{\ttUbar_{\nu,p}}$ and $A_{\ttt}$
be the matrix valued functions on 
$\nbigu^{\lambda}_{\nu,p}$
determined by
$\del_{E,\ttUbar_{\nu,p}}\vecu_{|\nbigu^{\lambda}_{\nu,p}}
=\vecu_{|\nbigu^{\lambda}_{\nu,p}}\cdot A_{\ttUbar_{\nu,p}}$
and 
$\del_{E,\ttt}\vecu_{|\nbigu^{\lambda}_{\nu,p}}
=\vecu_{|\nbigu^{\lambda}_{\nu,p}}\cdot A_{\ttt}$.
By the decay condition (\ref{eq;18.12.17.20}),
$A_{\ttUbar_{\nu,p}}$ and $A_{\ttt}$ 
extend to $C^{\infty}$-functions
on $\nbigubar^{\lambda}_{\nu,p}(I)$,
and 
$A_{\ttUbar_{\nu,p}|\pi^{-1}(I)\cap\Hhat^{\lambda}_{\nu,p}}
=A_{\ttt|\pi^{-1}(I)\cap\Hhat^{\lambda}_{\nu,p}}=0$.
Hence, $\del_{E,\ttUbar_{\nu,p}}$ and $\del_{E,\ttt}$
induce a mini-holomorphic structure on $V$.
There exists a mini-holomorphic frame 
$\vecw=(w_i)$ of $V$
on $\nbigubar^{\lambda}_{\nu,p}(I)$.
We have $|w_i|_h=O(|\ttU_{\nu,p}|^{-N})$ for some $N$.
Hence, we obtain that $w_i$
induce mini-holomorphic sections of
$\nbigp(E)$ on $\nbigubar^{\lambda}_{\nu,p}(I)$,
which are also denoted by the same notation.
Because $\vecw_{|\nbigu^{\lambda}_{\nu,p}(I)}$
is a holomorphic frame of
$E_{|\nbigu^{\lambda}_{\nu,p}(I)}$,
we obtain that
$\vecw$ is a frame of $\nbigp(E)$.
Then, the claim of the lemma follows.
\hfill\qed

\begin{lem}
For each $\ttt\in S^1_{\lambda}$,
$F$ induces an isomorphism of filtered bundles
\[
 \nbigp_{\ast}(E_{|\nbigu^{\lambda}_{\nu,p}(\ttt)})
\simeq
 \bigoplus
 \nbigp_{\ast}\bigl(
 \nbigl^{\lambda}_p(-\omega)_{|\nbigu^{\lambda}_{\nu,p}(\ttt)}
\bigr)
\otimes
 \nbigp_{\ast}\bigl(E_{\omega|\nbigu^{\lambda}_{\nu,p}(\ttt)}\bigr).
\]
\end{lem}
\pf
Take $\ttt\in S^{1}_{\lambda}$.
Let $\nbigl^{\lambda}_p(\omega)^{\ttt}$
denote the restriction of
$\nbigl^{\lambda}_p(\omega)$
to $\nbigu^{\lambda}_{\nu,p}(\ttt)$.
We set
$E_{\omega}^{\ttt}:=E_{\omega|\nbigu^{\lambda}_{\nu,p}(\ttt)}$.
Let $E^{\ttt}$ denote the restriction of $E$
to $\nbigu^{\lambda}_{\nu,p}(\ttt)$.

Let $s$ be a holomorphic section of
$\nbigp_{a}(
 \nbigl^{\lambda}_p(\omega)^{\ttt}
 \otimes
 E^{\ttt}_{\omega})$.
In particular,
$|s|_h=O(|\ttU_{\nu,p}|^{-a-\epsilon})$
for any $\epsilon>0$.
According to Lemma \ref{lem;18.12.17.30},
$s$ induces a section of
$\nbigp E^{\ttt}\otimes_{\nbigo_{\nbigubar^{\lambda}_{\nu,p}(\ttt)}}
 \nbigc^{\infty}_{\nbigubar^{\lambda}_{\nu,p}(\ttt)}$.
Note that
$\del_{E^{\ttt},\ttUbar_{\nu,p}}s=
 O\bigl(e^{-\epsilon_1(\log|\ttU_{\nu,p}|)^2}\bigr)$
for some $\epsilon_1>0$.
Hence, for any $N>0$,
there exists a $C^{\infty}$-section $b_N$ of  $E^{\ttt}$
such that
$|b_N|=O(|\ttU_{\nu,p}|^N)$
and 
$\del_{E^{\ttt},\ttUbar_{\nu,p}}(s-b_N)=0$.
Because $|s-b_N|_h=O(|\ttU_{\nu,p}|^{-a-\epsilon})$
for any $\epsilon>0$,
we obtain that 
$s-b_N$ is a section of $\nbigp_a(E^{\ttt})$.
Then, we obtain that
$s$ is a $C^{\infty}$-section of
$\nbigp_aE^{\ttt}\otimes_{\nbigo_{\nbigubar^{\lambda}_{\nu,p}(\ttt)}}
 \nbigc^{\infty}_{\nbigubar^{\lambda}_{\nu,p}(\ttt)}$.

We take a lift of $\ttt$ to $\real$,
and we set
\[
c(\omega,\ttt):=
 \left\{
 \begin{array}{ll}
 \omega \ttt/\gminit^{\lambda}
 & (\nu=0)
 \\
-\omega \ttt/\gminit^{\lambda}
 & (\nu=\infty).
 \end{array}
 \right.
\]
For $a\in\real$,
we take a holomorphic frame
$\vecv^{\ttt}_a$ of
$\nbigp_{a-c(\omega,\ttt)}(E_{\omega}^{\ttt})$,
which is compatible with the parabolic structure.
Let
$\sfv^{\lambda|\ttt}_{p,\omega}$
denote the restriction of
$\sfv^{\lambda}_{p,\omega}$
to $\nbigu^{\lambda}_{\nu,p}(\ttt)$.
We obtain a holomorphic frame
$\sfv^{\lambda|\ttt}_{p,\omega}
 \otimes\vecv^{\ttt}_{a}$
of $\nbigp_a\bigl(
 \nbigl^{\lambda}_{p}(\omega)^{\ttt}
\otimes
 E^{\ttt}_{\omega}
 \bigr)$.
We obtain an induced holomorphic frame
$\vecu^{\ttt}$ of 
$\bigoplus
 \nbigp_a\bigl(\nbigl^{\lambda}_p(\omega)^{\ttt}
 \otimes
 E^{\ttt}_{\omega}
 \bigr)$.
As observed above, $\vecu^{\ttt}$ induces
a tuple of $C^{\infty}$-sections of $\nbigp_a(E^{\ttt})$,
and $\vecu^{\ttt}_{|\widehat{0}}$ 
are tuples of holomorphic sections of 
$\nbigp_a(E^{\ttt})_{|\widehat{0}}$.
Hence, we obtain that
\[
 \bigoplus
 \nbigp_a\bigl(\nbigl^{\lambda}_p(\omega)^{\ttt}
 \otimes
 E^{\ttt}_{\omega}
 \bigr)_{|\widehat{0}}
\subset
 \nbigp_a(E^{\ttt})_{|\widehat{0}}.
\]

For each $u_i^{\ttt}$,
we have $\omega(i)$
such that $u_i^{\ttt}$
is a section of 
$\nbigp_a\bigl(
 \nbigl^{\lambda}_p(\omega(i))^{\ttt}\otimes
 E^{\ttt}_{\omega(i)}
 \bigr)$.
Moreover, we obtain $a-1<b(i)\leq a$
such that 
$u_i^{\ttt}$ is a section of 
$\nbigp_{b(i)}\bigl(
 \nbigl^{\lambda}_p(\omega)^{\ttt}\otimes
 E^{\ttt}_{\omega}
 \bigr)$,
and that 
the induced element in
$\Gr^{\nbigp}_{b(i)}\bigl(
 \nbigl^{\lambda}_p(\omega)^{\ttt}\otimes
 E^{\ttt}_{\omega}
\bigr)$
is non-zero.
We set
$u_i^{\prime\ttt}:=u_i^{\ttt}|\ttU_{\nu,p}|^{b(i)}$.
Then, for any $\epsilon>0$,
there exists $C(\epsilon)>1$ such that 
\[
 C(\epsilon)^{-1}
 |\ttU_{\nu,p}|^{\epsilon}
\leq
 \bigl|
 \bigwedge u_i^{\prime\ttt}
 \bigr|_h
\leq
 C(\epsilon)|\ttU_{\nu,p}|^{-\epsilon}.
\]

We can take a holomorphic section
$\utilde_i^{\ttt}$
of 
$\nbigp_{b(i)}\bigl(E \bigr)$
such that
$\utilde_i^{\ttt}-u_i^{\ttt}=O(|\ttU_{\nu,p}|^{-b(i)+N})$
for some $N>0$.
We set
$\utilde_i^{\prime\ttt}:=
 \utilde_i^{\ttt}|\ttU_{\nu,p}|^{b(i)}$.
Then, 
for any $\epsilon>0$,
there exists $C(\epsilon)>1$ such that 
\[
 C(\epsilon)^{-1}
 |\ttU_{\nu,p}|^{\epsilon}
\leq
 \bigl|
 \bigwedge \utilde_i^{\prime\ttt}
 \bigr|_h
\leq
 C(\epsilon)|\ttU_{\nu,p}|^{-\epsilon}.
\]
It implies that
$(\utilde_i)$ is a holomorphic frame of
$\nbigp_a(E^{\ttt})$
compatible with the parabolic structure.
Thus, we obtain
$\bigoplus
 \nbigp_a\bigl(\nbigl^{\lambda}_p(\omega)^{\ttt}
 \otimes
 E^{\ttt}_{\omega}
 \bigr)_{|\widehat{0}}
=\nbigp_a(E^{\ttt})_{|\widehat{0}}$.
\hfill\qed

\vspace{.1in}

Then, we obtain the claim of Proposition
\ref{prop;18.8.27.51}.
\hfill\qed

\subsection{Proof of Theorem \ref{thm;18.12.17.110}
and Theorem \ref{thm;19.1.6.10}}
\label{subsection;18.12.17.111}

We obtain that
$\nbigp_{\ast}E$ is a good filtered bundle
from Proposition \ref{prop;18.8.14.1},
Theorem \ref{thm;18.8.15.1}
and Proposition \ref{prop;18.8.27.51}.
Let $(E_{\omega},h_{\omega},\nabla_{\omega},\phi_{\omega})$
be as in (\ref{eq;19.1.6.11}).
Let $(V_{\omega},h_{V_{\omega}},\nabla_{V_{\omega}},\phi_{i,\omega})$
be as in Proposition \ref{prop;19.1.5.1}.
We set
$E_{\omega}^{\circledcirc}:=\Psi^{-1}(V_{\omega})$.
Let $h_{\omega}^{\circledcirc}:=\Psi^{-1}(h_{V_{\omega}})$
be the induced metric.
We set
$\nabla^{\circledcirc}_{\omega}:=
 \Psi^{\ast}(\nabla_{V_{\omega}})
+\sum_{i=1,2} \phi_{i,\omega}\,dy_i$.
We set $\phi^{\circledcirc}_{\omega}:=\Psi^{-1}(\phi_{3,\omega})$.
Note that 
$(E^{\circledcirc}_{\omega},h^{\circledcirc}_{\omega},
 \nabla^{\circledcirc}_{\omega},\phi^{\circledcirc}_{\omega})$
satisfies Condition \ref{condition;19.1.6.20}
and Condition \ref{condition;18.11.24.50}.
Moreover,
by modifying as in the proof of Lemma \ref{lem;18.12.16.10},
we may assume that
$\bigl(
 E^{\circledcirc}_{\omega},
 \del_{E^{\circledcirc}_{\omega},\ttubar},
 \del_{E^{\circledcirc}_{\omega},\ttt}
 \bigr)$
is a mini-holomorphic bundle
on $\nbigu^{\lambda}_{\nu,p}$.

\begin{lem}
\label{lem;19.1.7.20}
\mbox{{}}
\begin{itemize}
\item
The norm estimate holds for
$(\nbigp_{\ast}\bigl(E^{\circledcirc}
 _{\omega|\nbigu^{\lambda}_{\nu,p}(\ttt)}
\bigr),
 h^{\circledcirc}_{\omega})$.
\item
The isomorphism $F:E_{\omega}\simeq E^{\circledcirc}_{\omega}$ 
in Proposition 
{\rm\ref{prop;19.1.5.1}}
satisfies Condition {\rm\ref{condition;19.2.8.1}}.
\item
In particular,
there exists the isomorphism of local systems
$\Gr^{\nbigp}_a(E_{\omega})
\simeq
\Gr^{\nbigp}_a(E^{\circledcirc}_{\omega})$.
Moreover,
the norm estimate holds for 
$(\nbigp_{\ast}\bigl(E
 _{\omega|\nbigu^{\lambda}_{\nu,p}(\ttt)}
\bigr),
 h_{\omega})$.
\end{itemize}
\end{lem}
\pf
We can check the first claim by using
Lemma \ref{lem;19.1.7.10}
and Proposition \ref{prop;18.12.16.20}.
The second claim is clear.
The third claim follows from
Lemma \ref{lem;19.1.6.21}.
\hfill\qed

\vspace{.1in}

Take $\ttt_0\in S^1_{\lambda}$.
Let $I(\ttt_0)$ denote a small neighbourhood of $\ttt_0$.
We take a mini-holomorphic local frame
$\vecvhat$ of 
$\vecP^{(\ttt_0)}_a
 \sfp_p^{-1}E^{\lambda}_{|(\pihat^{\lambda}_{\nu,p})^{-1}(I(\ttt_0))}$
which is compatible with 
the slope decomposition
and the filtrations $\vecP^{(\ttt_0)}_{\ast}$ and $W$.
We may regard $\vecvhat$ as a mini-holomorphic frame of
$\Bigl(
\bigoplus
 \vecsfL_p(\omega)\otimes E_{\omega}
 \Bigr)_{|(\pihat^{\lambda}_{\nu,p})^{-1}(I(\ttt_0))}$,
which is compatible with 
the slope decomposition
and the filtrations $\vecP^{(\ttt_0)}_{\ast}$
and $W$.
There exists a $C^{\infty}$-frame 
$\vecv'$ of
$\vecP^{(\ttt_0)}_a
 \Bigl(
\bigoplus
 \vecsfL_p(\omega)\otimes E_{\omega}\Bigr)$
such that 
$\vecv'_{|(\pihat^{\lambda}_{\nu,p})^{-1}(I(\ttt_0))}
=\vecvhat$.
We may assume that $\vecv'$ is compatible
with the direct sum
$\bigoplus
 \vecsfL_p(\omega)\otimes E_{\omega}$,
i.e.,
$\vecv'=\bigcup_{\omega}\vecv'_{\omega}$,
where $\vecv'_{\omega}$ is a frame of
$\vecP^{(\ttt_0)}_a\bigl(
 \vecsfL_p(\omega)\otimes E_{\omega}\bigr)$.
Let $\vecv''_{\omega}$
be the frame of $E_{\omega}$
determined by
$\vecv'_{\omega}=
 \sfv^{\lambda}_{p,\omega}\otimes \vecv''_{\omega}$.
For each $v''_{\omega,i}$,
we have 
$k(\omega,i):=\deg^W(v''_{\omega,i})$
and 
$b(\omega,i):=\deg^{\vecP^{(\ttt_0)}}(v''_{\omega,i})$.
Let $h''_{0,\omega}$ be the metric
determined by
$h''_{0,\omega}(v''_{\omega,i},v''_{\omega,j})=0$ $(i\neq j)$
and 
$h''_{0,\omega}(v''_{\omega,i},v''_{\omega,i})=
 \bigl|\ttU_{p,\nu}\bigr|^{-2b(\omega,i)}
 \bigl|\log|\ttU_{p,\nu}|\bigr|^{k(\omega,i)}$.
Then, by Lemma \ref{lem;19.1.7.20},
we obtain that 
$h_{\omega}$ and $h''_{0,\omega}$ are mutually bounded.
There exists a mini-holomorphic local frame $\vecv$
of 
$\vecP^{(\ttt_0)}_a\sfp_p^{-1}(E^{\lambda})$
such that
$\vecv-\vecvhat=O(\ttU_{p,\nu}^N)$
for a sufficiently large $N$.
Then, by comparison of $\vecv$ and $\vecv'$,
we easily obtain that
the norm estimate holds for 
$(\nbigp_{\ast}E^{\lambda},h)$.

\vspace{.1in}
Let us prove Theorem \ref{thm;19.1.6.10}.
We have the induced $C^{\infty}$-isomorphism:
\begin{equation}
\label{eq;19.1.6.22}
 \sfp_{p}^{-1}\bigl(
 E^{\bullet}_{\omega}
 \bigr)
\simeq
 \vecsfL_p(\omega)\otimes
 E^{\circledcirc}_{\omega}.
\end{equation}
Note that 
$\vecsfL_p(\omega)\otimes
 E^{\circledcirc}_{\omega}$
is naturally equivariant with respect to
the action of 
$(\seisuu/p\seisuu)\cdot\tte_1$.
We obtain a tuple
$(E^{\bullet\prime}_{\omega},
 h^{\bullet\prime}_{\omega},
 \nabla^{\bullet\prime}_{\omega},
 \phi^{\bullet\prime}_{\omega})$
on $\nbigu_{1,\nu}^{\lambda}$
as the descent of 
$\vecsfL_p(\omega)\otimes
 \bigl(
 E^{\circledcirc}_{\omega},
 h^{\circledcirc}_{\omega},
 \nabla^{\circledcirc}_{\omega},
 \phi^{\circledcirc}_{\omega}
 \bigr)$.
Note that 
$E^{\bullet\prime\lambda}_{\omega}
=\bigl(
E^{\bullet\prime}_{\omega},
\del_{E^{\bullet\prime}_{\omega},\ttubar},
\del_{E^{\bullet\prime}_{\omega},\ttt}
\bigr)$
is a mini-holomorphic bundle.
We may assume that 
the isomorphism (\ref{eq;19.1.6.22})
is equivariant with respect to the action of
$(\seisuu/p\seisuu)\cdot \tte_1$.
We obtain the induced $C^{\infty}$-isomorphism
\begin{equation}
\label{eq;19.1.6.30}
 E_{\omega}^{\bullet}
\simeq
 E_{\omega}^{\bullet\prime}.
\end{equation}
We set
$E^{\prime}:=\bigoplus E^{\bullet\prime}_{\omega}$.
It is equipped with the induced metric $h^{\prime}$,
and the induced mini-holomorphic structure.
We obtain the mini-holomorphic bundle
$E^{\prime\lambda}$.
By Lemma \ref{lem;19.1.7.20},
(\ref{eq;19.1.6.30})
induces an isomorphism
$\ttG\bigl(
 \nbigp_{\ast}E^{\lambda}
 \bigr)
\simeq
 \ttG\bigl(
 \nbigp_{\ast}E^{\prime\lambda}
 \bigr)$.
Thus, we obtain 
the claim of Theorem \ref{thm;19.1.6.10}
from Proposition \ref{prop;19.2.8.2}.
(See also \S\ref{subsection;19.1.7.21}.)
\hfill\qed

\subsection{Initial metrics}

Let $\nbigubar^{\lambda}_{\nu,p}$ be a neighbourhood of
$H^{\lambda}_{\nu,p}$ in $\nbigmbar^{\lambda}_{p}$.
Let $\nbigp_{\ast}\gbigv$ be a good filtered bundle
on $(\nbigubar_{\nu,p},H^{\lambda}_{\nu,p})$.
Set $\nbigu^{\lambda}_{\nu,p}:=
 \nbigubar^{\lambda}_{\nu,p}\setminus H^{\lambda}_{\nu,p}$.
Let $V$ be the mini-holomorphic bundle on
$\nbigu^{\lambda}_{\nu,p}$
obtained as the restriction of $\nbigp\gbigv$.

\begin{prop}
\label{prop;18.9.2.3}
There exists a Hermitian metric $h_0$ of $V$
with the following property.
\begin{itemize}
\item
The norm estimate holds for
$(\nbigp_{\ast}\gbigv,h_0)$.
\item
$G(h_0)$
and its derivatives are
$O(e^{-\epsilon |y_0|})$.
\item
$F(h_0)$ is bounded.
\item
$\bigl[
 \del_{V,h_0,\ttu},
 \del_{V,\ttubar}
 \bigr]=O(y_0^{-2})$.
\end{itemize}
\end{prop}

\subsubsection{Approximation of regular filtered bundles}

Let $\nbigp_{\ast}\gbigv$
be a regular filtered bundle over
$(\nbigu^{\lambda}_{\nu,p},H^{\lambda}_{\nu,p})$.
We have the monodromy $F_{a}$ of
the local system $\Gr^{\nbigp}_a(\gbigv)$.
For each $a\in\Par(\nbigp_0\gbigv)$,
by using the results in 
\S\ref{subsection;18.12.21.22}--\ref{subsection;19.1.7.21},
we can construct a monopole
$(V_{0,a},\delbar_{V_{0,a}},h_{0,a})$
with the following property:
\begin{itemize}
\item
$F(h_{0,a})=O(y_0^{-2})$.
The associated Higgs field $\phi_{0,a}$ is bounded.
\item
$\Gr^{\nbigp}_b(V_{0,a})=0$
unless $b-a\in\seisuu$.
\item
We have an isomorphism of local systems
$\Gr^{\nbigp}_{a}(V_{0,a})
\simeq
 \Gr^{\nbigp}_a(V)$.
\end{itemize}
We set
$V_0:=\bigoplus V_{0,a}$.
We obtain the metric $h_0=\bigoplus h_{0,a}$.

\begin{lem}
\label{lem;18.9.2.11}
We have a $C^{\infty}$-isomorphism
$g:\nbigp_0V_0\simeq \nbigp_0V$
with the following property.
\begin{itemize}
\item
The induced isomorphism
$\nbigp_0V_{0|H^{\lambda}_{\nu,p}}
 \simeq
 \nbigp_0V_{|H^{\lambda}_{\nu,p}}$
preserves the parabolic filtrations.
\item
The induced morphism
$\Gr^{\nbigp}_aV
 \simeq
 \Gr^{\nbigp}_aV_0$
is an isomorphism of local systems.
\item
$g_{|\pi^{-1}(\ttt)}$ are holomorphic.
\item
Let $B$ be determined by
$B\,dt=\delbar_{V}-g^{\ast}\delbar_{V_0}$.
Then,
$B$ and its derivatives are $O(e^{-\epsilon|y_0|})$
with respect to $g^{\ast}(h_0)$.
\end{itemize}
\end{lem}
\pf
For each $-1<a\leq 0$,
we have the decomposition
$\Gr^{\nbigp}_a(V)=
 \bigoplus_{\alpha\in\cnum^{\ast}}
 \EE_{\alpha}\Gr^{\nbigp}_a(V)$
obtained as the generalized eigen decomposition
of the monodromy.
For each $\alpha\in\cnum^{\ast}$,
we take $\log\alpha\in\cnum$.
We take a $C^{\infty}$-frame $\vecu_{a,\alpha}$
of $\EE_{\alpha}\Gr^{\nbigp}_a(V)$
such that
$\del_t\vecu_{a,\alpha}=\vecu_a\cdot A_{a,\alpha}$,
where $A_{a,\alpha}$ is a constant matrix
with eigenvalues $(\gminit^{\lambda})^{-1}\log\alpha$.
We obtain a frame $\vecu_a$ of 
$\Gr^{\nbigp}_a(V)$.
We obtain a frame a frame $\vecu$ of
$\bigoplus_{-1<a\leq 0}\Gr^{\nbigp}_a(V)$.

Take $\ttt_0$.
Let $I(\ttt_0,\epsilon)$ be a small neighbourhood of $\ttt_0$
in $S^1_{\lambda}$.
We take a $C^{\infty}$-frame $\vecv^{(\ttt_0)}$
of $\nbigp_0\gbigv$ on $\pi^{-1}(I(\ttt_0,\epsilon))$
with the following property.
\begin{itemize}
\item
$\vecv^{(\ttt_0)}$ induces $\vecu_{|I(\ttt_0,\epsilon)}$.
\item
$\vecv^{(\ttt_0)}_{|\pi^{-1}(\ttt)}$ are holomorphic.
\end{itemize}
By using the partition of unity on $S^1$,
we can construct a $C^{\infty}$-frame $\vecv$
of $\nbigp_0\gbigv$ with the following property.
\begin{itemize}
\item
$\vecv$ induces $\vecu$.
\item
$\vecv_{|\pi^{-1}(\ttt)}$ are holomorphic.
\end{itemize}
We have $\nbigb$ determined by
$\delbar_{V}\vecv=\vecv\cdot\nbigb\,d\ttt$.
We have
$\del_{\ubar}\nbigb=0$.
Let $\nbigb^0$ be determined by
$\nbigb^{0}_{i,j}=\nbigb_{i,j|H^{\lambda}_{\nu,p}}$
if $\deg^{\nbigp}(v_i)=\deg^{\nbigp}(v_j)$,
and $\nbigb^0_{i,j}=0$
if $\deg^{\nbigp}(v_i)\neq \deg^{\nbigp}(v_j)$.
Then, the matrix $\nbigb^0$ 
represents the monodromy of
$\bigoplus_a\Gr^{\nbigp}_a(\gbigv)$
with the frame $\vecu$.

We take a $C^{\infty}$-frame $\vecv_0$ of 
$\nbigp_0V_0$ with similar properties.
We define 
$g:\nbigp_0\gbigv_0\lrarr\nbigp_0\gbigv_0$
by $g(\vecv_0)=\vecv$.
It has the desired property.
\hfill\qed

\subsubsection{Approximation of good filtered bundles}

Let $\nbigp_{\ast}\gbigv$ be a good filtered bundle.
Suppose that we are given 
$\nbigp_{\ast}(\gbigv_{0,\omega})$
for $\omega\in\nbigs(\gbigv)$,
and an isomorphism
\begin{equation}
\label{eq;18.9.2.1}
 \nbigp_{\ast}\gbigv_{|\Hhat^{\lambda}_{\nu,p}}
\simeq
\bigoplus_{\omega\in\nbigs(\gbigv)}
  \nbigp_{\ast}\gbigv_{0,\omega|\Hhat^{\lambda}_{\nu,p}}.
\end{equation}
We set
$\nbigp_{\ast}\gbigv_0:=\bigoplus\nbigp_{\ast}\gbigv_{0,\omega}$.
The following is easy to see.
\begin{lem}
\label{lem;18.9.2.10}
We have a $C^{\infty}$-isomorphism
$\nbigp_0\gbigv\simeq\nbigp_0\gbigv_0$
which induces {\rm(\ref{eq;18.9.2.1})}.
\hfill\qed
\end{lem}

Set $V_0:=
 \gbigv_{0|\nbigu^{\lambda}_{\nu,p}\setminus H^{\lambda}_{\nu,p}}$.
Let $h_0$ be a Hermitian metric of $V_0$
which has the properties in Proposition \ref{prop;18.9.2.3}
for $\nbigp_{\ast}\gbigv_0$.
By the isomorphism in Lemma \ref{lem;18.9.2.10},
we may regard $h_0$ as a Hermitian metric
of $V:=
 \gbigv_{|\nbigu^{\lambda}_{\nu,p}
 \setminus H^{\lambda}_{\nu,p}}$.
Then, $h_0$  also has the properties 
in Proposition \ref{prop;18.9.2.3}
for $\nbigp_{\ast}\gbigv$.

\subsubsection{Proof of Proposition \ref{prop;18.9.2.3}}

We can construct the desired metric
by using the approximations in
Lemma \ref{lem;18.9.2.11}
and Lemma \ref{lem;18.9.2.10},
and monopoles 
as in \S\ref{subsection;19.1.5.10}.
\hfill\qed

\subsection{Boundedness of curvature and adaptedness}

Let $\nbigp_{\ast}\gbigv$ be a good filtered bundle
on $(\nbigu^{\lambda}_{\nu,p},H^{\lambda}_{\nu,p})$.
Let $V$ be the mini-holomorphic bundle 
on $\nbigu^{\lambda}_{\nu,p}\setminus H^{\lambda}_{\nu,p}$
obtained as the restriction of $\gbigv$.
Let $h$ be a Hermitian metric of $V$
with the following property.
\begin{itemize}
\item
 $G(h)=0$.
\item
 $h$ is adapted to $\nbigp_{\ast}\gbigv$.
\end{itemize}

\begin{prop}
\label{prop;19.2.8.100}
$F(h)$ is bounded.
Moreover,
the norm estimate holds for
$(\nbigp_{\ast}\gbigv,h)$.
\end{prop}
\pf
We have a Hermitian metric $h_0$
as in Proposition \ref{prop;18.9.2.3}.
Let $s$ be the automorphism of $V$
determined by $h=h_0s$.
We obtain
\[
 \Delta\log\Tr(s)
\leq |G(h_0)|_{h_0}
\leq
 Ce^{-\epsilon |y_0|}.
\]
We obtain
\[
 \Delta\bigl(
 \log\Tr(s)-C_1e^{-\epsilon_1 |y_0|}
 \bigr)
\leq 0
\]
for some $C_1>0$ and $\epsilon_1>0$.
By the assumption, $\log\Tr(s)=O(\log |y_0|)$ holds.
We take $C_2>0$ 
such that 
$\log\Tr(s)< C_2$
on $\{y_0=R\}$.
Note that
$\Delta(\delta |y_0|)=0$ for any $\delta>0$.
We obtain
\[
 \Delta\bigl(
 \log\Tr(s)-C_1e^{-\epsilon_1|y_0|}
-\delta |y_0|-C_2
 \bigr)
\leq 0.
\]
Then, by a standard argument,
we obtain that 
$\log\Tr(s)-C_1e^{-\epsilon_1|y_0|}
\leq C_2+\delta |y_0|$
for any $\delta>0$.
Hence, 
$\log\Tr(s)-C_1e^{-\epsilon_1|y_0|}
\leq C_2$.
Thus, we obtain the boundedness of $s$.
Similarly, we obtain the boundedness of $s^{-1}$.
It implies that the norm estimate for
$(\nbigp_{\ast}\gbigv,h)$.

\begin{lem}
\label{lem;19.2.8.10}
$\int \bigl|
 \del_{E,h_0,\alphabar}s
 \bigr|^2
+\int\bigl|
 \del_{E,h_0,\tau}s
 \bigr|^2
<\infty$
and
$\int \bigl|
 \del_{E,h_0,\alpha}s
 \bigr|^2
+\int\bigl|
 \del'_{E,h_0,\tau}s
 \bigr|^2
<\infty$
hold.
\end{lem}
\pf
The following holds:
\[
 -\Bigl(
 \del_{\alpha}\del_{\alphabar}
+\frac{1}{4}
 \del_{\tau}\del_{\tau}
 \Bigr)
 \Tr(s)
=-\Tr\bigl( 
 sG(h_0)
\bigr)
-\bigl|
 s^{-1/2}\del_{E,h_0,\alpha}s
 \bigr|^2
-\frac{1}{4}\bigl|
 s^{-1/2}\del'_{E,h_0,\tau}s
 \bigr|^2.
\]
We set
\[
 b_1:=\int_{T^2}
 \Tr(s),
\quad
 b_2:=\int_{T^2}
 \Tr\bigl( 
 sG(h_0)
\bigr),
\quad
 b_3:=\int_{T^2}
\bigl|
 s^{-1/2}\del_{E,h_0,\alpha}s
 \bigr|^2
+\int_{T^2}
 \frac{1}{4}\bigl|
 s^{-1/2}\del'_{E,h_0,\tau}s
 \bigr|^2.
\]
Note that
$\del_{\alpha}\del_{\alphabar}
+\frac{1}{4}\del_{\tau}\del_{\tau}
=\frac{1}{4}
 \bigl(
 \del_{y_0}^2+\del_{y_1}^2+\del_{y_2}^2
 \bigr)$.
We obtain
\[
 -\del_{y_0}^2b_1
=-4b_2-4b_3.
\]
Note that
$|b_2|=O(e^{-\epsilon |y_0|})$.
Hence, there exists $c_2$ such that
$|c_2|=O(e^{-\epsilon |y_0|})$
and 
$-\del_{y_0}^2(b_1-c_2)=-4b_3$.
Note that $b_3\geq 0$.
Because $b_3-c_2$ is bounded and subharmonic,
we obtain that there exists
$\lim_{y_0\to\infty}(\del_{y_0}(b_3-c_2))$.
Then, we obtain the existence of
$\lim_{R\to\infty}\int_{C}^Rb_3$.
It implies the claim of the lemma.
\hfill\qed

\vspace{.1in}

As in \cite{Mochizuki-KH-infinite},
there exists $C>0$ such that
\[
 \Delta
 \bigl|
 s^{-1}\del_{E,h_0}s
 \bigr|^2_{h_0}
\leq
 C\Bigl(
 1+\bigl|s^{-1}\del_{E,h_0}s\bigr|_{h_0}^2
 \Bigr).
\]
By using \cite[Theorem 9.20]{Gilbarg-Trudinger}
and Lemma \ref{lem;19.2.8.10},
we obtain the boundedness of
$s^{-1}\del_{E,h_0}s$.
By using the equation for the monopole,
we also obtain that
$s$ and its derivatives are bounded.
\hfill\qed

\section{Rank one monopoles}

\subsection{Preliminary}

\subsubsection{Ahlfors type lemma}

Let $R>0$.
Let $g$ be a $C^{\infty}$-function
$\{t\geq R\}\lrarr\real_{\geq 0}$
such that 
$g=O(t^N)$ for some $N>0$.
Suppose
$-\del_t^2g\leq -C_0g+C_1e^{-a t}$
for some $C_i>0$ and $a>0$.

\begin{lem}
\label{lem;19.2.8.20}
We obtain
$g=O\bigl(\exp(-\epsilon t)\bigr)$ for some $\epsilon>0$.
\end{lem}
\pf
By making $C_0$ smaller,
we may assume that $C_0<a^2$.
We set
$C_2:=C_1(a^2-C_0)^{-1}$.
The following holds:
\[
 -\del_t^2C_2e^{-at}
=-\bigl(a^2C_2-C_1\bigr)e^{-at}-C_1e^{-at}
=-C_2e^{-at}-C_1e^{-at}.
\]
We obtain
\[
 -\del_t^2\bigl(
 g+C_2e^{-at}
 \bigr)
\leq
 -C_0(g+C_2e^{-at}).
\]

For $C_3>0$ and $\delta>0$,
we set
$F_{C_3,\delta}(t):=C_3\exp(-\epsilon t)+\delta\exp(\epsilon t)$.
There exists $C_3>0$ 
such that
$F_{C_3,\delta}(R)>(g+C_2e^{-at})_{|t=R}$ for any $\delta>0$.
Then, the set
$\{t\,|\,F_{C_3,\delta}(t)<g(t)\}$ is relatively compact
in $\{t>R\}$.
Set $\epsilon:=C_0^{1/2}$.
By using 
$-\del_t^2F_{C_3,\delta}=-C_0F_{C_3,\delta}$
with a standard argument,
we obtain that
$F_{C_3,\delta}>g+C_2e^{-at}$ on $\{t\geq R\}$ 
for any $\delta>0$.
By taking the limit $\delta\to 0$,
we obtain the desired estimate.
\hfill\qed

\subsubsection{Global subharmonic functions on $X\times \real$}

Let $(X,g_X)$ be a compact Riemannian manifold.
The Riemannian metric $g_X+dt\,dt$ on $X\times\real$
is induced.
\begin{lem}
Let $f$ be a bounded function
$X\times\real\lrarr\real_{\geq 0}$
such that
$\Delta f\leq 0$.
Then, $f$ is constant.
In particular, $\Delta f=0$.
\end{lem}
\pf
We obtain the decomposition $f=f_0+f_1$,
where $f_0$ is constant on $X\times\{t\}$,
and $\int_{X\times\{t\}}f_1=0$ holds for any $t$.
We obtain
$-\del_t^2f_0\leq 0$.
Because $f_0$ is bounded,
we obtain that $f_0$ is constant.
Let $d_X$ denote the exterior derivative
in the $X$-direction.
We obtain 
\[
 \Delta|f|^2\leq -\bigl|d_Xf\bigr|^2
=-\bigl|d_Xf_1\bigr|^2.
\]
We obtain
\[
 -\del_t^2\int_{X\times\{t\}}|f_1|^2=
 -\del_t^2\int_{X\times\{t\}}|f|^2\leq
 -\int_{X\times\{t\}}|d_Xf_1|^2
\leq
 -C_1\int_{X\times\{t\}}|f_1|^2.
\]
By Lemma \ref{lem;19.2.8.20},
we obtain
$\int_{X\times\{t\}}|f_1|^2
=O\bigl(\exp(-\epsilon|t|)\bigr)$
for some $\epsilon>0$.
Because
$\int_{X\times\{t\}}|f_1|^2\geq 0$ is subharmonic,
we obtain 
$\int_{X\times\{t\}}|f_1|^2$ is constantly $0$.
It implies $f_1=0$.
\hfill\qed

\subsubsection{Poisson equation on $X\times\real$}

Let $a$ be a $C^{\infty}$-function on $X\times\real$
such that
$a=O(\exp(-\epsilon|t|))$,
and that $\int_{X\times\real}a=0$.
For any $t\in\real$,
we set $X_t:=X\times\{t\}$.

\begin{lem}
\label{lem;18.11.20.1}
There exists a $C^{\infty}$-function $b$ on $X\times\real$
such that
(i) $\Delta b=a$
(ii) $|b|=O\bigl(\exp(\epsilon_1 t)\bigr)$ as $t\to-\infty$,
(iii) there exists the limit $\lim_{t\to\infty}b=b_{\infty}$,
and $|b-b_{\infty}|=O\bigl(\exp(-\epsilon_1 t)\bigr)$ 
as $t\to\infty$.
\end{lem}
\pf
Let $a=a_0+a_1$ be the decomposition such that
(i) $a_0$ is constant on $X_t$ for any $t$,
(ii) $\int_{X_t}a_1=0$ for any $t$.
We may regard $a_0$ as a $C^{\infty}$-function on $\real$
such that $a_0=O(\exp(-\epsilon|t|))$.
It is easy to see that there exists a function
$b_0$ on $\real$ such that 
(i) $-\del_t^2b_0=a_0$,
(ii) $b_0=O(\exp(\epsilon_1t))$ for some $\epsilon_1>0$
 as $t\to-\infty$,
(iii)  there exists $b_{\infty}:=\lim_{t\to\infty}b_0(t)$,
and  $b_0-b_{\infty}=O\bigl(\exp(-\epsilon_2t)\bigr)$
 for some $\epsilon_2>0$ as $t\to\infty$.

There exists a complete orthonormal set 
$\{\varphi\}$ in $C^{\infty}(X)$
such that $\Delta_X\varphi=\lambda(\varphi)\varphi$,
where $\lambda(\varphi)\in\real_{\geq 0}$.
Let $a_1=\sum_{\lambda(\varphi)>0}a_{1,\varphi}(t)\varphi$
be the expansion.
We set
\[
 b_{1,\varphi}(t):=
 e^{-\lambda(\varphi)^{1/2} t}\int_{-\infty}^t
 e^{2\lambda(\varphi)^{1/2} s}\,ds
\int_{s}^{\infty}e^{-\lambda(\varphi)^{1/2} u}a_{1,\varphi}(u)\,du.
\]
Then,
$(-\del_t^2+\lambda(\varphi))b_{1,\varphi}=a_{1,\varphi}$
holds.
Set $\|a_{1,\varphi}\|_{L^2}:=
\Bigl(
 \int_{\real}|a_{1,\varphi}(t)|^2\,dt
\Bigr)^{1/2}$.
We obtain
$|b_{1,\varphi}(t)|\leq 
 C\bigl\|a_{1,\varphi}\bigr\|_{L^2}$ for some $C>0$.
Because $\sum_{\varphi}\|a_{1,\varphi}\|_{L^2}^2<\infty$,
we obtain the locally $L^2$-function
$b_1:=\sum b_{1,\varphi}\varphi$ on $X\times\real$,
and 
$\Delta b_1=a_1$ holds 
in the sense of distributions.
By the elliptic regularity,
$b_1$ is $C^{\infty}$.
Set $f(t):=\int_{X_t}|b_1|^2$,
and then 
$|f(t)|\leq C\sum\|a_{1,\varphi}\|_{L^2}$.
The following holds:
\[
 \int_{X_t}a_1\overline{b_1}
=\int_{X_t}
 \bigl(-\del_t^2+\Delta_X
 \bigr)b_1\cdot \overline{b_1}
=-\del_t^2f+\int_{X_t}\|d_Xb_1\|^2.
\]
There exists $\epsilon_2>0$
such that
$\int_{X_t}\|d_Xb_1\|^2
\geq
 \epsilon_1 f$.
There exist $C_i$ $(i=1,2)$
and $\epsilon_i>0$ $(i=2,3,4)$ such that 
\[
 -\del_t^2f\leq 
 C_1e^{-\epsilon_2|t|}f^{1/2}
-\epsilon_1 f
\leq
 C_2e^{-\epsilon_3|t|}-\epsilon_4 f.
\]
on $\{|t|>R\}$ for some $R>0$.
Then, we obtain that
$|f|=O\bigl(\exp(-\epsilon_5|t|)\bigr)$
for some $\epsilon_5>0$.
Thus, we are done.
\hfill\qed

\subsection{Examples of monopoles of rank $1$ with Dirac type singularity}

\subsubsection{Filtered bundles of rank $1$}

Suppose that $\gminig^{\lambda}>0$.
Take a small $\epsilon>0$.
We set
$\nbigw:=\proj^1\times\openopen{-\epsilon}{1}$.
We have the open embedding
$\nbigw\lrarr\nbigmbar^{\lambda\cov}$
induced by 
$(\ttU,t)\longmapsto (\ttU,\gminit^{\lambda}t)$.
It induces the surjection
$\nbigw\lrarr \nbigmbar^{\lambda}$.
We have the isomorphism
$\Phi:\proj^1\times\openopen{-\epsilon}{0}
\simeq
\proj^1\times
 \openopen{1-\epsilon}{1}$
given by
$\Phi(\ttU,t)=(\gminiq^{\lambda}_p\ttU,t+1)$.
We regard $\nbigmbar^{\lambda}$
as the quotient space of $\nbigw$
by identifying 
$\proj^1\times\openopen{-\epsilon}{0}$
and 
$\proj^1\times\openopen{1-\epsilon}{1}$.

Let
$(\sfA_0,t_0)\in
 \cnum^{\ast}\times\closedopen{0}{1}$.
We set
$\nbigz_{\nbigw,(\sfA_0,t_0)}:=
 \{\sfA_0\}\times 
 \openopen{t_0}{1}$.
We set
\[
 \nbigv_{(\sfA_0,t_0),n}:=
 \nbigo_{\nbigw\setminus\{(\sfA_0,t_0)\}}
 \bigl(-n\nbigz_{\nbigw,(\sfA_0,t_0)}\bigr)
 \bigl(\ast(\{0,\infty\}\times\openopen{-\epsilon}{1})\bigr)
 \cdot v.
\]
Let $\pi:\nbigw\lrarr \openopen{-\epsilon}{1}$ 
denote the projection.
We define the filtered bundles
$\nbigp^{(a)}_{\ast}\bigl(\nbigv_{|\pi^{-1}(t)}\bigr)$
by the following conditions:
\begin{itemize}
\item
 The parabolic degree of 
 $v_{|\pi^{-1}(t)}$ at $\infty$ is constantly $0$.
\item
 The parabolic degree of
 $v_{|\pi^{-1}(t)}$ at $0$ is $a-nt$.
\end{itemize}

We define the isomorphism
$\Phi^{\ast}\bigl(
 \nbigv_{|\proj^1\times\openopen{1-\epsilon}{1}}
 \bigr)
\simeq
 \nbigv_{|\proj^1\times\openopen{-\epsilon}{0}}$
by
\[
 \Phi^{\ast}\bigl(
 (\ttU-\sfA_0)^n\ttU^{-n}v_{|\proj^1\times
 \openopen{1-\epsilon}{1}}
 \bigr)
=v_{|\proj^1\times\openopen{-\epsilon}{0}}.
\]
It induces an isomorphism of filtered bundles
for $t\in\openopen{-\epsilon}{0}$:
\[
 \Phi^{\ast}\bigl(
 \nbigp^{(a)}_{\ast}\nbigv_{|\pi^{-1}(t+1)}
 \bigr)
\simeq
 \nbigp^{(a)}_{\ast}\nbigv_{|\pi^{-1}(t)}.
\]
We set $\ttt_0:=\gminit^{\lambda}t_0$.
We obtain an induced
$\nbigo_{\nbigmbar^{\lambda}\setminus\{(\sfA_0,\ttt_0)\}}
 (\ast H^{\lambda}_p)$-module
$\nbigl(\sfA_0,\ttt_0,n)$,
and a filtered bundle
$\nbigp^{(a)}_{\ast}\nbigl(\sfA_0,\ttt_0,n)$
over $\nbigl(\sfA_0,\ttt_0,n)$.
We obtain the following lemma by a direct computation.
\begin{lem}
$\deg(\nbigp^{(a)}_{\ast}\nbigl(\sfA_0,\ttt_0,n))=
 \bigl|\gminit^{\lambda}\bigr|(-a-n/2+n\ttt_0/\gminit^{\lambda})$ holds.
In particular,
\[
 a(\ttt_0,n)=n(-1/2+\ttt_0/\gminit^{\lambda}),
\]
we obtain $\deg\nbigp^{(a(\ttt_0,n))}_{\ast}\nbigl(\sfA_0,\ttt_0,n)=0$.
\hfill\qed
\end{lem}

\subsubsection{Monopoles}

Set $\nbigu(\sfA_0,\ttt_0):=
 \nbigm^{\lambda}\setminus\{(\sfA_0,\ttt_0)\}$.
\begin{prop}
There exists a Hermitian metric $h$ of
$\nbigl(\sfA_0,\ttt_0,n)_{|\nbigu(\sfA_0,\ttt_0)}$
such that the following holds.
\begin{itemize}
\item
$\bigl(
\nbigl(\sfA_0,\ttt_0,n)_{|\nbigu(\sfA_0,\ttt_0)},
 h\bigr)$
is a monopole with Dirac type singularity
on $\nbigu(\sfA_0,\ttt_0)$.
\item
The norm estimate holds for 
$\nbigp^{(a(\ttt_0,n))}_{\ast}\nbigl(\sfA_0,\ttt_0,n)$
with $h$.
\end{itemize}
Such $h$ is unique up to the positive constant multiplications.
\end{prop}
\pf
Set $\nbigl:=\nbigl(\sfA_0,\ttt_0,n)$.
There exists a Hermitian metric $h_0$ of 
such that 
(i) $G(h_0)$ and its derivatives are $O(e^{-\epsilon|y_0|})$,
(ii) $(\nbigl,h_0)$ is a monopole with Dirac type singularity
on $U_{\sfA_0,\ttt_0}\setminus\{(\sfA_0,\ttt_0)\}$,
where $U_{\sfA_0,\ttt_0}$ denotes
a neighbourhood of $(\sfA_0,\ttt_0)$,
(iii) the norm estimate holds for $h_0$.
For another metric $h_0e^{\varphi}$,
we have
$G(h_0e^{\varphi})=G(h_0)+4^{-1}\Delta\varphi$.
Because $\deg(\nbigp^{(a(\ttt_0,n))}_{\ast}\nbigl)=0$,
we obtain $\int G(h_0)=0$,
and hence 
there exists a bounded $C^{\infty}$-function
$\varphi$
such that $G(h_0e^{\varphi})=0$
according to Lemma \ref{lem;18.11.20.1}.
The uniqueness is clear.
\hfill\qed

\subsection{Classification of rank one monopoles}

Let $Z_0=\{(\sfA_i,\ttt_i)\,|\,i=1,\ldots,m\}
 \subset \nbigm^{\lambda\cov}$ be a finite subset
such that $0\leq \ttt_i/\gminit^{\lambda}<1$.
Let $Z\subset\nbigm^{\lambda}$
be the induced subset.
For each $i$,
we set $a_i:=-1/2+\ttt_i/\gminit^{\lambda}$.
The following lemma is clear.
\begin{lem}
Let $\nbigp_{\ast}\nbigl$ be a good filtered bundle
with Dirac type singularity of degree $0$
on $(\nbigmbarlambda;Z\cup H^{\lambda})$.
Then, there exist
$\ell\in\seisuu$,
$(\sfa,\sfb)\in\real\times\cnum$
and an isomorphism
\[
 \nbigp_{\ast}\nbigl
\simeq
 \nbigp_{\ast}\bigl(
 \nbigl_1(\ell)
 \bigr)
 \otimes
 \nbigp_{\ast}\bigl(
 L_1(\lambda,\sfa,\sfb)
 \bigr)
\otimes
 \bigotimes_{i=1}^{m}
 \nbigp_{\ast}^{(a_i)}\nbigl(\sfA_i,\ttt_i,1).
\]
Here,
see {\rm\S\ref{subsection;19.2.8.30}}
for $\nbigp_{\ast}\nbigl_1(\ell)$,
and {\rm\S\ref{subsection;19.2.8.31}}
for $\nbigp_{\ast}L_1(\lambda,\sfa,\sfb)$.
\hfill\qed
\end{lem}

\begin{prop}
\label{prop;19.2.8.110}
There exists an equivalence between the following objects:
\begin{itemize}
\item
Monopoles of rank one $(E,h,\nabla,\phi)$ on 
$\nbigm^{\lambda}\setminus Z$
such that
(i) each point of $Z$ is Dirac type singularity,
(ii) $F(\nabla)$ is bounded.
\item
Filtered bundles with Dirac type singularity
of rank one with degree $0$
on $(\nbigmbar^{\lambda};H^{\lambda},Z)$.
\end{itemize}
The correspondence is induced by
$(E,h,\nabla,\phi)\longmapsto
 \nbigp^h_{\ast}E$.
\hfill\qed
\end{prop}

\section{Kobayashi-Hitchin correspondence for doubly-periodic monopoles}

\subsection{Main statement}

Let $Z$ be a finite subset of $\nbigm^{\lambda}$.
\begin{df}
A monopole $(E,h,\nabla,\phi)$  on $\nbigm^{\lambda}\setminus Z$
is called meromorphic if the following holds.
\begin{itemize}
\item
Any points of $Z$ are Dirac type singularity
of $(E,h,\nabla,\phi)$.
\item
There exists a compact subset $C$ of $\nbigm^{\lambda}$
such that (i) $Z\subset C$,
(ii) $F(\nabla)$ is bounded on $\nbigm^{\lambda}\setminus C$.
\hfill\qed
\end{itemize}
\end{df}

For any meromorphic monopole
$(E,h,\nabla,\phi)$,
we have the associated good filtered bundle with Dirac type singularity
$\nbigp_{\ast}E^{\lambda}$
on $(\nbigmbar^{\lambda};H^{\lambda},Z)$,
as explained in \S\ref{subsection;19.1.27.2}.
We shall prove the following theorem
in \S\ref{subsection;19.1.27.3}--\ref{subsection;19.1.27.4}.
\begin{thm}
\label{thm;18.11.21.20}
The above procedure induces the bijection of 
the isomorphism classes of the following objects:
\begin{itemize}
\item
Meromorphic monopoles
on $\nbigm^{\lambda}\setminus Z$.
\item
Polystable good filtered bundle with Dirac type singularity
of degree $0$
on $(\nbigmbar^{\lambda};H^{\lambda},Z)$.
\end{itemize}
\end{thm}

\subsection{Preliminary}
\label{subsection;19.1.27.3}

\subsubsection{Ambient good filtered bundles with appropriate metric}
\label{subsection;18.11.21.5}

Let $Z$ be a finite subset in $\nbigm^{\lambda}$.
Let $\nbigp_{\ast}\nbige^{\lambda}$
be a good filtered bundle with Dirac type singularity
on $(\nbigmbar^{\lambda};H^{\lambda},Z)$.
Let $(E,\delbar_E)$ denote the mini-holomorphic bundle
with Dirac type singularity on $\nbigmlambda\setminus Z$
obtained as the restriction of $\nbigp\nbige^{\lambda}$.

Let $h_1$ be a Hermitian metric
of $E$ adapted to $\nbigp_{\ast}\nbige$
such that the following holds.
\begin{description}
\item[(A1)]
 Around $H^{\lambda}$,
we have $G(h_1)=O\bigl(e^{-\epsilon|y_0|}\bigr)$
for some $\epsilon>0$,
and 
$(E,\delbar_E,h_1)$ satisfies the norm estimate
with respect to $\nbigp_{\ast}\nbige$.
Moreover,
we have
\begin{equation}
\label{eq;18.11.21.1}
 \bigl[
 \del_{E,\ttubar},\del_{E,h_1,\ttu}
 \bigr]=O\bigl(y_0^{-2}\bigr).
\end{equation}
\item[(A2)]
Around each point of $Z$,
$(E,\delbar_E,h_1)$ is a monopole with Dirac type singularity.
In particular,
it induces a $C^{\infty}$-metric 
of the Kronheimer resolution of $E$.
\end{description}

\subsubsection{Degree of filtered subbundles}

Let $\nbigp_{\ast}\nbige_1\subset
\nbigp_{\ast}\nbige$
be a filtered subbundle on 
$\bigl(
 \nbigmbar^{\lambda};H^{\lambda},Z
 \bigr)$.
Let $E_1$ be the mini-holomorphic bundle
with Dirac type singularity
on $(\nbigm^{\lambda},Z)$.
Let $h_{1,E_1}$ denote the metric of $E_1$
induced by $h_1$.
By the Chern-Weil formula,
the analytic degree
$\deg(E_1,h_{1,E_1})\in \real\cup\{-\infty\}$ makes sense.

\begin{prop}
\label{prop;19.2.8.50}
There exists $C>0$
such that 
$C\deg(\nbigp_{\ast}\nbige_1)
=\deg(E_1,h_{1,E_1})$
for any $\nbigp_{\ast}\nbige_1$.
\end{prop}
\pf
Because the argument is essentially the same
as the proof of \cite[Proposition 9.4]{Mochizuki-difference-modules},
we give only an outline.
We take a metric $h_{0,E_1}$ of $E_1$
which satisfies the conditions {\bf (A1,2)}
for $\nbigp_{\ast}\nbige_1$.
Because $G(h_{0,E_1})=O(e^{-\epsilon|y_0|})$ $(\epsilon>0)$ 
around $H^{\lambda}$,
and because $G(h_{0,E_1})=0$ around each point of $Z$,
$G(h_{0,E_1})$ is $L^1$.
Let $\nabla_{0}$ and $\phi_0$ be 
the Chern connection and the Higgs field
associated to $(E_1,\delbar_{E_1})$ with $h_{0,E_1}$.
Because $(E_1,\delbar_{E_1},h_{0,E_1})$ is a monopole
with Dirac type singularity around each point $P$ of $Z$,
we have $(\nabla_0\phi_0)_{|x}=O\bigl(d(x,P)^{-2}\bigr)$
around $P$, and hence $\nabla_0\phi_0$ is $L^1$ around $P$.
Let $\del_{E_1,\ttu}$ denote the operator
induced by $\del_{E_1,\ttubar}$ and $h_{0,E_1}$.
Because 
$[\del_{E_1,\ttu},\del_{E_1,\ttubar}]=
O\bigl(y_0^{-2}\bigr)$ around $H_{\infty}^{\lambda}$,
$[\del_{E_1,\ttu},\del_{E_1,\ttubar}]$
is $L^1$ around $H^{\lambda}$.
Hence, by Proposition \ref{prop;19.2.8.40},
we obtain the following equality:
\[
\int
 \Tr G(h_{0,E_1})\dvol
=C\int_{0}^{|\gminit^{\lambda}|}
 \pardeg(\nbigp_{\ast}\nbige_{1|\pi^{-1}(\ttt)})
 \,d\ttt
=C\deg\bigl(\nbigp_{\ast}\nbige_1\bigr).
\]

Let us prove the following equality:
\begin{equation}
 \int
 \Tr G(h_{1,E_1})\dvol
=\int
 \Tr G(h_{0,E_1})\dvol.
\end{equation}
By considering
$\det E_1\subset \bigwedge^{\rank E_1}E$,
it is enough to consider the case $\rank E_1=1$.
We have a Hermitian metric $h'_{E_1}$ of $E_1$
such that 
(i) $(E_1,\delbar_{E_1},h'_{E_1})$ is a meromorphic monopole,
(ii) the meromorphic extension
$\nbigp^{h'_{E_1}}E_1$ is equal to
$\nbigp\nbige_1$.
We have $\deg(\nbigp^{h'_{E_1}}_{\ast}E_1)=0$.
By considering
$(E,\delbar_E,h)\otimes(E_1,\delbar_{E_1},h'_{E_1})^{-1}$
and 
$\nbigp_{\ast}\nbige\otimes
 (\nbigp_{\ast}^{h'_{E_1}}E_1)^{\lor}$,
we may reduce the issue
to the case where
there exists an isomorphism
$\nbigp\nbige_1\simeq
\nbigo_{\nbigmbar^{\lambda}}(\ast H^{\lambda})$.
Let $g$ be a section of $\nbigp\nbige_1$
corresponding to
$1\in\nbigo_{\nbigmbar^{\lambda}}(\ast H^{\lambda})$
under the isomorphism.
We have the number $a_0$ 
such that
$g\in \nbigp_{a_0}\nbige_1$
and $f\not\in\nbigp_{<a_{0}}\nbige_1$
 around $H^{\lambda}_{0}$,
and the number $a_{\infty}$
such that
$g\in \nbigp_{a_{\infty}}\nbige_1$
and $g\not\in\nbigp_{<a_{\infty}}\nbige_1$
around $H^{\lambda}_{\infty}$.
By considering the metric
$h_1e^{-a_{0}y_0}$ on around $H^{\lambda}_{0}$
and 
$h_1e^{a_{\infty}y_0}$ on around $H^{\lambda}_{\infty}$,
it is enough to consider the case
$a_0=a_{\infty}=0$.

\begin{lem}
 Let $\nbigb^{\lambda}$ be a neighbourhood of
 $H^{\lambda}_{\infty}$ in $\nbigmbar^{\lambda}$.
 Let $E$ be a mini-holomorphic bundle
 on $\nbigb^{\lambda\ast}:=\nbigb^{\lambda}\setminus H^{\lambda}_{\infty}$
 with a metric $h$ such that 
 $G(h)$ is $L^1$.
Let $f$ be a mini-holomorphic section of $E$
such that 
\[
 C_1^{-1}
\leq
 |f|_hy_0^{-k}
\leq
 C_1
\]
for some $C_1>1$ and $k\in\real$.
Then,
$\bigl|\nabla_{\alpha}f\bigr|_h\cdot
 |f|_h^{-1}$
and
$\bigl|
 (\nabla_{\tau}+\sqrt{-1}\phi)f
 \bigr|_h\cdot
 |f|_h^{-1}$
 are $L^2$.

Similar claim holds on a neighbourhood of $H^{\lambda}_0$.
\end{lem}
\pf
It is enough to prove that 
$\bigl|\nabla_{\alpha}f\bigr|_hy^{-k}$
and 
$\bigl|(\nabla_{\tau}+\sqrt{-1}\phi)f
 \bigr|_hy_0^{-k}$
are $L^2$.
Because $f$ is mini-holomorphic,
we have
$\nabla_{\alphabar}f=0$
and $(\nabla_{\tau}-\sqrt{-1}\phi)f=0$.
We may assume that 
$\nbigb^{\lambda\ast}=\{y_0>R\}$.

We take a $C^{\infty}$-function
$\rho:\real\lrarr\{0\leq a\leq 1\}\subset\real_{\geq 0}$
such that,
(i) $\rho(t)=0$ $(t\geq 1)$,
(ii) $\rho(t)=1$ $(t\leq 1/2)$,
(iii) $\rho(t)^{1/2}$
and $\del_t\rho(t)\big/\rho(t)^{1/2}$
give $C^{\infty}$-functions.

For any large positive integer $N$,
we set
$\chi_N(y_0):=\rho\bigl(N^{-1}y_0\bigr)$.
We obtain $C^{\infty}$-functions
$\chi_N:\nbigb^{\lambda\ast}
\lrarr \real_{\geq 0}$
such that
$\chi_N(y_0)=0$ if $y_0>N$
and 
$\chi_N(y_0)=1$ if $y_0<N/2$.
Let $\mu:\nbigb^{\lambda\ast}\lrarr \real_{\geq 0}$
be a $C^{\infty}$-function such that
$\mu(y_0)=1-\rho\bigl(y_0-R\bigr)$.
We set $\chitilde_N:=\mu\cdot \chi_N$.
We have 
\[
 \del_{y_0}\chitilde_N(y_0)
=\del_{y_0}\mu(y_0)\chi_N(y_0)
+\mu(y_0)\rho'(N^{-1}y_0)N^{-1}.
\]
By the assumption on $\rho$,
$\del_{y_0}\chitilde_N(y_0)\big/\chitilde_N(y_0)^{1/2}$
naturally give $C^{\infty}$-functions
on $\nbigb^{\lambda\ast}$,
and there exists $C_2>0$, which is independent of $N$,
such that the following holds:
\[
 \bigl|
 \del_{y_0}\chitilde_N(y_0)\big/\chitilde_N(y_0)^{1/2}
 \bigr|
\leq
 C_2y_0^{-1}.
\]
Because $\del_{\alpha}y_0$ is constant,
we have $C_3>0$, which is independent of $N$,
such that the following holds:
\[
 \bigl|
 \del_{\alpha}\bigl(
 \chitilde_N(y_0)\bigr)\big/\chitilde_N(y_0)^{1/2}
 \bigr|
\leq
 C_3y_0^{-1}.
\]

We consider the following integral:
\begin{multline}
\int_{\nbigb^{\lambda\ast}}
 \chitilde_N(y_0)\cdot h(\nabla_{\alpha}f,\nabla_{\alpha}f)
 y_0^{-2k}\,\dvol
=
-\int_{\nbigb^{\lambda\ast}}
 \del_{\alpha}\bigl(\chitilde_N(y_0)\bigr)\cdot
 h(f,\nabla_{\alpha}f)
 y_0^{-2k}\,\dvol
 \\
-\int_{\nbigb^{\lambda\ast}}
 \chitilde_N(y_0)\cdot
 h(f,\nabla_{\alphabar}\nabla_{\alpha}f)
 y_0^{-2k}\,\dvol
+\int_{\nbigb^{\lambda\ast}}
 \chitilde_N(y_0)\cdot
 h(f,\nabla_{\alpha}f)
\cdot(-2k)y_0^{-2k-1}
 \del_{\alpha}y_0
 \,\dvol
\end{multline}
We have the following inequality:
\[
 \Bigl|
 \del_{\alpha}\chitilde_N\cdot
 h(f,\nabla_{\alpha}f)
 y_0^{-2k}
 \Bigr|
\leq
 C_3C_1y_0^{-1}
 \cdot
 \Bigl(
 \chitilde_N^{1/2}(y_0)\cdot
 \bigl|\nabla_{\alpha}f\bigr|_h
 y_0^{-k}
 \Bigr).
\]
We also have the following inequality:
\[
 \Bigl|
 \chitilde_N\cdot
 h(f,\nabla_{\alpha}f)
\cdot y_0^{-2k-1}
 \del_{\alpha}y_0
 \Bigr|
\leq \\
 2\Bigl(
 C_1
 \chitilde_N^{1/2}\cdot y_0^{-1}
 \Bigr)
 \cdot
 \Bigl(
 \chitilde_N^{1/2}
 \bigl|
 \nabla_{\alpha}f
 \bigr|_h 
 y_0^{-k}
 \Bigr).
\]
Note that
$\nabla_{\alphabar}\nabla_{\alpha}f
=(\nabla_{\alphabar}\nabla_{\alpha}
 -\nabla_{\alpha}\nabla_{\alphabar})f
=-F_{\alpha,\alphabar}(h)f$.
We have $C_4,C_5>0$
which are independent of $N$,
such that the following holds:
\begin{multline}
 \int_{\nbigb^{\lambda\ast}}
 \chitilde_N\cdot
 \bigl|\nabla_{\alpha}f\bigr|_h^2y_0^{-2k}\dvol
\leq
 C_4+
 C_5\Bigl(
 \int_{\nbigb^{\lambda\ast}}
 \chitilde_N\cdot
 \bigl|\nabla_{\alpha}f\bigr|_h^2y_0^{-2k}\dvol
 \Bigr)^{1/2}
+\int_{\nbigb^{\lambda\ast}}
 \chitilde_N\cdot
 h\bigl(
 f,F_{\alpha,\alphabar}f
 \bigr)y_0^{-2k}\dvol.
\end{multline}
Similarly, we have the following:
\begin{multline}
 \int_{\nbigb^{\lambda\ast}}
 \chitilde_N\cdot
 \bigl|(\nabla_{\tau}+\sqrt{-1}\ad\phi) f\bigr|_h^2y_0^{-2k}\dvol
\leq
 C_5+
 C_6\Bigl(
 \int_{\nbigb^{\lambda\ast}}
 \chitilde_N\cdot
 \bigl|(\nabla_{\tau}+\sqrt{-1}\phi)f\bigr|_h^2y_0^{-2k}\dvol
 \Bigr)^{1/2}
 \\
+\int_{\nbigb^{\lambda\ast}}
 \chitilde_N\cdot
 h\bigl(
 f,-2\sqrt{-1}\nabla_{\tau}\phi\cdot f
 \bigr)y_0^{-2k}\dvol.
\end{multline}
Here, $C_i$ $(i=5,6)$ are positive constants,
which are independent of $N$.
Because $G(h)$ is $L^1$,
we have a constant $C_7>0$,
which is independent of $N$,
such that the following holds:
\[
\int_{\nbigb^{\lambda\ast}}
 \chitilde_N\cdot
 h\bigl(
 f,F_{\alpha,\alphabar}f
 \bigr)y_0^{-2k}\dvol
+\frac{1}{4}
 \int_{\nbigb^{\lambda\ast}}
 \chitilde_N\cdot
 h\bigl(
 f,-2\sqrt{-1}\nabla_{\tau}\phi\cdot f
 \bigr)y_0^{-2k}\dvol
\leq C_7.
\]
We put
\[
 A_N:=
\int_{\nbigb^{\lambda\ast}}
 \chitilde_N\cdot
 \bigl|\nabla_{\alpha}f\bigr|_h^2y_0^{-2k}\dvol
+\frac{1}{4}\int_{\nbigb^{\lambda\ast}}
 \chitilde_N\cdot
 \bigl|(\nabla_{\tau}+\sqrt{-1}\ad\phi) f\bigr|_h^2y_0^{-2k}\dvol.
\]
We have constants $C_i>0$ $(i=8,9)$,
which are independent of $N$,
such that the following holds:
\[
A_N\leq C_8+C_9A_N^{1/2}.
\]
Hence, we obtain that $A_N$ are bounded.
By taking $N\to\infty$,
we obtain the claim of the lemma.
\hfill\qed

\vspace{.1in}

Let $h_{2,E_1}$ be a Hermitian metric of $E_1$
such that the following holds.
\begin{itemize}
\item
We have a neighbourhood $N_1$ of $Z$
and that 
$h_{2,E_1}=h_{0,E_1}$ on $\nbigmlambda\setminus N_1$.
\item
We have a neighbourhood $N_2$ of $Z$
contained in $N_1$ such that
$h_{2,E_1}=h_{1,E_1}$ on $N_2\setminus Z$.
\end{itemize}

We have the function $s$ determined by
$h_{1,E_1}=h_{2,E_1}\cdot s$.
We have the relation
$G(h_{1,E_1})-G(h_{2,E_1})=4^{-1}\Delta \log s$.
The support of $\log s$ is contained in 
$\nbigm^{\lambda}\setminus N_2$.
By using the previous lemma, 
we obtain 
$\int\Delta\log s=0$.
Hence, we have
$\int G(h_{1,E_1})=\int G(h_{2,E_1})$.
By using the argument
in the proof of \cite[Proposition 9.4]{Mochizuki-difference-modules},
we obtain 
$\int G(h_{0,E_1})=\int G(h_{2,E_1})$.
\hfill\qed

\subsubsection{Analytic degree of subbundles}

Let $E_2\subset E$ be a mini-holomorphic subbundle.
Let $h_{1,E_2}$ denote the metric of $E_2$
induced by $h_1$.
By the Chern-Weil formula,
$\deg(E_2,h_{1,E_2})\in\real\cup\{-\infty\}$ 
makes sense.

\begin{prop}
Suppose that
$\deg(E_2,h_{1,E_2})\neq-\infty$.
Then, there exists a good filtered subbundle
$\nbigp_{\ast}\nbige_2\subset
 \nbigp_{\ast}\nbige$ 
such that
$\nbigp_a\nbige_{2|\nbigmlambda\setminus Z}
=E_2$.
Moreover,
$\deg(E_2,h_{1,E_2})=C\deg(\nbigp_{\ast}\nbige_2)$
holds,
where $C$ is the constant in 
Proposition {\rm\ref{prop;19.2.8.50}}.
\end{prop}
\pf
By (\ref{eq;18.11.21.1}) and \cite[Lemma 10.6]{Simpson88},
$E_{2|\pi^{-1}(\ttt)\cap\nbigm^{\lambda}}$ are extended to 
a locally free $\nbigo_{\proj^1_{\beta_1}}(\ast\{0,\infty\})$-submodules
of $\nbigp\nbige_{|\pi^{-1}(t_1)}$.
We take $P\in H^{\lambda}_{\infty}$
and a small neighbourhood  $\nbigu_P$ of $P$ 
in $\nbigmbar^{\lambda}$.
On $\nbigu_P$,
we use a local mini-complex coordinate system $(\ttU^{-1},\ttt)$.
On $\nbigutilde_P:=\real_{\tts}\times \nbigu_P$,
we use the complex coordinate system
$(\ttU^{-1},\ttv)=(\ttU^{-1},\tts+\sqrt{-1}\ttt)$
as in \S\ref{subsection;18.11.21.2}.
We set $D:=\real_{\tts}\times (\nbigu_P\cap H^{\lambda}_{\infty})$.
Then, 
we have the locally free 
$\nbigo_{\nbigutilde_P}(\ast D)$-module
$\widetilde{\nbigp\nbige}$
induced by $\nbigp\nbige$.
We also have the holomorphic vector subbundle 
$\Etilde_2$ of
$\widetilde{\nbigp\nbige}_{|\nbigutilde_P\setminus D}$
induced by $E_2$.
Let $p:\nbigutilde_P\lrarr D$ be the projection
given by $p(\ttU^{-1},\ttv)=\ttv$.
By the above consideration,
$\Etilde_{2|p^{-1}(\ttv)}$
extends
$\nbigo_{p^{-1}(\ttv)}(\ast\infty)$-submodule of
$\widetilde{\nbigp\nbige}_{|p^{-1}(\ttv)}$.
By using \cite[Theorem 4.5]{Siu-extension},
we obtain that $\Etilde_2$
extends
$\nbigp_{\nbigutilde_P}(\ast D)$-submodule
$\widetilde{\nbigp\nbige_2}$
of $\widetilde{\nbigp\nbige}$.
By the construction,
$\widetilde{\nbigp\nbige_2}$
is naturally $\real$-equivariant,
we obtain that $E_{2|\nbigu_P\setminus H^{\lambda}_{\infty}}$ 
extends to a locally free
$\nbigo_{\nbigu_P}
 \bigl(\ast (H^{\lambda}_{\infty}\cap \nbigu_P)\bigr)$-submodule
of $\nbigp\nbige_{|\nbigu_P}$.
Hence, we obtain that $E_{2}$ is extended to 
a locally free 
$\nbigo_{\nbigmbar^{\lambda}}(\ast H^{\lambda}_{\infty})$-module
$\nbigp\nbige_2$.
We have the naturally induced good filtered bundle
$\nbigp_{\ast}\nbige_2$
over $\nbigp\nbige_2$.
The claim for the degree follows from 
Proposition \ref{prop;19.2.8.50}.
\hfill\qed

\vspace{.1in}

As a consequence,
we obtain the following.
\begin{cor}
\label{cor;18.11.21.4}
$\nbigp_{\ast}\nbige$ is stable if and only if
$(E,h_1)$ is analytic stable.
\hfill\qed
\end{cor}

\subsection{Proof of Theorem \ref{thm;18.11.21.20}}
\label{subsection;19.1.27.4}

\subsubsection{Associated filtered bundles}

Let $Z$ be a finite subset of $\nbigm^{\lambda}$.
Let $(E,\delbar_E,h)$ be a meromorphic monopole on
$\nbigm^{\lambda}\setminus Z$.
Let $\nbigp_{\ast}E$ be the associated filtered bundle
with Dirac type singularity on $(\nbigmbar^{\lambda};H^{\lambda},Z)$.
\begin{prop}
\label{prop;19.2.8.121}
The good filtered bundle $\nbigp_{\ast}E$ is 
polystable with $\deg(\nbigp_{\ast}E)=0$.
If the monopole $(E,\delbar_E,h)$ is irreducible,
$\nbigp_{\ast}E$ is stable.
\end{prop}
\pf
By Corollary \ref{cor;18.12.17.3},
$(E,\delbar_E,h)$ satisfies the condition 
in \S\ref{subsection;18.11.21.5}.
We obtain 
\[
 C\deg(\nbigp_{\ast}E)=\deg(E,h)=0.
\]
Let $\nbigp_{\ast}E_1$ be a good filtered subbundle of
$\nbigp_{\ast}E$.
We have 
$C\deg(\nbigp_{\ast}E_1)=
 \deg(E_1,h_{E_1})\leq 0$.
Moreover, if $\deg(\nbigp_{\ast}E_1)=0$,
$E_1$ is flat with respect to the Chern connection,
and the orthogonal decomposition
$E=E_1\oplus E_1^{\bot}$ is mini-holomorphic.
Hence, we have the decomposition
$\nbigp_{\ast}E=
 \nbigp_{\ast}E_1\oplus\nbigp_{\ast}E_1^{\bot}$.
We also have that
$E_1$ and $E_1^{\bot}$ with the induced metrics
are monopoles.
Hence, we obtain the poly-stability of
$\nbigp_{\ast}E$
by an easy induction on the rank of $E$.
\hfill\qed

\subsubsection{Uniqueness}

\begin{prop}
\label{prop;18.11.21.10}
Let $h'$ be another metric of $E$ such that
(i) $(E,\delbar_E,h')$ is a monopole,
(ii) any points of $Z$ are Dirac type singularity,
(iii) $h'$ is adapted to $\nbigp_{\ast}E$.
Then, the following holds.
\begin{itemize}
\item
There exists a mini-holomorphic decomposition
$(E,\delbar_E)
=\bigoplus_{i=1}^{m}
 (E_i,\delbar_{E_i})$,
which is orthogonal with respect to both
$h$ and $h'$.
\item
There exist positive numbers $a_i$ $(i=1,\ldots,m)$
such that $h_{E_i}=a_ih'_{E_i}$.
\end{itemize}
\end{prop}
\pf
By the norm estimate,
$h$ and $h'$ are mutually bounded.
Hence, we obtain the claim
from \cite[Proposition 2.4, Proposition 3.16]{Mochizuki-KH-infinite}.
\hfill\qed

\subsubsection{Construction of monopoles}

Let $Z$ be a finite subset.
Let $\nbigp_{\ast}\nbige$ be a  stable good filtered bundle
with Dirac type singularity on $(\nbigmbar^{\lambda};H_{\infty}^{\lambda},Z)$
with $\deg(\nbigp_{\ast}\nbige)=0$.
Set $E:=\nbigp_a\nbige_{|\nbigmlambda\setminus Z}$.

The following proposition is similar to
\cite[Proposition 9.10]{Mochizuki-difference-modules}.

\begin{prop}
\label{prop;19.2.8.120}
There exists a Hermitian metric $h$ such that
(i) $(E,\delbar_E,h)$ is a meromorphic monopole,
(ii) $(E,\delbar_E,h)$ satisfies the norm estimate
 with respect to $\nbigp_{\ast}\nbige$.
\end{prop}
\pf
We give only an outline.
By Proposition \ref{prop;18.9.2.3},
there exists a Hermitian metric $h_0$ of $E$
such that
(i) $(\nbigp_{\ast}\nbige,h_0)$ satisfies
 the norm estimate,
(ii) $(E,\delbar_E,h_0)$ is a monopole 
 with Dirac type singularity on a neighbourhood of
 each $P\in Z$,
(iii) $G(h_0)=O(e^{-|y_0|})$.
By Proposition \ref{prop;19.2.8.110},
we may assume that
$(\det(E),\delbar_{\det(E)},\det(h_0))$
is a monopole with Dirac type singularity
such that $\det(h_0)$ is adapted to
$\nbigp_{\ast}(\det(E))$.
By Corollary \ref{cor;18.11.21.4},
$(E,\delbar_E,h_0)$ is analytically stable.
By \cite[Theorem 2.5, Proposition 3.16]{Mochizuki-KH-infinite},
there exists a Hermitian metric $h$ of $E$
such that the following holds:
\begin{itemize}
\item
$\det(h)=\det(h_0)$.
\item
$G(h)=0$,
i.e., $(E,\delbar_E,h)$ is a monopole.
\item
Let $s$ be the automorphism of $E$
which is self-adjoint with respect to $h$ and $h_0$,
determined by $h=h_0s$.
Then, $s$ and $s^{-1}$ are bounded
with respect to $h_0$,
and $\delbar_Es$ is $L^2$
\end{itemize}
By Proposition \ref{prop;19.2.8.100},
there exists a compact subset $C\subset\nbigm^{\lambda}$
such that 
(i) $Z\subset C$,
(ii) $F(h)$ is bounded on $\nbigm^{\lambda}\setminus C$.
By \cite[Proposition 2.10]{Mochizuki-difference-modules},
each point of $Z$ is Dirac type singularity of
$(E,\delbar_E,h)$.
Because $s$ and $s^{-1}$ are bounded,
$(\nbigp_{\ast}\nbige,h)$
satisfies the norm estimate.
Thus, we obtain Proposition \ref{prop;19.2.8.120}.
\hfill\qed

\vspace{.1in}

The claim of Theorem \ref{thm;18.11.21.20}
follows from 
Proposition \ref{prop;19.2.8.121},
Proposition \ref{prop;18.11.21.10}
and Proposition \ref{prop;19.2.8.120}.
\hfill\qed

\section{Riemann-Hilbert correspondences
of filtered objects $(|\lambda|\neq 1)$}
\label{section;19.1.27.1}

We give a complement 
on the Riemann-Hilbert correspondence
for good filtered bundles with Dirac type singularity
on $(\nbigmbar^{\lambda};H^{\lambda},Z)$
for a finite subset $Z\subset\nbigm^{\lambda}$
in the case $|\lambda|\neq 1$.
It is a parabolic version of 
the Riemann-Hilbert correspondence
for local analytic $\gminiq$-difference modules,
due to Ramis, Sauloy and Zhang \cite{Ramis-Sauloy-Zhang}
and van der Put and Reversat \cite{van-der-Put-Reversat},
and for the global $\gminiq$-difference modules
due to Kontsevich and Soibelman,
where $|\gminiq|\neq 1$.

As a result,
from meromorphic doubly periodic monopoles,
for each $\lambda$ with $|\lambda|\neq 1$,
we obtain filtered objects
on the elliptic curve
$\cnum^{\ast}/(\gminiq^{\lambda})^{\seisuu}$.
They are constructed through
the associated good filtered bundles
on $(\nbigmbar^{\lambda};H^{\lambda},Z)$.
Recall that
$(\nbigmbar^{\lambda};H^{\lambda},Z)$
depends on the choice of $\tte_1$ and $\tts_1$.
However, the induced filtered objects
on $\cnum^{\ast}/(\gminiq^{\lambda})^{\seisuu}$
are essentially independent of
the choice of $\tte_1$ and $\tts_1$
(Theorem \ref{thm;19.2.8.130}).

\subsection{Analytic $\gminiq$-difference modules}

Let $\nbigk^{\an}$ denote the field of
the convergent Laurent power series $\cnum(\{y\})$.
Let $\nbigr^{\an}$ denote the ring of
the convergent power series $\cnum[\{y\}]$.
Let $\gminiq\in\cnum^{\ast}$.
Suppose that $|\gminiq|\neq 1$.
Let $\Phi^{\ast}:\nbigk^{\an}\lrarr\nbigk^{\an}$
be determined by
$\Phi^{\ast}(f)(y):=f(\gminiq y)$.
A $\gminiq$-difference $\nbigk^{\an}$-module
is a finite dimensional $\nbigk^{\an}$-vector space $\nbigv^{\an}$
equipped with a $\cnum$-linear automorphism
$\Phi^{\ast}$ such that
$\Phi^{\ast}(fs)=\Phi^{\ast}(f)\cdot\Phi^{\ast}(s)$
for any $f\in\nbigk^{\an}$ and $s\in\nbigv^{\an}$.
Let $\Diff(\nbigk^{\an},\gminiq)$
denote the category of
$\gminiq$-difference $\nbigk^{\an}$-modules.
By taking the formal completion
\[
 \gbigc(\nbigv^{\an},\Phi^{\ast}):=
 (\nbigv^{\an}\otimes_{\nbigk^{\an}}\nbigk,\Phi^{\ast}),
\]
we obtain the functor
$\gbigc:\Diff(\nbigk^{\an},\gminiq)
\lrarr
 \Diff(\nbigk,\gminiq)$.

\subsubsection{Pure isoclinic modules}

Let $\omega\in\rnum$.
A $\gminiq$-difference $\nbigk^{\an}$-module 
$(\nbigv^{\an},\Phi^{\ast})$
is called pure isoclinic of slope $\omega$
if 
$\gbigc(\nbigv^{\an},\Phi^{\ast})$
is pure isoclinic of slope $\omega$.
Let $\Diff(\nbigk^{\an},\gminiq;\omega)$
denote the full subcategory of
pure isoclinic $\gminiq$-difference
$\nbigk^{\an}$-modules of slope $\omega$.
It is known that $\gbigc$ induces an equivalence
\[
\gbigc:
 \Diff(\nbigk^{\an},\gminiq;\omega)
\simeq
 \Diff(\nbigk,\gminiq;\omega).
\]

\subsubsection{Slope filtrations}

Any $(\nbigv,\Phi^{\ast})\in\Diff(\nbigk,\gminiq)$
has a slope decomposition
$(\nbigv,\Phi^{\ast})
=\bigoplus_{\omega}
 (\nbigv_{\omega},\Phi^{\ast})$,
where
$(\nbigv_{\omega},\Phi^{\ast})
 \in\Diff(\nbigk,\gminiq)$.
We define the slope filtration $\gbigf$
of $(\nbigv,\Phi^{\ast})$ indexed by
$(\rnum,\leq)$ as follows:
\[
 \gbigf_{\mu}\nbigv:=
 \bigoplus_{\varrho(\gminiq)\omega\leq \mu}
 \nbigv_{\omega},
\]
where we put $\varrho(\gminiq):=1$ $(|\gminiq|>1)$
or $\varrho(\gminiq):=-1$ $(|\gminiq|<1)$.
We naturally have
$\Gr^{\gbigf}_{\mu}(\nbigv)
=\nbigv_{\varrho(\gminiq)\mu}$.

According to Sauloy \cite{Sauloy2004},
any $(\nbigv^{\an},\Phi^{\ast})
\in\Diff(\nbigk^{\an},\gminiq)$
has a unique filtration $\gbigf$
indexed by $(\rnum,\leq)$ 
such that
$\gbigc\gbigf_{\mu}(\nbigv^{\an})
=\gbigf_{\mu}(\gbigc(\nbigv^{\an}))$.
In particular,
$\Gr^{\gbigf}_{\mu}(\nbigv^{\an})$
is pure isoclinic of slope
$\varrho(\gminiq)\mu$.
The filtration is functorial,
i.e.,
for any morphism
$f:\nbigv^{\an}_1\lrarr\nbigv^{\an}_2$,
we have
$f(\gbigf_{\mu}\nbigv^{\an}_1)
\subset \gbigf_{\mu}\nbigv^{\an}_2$,
and more strongly
$f(\gbigf_{\mu}\nbigv^{\an}_1)
=\gbigf_{\mu}\nbigv^{\an}_2
 \cap f(\nbigv^{\an}_1)$.

\subsubsection{Equivalences}

We set $T:=\cnum^{\ast}/\gminiq^{\seisuu}$.
Let $\Vect(T)$ denote the category of
locally free $\nbigo_T$-modules of finite rank.
For any $\mu\in\rnum$,
let $\Vect^{ss}(T,\mu)\subset\Vect(T)$
denote the full subcategory of
semistable sheaves of slope $\mu$,
i.e., $\ttE\in\Vect(T)$ such that
$\deg(\ttE)/\rank(\ttE)=\mu$.

For $\ttE\in\Vect^{ss}(T)$,
a $\rnum$-anti-Harder-Narasimhan filtration of $\ttE$
is a filtration $\gbigf$ of $\ttE$ in $\Vect(T)$
indexed by $(\rnum,\leq)$
such that $\Gr^{\gbigf}_{\mu}(\ttE)\in \Vect^{ss}(T,\mu)$.
Let $\Vect^{\rnum\AHN}(T)$ denote 
the category of locally free $\nbigo_T$-modules
$\ttE$ equipped with a $\rnum$-anti-Harder-Narasimhan filtration $\gbigf$.

Let us recall that there exists a natural equivalence
\[
 \ttK:
\Vect^{\rnum\AHN}(T)
\simeq
 \Diff(\nbigk^{\an},\gminiq)
\]
due to van der Put, Reversat \cite{van-der-Put-Reversat}
and Ramis, Sauloy and Zhang \cite{Ramis-Sauloy-Zhang}.
Let $(\ttE,\gbigf)\in \Vect^{\rnum\AHN}(T)$.
We obtain the $\gminiq^{\seisuu}$-equivariant
locally $\nbigo_{\cnum^{\ast}}$-module $\nbige$
by the pull back
$\cnum^{\ast}\lrarr T$.
It is equipped with 
$\gminiq^{\seisuu}$-equivariant filtration $\gbigf$.
There exists a canonical extension of $\nbige$
to a $\gminiq^{\seisuu}$-equivariant
locally free $\nbigo_{\cnum}(\ast 0)$-module
$\nbigetilde$
equipped with a $\gminiq^{\seisuu}$-equivariant
filtration $\gbigf$
such that
the formal completion of
$\Gr^{\gbigf}_{\mu}(\nbigetilde)$
are pure isoclinic of slope $\varrho(\gminiq)\mu$.
By taking the stalk of $\nbigetilde$ at $0$,
we obtain $\ttK(E,\gbigf)\in \Diff(\nbigk^{\an},\gminiq)$.
The same procedure induces
$\ttK: \Vect^{ss}(T,\mu)
\simeq
\Diff(\nbigk^{\an},\gminiq;\varrho(\gminiq)\omega)$.

For any $\mu\in\rnum$,
we take
$\LLtilde_1(\mu)\in \Vect^{ss}(T,\mu)$
with an isomorphism
$\ttK(\LLtilde_1(\mu))\simeq
 \LL_1(\varrho(\gminiq)\mu)$
in $\Diff(\nbigk^{\an},\gminiq;\omega)$.
(See \S\ref{subsection;18.12.21.2}
for $\LL_m(\omega)$.)
For any $A\in \GL_r(\cnum)$,
we take $\VVtilde_1(A)\in\Vect^{ss}(T,0)$
with an isomorphism
$\ttK(\VVtilde_1(A))\simeq \VV_1(A)$
in $\Diff(\nbigk^{\an},\gminiq;0)$.
(See Example \ref{example;18.12.20.1}
for $\VV_m(A)$.)
Similarly,
for any finite dimensional $\cnum$-vector space
$V$ equipped with an automorphism $f$,
we take $\VVtilde_1(V,f)\in \Vect^{ss}(T,0)$
with an isomorphism
$\ttK(\VVtilde_1(V,f))\simeq\VV_1(V,f)$.

\subsection{Classification of good filtered formal
$\gminiq$-difference modules
in the case $|\gminiq|\neq 1$}
\label{subsection;19.2.4.30}

Let $\Vect^{ss}(T;\mu)^{\Par}$ denote
the category of $\ttE\in\Vect^{ss}(T,\mu)$
equipped with a filtration
$\nbigf_{\bullet}(\ttE)$ indexed by $(\rnum,\leq)$
such that
(i) $\nbigf_a(\ttE)=\bigcap_{a<b}\nbigf_b(\ttE)$,
(ii) $\Gr^{\nbigf}_a(\ttE):=\nbigf_a(\ttE)/\nbigf_{<a}(\ttE)
 \in\Vect^{ss}(T,\mu)$
for any $a\in\real$.
Note that
$\{a\in\real\,|\,\Gr^{\nbigf}_a(\ttE)\neq 0\}$ is finite.
For any $C>0$,
let us construct an equivalence
$\ttK^C:\Vect^{ss}(T;\mu)^{\Par}
\simeq
 \Diff(\nbigk,\gminiq;\varrho(\gminiq)\mu)^{\Par}$
depending on $C$.

\subsubsection{The case $\mu=0$}

Take $A_{\alpha}\in \GL_r(\cnum)$
which has a unique eigenvalue $\alpha$.
Let $\nbigf$ be a filtration of $\VVtilde(A_{\alpha})$
such that
$(\VVtilde(A_{\alpha}),\nbigf)\in \Vect^{ss}(T,0)^{\Par}$.
We obtain the induced filtration
$\nbigf$ on $\ttK(\VVtilde(A_{\alpha}))$
in $\Diff(\nbigk,\gminiq;0)$.
For $a\in\real$,
we set
\begin{equation}
\label{eq;19.1.25.10}
 b(\gminiq,\alpha,a):=
 C\cdot\Bigl(
 a+\frac{\log|\alpha|}{\log|\gminiq|}
\Bigr).
\end{equation}

We define the filtration $F$
of $\ttK(\VVtilde(A_{\alpha}))$
in $\Diff(\nbigk,\gminiq;0)$
indexed by $\real$
as follows:
\begin{equation}
\label{eq;19.1.25.11}
 F_{a}
 \ttK(\VVtilde(A_{\alpha}))
=\nbigf_{b(\gminiq,\alpha,a)}
 \ttK(\VVtilde(A_{\alpha})).
\end{equation}
There exists a frame $\vecv$ of $\ttK(\VVtilde(A_{\alpha}))$
such that 
(i) $\Phi^{\ast}\vecv=\vecv A_{\alpha}$,
(ii) $\vecv$ is compatible with $F$,
i.e., there exists a decomposition
$\vecv=\coprod_{c\in\real}\vecv_c$
such that
$\coprod_{c\leq a}\vecv_c$
is a frame of
$F_a\ttK(\VVtilde(A_{\alpha}))$.
For each $v_i$,
let $c(v_i)$ be determined by
$v_i\in \vecv_{c(v_i)}$.
We define
\[
 \nbigp_d\ttK(\VVtilde(A_{\alpha}))
=\bigoplus
 \nbigr\cdot y^{-[d-c(v_i)]}v_i.
\]
In this way,
we obtain the filtered bundle
$\ttK^C(\VVtilde(A_{\alpha}),\nbigf):=
\nbigp_{\ast}\ttK(\VVtilde(A_{\alpha}))$.

In general, 
for any $(\ttE,\nbigf)\in \Vect^{ss}(T,0)^{\Par}$,
there exist a partition $r=\sum r_i$,
matrices
$A_{\alpha_i}\in\GL_{r_i}(\cnum)$
with a unique eigenvalue $\alpha_i$,
objects
$(\VVtilde(A_{\alpha_i}),\nbigf)\in\Vect^{ss}(T;0)^{\Par}$,
and an isomorphism
\begin{equation}
\label{eq;19.1.18.10}
 (\ttE,\nbigf)
\simeq
 \bigoplus_{i=1}^{N}
 (\VVtilde(A_{\alpha_i}),\nbigf).
\end{equation}
We obtain the filtered bundle
$\nbigp_{\ast}\ttK(\ttE)$ over $\ttK(\ttE)$
induced by
the isomorphism
$\ttK(\ttE)\simeq
 \bigoplus \ttK(\VVtilde(A_{\alpha_i}))$
and the filtered bundle
$\bigoplus
 \ttK^C(\VVtilde(A_{\alpha_i}),\nbigf)$.
It is easy to check that
$\nbigp_{\ast}\ttK(\ttE)$
is independent of the choice of 
$A_{\alpha_i}$ and an isomorphism (\ref{eq;19.1.18.10}).
We define
$\ttK^C(\ttE,\nbigf):=\nbigp_{\ast}\ttK(\ttE)$.
Thus, we obtain a functor
\[
 \ttK^C:\Vect^{ss}(T;0)^{\Par}
\lrarr
 \Diff(\nbigk;0)^{\Par}.
\]

\begin{lem}
\label{lem;19.1.18.20}
$\ttK^C$ induces an equivalence
$\Vect^{ss}(T;0)^{\Par}
\simeq
 \Diff(\nbigk,\gminiq;0)$.
\end{lem}
\pf
Let $\nbigl$ be a lattice of $\nbigv=\VV(A_{\alpha})$
such that
$\Phi^{\ast}(\nbigl)=\nbigl$.
We obtain the automorphism
$\sigma(\Phi^{\ast};\nbigl)$ of $\nbigl_{|0}$,
and the generalized eigen decomposition
\[
 \nbigl_{|0}
=\bigoplus_{i\in\seisuu}
 \EE_{\alpha \gminiq^{-i}}(\nbigl_{|0}).
\]
We set
$i_0:=\max\{i\,|\,\EE_{\alpha\gminiq^{-i}}(\nbigl_{|0})\neq 0\}$
and 
$i_1:=\min\{i\,|\,\EE_{\alpha\gminiq^{-i}}(\nbigl_{|0})\neq 0\}$.

If $i_0>0$,
we define $\nbigl'$ 
as the kernel of 
$\nbigl\lrarr
 \EE_{\alpha \gminiq^{-i_0}}(\nbigl_{|0})$.
Then, it is easy to see that
$\Phi^{\ast}(\nbigl')=\nbigl'$.
We have
the natural inclusion
$\nbigl'\lrarr\nbigl$.
It induces
$\EE_{\alpha\gminiq^{-i}}(\nbigl'_{|0})
=\EE_{\alpha\gminiq^{-i}}(\nbigl_{|0})$
for $i<i_0-1$,
and the following exact sequence:
\begin{equation}
\label{eq;19.1.9.10}
 0\lrarr
 \EE_{\alpha\gminiq^{-i_0+1}}
 \bigl((y\nbigl)_{|0}\bigr)
\lrarr
 \EE_{\alpha\gminiq^{-i_0+1}}
 (\nbigl'_{|0})
\lrarr
 \EE_{\alpha\gminiq^{-i_0+1}}
 \bigl(\nbigl_{|0}\bigr)
\lrarr 0.
\end{equation}
Moreover, we have
$\EE_{\alpha\gminiq^{-i}}(\nbigl'_{|0})=0$
for $i\geq i_0$.

If $i_1<0$,
we define $\nbigl''$ as the kernel of the following:
\[
 y^{-1}\nbigl\lrarr 
 \bigoplus_{i>i_0+1}
 \EE_{\alpha\gminiq^{-i}}\bigl(
 (y^{-1}\nbigl)_{|0}
 \bigr).
\]
We have $\Phi^{\ast}(\nbigl'')=\nbigl''$.
We have the natural inclusion
$\nbigl\lrarr \nbigl''$.
It induces 
$\EE_{\alpha \gminiq^{-i}}(\nbigl_{|0})
\simeq
 \EE_{\alpha\gminiq^{-i}}(\nbigl''_{|0})$
for $i>i_1+1$,
and the following exact sequence:
\begin{equation}
 \label{eq;19.1.9.11}
 0\lrarr
 \EE_{\alpha \gminiq^{-i_1-1}}(\nbigl_{|0})
\lrarr
  \EE_{\alpha \gminiq^{-i_1-1}}(\nbigl''_{|0})
\lrarr
  \EE_{\alpha \gminiq^{-i_1-1}}((y^{-1}\nbigl)_{|0})
\lrarr 0.
\end{equation}

Suppose that 
each $\EE_{\alpha \gminiq^{-i}}(\nbigl_{|0})$
is equipped with a filtration
$F$ satisfying the following conditions.
\begin{itemize}
\item
 For $i_1<i<i_0$,
 $F_{\bullet}\EE_{\alpha \gminiq^{-i}}(\nbigl_{|0})$
is indexed by
 $\openclosed{-1}{0}$.
\item
 $F_{\bullet}\EE_{\alpha \gminiq^{-i_0}}(\nbigl_{|0})$
is indexed by $\real_{\leq 0}$.
\item
 $F_{\bullet}\EE_{\alpha \gminiq^{-i_1}}(\nbigl_{|0})$
is indexed by $\real_{>-1}$.
\end{itemize}
Note that 
 $F_{\bullet}\EE_{\alpha \gminiq^{-i_0}}(\nbigl_{|0})$
induces a filtration
 $F_{\bullet}\EE_{\alpha \gminiq^{-i_0+1}}((y\nbigl)_{|0})$
indexed by $\real_{\leq -1}$,
and that 
 $F_{\bullet}\EE_{\alpha \gminiq^{-i_1}}(\nbigl_{|0})$
induces a filtration
 $F_{\bullet}\EE_{\alpha \gminiq^{-i_1-1}}((y^{-1}\nbigl)_{|0})$
is indexed by $\real_{>0}$.

We obtain a filtration
$F_{\bullet}\EE_{\alpha \gminiq^{-i}}(\nbigl'_{|0})$
for $i<i_0-1$
by using the isomorphism
$\EE_{\alpha \gminiq^{-i}}(\nbigl'_{|0})
\simeq
\EE_{\alpha \gminiq^{-i}}(\nbigl_{|0})$.
We obtain the filtration 
$F_{\bullet}\EE_{\alpha \gminiq^{-i_0+1}}(\nbigl'_{|0})$
indexed by $\real_{\leq 0}$
by using exact sequence (\ref{eq;19.1.9.10}).

We obtain a filtration 
$F_{\bullet}\EE_{\alpha \gminiq^{-i}}(\nbigl''_{|0})$
for $i>i_1+1$
by using the isomorphism
$\EE_{\alpha \gminiq^{-i}}(\nbigl''_{|0})
\simeq
\EE_{\alpha \gminiq^{-i}}(\nbigl_{|0})$.
We obtain the filtration 
$F_{\bullet}\EE_{\alpha \gminiq^{-i_1-1}}(\nbigl''_{|0})$
indexed by $\real_{>-1}$
by using the exact sequence (\ref{eq;19.1.9.11}).

Let $\vece$ be a frame of $\nbigv$
such that $\Phi^{\ast}(\vece)=\vece\cdot A_{\alpha}$.
Let $\nbigl(A_{\alpha})$  be the lattice of $\nbigv$
generated by $\vece$.
Note that such $\nbigl(A_{\alpha})$ is 
independent of a choice of $\vece$.
Starting from a regular filtered bundle
$\nbigp_{\ast}\nbigv$,
by applying the above procedure inductively,
we obtain a filtration 
$F$ on $\nbigl(A_{\alpha})_{|0}$.
There exists a unique filtration $F$
of $\nbigl(A_{\alpha})$
such that
(i) it is preserved by $\Phi^{\ast}$,
(ii) it induces $F(\nbigl(A_{\alpha})_{|0})$.
We define $\nbigf(\nbigv)$ from $F$
by using (\ref{eq;19.1.25.10}) and (\ref{eq;19.1.25.11}).
It is easy to see that this gives a quasi-inverse
of $\ttK$.
\hfill\qed

\subsubsection{The case of general $\mu$}

Let $\mu\in\rnum$.
Let $(\ttE_{\mu},\nbigf)\in\Vect^{ss}(T;\mu)^{\Par}$.
There exists $(\ttE'_0,\nbigf)\in\Vect^{ss}(T;0)^{\Par}$
with an isomorphism
$\nbigf_{\bullet}\ttE_{\mu}\simeq
 \LLtilde_1(\mu)\otimes
 \nbigf_{\bullet}\ttE'_0$.
We define
\[
 \ttK^C(\ttE_{\mu},\nbigf)
:=
\nbigp_{\ast}^{(-\varrho(\gminiq)\mu/2)}\LL_1(\varrho(\gminiq)\mu)
\otimes
 \ttK^C(\ttE'_0,\nbigf)
\in 
 \Diff(\nbigk,\gminiq;\varrho(\gminiq)\mu)^{\Par}.
\]
Thus, we obtain a functor
$\ttK^C:\Vect^{ss}(T;\mu)^{\Par}
\lrarr
 \Diff(\nbigk,\gminiq;\varrho(\gminiq)\mu)^{\Par}$.
As a consequence of Lemma \ref{lem;19.1.18.20},
we obtain the following lemma.
\begin{lem}
$\ttK^C$ induces an equivalence
$\Vect^{ss}(T;\mu)^{\Par}
\simeq
 \Diff(\nbigk,\gminiq;\varrho(\gminiq)\mu)$.
\hfill\qed
\end{lem}

\begin{rem}
Let $(\ttE,\nbigf)\in \Vect^{ss}(T;\mu)^{\Par}$.
We define a new filtration
$\nbigf^{(C)}$ on $\ttE$
by 
$\nbigf^{(C)}_a(\ttE):=\nbigf_{Ca}(\ttE)$.
The correspondence
$(\ttE,\nbigf)\longmapsto
 (\ttE,\nbigf^{(C)})$
induces an equivalence
$\ttH^{C}:
 \Vect^{ss}(T;\mu)^{\Par}
\lrarr
 \Vect^{ss}(T;\mu)^{\Par}$.
The following is commutative
by the construction.
\[
 \begin{CD}
 \Vect^{ss}(T;\mu)^{\Par}
 @>{\ttK^1}>>
   \Diff(\nbigk,\gminiq;\varrho(\gminiq)\mu)\\
 @V{\ttH^C}VV @V{\id}VV \\
 \Vect^{ss}(T;\mu)^{\Par}
 @>{\ttK^C}>>
   \Diff(\nbigk,\gminiq;\varrho(\gminiq)\mu).
 \end{CD}
\]
\hfill\qed
\end{rem}

\subsubsection{Graded objects}

A $(\rnum,\real)$-grading of 
$\ttE\in\Vect(T)$
is a decomposition
\[
 \ttE=\bigoplus_{\mu\in\rnum}
 \bigoplus_{a\in\real}
  \ttE_{\mu,a}
\]
such that
$\ttE_{\mu,a}\in\Vect^{ss}(T,\mu)$.
Let $\Vect(T)_{(\rnum,\real)}$
is a category of
$\ttE\in\Vect(T)$ 
with a $(\rnum,\real)$-grading.
Let us construct a functor
$\Vect(T)_{(\rnum,\real)}
\lrarr
 \Diff(\cnum[y,y^{-1}],\gminiq)_{(\rnum,\real)}$.

For any $\ttE\in \Vect^{ss}(T;\mu)$,
we obtain 
$\gminiq^{\seisuu}$-equivariant
$\nbigo_{\cnum_y^{\ast}}$-module $\nbige$
as the pull back of $\ttE$
by $\cnum^{\ast}\lrarr T$.
It is extended to a locally free
$\gminiq^{\seisuu}$-equivariant
$\nbigo_{\cnum_y}(\ast 0)$-module $\nbigetilde_0$
such that the formal completion
$\nbigetilde_0\otimes\cnum(\!(y)\!)$
is naturally an isoclinic  $\gminiq$-difference 
$\cnum(\!(y)\!)$-module
of slope $\varrho(\gminiq)\mu$.
Similarly,
$\nbige$ is extended to
a locally free
$\gminiq^{\seisuu}$-equivariant
$\nbigo_{\cnum_{y^{-1}}}(\ast 0)$-module
$\nbigetilde_{\infty}$
such that 
$\nbigetilde_{\infty}\otimes\cnum(\!(y^{-1})\!)$
is naturally an isoclinic  $\gminiq^{-1}$-difference
$\cnum(\!(y^{-1})\!)$-module
of slope $-\varrho(\gminiq)\mu$.
By gluing $\nbigetilde_0$ and $\nbigetilde_{\infty}$,
we obtain a locally free
$\gminiq^{\seisuu}$-equivariant
$\nbigo_{\proj^1}(\ast\{0,\infty\})$-module
$\nbigetilde$.
Note that
$\nbigetilde$ has an $\nbigo_{\proj^1}$-lattice.
Hence, $\nbigetilde$ is an algebraic.
By taking the global section on $\proj^1$,
we obtain a $\gminiq$-difference
$\cnum[y,y^{-1}]$-module $\ttKtilde(E)$.

Suppose that $\ttE=\VVtilde(A_{\alpha})$,
where $A_{\alpha}$ has a unique eigenvalue $\alpha$.
Let $b\in\real$.
Let $\nbigf^{(b)}$ be the filtration
on $\ttE$ determined by
$\nbigf_b(\ttE)=\ttE$ and $\nbigf_{<b}(\ttE)=0$.
We obtain the filtered bundle
$\nbigp_{\ast}\bigl(
 \nbigetilde_{0}\otimes\cnum(\!(y)\!)
 \bigr)$ over 
$\nbigetilde_{0}\otimes\cnum(\!(y)\!)$
induced by $\ttK^C(\ttE,\nbigf^{(b)})$.
We set
\[
 a_0:=C^{-1}b-\frac{\log|\alpha|}{\log|\gminiq|}
\]
Then,
$\Gr^{\nbigp}_c(\nbigetilde_0)=0$
unless $c\in a_0-\varrho(\gminiq)\mu/2+\seisuu$.
We also obtain the filtered bundle
$\nbigp_{\ast}\bigl(
 \nbigetilde_{\infty}\otimes\cnum(\!(y^{-1})\!)
 \bigr)$ over 
$\nbigetilde_{\infty}\otimes\cnum(\!(y^{-1})\!)$.
We set
\[
 a_{\infty}:=C^{-1}b+\frac{\log|\alpha|}{\log|\gminiq|}.
\]
Then,
$\Gr^{\nbigp}_c(\nbigetilde_0)=0$
unless $c\in a_{\infty}+\varrho(\gminiq)\mu/2+\seisuu$.
Set $a:=a_0-\varrho(\gminiq)\mu/2$.
For any $n\in\seisuu$,
let $\nbigl_{(\mu,a+n)}
 \subset\nbigetilde$
be the lattice determined by
$\nbigp_{a_0+n-\varrho(\gminiq)\mu/2}\nbigetilde_0$
and 
$\nbigp_{a_{\infty}-n+\varrho(\gminiq)\mu/2}\nbigetilde_{\infty}$.
Then, it turns out that
$\nbigl_{(\mu,a+n)}$
is isomorphic to $\nbigo^{\rank E}$.
We set 
$\ttKtilde^C(\ttE)_{\mu,a+n}:=H^0(\proj^1,\nbigl_{(\mu,a+n)})$
for $n\in\seisuu$.
We also set
$\ttKtilde^C(\ttE)_{\mu,c}:=0$ unless $c-a\in \seisuu$.

Let $\ttE=\bigoplus\ttE_{\mu,b}\in \Vect(T)_{(\rnum,\real)}$.
For each $\ttE_{\mu,b}$, we apply the above construction,
and we obtain 
$\ttKtilde^C(\ttE_{\mu,b}) \in
 \Diff(\cnum[y,y^{-1}],\gminiq)_{(\rnum,\real)}$.
We define
$\ttKtilde^C(\ttE):=\bigoplus \ttKtilde^C(\ttE_{\mu,b})
 \in \Diff(\cnum[y,y^{-1}],\gminiq)_{(\rnum,\real)}$.

For each $(\ttE_{\mu},\nbigf)\in\Vect^{ss}(T;\mu)^{\Par}$,
we obtain 
$\Gr_{\bullet}^{\nbigf}(\ttE_{\mu})
=\bigoplus \Gr^{\nbigf}_a(\ttE_{\mu})
\in \Vect(T)_{(\rnum,\real)}$.
The following is easy to see
by the construction.
\begin{lem}
For $(\ttE_{\mu},\nbigf)\in\Vect^{ss}(T;\mu)^{\Par}$,
we have the natural isomorphism
$\ttG(\ttK^C(\ttE_{\mu},\nbigf))
\simeq
 \ttKtilde^C\bigl(\Gr_{\bullet}^{\nbigf}(\ttE_{\mu})\bigr)$.

\hfill\qed
\end{lem}

\subsubsection{Weight filtration}

Let $E\in\Vect^{ss}(T;\mu)$.
There exists an isomorphism
$E\simeq \bigoplus_i\LLtilde(\mu)\otimes \VVtilde(A_{\alpha_i})$,
where each $A_{\alpha_i}$ has a unique eigenvalue $\alpha_i$.
We obtain the logarithm $N_{\alpha_i}$ of the unipotent part of
$A_{\alpha_i}$,
and the nilpotent endomorphism
$N:=\bigoplus N_{\alpha_i}$ of $E$.
It is independent of the choice of 
an isomorphism
$E\simeq \bigoplus_i\LLtilde(\mu)\otimes \VVtilde(A_{\alpha_i})$
and $\alpha_i$.
We obtain the weight filtration $W$ of $E$ with respect to $N$.

Let $(E,\nbigf)\in \Vect^{ss}(T;\mu)^{\Par}$.
Each $\Gr^{\nbigf}_{\mu}(E)$
is equipped with the nilpotent endomorphism $N$ and $W$.
The following is clear by the construction.
\begin{lem}
The functor $\ttKtilde^C$
preserves the nilpotent endomorphism
and the weight filtrations.
\hfill\qed
\end{lem}

\subsubsection{Analytic case}

Let $\Vect^{\rnum\AHN}(T)^{\Par}$ denote 
the category of 
$(E,\gbigf)\in \Vect^{\rnum\AHN}(T)$
equipped with 
filtrations $\nbigf$ of $\Gr^{\gbigf}_{\mu}(E)$ for any
$\mu\in\rnum$
such that
$(\Gr^{\gbigf}_{\mu}(E),\nbigf)
 \in \Vect^{ss}(T;\mu)$.

For any $(\nbigv^{\an},\Phi^{\ast})\in \Diff(\nbigk^{an},\gminiq)$,
a good filtered bundle over $(\nbigv^{\an},\Phi^{\ast})$
means a good filtered bundle over $\gbigc(\nbigv^{\an},\Phi^{\ast})$.
Note that $\nbigr^{\an}$-lattices of $\nbigv^{\an}$
are equivalent to
$\nbigr$-lattices of $\gbigc(\nbigv^{\an})$.
Let $\Diff(\nbigk^{\an},\gminiq)^{\Par}$
denote the category of good filtered
$\gminiq$-difference $\nbigk^{\an}$-modules.

We obtain an equivalence
\[
\ttK^C:\Vect^{\rnum\AHN}(T)^{\Par}
\simeq
 \Diff(\nbigk^{\an},\gminiq)^{\Par}
\]
from the equivalence
$\ttK:\Vect^{\rnum\AHN}(T)
\simeq
 \Diff(\nbigk^{\an},\gminiq)$
and 
$\ttK^C:\Vect^{ss}(T;\mu)^{\Par}
\simeq
 \Diff(\nbigk,\gminiq;\mu)^{\Par}$.

\subsection{$\gminiq$-difference parabolic structure 
of sheaves on elliptic curves}
\label{subsection;19.2.9.110}

Let $\ttD\subset T=\cnum^{\ast}/\gminiq^{\seisuu}$ be a finite subset.
\begin{df}
Let $\ttEtilde$ be a locally free $\nbigo_{T}(\ast \ttD)$-module
of finite rank.
A $\gminiq$-difference parabolic structure of $\ttEtilde$ 
is data as follows:
\begin{itemize}
\item
A sequence
$s_{P,1}<s_{P,2}<\cdots<s_{P,\ttm(P)}$
in $\real$ for each $P\in \ttD$.

We formally set
$s_{P,0}:=-\infty$
and $s_{P,\ttm(P)+1}:=\infty$.
\item
A tuple of lattices
$\vecnbigk_P=\bigl(\nbigk_{P,i}\,|\,i=0,\ldots,\ttm(P)+1\bigr)$
of $\ttEtilde_{|\Phat}$.

Note that we obtain 
the lattice $\ttE_{-}\subset\ttEtilde$
determined by  $\nbigk_{P,0}$ $(P\in \ttD)$
and the lattice
$\ttE_{+}\subset\ttEtilde$
determined by 
$\nbigk_{P,\ttm(P)+1}$ $(P\in \ttD)$.
\item
Objects
$(\ttE_{\pm},\gbigf_{\pm},\nbigf_{\pm})
 \in\Vect^{\rnum\AHN}(T)^{\Par}$.
\end{itemize}
When we fix $(\vecs_P)_{P\in\ttD}$,
it is called a $\gminiq$-difference parabolic structure
at $(\ttD,(\vecs_P)_{P\in \ttD})$.
\hfill\qed
\end{df}

Let 
$\ttEtilde^{(i)}_{\ast}
=\bigl(
 \ttEtilde^{(i)},(\vecs_P,\vecnbigk^{(i)}_P)_{P\in \ttD},
 (\gbigf^{(i)}_{\pm},\nbigf^{(i)}_{\pm})
 \bigr)$
be locally free $\nbigo_T(\ast \ttD)$-modules of finite rank
with $\gminiq$-difference parabolic structure
at $(\vecs_P)_{P\in\ttD}$.
A morphism
$\ttEtilde^{(1)}_{\ast}\lrarr\ttEtilde^{(2)}_{\ast}$
is defined to be a morphism
$f:\ttEtilde^{(1)}\lrarr\ttEtilde^{(2)}$
of locally free $\nbigo_T(\ast \ttD)$-modules
such that the following holds:
\begin{itemize}
\item
$f(\nbigk^{(1)}_{P,i})\subset\nbigk^{(2)}_{P,i}$.
\item
The induced morphisms
$f:\ttE^{(1)}_{\pm}\lrarr \ttE^{(2)}_{\pm}$
are compatible with the filtrations
$(\gbigf_{\pm},\nbigf_{\pm})$,
i.e.,
they induce
$f:(\ttE^{(1)}_{\pm},\gbigf^{(1)},\nbigf^{(1)})
 \lrarr 
 (\ttE^{(2)}_{\pm},\gbigf^{(2)},\nbigf^{(2)})$
in $\Vect^{\rnum\AHN}(T)$.
\end{itemize}

Let $\Vect^{\gminiq}(T,(\vecs_P)_{P\in \ttD})$
denote the category of
locally free $\nbigo_T(\ast \ttD)$-modules of finite rank
with $\gminiq$-difference parabolic structure
at $(\ttD,(\vecs_P)_{P\in\ttD})$.

We define the degree of
$\ttEtilde_{\ast}=(\ttEtilde,(\vecs_P,\vecnbigk_P)_{P\in \ttD},
 (\gbigf_{\pm},\nbigf_{\pm}))$
as follows:
\begin{multline}
 \deg(\ttEtilde_{\ast}):=
 -\sum_{P\in D}
 \sum_{i=1}^{\ttm(P)}
 s_{P,i}\deg(\nbigk_{P,i},\nbigk_{P,i-1})
 \\
-\sum_{\omega\in\rnum}
 \sum_{b\in\real}
 b\rank\Gr^{\nbigf_-}_b\Gr^{\gbigf_-}_{\omega}(\ttE_-)
-\sum_{\omega\in\rnum}
 \sum_{b\in\real}
 b\rank\Gr^{\nbigf_+}_b\Gr^{\gbigf_+}_{\omega}(\ttE_+).
\end{multline}

\subsubsection{Rescaling of parabolic structure}

Let $\ttD\subset T$ be a finite subset.
Let 
$\ttEtilde_{\ast}=
 \bigl(
 \ttEtilde,(\vecs_{P},\vecnbigk_P)_{P\in\ttD},
 (\gbigf_{\pm},\nbigf_{\pm})
 \bigr)$
be a locally free $\nbigo_T(\ast \ttD)$-module
with $\gminiq$-difference parabolic structure.

Let $\gminit>0$.
We obtain a sequence
$\vecs^{(\gminit)}_P:=(\gminit s_{P,i})$.
We set
$\vecnbigk^{(\gminit)}_P:=\vecnbigk_P$
and 
$\gbigf^{(\gminit)}_{\pm}:=\gbigf_{\pm}$.
By setting
$\bigl(
 \nbigf^{(\gminit)}_{\pm}\bigr)_{\gminit a}
 \Gr^{\gbigf_{\pm}}(\ttE_{\pm}):=
 \bigl(\nbigf_{\pm}\bigr)_{a}
 \Gr^{\gbigf_{\pm}}(\ttE_{\pm})$,
we obtain filtrations
$\nbigf^{(\gminit)}_{\pm}$.
We set
\[
\ttH^{(\gminit)}(\ttEtilde_{\ast}):=
\bigl(
 \ttEtilde,
 (\vecs^{(\gminit)}_P,\vecnbigk^{(\gminit)}_{P})_{P\in \ttD},\,
 (\gbigf^{(\gminit)}_{\pm},\nbigf_{\pm}^{(\gminit)})
\bigr).
\]

Let $\gminit<0$.
We set
$s^{(\gminit)}_{P,i}:=\gminit s_{P,\ttm(P)-i+1}$.
We obtain a sequence $\vecs^{(\gminit)}_P$.
We set
$\nbigk^{(\gminit)}_{P,i}:=
 \nbigk_{P,\ttm(P)+1-i}$,
and we obtain a sequence
of lattices
$\vecnbigk^{(\gminit)}_{P}$.
We set
$\ttE^{(\gminit)}_{\pm}:=\ttE_{\mp}$.
Let
$\gbigf^{(\gminit)}_{\pm}(\ttE^{(\gminit)}_{\pm})$
denote the filtration induced by
$\gbigf_{\mp}(\ttE_{\mp})$.
We set
$\bigl(
 \nbigf^{(\gminit)}_{\pm}\bigr)_{|\gminit|a}
 \Gr^{\gbigf^{(\gminit)}_{\pm}}(\ttE^{(\gminit)}_{\pm})
=\bigl(
 \nbigf_{\mp}\bigr)_{a}
 \Gr^{\gbigf_{\mp}}(\ttE_{\mp})$.
Thus, we obtain 
\[
 \ttH^{(\gminit)}(\ttEtilde_{\ast}):=
 \bigl(
 \ttEtilde,
 (\vecs^{(\gminit)}_P,\vecnbigk^{(\gminit)}_P)_{P\in \ttD},\,\,
 (\gbigf^{(\gminit)}_{\pm},\nbigf^{(\gminit)})
 \bigr).
\]
The following is easy to check.
\begin{lem}
$\deg\bigl(\ttH^{(\gminit)}(\ttEtilde_{\ast})\bigr)
=|\gminit|\deg(\ttEtilde_{\ast})$.
\hfill\qed
\end{lem}

\subsection{Global correspondence for
parabolic $\gminiq$-difference modules}

\subsubsection{Parabolic $\gminiq$-difference modules}

Let $D\subset\cnum^{\ast}$ be a finite subset.
A parabolic structure of $\gminiq$-difference 
$\cnum[y,y^{-1}]$-module
is defined as in \S\ref{subsection;19.1.26.1}.
Let 
$\vecV^{(i)}_{\ast}=
 \bigl(\vecV^{(i)},V^{(i)},
 (\vect_{\alpha},\vecnbigl^{(i)}_{\alpha})_{\alpha\in D},
 \nbigp_{\ast}\vecV^{(i)}_{|\zerohat},
\nbigp_{\ast}\vecV^{(i)}_{|\inftyhat}\bigr)$ $(i=1,2)$
be $\gminiq$-difference
$\cnum[y,y^{-1}]$-modules
with good parabolic structure at infinity
and parabolic structure at 
$(D,(\vect_{\alpha})_{\alpha\in D})$.
A morphism
$\vecV^{(1)}_{\ast}\lrarr \vecV^{(2)}_{\ast}$
is defined to be 
a morphism of
$\gminiq$-difference $\cnum[y,y^{-1}]$-module
$f:\vecV^{(1)}\lrarr \vecV^{(2)}$
such that the following holds:
\begin{itemize}
\item
$f(V^{(1)})\subset V^{(2)}$.
\item
$f(\nbigl^{(1)}_{Q,i})\subset
 \nbigl^{(2)}_{Q,i}$.
\item
$f:\nbigp_{\ast}\vecVhat^{(1)}_{|\nuhat}
\lrarr
\nbigp_{\ast}\vecVhat^{(2)}_{|\nuhat}$
are induced for $\nu=0,\infty$.
\end{itemize}
Let $\Diff\bigl(\cnum[y,y^{-1}],
 \gminiq,(\vect_{\alpha})_{\alpha\in D}\bigr)^{\Par}$
be the category of
$\gminiq$-difference $\cnum[y,y^{-1}]$-modules
with good parabolic structure at infinity
and parabolic structure at $(\vect_{\alpha})_{\alpha\in D}$.

\subsubsection{An equivalence}

Let $\pi:\cnum^{\ast}\lrarr T:=\cnum^{\ast}/\gminiq^{\seisuu}$
denote the projection.
Let $\ttD\subset T$ be a finite subset.
For each $P\in \ttD$,
let $\vecs_P=(s_{P,1}<\cdots<s_{P,\ttm(P)})$ be a sequence
in $\real$.
For each $s_{P,i}$,
there exists $\alpha_{P,i}\in\pi^{-1}(P)\subset\cnum^{\ast}$
determined by the following conditions:
\[
0\leq
 s_{P,i}+\frac{\log|\alpha_{P,i}|}{\log|\gminiq|}
<1.
\]
We set
$u_{P,i}:=s_{P,i}+\frac{\log|\alpha_{P,i}|}{\log|\gminiq|}$.
We set
$D:=\coprod_{P\in \ttD}\{\alpha_{P,i}\,|\,i=1,\ldots,m(P)\}
\subset\cnum^{\ast}$.
For each $\alpha\in \pi^{-1}(P)\cap D$,
we set
$Z(\alpha):=
\{u_{P,i}\,|\,\alpha_{P,i}=\alpha\}
 \subset \closedopen{0}{1}$.
We obtain the sequence
$\vect_{\alpha}=\bigl(
0\leq t_{\alpha,0}<t_{\alpha,1}
<\cdots<t_{\alpha,m(\alpha)}<1\bigr)$
by ordering the elements of $Z(\alpha)$.
Let $i(\alpha)$
be determined by
$u_{P,i(\alpha)}=t_{\alpha,0}$.

Let us construct an equivalence
$\ttK:\Vect^{\gminiq}(T,(\vecs_{P})_{P\in \ttD})
\simeq
 \Diff\bigl(\cnum[y,y^{-1}],
 \gminiq,(\vect_{\alpha})_{\alpha\in D}\bigr)^{\Par}$.

Let
$\ttEtilde_{\ast}=(\ttEtilde,(\vect_P,\vecnbigk_P)_{P\in \ttD},
 (\gbigf_{\pm},\nbigf_{\pm}))$.
Let $\nbige$ be the locally free
$\nbigo_{\cnum^{\ast}}(\ast \pi^{-1}(D))$-module
obtained as the pull back of $\ttEtilde$.
For each $\alpha\in D\cap\pi^{-1}(P)$,
we obtain a lattice
$\nbigl_{\alpha}$
of $\nbige_{|\widehat{\alpha}}$
induced by $\nbigk_{P,i(\alpha)-1}$.
We obtain a locally free
$\nbigo_{\cnum^{\ast}}$-submodule 
$\nbigv\subset\nbige$
determined by
$\nbigl_{\alpha}$ $(\alpha\in \pi^{-1}(\ttD))$.
It is extended 
to a filtered bundle $\nbigp_{\ast}\nbigv$
on $(\proj^1,\{0,\infty\})$
by $(\gbigf_{\pm},\nbigf_{\pm})$
by using the functors $\ttK^1$.
We set $V:=H^0(\proj^1,\nbigp\nbigv)$,
which is $\cnum[y,y^{-1}]$-free module
of finite rank.
We set
$\vecVtilde:=V\otimes\cnum(y)$,
which is naturally a $\gminiq$-difference 
$\cnum(y)$-module.
Let $\vecV$ be the $\gminiq$-difference
$\cnum[y,y^{-1}]$-submodule
of $\vecVtilde$ generated by $V$.
For each $\alpha\in D\cap \pi^{-1}(P)$,
we obtain the lattices
$\nbigl_{\alpha,j}$ $(1\leq j\leq m(\alpha)-1)$
of $\nbigv(\ast\alpha)_{|\widehat{\alpha}}$
induced by
$\nbigk_{P,i(\alpha)+j-1}$.
We also obtain 
good filtered bundles
$\nbigp_{\ast}\vecV_{|\nuhat}$ $(\nu=0,\infty)$
over $\vecV_{\nuhat}$ from $\nbigp_{\ast}\nbigv$.
Thus, we obtain 
\[
 \vecV_{\ast}=
 \bigl(
 \vecV,V,(\vect_{\alpha},\vecnbigl_{\alpha})_{\alpha\in D},
 (\nbigp_{\ast}\vecV_{|\widehat{0}},\nbigp_{\ast}\vecV_{|\widehat{\infty}})
 \bigr)
 \in 
\Diff\bigl(\cnum[y,y^{-1}],\gminiq,
 (\vect_{\alpha})_{\alpha\in D}\bigr)^{\Par}.
\]
The following is clear by the construction.
\begin{prop}
\label{prop;19.2.9.100}
$\ttK$ induces an equivalence
$\Vect^{\gminiq}(T,(\vecs_{P})_{P\in\ttD})
\simeq
 \Diff\bigl(\cnum[y,y^{-1}],
 \gminiq,(\vect_{\alpha})_{\alpha\in D}\bigr)^{\Par}$.
Moreover,
$\deg(\ttK(\ttEtilde_{\ast}))
=\deg(\ttEtilde_{\ast})$ holds.
As a result,
the equivalence preserves the stable objects,
semistable objects and polystable objects.
\hfill\qed
\end{prop}

\subsection{Filtrations and growth orders of norms}
\label{subsection;19.2.4.31}

Let us consider the action of
$\seisuu\tte_2$ 
on $\nbigm_{\gminiq}^{\cov}:=\cnum^{\ast}\times\real$
by $\tte_2(y,t)=(\gminiq y,t+1)$.
It is extended to the action of $\seisuu\tte_2$
on $\nbigmbar^{\cov}_{\gminiq}:=\proj^1\times\real$.
Let $\nbigm_{\gminiq}$
and $\nbigmbar_{\gminiq}$
denote the quotient spaces of 
$\nbigm_{\gminiq}^{\cov}$
and $\nbigmbar^{\cov}_{\gminiq}$
by the action of $\seisuu\tte_2$,
respectively.
For $\nu=0,\infty$,
let $H_{\nu}$ denote the image of
$H_{\nu}^{\cov}:=
 \{\nu\}\times \real\lrarr \nbigmbar_{\gminiq}$.

Let $\nu$ denote $0$ or $\infty$.
Let $\nbigubar_{\nu}$
be a neighbourhood of
$H_{\nu}$ 
in $\nbigmbar_{\gminiq}$.
We set
$\nbigu_{\nu}:=
 \nbigubar_{\nu}\setminus
 H_{\nu}$.
Let $\nbigu^{\cov}_{\nu}$
denote the pull back of $\nbigu_{\nu}$
by $\nbigm_{\gminiq}^{\cov}\lrarr\nbigm_{\gminiq}$.
Similarly,
let $\nbigubar^{\cov}_{\nu}$
denote the pull back of $\nbigubar_{\nu}$
by $\nbigmbar^{\cov}_{\gminiq} \lrarr\nbigmbar_{\gminiq}$.
We set $y_{0}:=y$
and $y_{\infty}:=y^{-1}$.
We set
$\gminiq_{0}:=\gminiq$
and $\gminiq_{\infty}:=\gminiq^{-1}$.

\subsubsection{Equivalences}

Let $\LFM(\nbigubar_{\nu},H_{\nu})$
denote the category of locally free
$\nbigo_{\nbigubar_{\nu}}(\ast H_{\nu})$-modules.
We obtain an equivalence
$\Upsilon:
  \Diff(\nbigk^{\an},\gminiq_{\nu})
\simeq
\LFM(\nbigubar_{\nu},H_{\nu})$
as in the formal case.
(See \S\ref{subsection;19.1.19.1}.)
Hence, we obtain the following equivalence:
\[
 \ttK_{\nu}:\Vect^{\rnum\AHN}(T)
\simeq
 \LFM(\nbigubar_{\nu},H_{\nu}).
\]

Let $\LFM(\nbigubar_{\nu},H_{\nu})^{\Par}$
denote the category of 
good filtered bundles 
over $(\nbigubar_{\nu},H_{\nu})$.
By the definition of good filtered bundles,
we obtain the following equivalence:
\[
  \ttK_{\nu}:\Vect^{\rnum\AHN}(T)^{\Par}
\simeq
 \LFM(\nbigubar_{\nu},H_{\nu})^{\Par}.
\]

\subsubsection{Metrics and slope filtrations}

Let $(\ttE,\gbigf)\in\Vect^{\rnum\AHN}(T)$.
We obtain
$\gbigv_{\nu}:=\ttK_{\nu}(\ttE,\gbigf)\in \LFM(\nbigubar_{\nu},H_{\nu})$.
Let $h_{\nu}$ be a Hermitian metric of 
$\gbigv_{\nu|\nbigu_{\nu}}$
such that the following holds.
\begin{itemize}
\item
 Let $P$ be any point of $H_{\nu}$.
Let $\vecv$ be a frame of $\gbigv_{\nu}$
on a neighbourhood $U_P$ of $P$
 in $\nbigubar_{\nu}$.
Let $H(\vecv)$ be the Hermitian matrix valued
function on $U_P\setminus H_{\nu}$
determined by
$H(\vecv)_{i,j}=h_{\nu}(v_i,v_j)$.
Then, there exists $C>1$ and $N>0$
such that 
$C^{-1}|y_{\nu}|^{N} 
\leq H(\vecv)\leq
 C|y_{\nu}|^{-N}$.
\end{itemize}
It is easy to construct such a Hermitian metric $h_{\nu}$.

Let $\gbigv^{\cov}$ be the pull back of $\gbigv$
by $\nbigubar^{\cov}_{\nu}\lrarr\nbigubar_{\nu}$.
Let $h^{\cov}_{\nu}$ be the metric of
$\gbigv^{\cov}_{|\nbigu^{\cov}_{\gminiq}}$
induced by $h$.
Let $Q$ be any point of $T$.
We take $\alpha_0\in\cnum^{\ast}$,
which is mapped to $Q$
by $\pi:\cnum^{\ast}\lrarr T=\cnum^{\ast}/\gminiq^{\seisuu}$.
Set $\alpha_{\infty}=(\gminiq_0^n\alpha_0)^{-1}$
for an appropriate $n\in\seisuu$.
We may assume that the half line 
$\ttl_{\nu}:=\{(\alpha_{\nu},t)\,|\, \varrho(\gminiq_{\nu})t\geq 0\}$
is contained in
$\nbigu^{\cov}_{\nu}$.
For each $s\in \ttE_{|Q}$,
we obtain a flat section $\stilde_{\nu}$
of $\gbigv^{\cov}_{\nu}$ along $\ttl_{\nu}$.
The following is easy to see.

\begin{lem}
\label{lem;19.1.19.11}
$s$ is contained in $\gbigf_{\mu}(\ttE)_{|Q}$
if and only if the following holds for any $\epsilon>0$:
\[
 \log|\stilde_{\nu}|_{h^{\cov}_{\nu}}=O\left(
 \frac{\mu}{2}
 \Bigl|\log|\gminiq_{\nu}|\Bigr|
 \Bigl(
 t-\frac{\log|\alpha_{\nu}|}{\log|\gminiq_{\nu}|}
 \Bigr)^2
+\epsilon 
 \Bigl(
 t-\frac{\log|\alpha_{\nu}|}{\log|\gminiq_{\nu}|}
 \Bigr)^2
\right).
\]
More strongly,
for any $s\in\gbigf_{\mu}\setminus\gbigf_{<\mu}$,
the following holds:
\[
 \log|\stilde_{\nu}|_{h^{\cov}_{\nu}}=
 \frac{\mu}{2}
 \Bigl|
 \log|\gminiq_{\nu}|
 \Bigr|
 \Bigl(
 t-\frac{\log|\alpha_{\nu}|}{\log|\gminiq_{\nu}|}
 \Bigr)^2
+O\Bigl(
 t\log|\gminiq_{\nu}|-\log|\alpha_{\nu}|
 \Bigr).
\]
\hfill\qed
\end{lem}

\subsubsection{Refinement}

Let $(\ttE,\gbigf,\nbigf)\in\Vect^{\AHN}(T)^{\Par}$.
We set 
$\nbigp_{\ast}\gbigv_{\nu}:=
 \ttK^1_{\nu}(\ttE,\gbigf,\nbigf)$.
Suppose that 
$h_{\nu}$ 
is adapted to $\nbigp_{\ast}\gbigv_{\nu}$.
For $s\in \gbigf_{\mu}(\ttE)_{|Q}$,
let $[s]$ denote the induced element of
$\Gr^{\gbigf}_{\mu}(\ttE)_{|Q}$.

\begin{lem}
\label{lem;19.1.19.12}
$[s]\in \nbigf_b\Gr^{\gbigf}_{\mu}(\ttE)_{|Q}$
if and only if the following holds for any $\epsilon>0$:
\[
 \log|\stilde|_{h^{\cov}_{\nu}}=
 O\left(
 \frac{\mu}{2}
 \Bigl|\log|\gminiq_{\nu}|\Bigr|
 \Bigl(
 t-\frac{\log|\alpha_{\nu}|}{\log|\gminiq_{\nu}|}
 \Bigr)^2
+(b+\epsilon)
\Bigl(
t\log|\gminiq_{\nu}|-\log|\alpha_{\nu}|
\Bigr)
\right).
\]
\hfill\qed
\end{lem}

Let $W_k\nbigf_b\Gr^{\gbigf}_{\mu}(E)$
denote the inverse image of
$W_k\Gr^{\nbigf}_b\Gr^{\gbigf}_{\mu}(E)$
by the surjection
$\nbigf_b\Gr^{\gbigf}_{\mu}(E)
\lrarr
 \Gr^{\nbigf}_b\Gr^{\gbigf}_{\mu}(E)$.

\begin{lem}
Suppose moreover that
the norm estimate holds for
$(\nbigp_{\ast}\gbigv,h)$.
Then, $[s]\in W_k\Gr^{\nbigf}_b\Gr^{\gbigf}_{\mu}(E)$
if and only if
the following holds:
\[
 \log|\stilde|_{h^{\cov}_{\nu}}=
 O\left(
 \frac{\mu}{2}
 \Bigl|\log|\gminiq_{\nu}|\Bigr|
 \Bigl(
 t-\frac{\log|\alpha_{\nu}|}{\log|\gminiq_{\nu}|}
 \Bigr)^2
+b
\Bigl(
 t\log|\gminiq_{\nu}|
-\log|\alpha_{\nu}|
\Bigr)
+\frac{k}{2}
 \log
\Bigl(
 t\log|\gminiq_{\nu}|
-\log|\alpha_{\nu}|
\Bigr)
\right).
\]
\hfill\qed
\end{lem}

\subsection{Filtered objects on elliptic curves associated to monopoles}

\subsubsection{Induced filtered objects on the elliptic curve}

We use the notation in \S\ref{subsection;18.8.19.20}.
Suppose that $|\lambda|\neq 1$.
We set 
$T^{\lambda}:=\cnum_{\ttU}^{\ast}/(\gminiq^{\lambda})^{\seisuu}$.
Let $\pi:\nbigm^{\lambda}\lrarr T^{\lambda}$
denote the morphism
induced by 
$\cnum_{\ttU}^{\ast}\times
 \real_{\ttt}\lrarr\cnum_{\ttU}^{\ast}$.
Let $Z\subset\nbigm^{\lambda}$ be a finite subset.
We set $\ttD:=\pi(Z)\subset T^{\lambda}$.
Note that the function $\ttU$ on $\nbigm^{\lambda}$
is independent of the choice of $(\tte_1,\tts_1)$,
but $\ttt$ depends on $(\tte_1,\tts_1)$.
Hence, we use the notation
$\ttt(\tte_1,\tts_1)$ to emphasize 
the dependence on $(\tte_1,\tts_1)$.
Similarly, we use the notation $\nbigmbar^{\lambda}_{(\tte_1,\tts_1)}$ 
to denote $\nbigmbar^{\lambda}$ in \S\ref{subsection;18.8.19.20}
to emphasize the dependence on $(\tte_1,\tts_1)$.
The sets $H^{\lambda}$
are also denoted by $H^{\lambda}_{(\tte_1,\tts_1)}$.
The number $\gminit^{\lambda}$
is denoted by $\gminit^{\lambda}(\tte_1,\tts_1)$.
Let us denote $\gminiq^{\lambda}$ by $\gminiq^{\lambda}(\tte_1)$
to emphasize the dependence on $\tte_1$.

Let $(E,h,\nabla,\phi)$ be
a meromorphic monopole 
on $\nbigm^{\lambda}\setminus Z$.
We obtain a good filtered bundle
with Dirac type singularity
$\nbigp_{\ast}E_{(\tte_1,\tts_1)}$ 
on
$(\nbigmbar_{(\tte_1,\tts_1)}^{\lambda};H_{(\tte_1,\tts_1)}^{\lambda},Z)$.
It is equivalent to
a parabolic $\gminiq^{\lambda}(\tte_1)$-difference
$\cnum[\ttU,\ttU^{-1}]$-module.
Let $\ttEtilde_{(\tte_1,\tts_1)\,\ast}$
denote the corresponding
locally free $\nbigo_{T^{\lambda}}(\ast \ttD)$-module
with a $\gminiq^{\lambda}(\tte_1)$-difference parabolic structure.
(See Proposition \ref{prop;19.2.9.100}.)
By rescaling the parabolic structure,
we obtain a locally free $\nbigo_{T^{\lambda}}(\ast \ttD)$-module
with a $\gminiq^{\lambda}(\tte_1)$-difference parabolic structure
$\ttH^{(\gminit^{\lambda}(\tte_1,\tts_1))}(\ttEtilde_{(\tte_1,\tts_1)\,\ast})$.

\begin{thm}
\label{thm;19.2.8.130}
$\ttH^{(\gminit^{\lambda}(\tte_1,\tts_1))}(\ttEtilde_{(\tte_1,\tts_1)\,\ast})$
is independent of the choice of
$(\tte_1,\tts_1)$.
\end{thm}
\pf
Recall that
the filtered object
$\ttEtilde_{(\tte_1,\tts_1)\ast}$
consists of
\begin{itemize}
\item
a locally free
$\nbigo_{T^{\lambda}}(\ast\ttD)$-module
$\ttEtilde_{(\tte_1,\tts_1)}$,
\item
a tuple
$(\vecs_{P}(\tte_1,\tts_1),\vecnbigl_P(\tte_1,\tts_1))
 _{P\in\ttD}$,
\item
filtrations
$\gbigf_{\pm}(\tte_1,\tts_1)$
on $\ttEtilde_{(\tte_1,\tts_1),\pm}$,
\item
filtrations
$\nbigf_{\pm}(\tte_1,\tts_1)$
on $\Gr^{\gbigf_{\pm}(\tte_1,\tts_1)}(\ttEtilde_{(\tte_1,\tts_1),\pm})$.
\end{itemize}
(See \S\ref{subsection;19.2.9.110}.)
We have the isomorphism
$f_{(\tte_1,\tts_1)}:\nbigm^{\lambda}
\simeq
\nbigm_{\gminiq^{\lambda}(\tte_1)}$
induced by
\[
 (\ttU,\ttt(\tte_1,\tts_1))\longmapsto 
(\ttU,\ttt(\tte_1,\tts_1)/\gminit^{\lambda}(\tte_1,\tts_1)).
\]
Note that 
$\ttEtilde_{(\tte_1,\tts_1)}$
depend only on 
the mini-holomorphic bundle
$(E,\delbar_{E})$ on 
$\nbigm^{\lambda}\setminus Z\simeq
 \nbigm_{\gminiq}\setminus f_{(\tte_1,\tts_1)}(Z)$
underlying the monopole $(E,h,\nabla,\phi)$.
Hence, they are independent of $(\tte_1,\tts_1)$.
According to Lemma \ref{lem;19.2.8.200},
\[
\ttt(\tte_1,\tts_1)-
\gminit^{\lambda}(\tte_1,\tts_1)
\frac{\log|\ttU|}{\log|\gminiq^{\lambda}(\tte_1,\tts_1)|}
\]
is independent of $(\tte_1,\tts_1)$.
Therefore, we obtain that
the sequence
$\bigl(
 \vecs(\tte_1,\tts_1)_P
 \bigr)^{(\gminit^{\lambda}(\tte_1,\tts_1))}$
and $\vecnbigl_P(\tte_1,\tts_1)^{(\gminit^{\lambda})}$
are independent of 
$(\tte_1,\tts_1)$.
According to 
Lemma \ref{lem;19.1.19.11}
and 
Lemma \ref{lem;19.1.19.12},
the filtrations $\gbigf_{\pm}(\tte_1,\tts_1)$
and $\nbigf_{\pm}(\tte_1,\tts_1)$
are characterized by the growth order of
the norms of the $\del_{\ttt}$-flat sections
with respect to $h$.
Then, 
Lemma \ref{lem;19.2.8.201}
and Lemma \ref{lem;19.2.8.200}
imply that 
the filtrations 
$(\gbigf^{(\gminit)}_{\pm},\nbigf^{(\gminit)}_{\pm})$
are independent of
$(\tte_1,\tts_1)$.
Thus, we obtain Theorem \ref{thm;19.2.8.130}.
\hfill\qed